\documentclass{memo-l}

\usepackage{amsthm}
\usepackage{amssymb}
\usepackage{amsfonts}
\usepackage{amsxtra}     
\usepackage{epsfig}
\usepackage{verbatim}
\usepackage{xypic}
\usepackage{xspace}
\usepackage{layout}
\usepackage{longtable}
\usepackage[all]{xy}
\usepackage{url}    
\input xypic
\input xy
\xyoption{all}

\newcommand{\ie}{\textit{i.e.\xspace}}
\newcommand{\eg}{\textit{e.g.\xspace}}
\newcommand{\etc}{\textit{etc.}\xspace}
\newcommand{\fig}[2]{
\begin{figure}[#1]
#2
\end{figure}}

\makeatletter
\renewcommand{\maketag@@@}[1]{\hbox{\m@th\normalsize\normalfont#1}}
\makeatother

\makeatletter
\def\@strippedMR{}
\def\@scanforMR#1#2#3\endscan{%
  \ifx#1M\ifx#2R\def\@strippedMR{#3}%
  \else\def\@strippedMR{#1#2#3}%
  \fi\fi}
\renewcommand\MR[1]{\relax\ifhmode\unskip\spacefactor3000 \space\fi
  \@scanforMR#1\endscan
  MR\MRhref{\@strippedMR}{\@strippedMR}}
\makeatother

\def\square{\hfill${\vcenter{\vbox{\hrule height.4pt \hbox{\vrule width.4pt
height7pt \kern7pt \vrule width.4pt} \hrule height.4pt}}}$}

\newenvironment{pf}{{\it Proof:}\quad}{\square \vskip 12pt}

\theoremstyle{plain}
\newtheorem{thm}{Theorem}[chapter]

\newtheorem{lem}[thm]{Lemma}

\theoremstyle{definition}
\newtheorem{defn}[thm]{Definition}
\newtheorem{notation}[thm]{Notation}

\newtheorem{examp}[thm]{Example}

\newtheorem{rem}[thm]{Remark}
\newtheorem{rmk}[thm]{Remark}

\theoremstyle{remark}

\numberwithin{section}{chapter}
\numberwithin{equation}{chapter}

\begin{document}
\frontmatter
\title{Pseudo Limits, Biadjoints, and Pseudo
Algebras: \linebreak Categorical Foundations of Conformal Field
Theory}

\author{Thomas M. Fiore}
\address{Department of Mathematics \\ University of Michigan \\ Ann Arbor, MI 48109-1109}
\curraddr{Department of Mathematics \\ University of Chicago \\ 5734
S. University Avenue \\ Chicago, IL 60637}

\email{fiore@math.uchicago.edu}
\date{April 2nd, 2004}



\begin{abstract}
In this paper we develop the categorical foundations needed for
working out completely the rigorous approach to the definition of
conformal field theory\index{conformal field theory} outlined by
Graeme Segal\index{Segal, Graeme}. We discuss pseudo
algebras\index{algebra!pseudo algebra} over theories\index{theory}
and 2-theories\index{2-theory}, their pseudo
morphisms\index{morphism}, bilimits\index{bilimit},
bicolimits\index{bicolimit}, biadjoints\index{biadjoint},
stacks\index{stack}, and related concepts.

These 2-categorical concepts are used to describe the algebraic
structure on the class of rigged surfaces. A {\it rigged
surface}\index{rigged surface} is a real, compact, not necessarily
connected, two dimensional manifold\index{manifold} with complex
structure\index{complex structure} and analytically parametrized
boundary components\index{boundary components}. This class admits
algebraic operations of {\it disjoint union}\index{disjoint union}
and {\it gluing}\index{gluing} as well as a {\it unit}\index{unit}.
These operations satisfy axioms such as unitality\index{unit} and
distributivity\index{distributivity} up to coherence
isomorphisms\index{coherence isomorphism} which satisfy coherence
diagrams\index{coherence diagram}. These operations, coherences, and
their diagrams are neatly encoded as a {\it pseudo algebra over the
2-theory of commutative monoids with
cancellation}\index{algebra!pseudo
algebra}\index{cancellation}\index{2-theory!2-theory of commutative
monoids with cancellation}. A {\it conformal field
theory}\index{conformal field theory} is a morphism of stacks of
such structures\index{stack}.

This paper begins with a review of 2-categorical concepts,
Lawvere\index{Lawvere} theories, and algebras\index{algebra} over
Lawvere theories. We prove that the 2-category of small pseudo
algebras\index{algebra!pseudo algebra} over a theory\index{theory}
admits weighted pseudo limits\index{limit!weighted pseudo limit} and
weighted\index{weighted} bicolimits\index{bicolimit!weighted
bicolimit}. This 2-category is biequivalent to the 2-category of
algebras over a 2-monad\index{2-monad} with pseudo morphisms. We
prove that a pseudo functor admits a left
biadjoint\index{biadjoint!left biadjoint} if and only if it admits
certain biuniversal arrows. An application of this theorem implies
that the forgetful 2-functor\index{forgetful
2-functor}\index{2-functor!forgetful 2-functor} for pseudo
algebras\index{algebra!pseudo algebra} admits a left biadjoint. We
introduce stacks\index{stack} for Grothendieck
topologies\index{Grothendieck topology} and prove that the
traditional definition of stacks\index{stack} in terms of descent
data\index{descent data} is equivalent to our definition via
bilimits\index{bilimit}. The paper ends with a proof that the
2-category of pseudo algebras\index{algebra!pseudo algebra} over a
2-theory admits weighted\index{weighted} pseudo
limits\index{limit!pseudo limit}. This result is relevant to the
definition of conformal field theory\index{conformal field theory}
because bilimits are necessary to speak of stacks.

\end{abstract}

\maketitle




\setcounter{page}{3} \cleardoublepage \thispagestyle{empty}
\vspace*{13.5pc}
\begin{center}
In memory of my Mother
\end{center}
\clearpage



\setcounter{page}{6}

\tableofcontents


\chapter*{Acknowledgements}

It is with great pleasure that I acknowledge the many people who
have aided me in the creation of this book. I am deeply grateful to
Igor Kriz for his careful guidance and encouragement. My gratitude
extends to Po Hu, F. W. Lawvere, Ross Street, Steve Lack, John Baez,
Tibor Beke, Bob Bruner, James McClure, Jeff Smith, Art Stone, Martin
Hyland, John Power, Michael Johnson, Mark Weber, Craig Westerland,
and Bart Kastermans for helpful comments. Barbara Beeton enhanced
the format of this document through her valuable typesetting advice.

This research was generously supported through a VIGRE grant of the
National Science Foundation. The Mathematics Department of the
University of Michigan and the Horace H. Rackham Graduate School
also provided assistance.

Special thanks go to my wife Eva Ackermann and to my parents.

\mainmatter

\chapter{Introduction}
\label{sec:introduction} The purpose of this paper is to work out
the categorical basis for the foundations of conformal field
theory\index{conformal field theory}. The definition of conformal
field theory was outlined in Segal\index{Segal, Graeme}
\cite{segal1} and recently given in \cite{hu} and \cite{hu1}.
Concepts of 2-category theory, such as versions of
algebra\index{algebra}, limit\index{limit}, colimit\index{colimit},
and adjunction\index{adjunction}, are necessary for this definition.

The structure present on the class $\mathcal{C}$ of rigged
surfaces\index{rigged surface} is captured by these concepts of
2-category theory. Here a {\it rigged surface}\index{rigged surface}
is a real, compact, not necessarily connected, two dimensional
manifold\index{manifold} with complex structure\index{complex
structure} and analytically parametrized boundary
components\index{boundary components}. Isomorphisms of such rigged
surfaces\index{rigged surface} are holomorphic\index{holomorphic}
diffeomorphisms\index{diffeomorphism} preserving the boundary
parametrizations\index{boundary
parametrization}\index{parametrization}. These rigged
surfaces\index{rigged surface} and isomorphisms form a groupoid and
are part of the structure present on $\mathcal{C}$. Concepts of
2-categories enter when we describe the operations of disjoint
union\index{disjoint union} of two rigged surfaces\index{rigged
surface} and gluing\index{gluing} along boundary
components\index{boundary components} of opposite
orientation\index{orientation}. We need a mathematical structure to
capture all of these features. This has been done in \cite{hu}.

One step in this direction is the notion of
algebra\index{algebra!algebra over a theory} over a
theory\index{theory} in the sense of Lawvere
\cite{lawvere}\index{Lawvere theory}. We need a weakened notion in
which relations are replaced by coherence isos\index{coherence
isomorphism}. This weakened notion is called a {\it pseudo algebra}
in this paper\index{algebra!pseudo algebra}. Coherence
diagrams\index{coherence diagram} are required in a pseudo
algebra\index{algebra!pseudo algebra}, but it was noticed in
\cite{hu} that Lawvere's\index{Lawvere theory} notion of a
theory\index{theory} allows us to write down all such diagrams
easily. See Chapter \ref{sec:laxTalgebras} below. A symmetric
monoidal category\index{symmetric monoidal category} as defined in
\cite{maclane1} provides us with a classical example of a pseudo
algebra\index{algebra!pseudo algebra} over the theory of commutative
monoids\index{commutative monoid!theory of commutative
monoids}\index{theory!theory of commutative monoids}. Theories,
duality, and related topics are discussed further in \cite{adamek},
\cite{adamek1}, \cite{adamek2}, \cite{lawvere1}, and
\cite{lawvere2}.

Unfortunately, pseudo algebras\index{algebra!pseudo algebra} over a
theory are not enough to capture the structure on $\mathcal{C}$. The
reason is that the operation of gluing\index{gluing} is
indexed\index{indexed} by the variable set of pairs of boundary
components\index{boundary components} of opposite
orientation\index{orientation}.  The operation of disjoint
union\index{disjoint union} also has an indexing\index{indexing}. We
need pseudo algebras\index{algebra!pseudo algebra} over a ``theory
indexed over another theory,''\index{theory!theory indexed over a
theory} which we call a 2-theory\index{2-theory}. More precisely,
the pseudo algebras we need are pseudo algebras over the 2-theory of
{\it commutative monoids with
cancellation}\index{cancellation}\index{2-theory!2-theory of
commutative monoids with cancellation}. See \cite{hu} and Chapter
\ref{sec:2-theories} below. The term 2-theory does {\it not} mean a
theory in 2-categories.

Nevertheless, 2-categories are relevant. This is because we want to
capture the behavior of holomorphic\index{holomorphic} families of
rigged surfaces\index{rigged surface!holomorphic families of rigged
surfaces}\index{holomorphic families of rigged surfaces} in our
description of the structure of $\mathcal{C}$. This amounts to
saying that $\mathcal{C}$ is a stack\index{stack} of pseudo
commutative monoids with
cancellation\index{cancellation}\index{commutative monoid with
cancellation!pseudo commutative monoid with cancellation}. To
consider this, we must remark that pseudo algebras over a
theory\index{theory} and pseudo algebras\index{algebra!pseudo
algebra} over a 2-theory form 2-categories. A {\it
stack}\index{stack} is a contravariant\index{contravariant} pseudo
functor\index{functor!pseudo functor} from a Grothendieck
site\index{Grothendieck site} into a 2-category which takes
Grothendieck covers\index{Grothendieck cover} into
limits\index{limit} of certain type, which are called
bilimits\index{bilimit}. They are defined below, in \cite{kelly2},
and \cite{street3}, while a slightly stronger notion is called
pseudo limit\index{limit!pseudo limit} in \cite{street3}. One needs
to understand such notions for the rigorous foundations of conformal
field theory\index{conformal field theory}. More elaborate notions,
such as analogous kinds of colimits are also needed in \cite{hu1}.

In this article we introduce the general concepts of weighted
bilimits\index{bilimit}, weighted\index{weighted}
bicolimits\index{bicolimit}, and biadjoints\index{biadjoint} for
pseudo functors between 2-categories in the sense below and prove
statements about their existence in certain cases. There are many
versions of such concepts and many (but not all) of the theorems we
give are in the literature, see \cite{bird}, \cite{borceux2},
\cite{carboni1}, \cite{gray1}, \cite{gray3}, \cite{gray4},
\cite{gray5}, \cite{kelly}, \cite{kelly3}, \cite{kelly2},
\cite{simmons}, \cite{street1}, \cite{street2}, \cite{street3}, and
\cite{street4}. Bicategories\index{bicategory} were first introduced
in \cite{benabou1} and \cite{ehresmann}. The circumstances of
conformal field theory\index{conformal field theory} suggest a
particular choice of concepts. To a topologist, the most natural and
naive choice of terminology may be to use the term
``lax''\index{lax|textbf} to mean ``up to coherence
isos''\index{coherence isomorphism} with these coherence isos
required to satisfy appropriate coherence diagrams\index{coherence
diagram}. ``Iso''\index{iso} seems to be the only natural concept in
the case of pseudo algebras over a theory\index{theory}: there seems
to be no reasonable notion where coherences would not be
iso\index{algebra!pseudo algebra}. For this reason, the authors of
\cite{hu}, \cite{hu1}, and \cite{hu2} use the ``lax=up to coherence
isos'' philosophy. This terminology however turns out to be
incorrect from the point of view of category theory (other ad hoc
terminology also appears in \cite{hu}, \cite{hu1}, and \cite{hu2}).
In this paper, we decided to follow established categorical
terminology while giving a precise translation of the notions in
\cite{hu}, \cite{hu1}, and \cite{hu2}. In the established
categorical terminology, what is called a lax
algebra\index{algebra!lax algebra} in \cite{hu}, \cite{hu1}, and
\cite{hu2} is called a {\it pseudo algebra}\index{algebra!pseudo
algebra}, what is called a lax morphism\index{morphism!lax morphism}
(morphism which commutes with operations up to coherence isos) in
\cite{hu}, \cite{hu1}, and \cite{hu2} is called a {\it pseudo
morphism} (or just {\it morphism}), and what is called a lax
functor\index{functor!lax functor} in \cite{hu}, \cite{hu1}, and
\cite{hu2} is called a {\it pseudo functor}. In addition, the
notions which the authors of \cite{hu}, \cite{hu1}, and \cite{hu2}
refer to as lax\index{lax|textbf} limit\index{limit!lax limit}, lax
colimit\index{colimit!lax colimit}, and lax
adjoint\index{adjoint!lax adjoint} are called {\it
bilimit}\index{bilimit}, {\it bicolimit}\index{bicolimit}, and {\it
biadjoint}\index{biadjoint} in established categorical terminology.
The stronger categorical notions of pseudo limit\index{limit!pseudo
limit}, pseudo colimit\index{colimit!pseudo colimit}, and pseudo
adjoint\index{adjoint!pseudo adjoint} are also sometimes relevant.

The term ``lax''\index{lax|textbf} in standard categorical
terminology is reserved for notions ``up to 2-cells\index{2-cell}
which are not necessarily iso''\index{iso}. However, such notions
will not play a central role in the present paper, as our motivation
here is the same as in \cite{hu}, \cite{hu1}, and \cite{hu2}, namely
conformal field theory\index{conformal field theory} and
stacks\index{stack}.

We show that every pseudo functor\index{functor!pseudo functor} from
a 1-category to the 2-category of small categories admits both a
pseudo limit\index{limit!pseudo limit} and a pseudo
colimit\index{colimit!pseudo colimit} by constructive proofs.
Furthermore, the 2-category of small categories admits
weighted\index{weighted} pseudo limits\index{limit!weighted pseudo
limit} and weighted pseudo colimits\index{colimit!weighted pseudo
colimit}. After that we introduce the notions of a
theory\index{theory}, an algebra over a theory\index{algebra!algebra
over a theory}, and a pseudo algebra\index{algebra!pseudo algebra}
over a theory. We then go on to show that any pseudo functor from a
1-category to the 2-category of pseudo $T$-algebras admits a pseudo
limit by an adaptation of the proof for small categories. After a
proof of the existence of cotensor products\index{cotensor product}
in the 2-category of pseudo $T$-algebras, we conclude from a theorem
of Street that this 2-category admits weighted\index{weighted}
pseudo limits\index{limit!weighted pseudo limit}.

We continue the study of weakened structures by turning to
biadjoints\index{biadjoint}. First we show that a pseudo functor
admits a left biadjoint\index{biadjoint!left biadjoint} if and only
if for each object of the source category we have an appropriate
biuniversal arrow\index{arrow!biuniversal arrow}\index{biuniversal
arrow} in analogy to the standard result in 1-category theory. By
means of this description we show that for any morphism of
theories\index{morphism!morphism of theories} $\phi:S \rightarrow T$
the associated forgetful 2-functor\index{forgetful
2-functor}\index{2-functor!forgetful 2-functor} from the 2-category
of pseudo $T$-algebras to the 2-category of pseudo $S$-algebras
admits a left biadjoint. The formalism developed for biadjoints is
then adapted to treat bicolimits of pseudo $T$-algebras. Moreover,
the universal property of these bicolimits\index{bicolimit} is
slightly weaker than the universal property of the pseudo
limits\index{limit!pseudo limit}. Similarly, the 2-category of
pseudo $T$-algebras admits bitensor products\index{bitensor
product}, and hence also weighted\index{weighted}
bicolimits\index{bicolimit!weighted bicolimit}.

Lastly, we construct pseudo limits\index{limit!pseudo limit} of
pseudo algebras\index{algebra!pseudo algebra} over a 2-theory.
Again, a theorem of Street and the existence of cotensor
products\index{cotensor product} imply that the 2-category of pseudo
algebras over a 2-theory\index{algebra!pseudo algebra over a
2-theory} admits weighted\index{weighted} pseudo limits. An example
of a pseudo algebra over a 2-theory comes from the category of
rigged surfaces\index{rigged surface} in \cite{hu}.

Some of these results may be found in some form in the literature.
There are many different ways to weaken 1-categorical concepts. This
study only sets up the weakened notions needed for utilizing
stacks\index{stack} to rigorously define conformal field
theory\index{conformal field theory} as in \cite{hu}. Results about
bilimits\index{bilimit} can be found in the references mentioned
above. In particular, Gray explicitly describes
quasilimits\index{quasilimit} and quasicolimits\index{quasicolimit}
of {\it strict} 2-functors from an arbitrary small 2-category to the
2-category $Cat$ of small categories on pages 201 and 219 of
\cite{gray1}, although his quasilimit is defined in terms of
quasiadjunction\index{quasiadjunction} rather than cones. In any
case, he does not have formulas for pseudo limits\index{limit!pseudo
limit} of {\it pseudo} functors. Street has the most general result
in this context. In \cite{street3}, he states that $Cat$ admits
indexed pseudo limits\index{limit!indexed pseudo limit} of pseudo
functors and writes down the indexed\index{indexed} pseudo limit.
His indexed pseudo limit is the same as the weighted\index{weighted}
pseudo limit\index{limit!weighted pseudo limit} in this paper.
Results about notions similar to the notion of
biadjoint\index{biadjoint} can be found in \cite{gray1},
\cite{gray2}, \cite{kelly2}, and \cite{street3}. These similarities
are discussed in the introduction to Chapter \ref{sec:laxadjoints}.
Blackwell, Kelly, and Power have limit and adjoint results similar
to ours for strict 2-functors into 2-categories of strict
algebras\index{algebra!strict algebra} and pseudo morphisms over a
2-monad\index{2-monad} in \cite{blackwell}. In fact, we prove below
that pseudo algebras over a theory are the strict algebras for a
2-monad\index{2-monad} in Chapter
\ref{sec:laxTalgebras}\index{algebra!pseudo algebra}.

Any discussion of weakened algebraic structures must involve
coherence questions, one of which was first treated in the classic
paper \cite{maclanecoherence1963} of Mac Lane. Many authors,
including Laplaza, Kelly, Mac Lane, and Par\'e, have contributed to
the theory of coherence as evidenced by the bibliographies of
\cite{maclane1991} and \cite{maclane1}. Some recent treatments in
the context of $n$-categories and categorification are \cite{baez2},
\cite{baez}, and \cite{dunn}. See also \cite{yanofskyproposal},
\cite{yanofsky2000}, and \cite{yanofsky2001} for an approach to
coherence involving a notion of 2-theory distinct from the notion of
2-theory in \cite{hu}, \cite{hu1}, \cite{hu2}, and Chapter
\ref{sec:2-theories}.

We follow the usual convention that 2-categories\index{2-category}
are denoted by capital script letters $\mathcal{A,C,D,X}$, pseudo
functors\index{functor!pseudo functor} are denoted by capital
letters $F,G$, morphisms\index{morphism} are denoted by $e,f,g,h$,
and 2-cells\index{2-cell} are denoted by Greek letters $\alpha,
\beta,\gamma$. The identity 2-cell on a morphism $f$ is denoted
$i_f$. Natural transformations\index{natural transformation} and
pseudo natural transformations\index{natural transformation!pseudo
natural transformation}\index{natural transformation} are also
denoted by lowercase Greek letters. The double arrow $\Rightarrow$
is used to denote 2-cells, natural transformations, and pseudo
natural transformations, which in some cases are all the same thing.
The notation $A \in \mathcal{D}$ means that $A$ is an object of
$\mathcal{D}$.

We usually reserve the notation $\mathcal{C}$ for a 2-category in
which we are building various limits and colimits. For example, in
Chapters \ref{sec:laxcolimitsinCat} and \ref{sec:laxlimitsinCat} the
letter $\mathcal{C}$ denotes the 2-category of small categories,
while it stands for the 2-category of small pseudo $T$-algebras in
Chapters \ref{sec:laxlimitsinlaxalgebras} and \ref{sec:laxcolaxT}.
In Chapter 13, the notation $\mathcal{C}$ stands for the 2-category
of small pseudo $(\Theta, T)$-algebras. We use the same letter to
highlight the similarities of the various proofs. In this
introduction $\mathcal{C}$ stands for the category of rigged
surfaces\index{rigged surface}.

All sets, categories, and 2-categories appearing in this paper are
assumed to be small.

\chapter{Some Comments on Conformal Field Theory}
\label{sec:cftcomments} In this chapter we make some motivational
remarks about conformal field theory\index{conformal field
theory|(}. Most of these terms will not appear in the rest of the
paper, and are therefore only briefly discussed. More detail can be
found in the articles \cite{hu} and \cite{hu1}, which this paper
accompanies.

Conformal field theory has recently received considerable attention
from mathematicians and physicists. It offers one approach to string
theory\index{string theory}, which aims to unify the four
fundamental forces of nature. This is one reason why physicists are
interested in conformal field theory as in \cite{polchinski}. The
motivation for the axioms of conformal field theory comes from the
path integral\index{path integral} formalism of quantum field
theory\index{quantum field theory}. Mathematicians have become
interested in conformal field theory because it gives rise to a
geometric definition of elliptic cohomology\index{elliptic
cohomology}, which is related to Borcherds' proof \cite{borcherds}
of the Moonshine conjectures\index{Moonshine}.

The formalism necessary to rigorously define conformal field theory,
and to prove theorems about it, is called {\it stacks of lax
commutative monoids with cancellation (SLCMC's)} in
\cite{hu}\index{cancellation}\index{SLCMC|textbf}\index{stack}\index{stack!stack
of lax commutative monoids with cancellation}. These are the same as
{\it stacks of pseudo algebras\index{algebra!pseudo algebra} over
the 2-theory of commutative monoids with
cancellation}\index{cancellation}\index{2-theory!2-theory of
commutative monoids with cancellation} defined in Chapters
\ref{sec:stacks} and \ref{sec:2-theories}. Roughly speaking, a
strict commutative monoid with
cancellation\index{cancellation}\index{commutative monoid with
cancellation|textbf} consists of a commutative
monoid\index{commutative monoid} $I$ and a function $X:I^2
\rightarrow Sets$ equipped with operations
$$+_{a,b,c,d}: X_{a,b} \times X_{c,d} \rightarrow X_{a+c, b+d}$$
$$\check{?}_{a,b,c}:X_{a+c,b+c} \rightarrow X_{a,b}$$
$$0 \in X_{0,0}$$
for all $a,b,c,d \in I$. These operations, called {\it disjoint
union}\index{disjoint union}, {\it cancellation\index{cancellation}
(gluing)}\index{gluing}, and {\it unit}\index{unit} must be
commutative, associative, unital, and distributive in the
appropriate senses. Whenever we add the adjective
``pseudo''\index{pseudo|textbf} (or ``lax''\index{lax|textbf} in
\cite{hu}, \cite{hu1}, \cite{hu2}), it means that we replace sets by
categories, functions by functors, and axioms by coherence
isos\index{coherence isomorphism} that satisfy coherence
diagrams\index{coherence diagram}. The theory and 2-theory apparatus
gives us a concise way to list the necessary coherence
isos\index{coherence isomorphism} and coherence
diagrams\index{coherence diagram}. A thorough treatment of theories,
2-theories, their pseudo algebras\index{algebra!pseudo algebra}, and
their relevant diagrams are part of this paper. This formalism
allows the authors of \cite{hu} and \cite{hu1} to rigorously define
conformal field theory in the sense of Segal\index{Segal, Graeme},
in particular all of the coherence isos\index{coherence isomorphism}
and coherence diagrams\index{coherence diagram} are neatly encoded.

The first example of a pseudo commutative monoid with
cancellation\index{cancellation}\index{commutative monoid with
cancellation!pseudo commutative monoid with cancellation!examples of
pseudo commutative monoids with cancellation} is the category of
rigged surfaces\index{rigged surface}. In this example the pseudo
commutative monoid\index{commutative monoid!pseudo commutative
monoid!examples of pseudo commutative monoids} $I$ is the category
of finite sets and bijections equipped with disjoint
union\index{disjoint union}. The 2-functor $X:I^2 \rightarrow Cat$
from $I^2$ to the 2-category of small categories is given by
defining $X_{a,b}$ to be the category of rigged
surfaces\index{rigged surface} with inbound\index{inbound}
components labelled\index{label} by $a$ and outbound\index{outbound}
components labelled\index{label} by $b$. The operation $+$ is
disjoint union\index{disjoint union} of labelled rigged
surfaces\index{rigged surface} (this is why the
indices\index{indices} are added). The stack structure for this
example is described in Section \ref{sec:rigged}\index{stack}.

There are two other examples of SLCMC's that we need before defining
conformal field theory and modular functor\index{algebra!pseudo
algebra!examples of pseudo
algebras}\index{SLCMC}\index{SLCMC!examples of
SLCMC's}\index{commutative monoid with cancellation!pseudo
commutative monoid with cancellation!examples of pseudo commutative
monoids with cancellation}. These are $C(\mathcal{M})$ and
$C(\mathcal{M},H)$ from page 235 of \cite{hu1}. The notation
$\mathbb{C}_2$ denotes the pseudo commutative
semi-ring\index{commutative semi-ring!pseudo commutative
semi-ring}\index{semi-ring} of finite dimensional complex vector
spaces\index{vector space}, $\mathbb{C}_2^{Hilb}$ denotes the pseudo
$\mathbb{C}_2$-algebra of complex Hilbert spaces\index{Hilbert
space} equipped with the operation $\hat{\otimes}$ of Hilbert tensor
product, $\mathcal{M}$ is a pseudo module\index{module!pseudo
module} over $\mathbb{C}_2$, $\mathcal{M}^{Hilb}$ denotes
$\mathcal{M} \otimes \mathbb{C}_2^{Hilb}$, and $H$ is an object of
$\mathcal{M}^{Hilb}$. If $\mathcal{M}$ has only one object, then $H$
is a Hilbert space\index{Hilbert space}, otherwise $H$ is a
collection of Hilbert spaces\index{Hilbert space}
indexed\index{indexed} by the objects of $\mathcal{M}$. For finite
sets $a,b$ the category $C(\mathcal{M})_{a,b}$ is
$\mathcal{M}^{\otimes a} \otimes \mathcal{M}^{* \otimes b}$ where
$\mathcal{M}^*:=Hom_{pseudo}(\mathcal{M}, \mathbb{C}_2)$. The
operation $+$ is given by $\otimes$ and gluing\index{gluing} is
given by evaluation $tr:\mathcal{M} \otimes \mathcal{M}^*
\rightarrow \mathbb{C}_2$. The pseudo commutative monoid with
cancellation\index{cancellation}\index{commutative monoid with
cancellation!pseudo commutative monoid with cancellation}
$C(\mathcal{M},H)$ is defined similarly, except that an object of
$C(\mathcal{M},H)_{a,b}$ is an object $M$ of $C(\mathcal{M})_{a,b}$
equipped with a morphism $M \rightarrow H^{\hat{\otimes} a}
\hat{\otimes} H^{*\hat{\otimes} b}$ in $C(\mathcal{M})_{a,b}$ whose
image consists of trace class\index{trace class} elements. The
morphisms of $C(\mathcal{M},H)_{a,b}$ are the appropriate
commutative triangles in $C(\mathcal{M})_{a,b}$. These two
LCMC's\index{LCMC} can be made into stacks\index{stack}
appropriately. Finally we are ready to give the rigorous definition
of modular functor and conformal field theory.

\begin{defn} Let $C$ be a stack\index{stack} of pseudo commutative monoids
with cancellation\index{cancellation} (SLCMC)\index{SLCMC}. A {\it
modular functor\index{modular functor|textbf} on $C$ with
labels\index{label|textbf} $\mathcal{M}$} is a (pseudo) morphism
$\phi:C \rightarrow C(\mathcal{M})$ of stacks of pseudo commutative
monoids with cancellation\index{cancellation}. A {\it conformal
field theory\index{conformal field theory|textbf} on $C$ with
modular functor\index{modular functor|textbf} on labels
$\mathcal{M}$ with state space\index{state space|textbf} $H$} is a
(pseudo) morphism $\Phi:C \rightarrow C(\mathcal{M}, H)$ of
stacks\index{stack} of pseudo commutative monoids with
cancellation\index{cancellation}\index{SLCMC!morphism of SLCMC's}.
\end{defn}

If we take $C$ to be the SLCMC\index{SLCMC} of rigged
surfaces\index{rigged surface}, then we recover the usual definition
of conformal field theory which assigns (up to a finite dimensional
vector space\index{vector space}) a trace class\index{trace class}
operator to a rigged surface\index{rigged surface} in such a way
that gluing\index{gluing} surfaces corresponds to composing
operators. Notice that modular functor and conformal field
theory\index{conformal field theory} are both morphisms of the same
algebraic structure. This was first noted by the authors of
\cite{hu} and \cite{hu1}.

It is also possible to define one dimensional modular
functors\index{modular functor!one dimensional modular functor} (\ie
those with one object in $\mathcal{M}$) in terms of {\it
$\mathbb{C}^{\times}$-central extensions}\index{central extension}
of SLCMC's\index{SLCMC}\index{SLCMC!central extension of SLCMC's}. A
$\mathbb{C}^{\times}$-central extension of an SLCMC $\mathcal{D}$ is
a strict morphism $\psi:\tilde{\mathcal{D}} \rightarrow \mathcal{D}$
of SLCMC's such that for fixed finite sets $s,t$, a fixed finite
dimensional complex manifold\index{manifold!complex manifold} $B$,
and fixed $\alpha \in \mathcal{D}(B)_{s,t}$, the pre-images
$\psi^{-1}(\alpha|_{B'})$ patch together for varying $B' \rightarrow
B$ to form the sheaf of sections of a complex
holomorphic\index{holomorphic} line bundle\index{line bundle} over
$B$. The maps on these sections\index{section} induced by disjoint
union\index{disjoint union} and gluing\index{gluing} are required to
be isomorphisms of sheaves\index{sheaf} of vector
spaces\index{vector space}. If $\mathcal{H}$ is a Hilbert
space\index{Hilbert space}, then there is an SLCMC\index{SLCMC}
$\underline{\mathcal{H}}$ in which $+$ is the operation of taking
the Hilbert space\index{Hilbert space} tensor product\index{Hilbert
tensor product}\index{tensor product} and then the subset of trace
class\index{trace class} elements and $\check{?}$ is the trace
map\index{trace map|textbf}. Then a {\it chiral\index{chiral}
conformal field theory with one dimensional modular
functor\index{modular functor!one dimensional modular functor} over
$\mathcal{D}$} is a morphism of SLCMC's\index{SLCMC}
$\phi:\tilde{\mathcal{D}} \rightarrow \underline{\mathcal{H}}$ which
is linear on the spaces of sections $\psi^{-1}(\alpha|_{B'})$.

The present paper deals with the 2-categorical foundations of the
above project. We begin by introducing 2-categories and proving the
existence of various types of limits in various 2-categories in
Chapters \ref{sec:laxcolimits}, \ref{sec:laxcolimitsinCat},
\ref{sec:theories}, \ref{sec:laxlimitsinlaxalgebras} ,
\ref{sec:laxcolaxT}, and \ref{sec:2-theories}. We need the existence
of certain limits in the above project because a stack\index{stack}
is a contravariant pseudo
functor\index{contravariant}\index{functor!pseudo functor} that
takes Grothendieck covers\index{Grothendieck cover} to
bilimits\index{bilimit}. Grothendieck topologies and stacks are
discussed in Chapter \ref{sec:stacks}. The fundamentals of Lawvere
theories\index{Lawvere theory} and algebras are treated in Chapter
\ref{sec:theories}. The passage from strict
algebras\index{algebra!strict algebra} to pseudo
algebras\index{algebra!pseudo algebra}, which is so important for
the definition of conformal field theory, is discussed in Chapter
\ref{sec:laxTalgebras}. The biadjoints of Chapters
\ref{sec:laxadjoints} and \ref{sec:forgetfulfunctor} allow a
universal description of the stack\index{stack} of covering
spaces\index{covering space}\index{stack!stack of covering
spaces}\index{covering space!stack of covering spaces} on page 337
of \cite{hu}. Lastly, the 2-theory of commutative monoids with
cancellation\index{cancellation}\index{2-theory!2-theory of
commutative monoids with cancellation} is presented in Chapter
\ref{sec:2-theories} along with the example of rigged
surfaces.\index{conformal field theory|)}\index{rigged surface}

\chapter{Weighted Pseudo Limits in a 2-Category}
\label{sec:laxcolimits} In this chapter we introduce the notion of a
weighted\index{weighted} pseudo limit and related concepts. The most
important examples of 2-categories to keep in mind are the
following.

\begin{examp}
The 2-category\index{2-category} of small categories\index{category
of small categories} is formed by taking the objects
(0-cells)\index{0-cell} to be small categories, the morphisms
(1-cells)\index{1-cell} to be functors, and the
2-cells\index{2-cell} to be natural transformations. This 2-category
is denoted $Cat$\index{$Cat$}.
\end{examp}

\begin{examp}
A full sub-2-category of the previous example is the 2-category with
objects groupoids\index{groupoid} and 1-cells and 2-cells the same
as above.
\end{examp}

\begin{examp}
An example of a different sort is the 2-category with objects
topological spaces\index{topological space}, morphisms continuous
maps\index{continuous map}, and 2-cells homotopy classes of
homotopies\index{homotopy}. The 2-cells must be homotopy classes of
homotopies in order to make the various compositions associative and
unital.
\end{examp}

\begin{examp}
Let $\mathcal{J}$ be a small 1-category. Then $\mathcal{J}$ has the
structure of a 2-category if we regard $Mor_{\mathcal{J}}(i,j)$ as a
discrete\index{discrete category} category for all $i,j \in Obj
\hspace{1mm} \mathcal{J}$.
\end{examp}

These examples show that there are two ways of composing the
2-cells: vertically and horizontally. Natural transformations can be
composed in two ways. Homotopy classes of homotopies can also be
composed in two ways. To clarify which composition we mean, we
follow Borceux's notation. See \cite{borceux} for a more thorough
discussion.

\begin{defn}
Let $\mathcal{C}$ be a 2-category. If $A,B \in Obj \hspace{1mm}
\mathcal{C}$ and $f,g,h:A \rightarrow B$ are objects of the category
$Mor(A,B)$ with 2-cells $\alpha:f \Rightarrow g$ and $\beta: g
\Rightarrow h$ then the composition
$$\xymatrix@R=3pc@C=3pc{A \ar[rr]^(.3)f & \ar@{=>}[d]^{\alpha} & B \\ A \ar[rr]^(.3)g &
\ar@{=>}[d]^{\beta} & B \\ A \ar[rr]^(.3)h & & B}$$ in the category
$Mor(A,B)$ is called the {\it vertical composition}\index{vertical
composition|textbf} of $\alpha$ and $\beta$. This composition is
denoted $\beta \odot \alpha$. The {\it identity}\index{vertical
identity|textbf} on $f$ with respect to vertical composition is
denoted $i_f$.
\end{defn}

\begin{defn}
Let $\mathcal{C}$ be a 2-category and $A,B,C \in Obj \hspace{1mm}
\mathcal{C}$. Let $c:Mor(A,B) \times Mor(B,C) \rightarrow Mor(A,C)$
denote the functor of composition\index{composition} in the
2-category $\mathcal{C}$. If $f,g:A \rightarrow B$ and $m,n:B
\rightarrow C$ are objects of the respective categories $Mor(A,B)$
and $Mor(B,C)$ and $\alpha:f \Rightarrow g$, $\beta:m \Rightarrow n$
are 2-cells, then the composite 2-cell $c(\alpha, \beta): c(f,m)
\Rightarrow c(g,n)$ is called the {\it horizontal
composition}\index{horizontal composition|textbf} of $\alpha$ and
$\beta$. It is a morphism in the category $Mor(A,C)$ and is denoted
$\beta \ast \alpha$.
$$\xymatrix@R=3pc@C=3pc{A \ar[rr]^f & \ar@{=>}[d]^{\alpha} & \ar[rr]^m B &
\ar@{=>}[d]^{\beta} & C \\ A \ar[rr] _g & & B \ar[rr]_n & & C}$$
\end{defn}

To define the concept of weighted pseudo limit, we need to discuss
pseudo functors and pseudo natural transformations. A pseudo functor
is like a 2-functor except that it preserves composition and
identity only up to iso coherence 2-cells which satisfy coherence
diagrams. A pseudo natural transformation\index{natural
transformation!pseudo natural transformation} is like a 2-natural
transformation except that it is natural only up to an iso coherence
2-cell which satisfies coherence diagrams. We define these notions
more carefully to fix some notation. We reproduce Borceux's
treatment in \cite{borceux}. The coherence 2-cells for pseudo
functors and pseudo natural transformations in this paper are always
assumed to be iso. Recall again that a pseudo functor in this paper
is a lax\index{lax}\index{functor!lax functor|textbf} functor in
\cite{hu}, \cite{hu1}, and \cite{hu2} as well as in other previous
papers.

\begin{defn}
Let $\mathcal{C},\mathcal{D}$ be 2-categories. A {\it pseudo
functor}\index{functor}\index{functor!pseudo functor|textbf}
$F:\mathcal{C} \rightarrow \mathcal{D}$ consists of the following
assignments and iso coherence 2-cells:

\begin{itemize}
\item
For every object $A \in Obj \hspace{1mm}  \mathcal{C}$ an object $FA
\in Obj \hspace{1mm}  \mathcal{D}$
\item
For all objects $A,B \in Obj \hspace{1mm}  \mathcal{C}$ a functor
$F:Mor_{\mathcal{C}}(A,B) \rightarrow Mor_{\mathcal{D}}(FA,FB)$
\item
For all objects $A,B,C \in Obj \hspace{1mm} \mathcal{C}$ a natural
isomorphism $\gamma$ between the composed functors
$$\xymatrix@C=3pc@R=3pc{Mor_{\mathcal{C}}(A,B) \times Mor_{\mathcal{C}}(B,C)
\ar[r]^-c \ar[d]_{F \times F} & Mor_{\mathcal{C}}(A,C) \ar[d]^F \\
Mor_{\mathcal{D}}(FA,FB) \times Mor_{\mathcal{D}}(FB,FC) \ar[r]_-c
 \ar@{=>}[ur]^{\gamma} & Mor_{\mathcal{D}}(FA,FC)}$$
\item
For every object $A \in \mathcal{C}$ a natural isomorphism $\delta$
between the following composed functors.
$$\xymatrix@R=3pc@C=3pc{\mathbf{1} \ar[r]^-u \ar@{=}[d] & Mor_{\mathcal{C}}(A,A)
\ar[d]^F \\ \mathbf{1} \ar@{=>}[ur]^{\delta} \ar[r]_-u &
Mor_{\mathcal{D}}(FA,FA)}$$ where the functor $u:\mathbf{1}
\rightarrow Mor_{\mathcal{C}}(A,A)$ from the terminal
object\index{terminal object} $\mathbf{1}$ in the category of small
categories to the category $Mor_{\mathcal{C}}(A,A)$ takes the unique
object $*$ of $\mathbf{1}$ to the identity morphism on $A$.
\end{itemize}
These coherence 2-cells\index{coherence 2-cell|textbf} must satisfy
the following coherence diagrams\index{coherence diagram|textbf}.
\begin{itemize}
\item
For every morphism $f:A \rightarrow B$ in $\mathcal{C}$ we require
$$\xymatrix@C=4pc@R=3pc{Ff \circ 1_{FA} \ar@{=>}[r]^-{i_{Ff} \ast \delta_{A \ast}}
\ar@{=>}[d]_{i_{Ff}} & Ff \circ F1_A \ar@{=>}[d]^{\gamma_{1_A,f}} &
1_{FB} \circ Ff \ar@{=>}[r]^-{\delta_{B \ast} \ast i_{Ff}}
\ar@{=>}[d]_{i_{Ff}} & F(1_B) \circ Ff \ar@{=>}[d]^{\gamma_{f,1_B}}
\\ Ff \ar@{=>}[r]_-{i_{Ff}}
& F(f \circ 1_A) & Ff \ar@{=>}[r]_-{i_{Ff}} & F(1_B \circ f)}$$ to
commute. Here $\delta_{A \ast}$ means the natural transformation
$\delta_A$ evaluated at the unique object $\ast$ of $\mathbf{1}$.
This is called the {\it unit axiom}\index{unit axiom} for the pseudo
functor $F$.
\item
For all morphisms $f,g,h$ of $\mathcal{C}$ such that $h \circ g
\circ f$ exists we require that
$$\xymatrix@C=4pc@R=3pc{Fh \circ Fg \circ Ff \ar@{=>}[r]^{i_{Fh} \ast \gamma_{f,g}}
\ar@{=>}[d]_{\gamma_{g,h} \ast i_{Ff}} & Fh \circ F(g \circ f)
 \ar@{=>}[d]^{\gamma_{g \circ f, h}}
\\ F(h \circ g) \circ Ff \ar@{=>}[r]_{\gamma_{f, h \circ g}} &
F(h \circ g \circ f)}$$ commutes. This is called the {\it
composition axiom}\index{composition axiom} for the pseudo functor
$F$.
\end{itemize}
\end{defn}

Each of these functors and natural transformations of course depends
on the objects, so they really need indices, \eg $c_{A,B,C},
F_{A,B}, \gamma_{A,B,C}, u_A,u_{FA},$ and $\delta_A$. Often we leave
the indices off for more convenient notation. Note that the first
diagram in the definition says that the pseudo functor preserves
composition of morphisms up to coherence 2-cell\index{coherence
2-cell} because for morphisms $\xymatrix@1{A \ar[r]^f & B \ar[r]^g &
C}$ in $\mathcal{C}$ we have $\gamma_{f,g}:F(g) \circ F(f)
\Rightarrow F(g \circ f)$ and $\gamma$ is natural in $f$ and $g$.
The second diagram in the definition says that the pseudo functor
preserves identity up to coherence 2-cell because $\delta_{A
\ast}:1_{FA} \Rightarrow F(1_A)$ for all $A \in Obj \hspace{1mm}
\mathcal{C}$.

\begin{defn}
Let $\xymatrix@1{\mathcal{C} \ar[r]^F & \mathcal{D} \ar[r]^G &
\mathcal{E}}$ be pseudo functors. Then the {\it composition $G \circ
F$ of pseudo functors}\index{functor!pseudo functor!composition of
pseudo functors|textbf}\index{composition of pseudo functors|textbf}
is the composition of the underlying maps of objects and the
composition of the underlying functors on the morphism categories.
The coherence 2-cells are as follows.
\begin{itemize}
\item
For morphisms $\xymatrix@1{A \ar[r]^f & B \ar[r]^g & C}$ in
$\mathcal{C}$ the 2-cell $\gamma^{GF}_{f,g}$ is the composition
$$\xymatrix@1@C=4pc{GF(g) \circ GF(f) \ar@{=>}[r]^-{\gamma^G_{Ff,Fg}}
& G(Fg \circ Ff) \ar@{=>}[r]^-{G(\gamma^F_{f,g})} & GF(g \circ
f)}.$$
\item
For each object $A \in Obj \hspace{1mm}  \mathcal{C}$ the 2-cell
$\delta^{GF}_{A \ast}$ is the composition
$$\xymatrix@1@C=4pc{1_{GFA} \ar@{=>}[r]^-{\delta^G_{FA \ast}}
& G(1_{FA}) \ar@{=>}[r]^-{G(\delta^F_{A \ast})} & GF(1_A)}.$$
\end{itemize}
\end{defn}
Then the assignment $(f,g) \mapsto \gamma^{GF}_{f,g}$ is natural and
$\gamma^{GF}$ and $\delta_A^{GF}$ satisfy the coherences to make
$GF$ a pseudo functor.

\begin{defn}
A {\it pseudo natural transformation}\index{natural
transformation!pseudo natural transformation|textbf} $\alpha:F
\Rightarrow G$ from the pseudo functor $F: \mathcal{C} \rightarrow
\mathcal{D}$ to the pseudo functor $G: \mathcal{C} \rightarrow
\mathcal{D}$ consists of the following assignments:

\begin{itemize}
\item
For each $A \in Obj \hspace{1mm} \mathcal{C}$ a morphism
$\alpha_A:FA \rightarrow GA$ in the category $\mathcal{D}$
\item
For all objects $A,B \in Obj \hspace{1mm} \mathcal{C}$ a natural
isomorphism $\tau$ between the following functors.
$$\xymatrix@R=3pc@C=3pc{Mor_{\mathcal{C}}(A,B) \ar[r]^-F \ar[d]_G &
Mor_{\mathcal{D}}(FA,FB) \ar[d]^{\alpha_B \circ}
\\ Mor_{\mathcal{D}}(GA,GB) \ar[r]_{\circ \alpha_A} \ar@{=>}[ur]^{\tau} &
Mor_{\mathcal{D}}(FA,GB)}$$
\end{itemize}
The natural transformations $\tau$ must satisfy the following
coherence diagrams\index{coherence diagram|textbf} involving
$\delta$ and $\gamma$.
\begin{itemize}
\item
For every $A \in Obj \hspace{1mm}  \mathcal{C}$ we require
$$\xymatrix@C=4pc@R=3pc{\alpha_A \ar@{=>}[r]^{i_{\alpha_A}}
\ar@{=>}[d]_{i_{\alpha_A}} & 1_{GA} \circ \alpha_A
\ar@{=>}[r]^{\delta^G_{A \ast} \ast i_{\alpha_A}} & G(1_A) \circ
\alpha_A \ar@{=>}[d]^{\tau_{1_A}} \\ \alpha_A \circ 1_{FA}
\ar@{=>}[rr]_{i_{\alpha_A} \ast \delta^F_{A \ast}} & & \alpha_A
\circ F(1_A)}$$ to commute. This is called the {\it unit
axiom}\index{unit axiom} for the pseudo natural transformation
$\alpha$.
\item
For all morphisms $\xymatrix@1{A \ar[r]^f & B \ar[r]^g & C}$ in
$\mathcal{C}$ we require
$$\xymatrix@C=3pc@R=3pc{Gg \circ Gf \circ \alpha_A \ar@{=>}[r]^{i_{Gg} \ast \tau_f}
\ar@{=>}[d]_{\gamma_{f,g}^G \ast i_{\alpha_A}} & Gg \circ \alpha_B
\circ Ff \ar@{=>}[r]^{\tau_g \ast i_{Ff}} & \alpha_C \circ Fg \circ
Ff \ar@{=>}[d]^{i_{\alpha_C} \ast \gamma_{f,g}^F}
\\ G(g \circ f) \circ \alpha_A \ar@{=>}[rr]_{\tau_{g \circ f}}
& & \alpha_C \circ F(g \circ f)}$$ to commute. This is called the
{\it composition axiom}\index{composition axiom} for the pseudo
natural transformation $\alpha$.
\end{itemize}
\end{defn}

Here $\tau$ should of course also be indexed by the objects $A,B$
\etc, but we leave off these indices for convenience. The coherence
required on $\gamma$ and $\tau$ is the commutivity of the 2-cells
(from $\tau$ and $\gamma$) written on the faces of the prism with
edges $Ff,\hspace{1mm} Fg, \hspace{1mm} F(g \circ f), \hspace{1mm}
Gf, \hspace{1mm}Gg,$ and $G(g \circ f)$ where $f$ and $g$ are
composable morphisms in the 2-category $\mathcal{C}$. There are
several ways to compose these 2-cells, but they are related by the
interchange law. Here we must sometimes horizontally precompose or
postcompose a 2-cell with identity 2-cells in order to horizontally
compose. Note the diagram for $\tau$ drawn in the definition says
that the assignment of $A \mapsto \alpha_A$ is natural up to
coherence 2-cell\index{coherence 2-cell} because for $f \in
Mor_{\mathcal{C}}(A,B)$ we have the diagram
$$\xymatrix@C=3pc@R=3pc{FA \ar[r]^{\alpha_A} \ar[d]_{Ff}
& GA \ar[d]^{Gf} \ar@{=>}[ld]_{\tau_f}
\\ FB \ar[r]_{\alpha_B} & GB}$$
in $\mathcal{D}$. The assignment $f \mapsto \tau_f$ is natural in
$f$, \ie $\tau_{A,B}$ is a natural transformation.

Some authors prefer to denote the coherence 2-cells\index{coherence
2-cell} of $\alpha$ by $\alpha_f$ instead of $\tau_f$. However we
follow Borceux's notation in \cite{borceux} and use the
distinguished notation $\tau$ in order to navigate complicated
diagrams with less effort.

Pseudo natural transformations can also be
horizontally\index{horizontal composition} and vertically
composed\index{vertical composition}. For example, if $\xymatrix@1{F
\ar@{=>}[r]^{\alpha} & G \ar@{=>}[r]^{\beta} & H}$ are pseudo
natural transformations, the vertical composition $\beta \odot
\alpha$ has coherence 2-cells\index{coherence 2-cell} $\tau^{\beta
\odot \alpha}_f=(i_{\beta_B}*\tau^{\alpha}_f) \odot (\tau^{\beta}_f
* i_{\alpha_A})$ for $f:A \rightarrow B$ as in the following
diagram. \label{compositionoflaxnaturals}
$$\xymatrix@C=4pc@R=4pc{FA \ar[r]^{\alpha_A} \ar[d]_{Ff} &
GA \ar@{=>}[ld]_{\tau^{\alpha}_f}
 \ar[r]^{\beta_A} \ar[d]^{Gf} & HA \ar[d]^{Hf} \ar@{=>}[ld]_{\tau^{\beta}_f}
  \\ FB \ar[r]_{\alpha_B} & GB \ar[r]_{\beta_B} & HB}$$
Natural transformations can be seen as morphisms between functors.
In the context of 2-categories there is a similar notion of a
modification between pseudo natural transformations.

\begin{defn}
Let $F,G: \mathcal{C} \rightarrow \mathcal{D}$ be pseudo functors
and $\alpha, \beta: F \Rightarrow G$ pseudo natural transformations.
A {\it modification}\index{modification|textbf} $\Xi:\alpha
\rightsquigarrow \beta$ is a function which assigns to every $A \in
Obj \hspace{1mm} \mathcal{C}$ a 2-cell $\Xi_A: \alpha_A \Rightarrow
\beta_A$ in $\mathcal{D}$ in such a way that $\tau_{A,B}^{\beta} (g)
\odot (G\gamma \ast \Xi_A) = (\Xi_{B} \ast F\gamma) \odot
\tau_{A,B}^{\alpha}(f)$ for all $A, B \in Obj \hspace{1mm}
\mathcal{C}$ and all morphisms $f,g:A \rightarrow B$ and all 2-cells
$\gamma:f \Rightarrow g$. Here $\tau^{\alpha}$ and $\tau^{\beta}$
denote the natural transformations belonging to the pseudo natural
transformations $\alpha$ and $\beta$ respectively, while $\gamma$ is
an arbitrary 2-cell in $\mathcal{C}$.  This means that the following
two compositions of 2-cells are the same.
\begin{equation} \label{defnmodification1}
\xymatrix@R=3pc@C=3pc{FA \ar[rr]^{\alpha_A} & \ar@{=>}[d]^{\Xi_A} &
GA \ar[rr]^{Gf} &
 \ar@{=>}[d]^{G \gamma}  & GB\\
FA \ar[rr]_{\beta_A} & & GA \ar[rr]_{Gg}
\ar@{=>}[d]^{\tau_{A,B}^{\beta}(g)} & & GB \\ FA \ar[rr]_{Fg} & & FB
\ar[rr]_{\beta_{B}} & & GB}
\end{equation}
\begin{equation} \label{defnmodification2}
\xymatrix@R=3pc@C=3pc{FA \ar[rr]^{\alpha_A} & & GA \ar[rr]^{Gf}
\ar@{=>}[d]^{\tau_{A,B}^{\alpha}(f)} & & GB\\
FA \ar[rr]^{Ff} & \ar@{=>}[d]^{F \gamma} & FB \ar[rr]^{\alpha_{B}}
& \ar@{=>}[d]^{\Xi_{B}} & GB\\
FA \ar[rr]_{Fg} & & FB \ar[rr]_{\beta_{B}} & & GB}
\end{equation}
These two diagrams can be combined to make a cube whose faces have
2-cells inscribed in them. In this definition $\gamma$ is not to be
confused with the required coherence 2-cell in the definition of
pseudo functor.
\end{defn}

\begin{defn}
If $F:\mathcal{D} \rightarrow \mathcal{C}$ is a pseudo functor, then
a {\it pseudo limit}\index{limit}\index{limit!pseudo limit|textbf}
of $F$ consists of an object $W \in Obj \hspace{1mm} \mathcal{C}$
and a pseudo natural transformation $\pi:\Delta_W \Rightarrow F$
from the constant 2-functor $W$ to the pseudo functor $F$ which is
universal in the following sense: the functor $(\pi \circ):
Mor_{\mathcal{C}}(C,W) \rightarrow PseudoCone(C, F)$ is an {\it
isomorphism} of categories for every object $C \in Obj \hspace{1mm}
{\mathcal{C}}$.
\end{defn}

$PseudoCone(C,F)$\index{cone!pseudo cone|textbf}\index{$PseudoCone$}
denotes here the category with objects taken to be the pseudo
natural transformations $\Delta_C \Rightarrow F$ and with morphisms
taken to be the modifications. Pseudo colimits\index{colimit!pseudo
colimit|textbf}\index{colimit} can be defined in terms of
$PseudoCone(F,C)$ and $( \circ \pi): Mor_{\mathcal{C}}(W,C)
\rightarrow PseudoCone(F,C)$ similarly.

\begin{thm}
Any two pseudo limits of a pseudo functor are isomorphic.
\end{thm}

\begin{defn}
If $F:\mathcal{D} \rightarrow \mathcal{C}$ is a pseudo functor, then
a {\it bilimit}\index{bilimit|textbf} of $F$ consists of an object
$W \in Obj \hspace{1mm} \mathcal{C}$ and a pseudo natural
transformation $\pi:\Delta_W \Rightarrow F$ from the constant
2-functor $W$ to the pseudo functor $F$ which is universal in the
following sense: the functor $(\pi \circ): Mor_{\mathcal{C}}(C,W)
\rightarrow PseudoCone(C, F)$ is an {\it equivalence} of categories
for every object $C \in Obj \hspace{1mm} {\mathcal{C}}$.
\end{defn}

Some authors would call this bilimit a {\it conical
bilimit}\index{conical}, see \cite{kelly2} and \cite{street3} for
example. They discuss the more general notion\index{weighted} of
{\it weighted bilimit}\index{bilimit!weighted bilimit} or {\it
indexed bilimit}\index{bilimit!indexed bilimit}, which is defined
below. Limits defined in terms of cones\index{cone}, such as this
bilimit, have constant weight\index{weighted}\index{weight} or
constant index\index{index}. For our applications to conformal field
theory\index{conformal field theory}, it is sufficient to consider
only conical bilimits\index{bilimit!conical bilimit} although we
prove results for more general weighted bilimits in this paper. The
existence of conical\index{conical} bilimits is sufficient to speak
of stacks\index{stack}. The term {\it lax
limit}\index{lax}\index{limit!lax limit|textbf} in \cite{hu},
\cite{hu1}, and \cite{hu2} is synonymous with the term {\it
bilimit}\index{bilimit} defined above.

Every pseudo limit for a fixed pseudo functor is obviously a
bilimit\index{bilimit} of that pseudo functor. One can ask whether
or not bilimits and pseudo limits are the same. The following
trivial example shows that bilimits and pseudo limits are not the
same.

\begin{examp} \label{stronglaxneqlax}
Let $\mathbf{1}$ denote the terminal object\index{terminal object}
in the category of small categories, in other words $\mathbf{1}$ is
the category with one object $*$ and one morphism, namely the
identity morphism. This category can be viewed as a 2-category with
no nontrivial 2-cells. Suppose $\mathcal{C}$ is a 2-category with at
least two objects $W, W'$ such that we have a morphism $\pi':W'
\rightarrow W$ which is a pseudo isomorphism. This means that there
exists a morphism $\theta:W \rightarrow W'$ and iso 2-cells $\theta
\circ \pi' \Rightarrow 1_{W'}$ and $\pi' \circ \theta \Rightarrow
1_W$. Suppose further that $\pi'$ is not monic. This means there
exists an object $C \in Obj \hspace{1mm}  \mathcal{C}$ and distinct
morphisms $f_1,f_2:C \rightarrow W'$ such that $\pi' \circ f_1 =
\pi' \circ f_2$. Let $F:\mathbf{1} \rightarrow \mathcal{C}$ be the
constant functor $\Delta_W$, \ie $F(*)=W$ and the identity gets
mapped to $1_W$. Then $PseudoCone(C,F)$ is isomorphic to
$Mor_{\mathcal{C}}(C,W)$. We identify these two categories.
Obviously $W$ and the pseudo natural transformation $\pi=1_W$ (under
the identification) form a pseudo limit\index{limit!pseudo
limit!example of pseudo limit}, while $W'$ and $\pi'$ form a
bilimit\index{bilimit!example of bilimit}. However, $W'$ and $\pi'$
do not form a pseudo limit because $(\pi' \circ):
Mor_{\mathcal{C}}(C,W') \rightarrow Mor_{\mathcal{C}}(C,W)$ is not
an isomorphism of categories, since $\pi' \circ f_1 = \pi' \circ
f_2$ although $f_1 \neq f_2$.
\end{examp}

\begin{examp}
There are also examples where a bicolimit\index{bicolimit!example of
bicolimit} exists but not a pseudo colimit\index{colimit!pseudo
colimit}. This example goes back to \cite{blackwell}. Let
$Lex$\index{$Lex$} denote the 2-category of small finitely complete
categories, left exact functors, and natural transformations. A
functor is called {\it left exact}\index{left
exact}\index{exact!left exact} if it preserves all finite limits. An
initial object\index{initial object} is a colimit of the empty
2-functor. A pseudo colimit\index{colimit}\index{colimit!pseudo
colimit} and a 2-colimit\index{2-colimit} of the empty 2-functor are
the same thing. The 2-category $Lex$ does not admit an initial
object because there are always two distinct functors $A \rightarrow
I$ where $I$ is the category with only two isomorphic objects and no
nontrivial morphisms besides the isomorphism and its inverse. The
two constant functors provide us with two distinct functors $A
\rightarrow I$ for each $A \in Obj \hspace{1mm}  Lex$. The empty
functor does however admit a bicolimit\index{bicolimit} because
$Lex$ is the 2-category of strict algebras, pseudo algebra
morphisms, and 2-cells for some finitary 2-monad on $Cat$.
Blackwell, Kelly, and Power prove in \cite{blackwell} that such
algebra categories admit bicolimits.

Many pseudo algebra\index{algebra!pseudo algebra} categories do not
admit pseudo colimits\index{colimit!pseudo colimit} because the
morphisms are not strict. Another example can be obtained by
adapting Example \ref{nostrictification} on page
\pageref{nostrictification} to colimits.
\end{examp}

After Example \ref{stronglaxneqlax}, one might wonder whether or not
the equivalences of categories in the definition of bilimit can be
chosen in some natural way. They can in fact be chosen pseudo
naturally as follows. We write it explicitly only for the bicolimit,
although a completely analogous statement holds for the bilimit.

\begin{rem}
Let $\mathcal{C}, \mathcal{D}$ be 2-categories. Let
$\mathcal{F}:\mathcal{D} \rightarrow \mathcal{C}$ be a pseudo
functor. Suppose $W \in Obj \hspace{1mm}  \mathcal{C}$ is a
bicolimit with universal pseudo cone $\pi:\mathcal{F} \Rightarrow
\Delta_W$. Let $\phi_C$ denote the equivalence of categories $(\circ
\pi):Mor_{\mathcal{C}} (W,C) \rightarrow PseudoCone(\mathcal{F},C)$.
Let $G(C):=Mor_{\mathcal{C}}(W,C)$ and
$F(C):=PseudoCone(\mathcal{F},C)$. Then $G$ and $F$ are strict
2-functors and $C \mapsto \phi_C$ is a 2-natural transformation $G
\Rightarrow F$.
\end{rem}
\begin{pf}
This follows from the definitions.
\end{pf}

\begin{rem} \label{laxcolimitnatural}
Let the notation be the same as in the previous remark. For $C \in
Obj \hspace{1mm}  \mathcal{C}$ let $\psi_C:FC \rightarrow GC$ be a
right adjoint to $\phi_C$ such that the unit $\eta_C:1_{GC}
\Rightarrow \psi_C \circ \phi_C$ and counit $\varepsilon_C: \phi_C
\circ \psi_C \Rightarrow 1_{FC}$ are natural isomorphisms. Then $C
\mapsto \psi_C$ is a pseudo natural transformation from $F$ to $G$
and there exist iso modifications $\eta:i_G \rightsquigarrow \psi
\odot \phi$ and $\varepsilon: \phi \odot \psi \rightsquigarrow i_F$
which satisfy the triangle identities, namely $C \mapsto \eta_C$ and
$C \mapsto \varepsilon_C$. In the terminology of \cite{street3},
this means that $F$ and $G$ are equivalent in the 2-category
$Hom[\mathcal{C},Cat]$ of pseudo functors, pseudo natural
transformations, and modifications. The equivalences in
$Hom[\mathcal{C},Cat]$ are precisely the pseudo natural
transformations whose components are equivalences of categories.
\end{rem}
\begin{pf}
Since $\phi_C$ is an equivalence of categories, there exists such a
functor $\psi_C$ with unit and counit as above. For $f:A \rightarrow
B$ in $\mathcal{C}$ define the coherence iso\index{coherence
isomorphism} $\tau_f^{\psi}:Gf \circ \psi_A \Rightarrow \psi_B \circ
Ff$ to be the composition of 2-cells in the following diagram.
$$\xymatrix@R=4pc@C=4pc{FA \ar[d]_{1_{FA}} \ar[r]^{\psi_A} & GA \ar[d]^{1_{GA}}
\ar@{=>}[ld]_{\varepsilon_A}
\\ FA \ar[d]_{Ff} \ar@{=}[rd] & \ar[l]^{\phi_A} GA \ar[d]^{Gf}
\\ FB \ar[d]_{1_{FB}} & \ar[l]_{\phi_B} GB \ar@{=>}[ld]_{\eta_B}  \ar[d]^{1_{GB}}
\\ FB \ar[r]_{\psi_B} &  GB}$$
The middle square commutes because $\phi$ is a 2-natural
transformation. We can see that the assignment $f \mapsto
\tau_f^{\psi}$ is natural after segmenting the naturality diagram
into three inner squares and using the fact that $\phi$ is a
2-natural transformation as follows. Let $f,g:A \rightarrow B$ and
$\mu:f \rightarrow g$ in $\mathcal{C}$.

\begingroup
\vspace{-2\abovedisplayskip} \small
$$\xymatrix@R=4pc{1_{GB} \circ Gf \circ \psi_A
\ar@{=>}[r]^{\underset{\phantom{\eta_B * i_{Gf}*i_{\psi_A}}} {\eta_B
* i_{Gf}*i_{\psi_A}}} \ar@{=>}[d]|{i_{1_{GB}}*G \mu
*i_{\psi_A}} & \psi_B \circ \phi_B \circ  Gf \circ \psi_A
\ar@{=}[r] \ar@{=>}[d]|{i_{\psi_B \circ \phi_B}*G \mu * i_{\psi_A}}
& \psi_B \circ Ff \circ \phi_A \circ \psi_A
\ar@{=>}[r]^{\underset{\phantom{\eta_B * i_{Gf}*i_{\psi_A}}}
{i_{\psi_B}*i_{Ff}*\varepsilon_A}}
\ar@{=>}[d]|{i_{\psi_B}*F\mu*i_{\phi_A \circ \psi_A}} & \psi_B \circ
Ff \circ 1_{FA} \ar@{=>}[d]|{i_{\psi_B}*F \mu *i_{1_{FA}}}
\\ 1_{GB} \circ Gg \circ \psi_A
\ar@{=>}[r]_{\overset{\phantom{\eta_B * i_{Gg}*i_{\psi_A}}} {\eta_B
* i_{Gg}*i_{\psi_A}}} & \psi_B \circ \phi_B \circ  Gg \circ \psi_A
\ar@{=}[r] & \psi_B \circ Fg \circ \phi_A \circ \psi_A
\ar@{=>}[r]_{\overset{\phantom{\eta_B * i_{Gf}*i_{\psi_A}}}
{i_{\psi_B}*i_{Fg}*\varepsilon_A}} & \psi_B \circ Fg \circ 1_{FA}}$$
\endgroup
\noindent The left square and the right square commute because of
the interchange law and the defining property of identity 2-cells.
The middle square commutes because $\phi$ is a 2-natural
transformation. Hence the outermost rectangle commutes and $f
\mapsto \tau_f^{\psi}$ is natural.

Since $F$ and $G$ are strict 2-functors, verifying the unit axiom
for $\psi$ reduces to proving that $\tau_{1_C}^{\psi}$ is
$i_{\psi_C}$ for all $C \in Obj \hspace{1mm}  \mathcal{C}$. That
follows from the definition of $\tau_{1_C}$ and one of the triangle
identities.

Since $F$ and $G$ are strict 2-functors, verifying the composition
axiom for $\psi$ amounts to proving for $\xymatrix{A \ar[r]^f & B
\ar[r]^g & C}$ in $\mathcal{C}$ that the composition $(\tau_g^{\psi}
* i_{Ff}) \odot (i_{Gg} * \tau_f^{\psi})$ in
(\ref{laxcolimitlaxnatural1}) is the same as $\tau_{g \circ
f}^{\psi}$ in (\ref{laxcolimitlaxnatural2}). That follows since the
middle parallelogram in (\ref{laxcolimitlaxnatural2}) is
$i_{\phi_B}$ by the triangle identity. Hence $\psi$ with
$\tau^{\psi}$ satisfies the composition axiom and we conclude that
$C \mapsto \psi_C$ is a pseudo natural transformation $F \Rightarrow
G$.
\begin{equation} \label{laxcolimitlaxnatural1}
\xymatrix@C=4pc@R=4pc{ FA \ar[d]_{1_{FA}} \ar[r]^{\psi_A}
\ar@{}[dr]|(.25){\overset{\varepsilon_A}{\Leftarrow}} & GA
\ar[ld]^{\phi_A} \ar[d]^{Gf}
\\  FA \ar@{}[rd]|(.75){\overset{\eta_B}{\Leftarrow}}
\ar[d]_{Ff} & GB \ar[d]^{1_{GB}} \ar[ld]_{\phi_B}
\\  FB \ar[r]_{\psi_B} \ar[d]_{1_{FB}}
\ar@{}[rd]|(.25){ \overset{\varepsilon_B}{\Leftarrow}} & GB
\ar[d]^{Gg} \ar[ld]^{\phi_B}
\\  FB \ar[d]_{Fg}
\ar@{}[rd]|(.75){ \overset{\eta_C}{\Leftarrow}} & GC \ar[d]^{1_{GC}}
\ar[ld]_{\phi_C}
\\  FC \ar[r]_{\psi_C} & GC}
\end{equation}
\begin{equation} \label{laxcolimitlaxnatural2}
\xymatrix@R=4pc@C=4pc{FA \ar[d]_{1_{FA}} \ar[r]^{\psi_A} & GA
\ar[d]^{1_{GA}} \ar@{=>}[ld]_{\varepsilon_A}
\\ FA \ar[d]_{F(g \circ f)} & \ar[l]_{\phi_A} GA \ar[d]^{G(g \circ f)}
\\ FC \ar[d]_{1_{FC}} &
\ar[l]^{\phi_B} GC \ar[d]^{1_{GC}} \ar@{=}[ul] \ar@{=>}[ld]_{\eta_C}
\\ FC \ar[r]_{\psi_C} &  GC}
\end{equation}
Next we prove that $A \mapsto \eta_A$ is a modification $i_G
\rightsquigarrow \psi \odot \phi$. This requires a proof that
(\ref{defnmodification1}) is the same as (\ref{defnmodification2}).
Let $f,g:A \rightarrow B$ be morphisms in $\mathcal{C}$ and
$\gamma:f \Rightarrow g$ a 2-cell in $\mathcal{C}$. Since $\phi$ is
a 2-natural transformation, we see that (\ref{defnmodification2}) is
$\eta_B * G \gamma$. We proceed by showing that
(\ref{defnmodification1}) is $\eta_B * G \gamma$. Note that
$\tau^{\beta}_{A,B}(g)=\tau^{\psi \circ \phi}_g$ in
(\ref{defnmodification1}) is $(i_{\psi_B}*i_{\phi_B \circ Gg}) \odot
(\tau^{\psi}_g *i_{\phi_A})$ by the remarks on page
\pageref{compositionoflaxnaturals} about coherence
isos\index{coherence isomorphism} for a vertical composition of
pseudo natural transformations. Writing out
(\ref{defnmodification1}) with $\alpha= i_G,$ $\beta=\psi \odot
\phi$, $\Xi=\eta,$ and including many trivial arrows gives
(\ref{laxcolimitlaxnatural3}).

\begingroup
\vspace{-2\abovedisplayskip} \small
\begin{equation} \label{laxcolimitlaxnatural3}
\xymatrix@R=4pc@C=4pc{GA \ar[rr]^{1_{GA}} & \ar@{=>}[d]^{\eta_A}
 & GA \ar[r]^{Gf} \ar@{=>}@<8ex>[d]^{G \gamma}
& GB \ar[rr]^{1_{GB}} & \ar@{=>}[d]^{i_{1_{GB}}} & GB
\\ GA \ar[r]^{\phi_A} \ar@{=>}@<8ex>[d]^{i_{\phi_A}}
& FA \ar[r]^{\psi_A} \ar@{=>}@<8ex>[d]^{i_{\psi_A}} & GA
\ar[r]^(.3){Gg} \ar@{=>}@<8ex>[d]^{i_{Gg}} & GB \ar[rr]_(.3){1_{GB}}
& \ar@{=>}[d]^{\eta_B} & GB
\\ GA \ar[r]_(.3){\phi_A} \ar@{=>}@<8ex>[d]^{i_{\phi_A}}
& FA \ar[r]_(.3){\psi_A} \ar@{=>}@<8ex>[d]^{i_{\psi_A}} & GA
\ar[r]_{Gg} & GB \ar@{=>}[d]^{i_{Fg \circ \phi_A}} \ar[r]^{\phi_B} &
FB \ar[r]^{\psi_B} \ar@{=>}@<8ex>[d]^{i_{\psi_B}} & GB
\\ GA \ar[r]_(.3){\phi_A} \ar@{=>}@<8ex>[d]^{i_{\phi_A}}
& FA \ar[r]_{\psi_A} & GA \ar[r]_{\phi_A}
\ar@{=>}[d]^{\varepsilon_A} & FA \ar@{=>}@<8ex>[d]^{i_{Fg}}
\ar[r]^{Fg} & FB \ar[r]_(.3){\psi_B} \ar@{=>}@<8ex>[d]^{i_{\psi_B}}
& GB
\\ GA \ar[r]_{\phi_A} & FA \ar[rr]_(.3){1_{FA}}
& \ar@{=>}[d]^{i_{\phi_B \circ Gg}} & FA \ar[r]_{Fg} & FB
\ar[r]_(.3){\psi_B} \ar@{=>}@<8ex>[d]^{i_{\psi_B}} & GB
\\ GA \ar[rrr]_{Gg} & & & GB \ar[r]_{\phi_B} & FB \ar[r]_{\psi_B} & GB}
\end{equation}
\endgroup
\noindent Using a triangle identity and contracting all the trivial
identities, we see that the only thing that does not cancel is
$\eta_B * G \gamma$. Hence (\ref{defnmodification1}) is the same as
(\ref{defnmodification2}) and $A \mapsto \eta_A$ is a modification.

One can similarly show that $A \mapsto \varepsilon_A$ is a
modification. The modifications $\eta$ and $\varepsilon$ satisfy the
triangle identities because their constituent arrows do.
\end{pf}

\begin{defn}
A 2-category $\mathcal{C}$ {\it admits bilimits}\index{admits
bilimits|textbf} if every pseudo functor $F:\mathcal{J} \rightarrow
\mathcal{C}$ from a small 1-category $\mathcal{J}$ to $\mathcal{C}$
admits a bilimit in $\mathcal{C}$.
\end{defn}

There are analogous definitions for pseudo limits\index{admits
pseudo limits|textbf}, bicolimits\index{bicolimit}, and pseudo
colimits\index{colimit!pseudo colimit}.  If we view the category
$\mathcal{J}$ as an indexing category\index{indexing category}, then
we can speak of bilimits of diagrams\index{bilimit!bilimit of a
diagram}, \ie we can view a diagram in $\mathcal{C}$ as the image of
a pseudo functor from a source diagram $\mathcal{J}$ to the
2-category $\mathcal{C}$.

The concept of pseudo limit can be further generalized to weighted
pseudo limit. For any small 2-category $\mathcal{C}$ we denote the
small category $Mor_{\mathcal{C}}(A,B)$ by $\mathcal{C}(A,B)$ for
$A,B \in Obj \hspace{1mm}  \mathcal{C}$.

\begin{defn}

Let $\mathcal{C}, \mathcal{D}$ be 2-categories. Let $J:\mathcal{D}
\rightarrow Cat$ and $F:\mathcal{D} \rightarrow \mathcal{C}$ be
pseudo functors. Let $Hom[\mathcal{D}, Cat]$\index{$Hom$} denote the
2-category with pseudo functors \newline $\mathcal{D} \rightarrow
Cat$ as objects, pseudo natural transformations as morphisms, and
modifications as 2-cells. Then $\{J,F\}_{p} \in Obj \hspace{1mm}
\mathcal{C}$ is called a {\it $J$-weighted pseudo limit of
$F$}\index{limit!weighted pseudo
limit|textbf}\index{weighted|textbf} if the strict 2-functors
$\mathcal{C}^{op} \rightarrow Cat$
$$ C \mapsto \mathcal{C}(C, \{J,F\}_{p})$$
$$ C \mapsto Hom[\mathcal{D},Cat](J,\mathcal{C}(C,F-))$$
are 2-isomorphic. The image $\xi:J \Rightarrow
\mathcal{C}(\{J,F\}_{p},F-)$ of $1_{\{J,F\}_{p}}$ under this
2-representation is called the {\it unit}\index{unit|textbf}.
\end{defn}

Street refers to this as the {\it $J$-indexed pseudo limit of
F}\index{limit!indexed pseudo limit|textbf} in \cite{street3},
although now the term weighted\index{weighted|textbf} is used
instead of indexed. This is similar to Kelly's definition in
\cite{kelly2}, except that his definition is for strict 2-functors
$J,F$ and he uses the full sub-2-category $Psd[\mathcal{D},Cat]$ of
$Hom[\mathcal{D},Cat]$ in place of $Hom[\mathcal{D},Cat]$. The
2-category $Psd[\mathcal{D},Cat]$\index{$Psd$} consists of strict
2-functors, pseudo natural transformations, and modifications.

We recover the usual definition of pseudo limit whenever $J$ is the
constant functor which takes everything to the trivial category with
one object. A weighted pseudo limit is said to be {\it
conical}\index{conical|textbf} whenever $J$ is this constant
functor. Another special type of weighted limit called {\it cotensor
product} occurs when $\mathcal{D}$ is the trivial 2-category with
one object and $J$ and $F$ are strict 2-functors. In this case $J$
and $F$ can be identified with objects of $Cat$ and $\mathcal{C}$
respectively. Tensor products\index{tensor product} can be defined
similarly.

\begin{defn}
Let $J \in Obj \hspace{1mm}  Cat$ and $F \in Obj \hspace{1mm}
\mathcal{C}$. Then $\{J,F\} \in Obj \hspace{1mm}  \mathcal{C}$ is
called a {\it cotensor product}\index{cotensor product|textbf} of
$J$ and $F$ if the strict 2-functors $\mathcal{C}^{op} \rightarrow
Cat$
$$C \mapsto \mathcal{C}(C,\{J,F\})$$
$$C \mapsto Cat(J, \mathcal{C}(C,F))$$ are 2-naturally isomorphic.
\end{defn}

\begin{rmk} \label{cotensorunit}
(Kelly) We can rephrase the definition of cotensor product entirely
in terms of the unit\index{unit|textbf} $\pi:J \rightarrow
\mathcal{C}(\{J,F\},F)$. The object $\{J,F\}$ of $\mathcal{C}$ is a
cotensor product of $J$ and $F$ with unit $\pi:J \rightarrow
\mathcal{C}(\{J,F\},F)$ if and only if the functor
$\mathcal{C}(C,\{J,F\}) \rightarrow Cat(J, \mathcal{C}(C,F))$
defined by composition with $\pi$
$$b \mapsto \mathcal{C}(b,F) \circ \pi$$
$$\alpha \mapsto \mathcal{C}(\alpha,F)*i_{\pi}$$ for arrows
$b:C \rightarrow \{J,F\}$ and 2-cells $\alpha:b \rightarrow b'$ in
$\mathcal{C}$ is an isomorphism of categories for all $C \in Obj
\hspace{1mm} \mathcal{C}$. More specifically:
\begin{enumerate}
\item
For every functor $\sigma:J \rightarrow \mathcal{C}(C,F)$ there is a
unique arrow $b:C \rightarrow \{J,F\}$ in $\mathcal{C}$ such that
$\mathcal{C}(b,F) \circ \pi = \sigma$.
\item
For every natural transformation $\Xi:\sigma \Rightarrow \sigma'$
there is a unique 2-cell \newline $\alpha:b \Rightarrow b'$ in
$\mathcal{C}$ such that $\mathcal{C}(\alpha,F) * i_{\pi} =\Xi$.
\end{enumerate}
\end{rmk}

A useful reformulation of an observation by Street on page 120 of
\cite{street3} illustrates the importance of cotensor products in
the context of weighted\index{weighted|textbf} pseudo limits.

\begin{thm} \label{streetpseudo}
(Street) A 2-category $\mathcal{C}$ admits weighted pseudo limits if
and only if it admits 2-products\index{2-product}, cotensor
products, and pseudo equalizers\index{equalizer!pseudo equalizer}.
\end{thm}

\begin{rmk} (Street) Pseudo equalizers\index{equalizer!pseudo equalizer}
can be constructed from cotensor products\index{cotensor product}
and 2-pullbacks\index{2-pullback}, while 2-pullbacks can be
constructed from 2-products\index{2-product} and
2-equalizers\index{2-equalizer}. Thus it is sufficient to require
2-equalizers instead of pseudo equalizers\index{equalizer!pseudo
equalizer} in the previous theorem.
\end{rmk}

\begin{defn}
Let $\mathcal{C}, \mathcal{D}$ be 2-categories. Let $J:\mathcal{D}
\rightarrow Cat$ and $F:\mathcal{D} \rightarrow \mathcal{C}$ be
pseudo functors. As above, let $Hom[\mathcal{D}, Cat]$ denote the
2-category with pseudo functors $\mathcal{D} \rightarrow Cat$ as
objects, pseudo natural transformations as morphisms, and
modifications as 2-cells. Then $\{J,F\}_{b} \in Obj \hspace{1mm}
\mathcal{C}$ is called a {\it $J$-weighted bilimit of
$F$}\index{bilimit!weighted bilimit}\index{weighted|textbf} if the
strict 2-functors $\mathcal{C}^{op} \rightarrow Cat$
$$ C \mapsto \mathcal{C}(C, \{J,F\}_{b})$$
$$ C \mapsto Hom[\mathcal{D},Cat](J,\mathcal{C}(C,F-))$$
are equivalent in the 2-category $Hom[\mathcal{C}^{op}, Cat]$, \ie
there is a pseudo natural transformation going from one to the other
whose arrow components are equivalences of categories. The image
$\xi:J \Rightarrow \mathcal{C}(\{J,F\}_{b},F-)$ of $1_{\{J,F\}_{b}}$
under this birepresentation\index{birepresentation} is called the
{\it unit}\index{unit|textbf}.
\end{defn}

Kelly refers to this in \cite{kelly2} as the {\it J-indexed bilimit
of $F$}\index{bilimit!indexed bilimit}. The concepts
weighted\index{weighted|textbf} bicolimit\index{bicolimit!weighted
bicolimit|textbf} and bitensor product can be defined similarly.
Later we will need bitensor products, so we formulate this precisely
and describe it entirely in terms of the unit like Kelly in
\cite{kelly2}.

\begin{defn}
Let $J \in Obj \hspace{1mm}  Cat$ and $F \in Obj \hspace{1mm}
\mathcal{C}$. Then $J*F \in Obj \hspace{1mm}  \mathcal{C}$ is called
a {\it bitensor product}\index{bitensor product|textbf} of $J$ and
$F$ if the strict 2-functors $\mathcal{C}^{op} \rightarrow Cat$
$$C \mapsto \mathcal{C}(J*F,C)$$
$$C \mapsto Cat(J, \mathcal{C}(F,C))$$ are equivalent in the
2-category $Hom[\mathcal{C}^{op}, Cat]$.
\end{defn}

\begin{rmk} \label{laxtensorunit}
We can rephrase the definition of bitensor product entirely in terms
of the unit\index{unit|textbf} $\pi:J \rightarrow
\mathcal{C}(F,J*F)$. The object $J*F$ of $\mathcal{C}$ is a bitensor
product of $J$ and $F$ with unit $\pi:J \rightarrow
\mathcal{C}(F,J*F)$ if and only if the functor $\mathcal{C}(J*F,C)
\rightarrow Cat(J,\mathcal{C}(F,C))$ defined by $$b \mapsto
\mathcal{C}(F,b) \circ \pi$$ $$\alpha \mapsto \mathcal{C}(F,\alpha)
* i_{\pi}$$ for arrows $b:J*F \rightarrow C$ and 2-cells $\alpha:b
\rightarrow b'$ in $\mathcal{C}$ is an equivalence of categories for
all $C \in Obj \hspace{1mm} \mathcal{C}$.
\end{rmk}

Street points out the dual version of the following theorem on page
120 of \cite{street3}.

\begin{thm} \label{streetlaxcolimit}
 A 2-category $\mathcal{C}$ admits weighted bicolimits\index{bicolimit!weighted bicolimit|textbf} if and
only if it admits bicoproducts\index{bicoproduct|textbf}, bitensor
products, and bicoequalizers\index{bicoequalizer|textbf}.
\end{thm}

Cotensor products, bitensor products, and the theorems above will be
used later to show that the 2-categories of interest to us admit
weighted\index{weighted|textbf} pseudo limits as well as weighted
bicolimits\index{bicolimit!weighted bicolimit}.

\chapter{Weighted Pseudo Colimits in the 2-Category of Small Categories}
\label{sec:laxcolimitsinCat} In this chapter we show constructively
that the 2-category $\mathcal{C}$ of small categories admits pseudo
colimits. The dual version of Theorem \ref{streetpseudo} will imply
that this 2-category also admits weighted\index{weighted} pseudo
colimits. One of the concepts in the proof is the free category
generated by a directed graph.

\begin{defn}
A {\it directed graph}\index{graph}\index{graph!directed
graph|textbf}\index{directed graph|textbf} $G$ consists of a set $O$
of {\it objects}\index{object} and a set $A$ of {\it
arrows}\index{arrow} and two functions $S,T:A \rightarrow O$ called
{\it source}\index{source|textbf} and {\it
target}\index{target|textbf}.
\end{defn}

A directed graph is like a category except composition and identity
arrows are not necessarily defined. Any directed graph $G$ whose
sets of arrows and objects are both small generates a {\it free
category}\index{free category} on $G$, which is also called the {\it
path category}\index{path category} of $G$. Similarly $G$ generates
a {\it free groupoid}\index{free groupoid}. We can force commutivity
of certain diagrams by putting a congruence\index{congruence|textbf}
on the morphism sets of the free category or free groupoid and then
passing to the {\it quotient category}\index{quotient category}. We
use this construction in the proof below. The $S,T$ in the
definition of directed graph will also be used to denote the source
and target of a morphism in a category.

\begin{thm} \label{catlaxcolimits}
The 2-category $\mathcal{C}$ of small categories admits pseudo
colimits\index{colimit!pseudo colimit|textbf}.
\end{thm}

\begin{pf}
Let $\mathcal{J}$ be a small 1-category and $F: \mathcal{J}
\rightarrow \mathcal{C}$ a pseudo functor. Here we view
$\mathcal{J}$ as a 2-category which has no nontrivial 2-cells. The
category $\mathcal{J}$ plays the role of an indexing
category\index{indexing category}. For any $X \in Obj \hspace{1mm}
\mathcal{C}$ let $\Delta_X$ denote the constant 2-functor which
takes every object of $\mathcal{J}$ to $X$, every morphism to $1_X$,
and every 2-cell to the identity 2-cell $i_X:1_X \Rightarrow 1_X$.
Then a pseudo cone from $F$ to $X$ is a pseudo natural
transformation  $F \Rightarrow \Delta_X$. Recall $PseudoCone(F,X)$
denotes the category with objects the pseudo cones from $F$ to $X$
with morphisms the modifications between them. The pseudo colimit of
$F$ is an object $W \in \mathcal{C}$ with a pseudo cone $\pi:F
\Rightarrow \Delta_W$ which are universal in the sense that $(\circ
\pi): Mor_{\mathcal{C}}(W,V) \rightarrow PseudoCone(F,V)$ is an
isomorphism of categories for all small categories $V$.

First we define candidates $W \in Obj \hspace{1mm} \mathcal{C}$ and
$\pi:F \Rightarrow \Delta_W$. Then we show that they are universal.
For each $j \in Obj \hspace{1mm}  \mathcal{J}$ let $A_j$ denote the
small category $Fj$ and let $a_f$ denote the functor $Ff$ between
small categories. Since $F$ is a pseudo functor, for every pair
$f,g$ of morphisms of $\mathcal{J}$ such that $g \circ f$ exists we
have a natural transformation (a 2-cell in the 2-category of small
categories) $\gamma_{f,g}: Fg \circ Ff \Rightarrow F(g\circ f)$. We
define a directed graph with objects $O$ and arrows $A$ as follows.
Let $O = \coprod_{j \in \mathcal{J}} Obj \hspace{1mm} A_j$. There is
a well defined function $p: O \rightarrow Obj \hspace{1mm}
\mathcal{J}$ satisfying $p(Obj \hspace{1mm}  A_j) = \{j\}$ because
this union is disjoint, \ie even if the small categories $A_i$ and
$A_j$ are the same, we distinguish them in the disjoint
union\index{disjoint union} by their indices. Let the collection of
arrows be $A = (\coprod_{j \in \mathcal{J}} Mor \hspace{1mm}  A_j)
\coprod \{h_{(x,f)}, h_{(x,f)}^{-1}:(x,f) \in O \times Mor
\hspace{1mm}  \mathcal{J}$ such that $p(x)=Sf\}$ where the elements
of $\coprod_{j \in \mathcal{J}} Mor \hspace{1mm}  A_j$ have the
obvious source and target while $Sh_{(x,f)} = x$ and $Th_{(x,f)} =
a_f(x)$.  Let $W'$ be the free category generated by this graph. We
put the smallest congruence\index{congruence} $\sim$ on $Mor
\hspace{1mm}  W'$ such that:
\begin{itemize}
\item
All of the relations in each $A_i$ are contained in $\sim$, \ie for
$m,n \in Mor \hspace{1mm}  A_i \subseteq Mor \hspace{1mm} W'$ with
$Sn=Tm$ we have $n \circ_{W'} m \sim n \circ_{A_i} m$ where the
composition on the left is the composition in the free category $W'$
and the composition on the right is the composition in the small
category $A_i$.
\item
For all $f,g \in Mor \hspace{1mm}  \mathcal{J}$ with $Sg=Tf$ and all
$x \in Obj \hspace{1mm} A_{Sf}$ we have \linebreak $\gamma_{f,g}(x)
\circ_{W'} h_{(a_f(x),g)} \circ_{W'} h_{(x,f)} \sim h_{(x, g \circ
f)}$ and also every identity $1_x \in A_i$ is congruent to the
identity in the free category on the object $x$.
\item
For all $i,j \in Obj \hspace{1mm}  \mathcal{J}$ and all $f \in
Mor_{\mathcal{J}}(i,j)$ and all morphisms $m: x \rightarrow y$ of
$A_i$ we have $h_{(y,f)} \circ_{W'} m \sim a_f(m) \circ_{W'}
h_{(x,f)}$.
\item
For all $j \in Obj \hspace{1mm}  \mathcal{J}$ and all $x \in Obj
\hspace{1mm} A_j$ we have $(\delta_{j\ast}^F)_x \sim h_{(x,1_j)}$
where $\ast$ denotes the unique object of the terminal
object\index{terminal object} $\mathbf{1}$ in the category of small
categories and $\delta_{j\ast}^F$ is the natural transformation
$\delta_{j}^F$ evaluated at $\ast$.
\item
For all $h_{(x,f)}$ from above we have $h_{(x,f)}^{-1} \circ_{W'}
h_{(x,f)} \sim 1_x$ and  $h_{(x,f)} \circ_{W'} h_{(x,f)}^{-1} \sim
1_{a_fx}$.
\end{itemize}
Define $W$ to be the quotient category of the free category $W'$ by
the congruence\index{congruence} $\sim$.  This is the candidate for
the pseudo colimit.

Now we define a pseudo natural transformation $\pi:F \Rightarrow
\Delta_W$ and its coherence 2-cells $\tau$, \ie we define an element
of $PseudoCone(F,W)$.
 For each object $j
\in Obj \hspace{1mm}  \mathcal{J}$ we need a  morphism in
$\mathcal{C}$ (\ie a functor) $\pi_j:Fj=A_j \rightarrow
W=\Delta_W(j)$. Define $\pi_j:A_j \rightarrow W$ to be the inclusion
functors $A_j \hookrightarrow W$.  In order for $\pi$ to be a pseudo
natural transformation, this assignment must be natural up to
coherence 2-cell, \ie for all $i,j \in Obj \hspace{1mm} \mathcal{J}$
we should have a natural isomorphism $\tau_{i,j}$ of the following
sort.
$$\xymatrix@R=3pc@C=3pc{Mor_{\mathcal{J}}(i,j) \ar[r]^F \ar[d]_{\Delta_W} &
Mor_{\mathcal{C}}(A_i,A_j) \ar[d]^{\pi_j \circ}
\\ Mor_{\mathcal{C}}(W,W) \ar[r]_{\circ \pi_i} \ar@{=>}[ur]^{\tau_{i,j}} &
Mor_{\mathcal{C}}(A_i,W)}$$ Evaluating this diagram at a morphism
$f:i \rightarrow j$ of $\mathcal{J}$ we should have a natural
isomorphism between functors $\tau_{i,j}(f):\pi_i \Rightarrow \pi_j
\circ a_f$.  In other words, $\tau_{i,j}(f)$ should be a 2-cell in
the 2-category $\mathcal{C}$ of small categories. For each $x \in
Obj \hspace{1mm}  A_i$ define $\tau_{i,j}(f)_x: \pi_i(x)=x
\rightarrow a_f(x)= \pi_j \circ a_f(x)$ to be the isomorphism
$h_{(x,f)}$.

\begin{lem} \label{pilaxnatural}
The map $\pi:F \Rightarrow \Delta_W$ is a pseudo natural
transformation with coherence 2-cells given by the natural
isomorphisms $\tau$.
\end{lem}
\begin{pf}
First we show for fixed $f:i \rightarrow j$ that the assignment $Obj
\hspace{1mm}  A_i \ni x \mapsto \tau_{i,j}(f)_x \in Mor_W(\pi_i(x),
\pi_j \circ a_f(x))$ is a natural transformation. To this end, let
$m:x \rightarrow y$ be a morphism in the small category $A_i$. By
definition, $\tau_{i,j}(f)_x = h_{(x,f)}$, $\tau_{i,j}(f)_y =
h_{(y,f)}$, $\pi_i(m)=m$, $\pi_i(x)=x$, $\pi_j \circ a_f (x) = a_f
(x)$, and $\pi_j \circ a_f (m) = a_f(m)$. Some similar statements
hold for the object $y$. The third requirement on the
congruence\index{congruence} in $W'$ gives us the following
commutative diagram in the small category $W$.
$$\xymatrix@R=3pc@C=3pc{x \ar[r]^-{h_{(x,f)}} \ar[d]_m & a_f(x) \ar[d]^{a_f(m)} \\ y
\ar[r]_-{h_{(y,f)}} & a_f(y)}$$ Using the identities just mentioned,
the commutivity of this diagram says precisely that $x \mapsto
\tau_{i,j}(f)_x$ is a natural transformation. Thus
$\tau_{i,j}(f):\pi_i \Rightarrow \pi_j \circ a_f$ is a natural
transformation between functors, \ie a 2-cell in the 2-category
$\mathcal{C}$ of small categories.

The assignment $f \mapsto \tau_{i,j}(f)$ for fixed $i,j$ is natural
because the category $Mor_{\mathcal{J}}(i,j)$ has no nontrivial
morphisms. Thus $\tau_{i,j}$ is a natural transformation between the
indicated functors.

Next we verify the composition axiom for pseudo natural
transformations which involves $\tau$ and $\gamma$. The diagram
states that $\tau$ must satisfy for all $\xymatrix@1{i \ar[r]^f & j
\ar[r]^g & k}$ in $\mathcal{J}$ the coherence axiom $(i_{\pi_k} \ast
\gamma_{f,g}) \odot (\tau_{j,k}(g) \ast i_{a_f}) \odot (i_{1_W} \ast
\tau_{i,j}(f))= \tau_{i,k}(g \circ f) \odot (i_{1_W} \ast
i_{\pi_i})$ as natural transformations. This coherence is satisfied
because of the second requirement on the relation in $W'$ for each
$x \in Obj \hspace{1mm}  A_i$ which states $\gamma_{f,g}(x) \circ
\tau_{j,k}(g)_{a_f(x)} \circ \tau_{i,j}(f)_x = \tau_{i,k}(g \circ
f)_x$.  Note that $(i_{\pi_k} \ast
\gamma_{f,g})(x)=\pi_k(\gamma_{f,g}(x))= \gamma_{f,g}(x)$.

Lastly we verify the unit axiom for pseudo natural transformations
which involves $\tau$ and $\delta$. This coherence requires the
commutivity of the following diagram for all $j \in Obj \hspace{1mm}
\mathcal{J}$.
$$\xymatrix@R=3pc@C=4pc{\pi_j \ar@{=>}[r]^{i_{\pi_j}} \ar@{=>}[d]_{i_{\pi_j}} &
1_W \circ \pi_j \ar@{=>}[r]^-{\delta_{j\ast}^{\Delta_W} \ast
i_{\pi_j}} & \Delta_W(1_j) \circ \pi_j
\ar@{=>}[d]^{\tau_{1_j}=\tau_{j,j}(1_j)}
\\ \pi_j \circ 1_{Fj} \ar@{=>}[rr]_{i_{\pi_j}\ast \delta_{j\ast}^F} &
& \pi_j \circ F(1_j)}$$ Here $\delta_j^{\Delta_W}$ and $\delta_j^F$
are the natural transformations associated to the pseudo functors
$\Delta_W$ and $F$ which make them preserve the identity morphisms
$1_j$ up to coherence 2-cell. In fact, $\delta_{j\ast}^{\Delta_W}$
is trivial. The coherences $\delta_j^{\Delta_W}$ and $\delta_j^F$
fill in the following diagrams for all objects $j$ of $\mathcal{J}$.
$$\xymatrix@R=3pc@C=3pc{\mathbf{1} \ar[r]^-{u_j} \ar@{=}[d] & Mor_{\mathcal{J}}(j,j)
\ar[d]^{\Delta_W} \\ \mathbf{1} \ar[r]_-{u_W}
\ar@{=>}[ur]^{\delta_j^{\Delta_W}} & Mor_{\mathcal{C}}(W,W)}$$
$$\xymatrix@R=3pc@C=3pc{\mathbf{1} \ar[r]^-{u_j} \ar@{=}[d] & Mor_{\mathcal{J}}(j,j)
\ar[d]^{F} \\ \mathbf{1} \ar[r]_-{u_{Fj}}
\ar@{=>}[ur]^{\delta_j^{F}} & Mor_{\mathcal{C}}(Fj,Fj)}$$

Using the fact that $\delta_j^{\Delta_W}$ evaluated on the unique
object $*$ of $\mathbf{1}$ gives the identity 2-cell $i_W:1_W
\Rightarrow 1_W$, the desired coherence diagram simplifies to the
following.
$$\xymatrix@R=3pc@C=3pc{\Delta_W(1_j) \circ \pi_j \ar@{=>}[dr]^{\hspace{2mm}\tau_{j,j}(1_j)}
\ar@{=>}[d]_{i_{\pi_j}}  \\ \pi_j \circ 1_{Fj}
\ar@{=>}[r]_-{i_{\pi_j} \ast \delta_{j\ast}^F} & \pi_j \circ
F(1_j)}$$

Recall that $(\delta_{j\ast}^F)_x = h_{(x,1_j)}$ in $W$ by the
fourth requirement on the congruence\index{congruence} in $W'$. By
definition we also have $h_{(x,1_j)}=\tau_{j,j}(1_j)_x$. This
implies $(\delta_{j\ast}^F)_x = h_{(x,1_j)}=\tau_{j,j}(1_j)_x$ and
the simplified diagram commutes because $\pi_j$ is the inclusion
functor. Hence the required coherence diagram involving $\tau$ and
$\delta$ is actually satisfied.

Thus $\pi:F \Rightarrow \Delta_W$ is a pseudo natural transformation
with the indicated coherence 2-cells.
\end{pf}

Now we must show that the small category $W$ and the pseudo natural
transformation $\pi:F \Rightarrow \Delta_W$ are universal in the
sense that the functor $\phi:Mor_{\mathcal{C}}(W,V) \rightarrow
PseudoCone(F,V)$ defined by $\phi (b) = b \circ \pi$ for objects $b$
is an isomorphism of categories for all objects $V$ of
$\mathcal{C}$. More precisely, $\phi$ is defined for $b \in Obj
\hspace{1mm} Mor_{\mathcal{C}}(W,V)$ and $j \in Obj \hspace{1mm}
\mathcal{J}$ as $\phi(b)(j)=b \circ \pi_j$. The coherence 2-cells
for the pseudo cone $\phi (b)$ are $i_b \ast \tau_{i,j}(f)$ for all
$f:i \rightarrow j$ in $\mathcal{J}$.
 For morphisms $\gamma:b
\Rightarrow b'$ in $Mor \hspace{1mm}  Mor_{C}(W,V)$ we define
$\phi(\gamma):b\circ \pi \rightsquigarrow b' \circ \pi$ to be the
modification which takes $j \in Obj \hspace{1mm}  \mathcal{J}$ to
$\phi(\gamma)(j)=\gamma \ast i_{\pi_j}$.  In the following, $V$ is a
fixed object of the 2-category $\mathcal{C}$ of small categories.

\begin{lem} \label{phifunctor}
The map $\phi:Mor_{\mathcal{C}}(W,V) \rightarrow PseudoCone(F,V)$ is
a functor.
\end{lem}
\begin{pf}
Let $b \in Obj \hspace{1mm}  Mor_{\mathcal{C}}(W,V)$ be a functor
and $i_b:b \Rightarrow b$ its identity natural transformation. Then
obviously $\phi(i_b)(j)=i_b \ast i_{\pi_j}:b \circ \pi_j \Rightarrow
b \circ \pi_j$ is the identity natural transformation $i_{b\circ
\pi_j}$ for all $j \in Obj \hspace{1mm} \mathcal{J}$ and thus
$\phi(i_b)$ is the identity modification.  Hence $\phi$ preserves
identities.

To verify that $\phi$ preserves compositions, let $\gamma:b
\Rightarrow b'$ and $\gamma':b' \Rightarrow b''$ be natural
transformations. Then for each $j \in Obj \hspace{1mm} \mathcal{J}$
we have $\phi(\gamma' \odot \gamma)(j)=(\gamma' \odot \gamma) \ast
i_{\pi_j}=(\gamma' \odot \gamma) \ast (i_{\pi_j} \odot i_{\pi_j})$.
By the interchange law we have $(\gamma' \odot \gamma) \ast
(i_{\pi_j} \odot i_{\pi_j})=$ $(\gamma' \ast i_{\pi_j}) \odot
(\gamma \ast i_{\pi_j})=(\phi(\gamma')(j)) \odot
(\phi(\gamma)(j))=(\phi(\gamma') \diamond \phi(\gamma))_j$ where the
last equality follows from the definition of vertical composition of
modifications. Thus $\phi(\gamma' \odot \gamma)=\phi(\gamma')
\diamond \phi(\gamma)$ and $\phi$ preserves compositions. Thus
$\phi$ is a functor.

\end{pf}

The purpose of the next few lemmas is to exhibit an inverse functor
$\psi$ for $\phi$.

\begin{lem} \label{psifunctor}
There is a functor $\psi: PseudoCone(F,V) \rightarrow
Mor_{\mathcal{C}}(W,V)$.
\end{lem}
\begin{pf}
First we define $\psi$ for objects. Then we define $\psi$ for
morphisms. Finally we verify that $\psi$ is a functor.

Let $\pi'$ be an object of $PseudoCone(F,V)$, \ie $\pi': F
\Rightarrow \Delta_V$ is a pseudo natural transformation with
coherence 2-cells $\tau'$ up to which $\pi'$ is natural. To define a
functor $\psi \pi'= b \in Obj \hspace{1mm} Mor_{\mathcal{C}}(W,V)$
we use the universal mapping property of the quotient category $W$
as follows.  Define an auxiliary functor $d:W' \rightarrow V$ as the
functor induced by the map of directed graphs below which is also
called $d$.

\begin{itemize}
\item
For all $i\in Obj \hspace{1mm}  \mathcal{J}$ and $x \in Obj
\hspace{1mm}  A_i \subseteq Obj \hspace{1mm}  W'$ let
$$dx:=\pi_i'x.$$
\item
For all $i\in Obj \hspace{1mm}  \mathcal{J}$, $x,y \in Obj
\hspace{1mm}  A_i$, and all $g \in Mor_{A_i}(x,y) \subseteq Mor_{W'}
(x,y)$ let $$dg:=\pi_i'g.$$
\item
For all $i,j \in Obj \hspace{1mm}  \mathcal{J}$, $f \in
Mor_{\mathcal{J}}(i,j)$, and all $x \in Obj \hspace{1mm}  A_i
\subseteq Obj \hspace{1mm}  W'$ define $$d(h_{(x,f)}):=
\tau_{i,j}'(f)_x:\pi_i'x \rightarrow\pi_j' \circ a_fx$$
$$d(h_{(x,f)}^{-1}):= \tau_{i,j}'(f)_x^{-1}: \pi_j' \circ a_fx
\rightarrow \pi_i'x.$$
\end{itemize}

We claim that $d$ preserves the congruence\index{congruence} placed
on the category $W'$. Following the order in the definition of
$\sim$ we have the verifications:
\begin{itemize}
\item
For $m,n \in Mor \hspace{1mm}  A_i \subseteq Mor \hspace{1mm}  W'$
with $Sn=Tm$ we have $d(n \circ_{W'} m)=dn \circ_V dm = \pi_in
\circ_V \pi_im=\pi_i(n \circ_{A_i} m)=d(n \circ_{A_i} m)$ and for
all $1_x \in A_i$ we have $d1_x=\pi_i'(1_x)=1_{\pi_i'x}$ because
$\pi_i'$ is a functor. But $1_{\pi_i'x}$ is also the same as $d$
applied to the identity on $x$ in the free category $W'$.
\item
Since $\pi'$ is a pseudo natural transformation, for all
$\xymatrix@1{i \ar[r]^f & j \ar[r]^g & k}$ in $\mathcal{J}$ we have
\newline $(i_{\pi_k'} \ast \gamma_{f,g}) \odot (\tau_{j,k}'(g) \ast
i_{a_f}) \odot (i_{1_V} \ast \tau_{i,j}'(f))= \tau_{i,k}'(g \circ f)
\odot (i_{1_V} \ast i_{\pi_i})$ as natural transformations.
Evaluating this at $x \in Obj \hspace{1mm}  A_i$ yields $$(\pi_k'
\gamma_{f,g}(x)) \circ \tau_{j,k}'(g)_{a_fx} \circ
\tau_{i,j}'(f)_x=\tau_{i,k}'(g \circ f)_x.$$  This says precisely
$d(\gamma_{f,g}(x) \circ_{W'} h_{(a_f(x),g)} \circ_{W'} h_{(x,f)}) =
d(h_{(x, g \circ f)})$.
\item
For all $i,j \in Obj \hspace{1mm}  \mathcal{J}$, all $f \in
Mor_{\mathcal{J}}(i,j)$, and all morphisms $m: x \rightarrow y$ of
$A_i$ we have to show $d(h_{(y,f)} \circ_{W'} m) = d(a_f(m)
\circ_{W'} h_{(x,f)})$. Writing out $d$, we see that this is the
same as verifying $\tau_{i,j}'(f)_y \circ_V \pi_i'm = (\pi_j' \circ
a_f)m \circ_V \tau_{i,j}'(f)_x$, which is true because the
assignment $x \mapsto \tau_{i,j}'(f)_x$ is a natural transformation
from $\pi_i'$ to $\pi_j' \circ a_f$.
\item
For all $j \in Obj \hspace{1mm}  \mathcal{J}$ and all $x \in Obj
\hspace{1mm} A_j$ we have to show $d(\delta_{j\ast}^F)_x =
dh_{(x,1_j)}$.  Writing out $d$ we see that this is the same as
verifying $\pi_j'(\delta_{j\ast}^F)_x=\tau_{j,j}'(1_j)_x$. Since
$\pi'$ is a pseudo natural transformation from $F$ to $\Delta_V$,
the natural transformation $\tau'$ must satisfy the coherence
$(i_{\pi_j'} \ast \delta_{j\ast}^F) \odot
i_{\pi_j'}=\tau_{j,j}'(1_j) \odot (i_{1_V} \ast i_{\pi_j'}) \odot
i_{\pi_j'}$ as natural transformations. Evaluating this coherence at
$x \in Obj \hspace{1mm}  A_j$ we get $\pi_j'(\delta_{j\ast}^F)_x
\circ 1_{\pi_j'x}=\tau_{j,j}'(1_j)_x \circ 1_{\pi_j'x} \circ
1_{\pi_j'x}$, which implies $d(\delta_{j\ast}^F)_x = dh_{(x,1_j)}$
by the remarks above.
\item
For all $i,j \in Obj \hspace{1mm}  \mathcal{J}$, $f \in
Mor_{\mathcal{J}}(i,j)$, and all $x \in Obj \hspace{1mm}  A_i
\subseteq Obj \hspace{1mm}  W'$ we have $d(h_{(x,f)}^{-1} \circ_{W'}
h_{(x,f)})=\tau_{i,j}'(f)^{-1}_x \circ
\tau_{i,j}'(f)_x=1_{\pi_j'x}=d(1_x)$ and similarly $d(h_{(x,f)}
\circ_{W'} h^{-1}_{(x,f)})=d(1_{a_fx})$.
\end{itemize}

Thus $d:W' \rightarrow V$ is a functor that preserves the
congruence\index{congruence} on $W'$. By the universal mapping
property of quotient category $W$ of $W'$, there exists a unique
functor $b:W \rightarrow V$ which factors $d$ via the projection.
Define $\psi(\pi'):=b \in Obj \hspace{1mm}  Mor_{\mathcal{C}}(W,V)$.
This is how $\psi$ is defined on the objects of the category
$PseudoCone(F,V)$.

Next we define $\psi$ on morphisms of the category
$PseudoCone(F,V)$. Let $\Xi:\sigma \rightsquigarrow \sigma'$ be a
morphism in $PseudoCone(F,V)$, \ie $\Xi$ is a modification from the
pseudo natural transformation $\sigma:F
 \Rightarrow \Delta_V$ to the pseudo natural transformation $\sigma':F
\Rightarrow \Delta_V$.  Let $\tau$ and $\tau'$ respectively denote
the natural transformations that make the pseudo natural
transformations $\sigma$ and $\sigma'$ natural up to cell.  We
define a morphism $\psi(\Xi)$ of $Mor_{\mathcal{C}}(W,V)$ as
follows. Note that such a morphism is by definition a natural
transformation between functors from the small category $W$ to the
small category $V$.  Since $\Xi$ is a modification, we have a 2-cell
$\Xi_i:\sigma_i \Rightarrow \sigma_i'$ in the category $\mathcal{C}$
for each $i \in Obj \hspace{1mm}  \mathcal{J}$. Let $b,b'$ denote
the respective functors $\psi(\sigma),\psi(\sigma'):W \rightarrow
V$. For $x \in Obj \hspace{1mm}  A_i \subseteq Obj \hspace{1mm}  W$
define $\psi(\Xi)_x:bx=\sigma_ix \rightarrow \sigma'_ix=b'x$ to be
 $\Xi_i(x):\sigma_ix \rightarrow \sigma_i'x$. The following two commutative
diagrams show that $\psi(\Xi)$ is a natural transformation. For $x,y
\in Obj \hspace{1mm}  A_i$ and $m \in Mor_{A_i}(x,y) \subseteq
Mor_W(x,y)$ the diagram
$$\xymatrix@R=3pc@C=3pc{bx \ar[r]^{\Xi_ix} \ar[d]_{\sigma_im=bm} & b'x
\ar[d]^{b'm=\sigma_i'm}\\ by \ar[r]_{\Xi_iy} & b'y}$$ in $V$
commutes because $\Xi_i:\sigma_i \Rightarrow \sigma_i'$ is a natural
transformation.  For a morphism $f:i \rightarrow j$ in $\mathcal{J}$
the diagram
$$\xymatrix@R=3pc@C=4pc{bx \ar[r]^{\Xi_ix} \ar[d]_{\tau_{i,j}(f)_x=bh_{(x,f)}} & b'x
\ar[d]^{b'h_{(x,f)}=\tau_{i,j}'(f)_x} \\ ba_f(x)
\ar[r]_{\Xi_ja_f(x)} & b'a_f(x)}$$ commutes because of the coherence
in the definition of modification and because of the definitions of
$b,b'$ on $h_{(x,f)}$. We see this by taking $\gamma=i_f$ in
diagrams (\ref{defnmodification1}) and (\ref{defnmodification2}) in
the definition of modification. An inductive argument shows that
$\psi(\Xi)$ is natural for all other arrows in $W$ as well. Hence
$\psi(\Xi):\psi(\sigma) \Rightarrow \psi(\sigma')$ is a morphism in
the category $Mor_{\mathcal{C}}(W,V)$.

Lastly we verify that $\psi$ is a functor, \ie that $\psi$ preserves
the identity modifications and the composition of modifications. Let
$\Xi:\sigma \rightsquigarrow \sigma$ be the identity modification
belonging to a pseudo natural transformation $\sigma:F \Rightarrow
\Delta_V$. This means that $\Xi_i:\sigma_i \Rightarrow \sigma_i$ is
the identity natural transformation for the functor $\sigma_i:A_i
\rightarrow V$. For all $i\in Obj \hspace{1mm} \mathcal{J}$ and all
$x \in Obj \hspace{1mm} A_i$ we have by definition of $\psi$ that
$\psi(\Xi)_x:\psi(\sigma)x=\sigma_ix \rightarrow
\sigma_ix=\psi(\sigma)x$ is $\Xi_i(x):\sigma_ix \rightarrow
\sigma_ix$, which is the identity morphism on the object $\sigma_ix$
of the small category $V$ by hypothesis.  Hence
$\psi(\Xi):\psi(\sigma) \rightarrow \psi(\sigma)$ is the identity
natural transformation and $\psi$ preserves identity modifications.

To verify that $\psi$ preserves compositions, let $\Xi:\sigma
\rightsquigarrow \sigma'$ and $\Xi':\sigma' \rightsquigarrow
\sigma''$ be modifications. Then the vertical composition of
modifications (which makes $PseudoCone(F,V)$ a category) is defined
as $(\Xi' \diamond \Xi)_i := \Xi'_i \odot \Xi_i$ where $\Xi'_i \odot
\Xi_i$ is the vertical composition of the natural transformations
$\Xi_i:\sigma_i \Rightarrow \sigma_i'$ and $\Xi_i':\sigma_i'
\Rightarrow \sigma_i''$ as usual. Then for all $i\in Obj
\hspace{1mm} \mathcal{J}$ and all $x \in Obj \hspace{1mm}  A_i
\subseteq Obj \hspace{1mm}  W$ we have $\psi(\Xi' \diamond \Xi)_x =
(\Xi' \diamond \Xi)_i(x) =(\Xi_i' \odot \Xi_i)_x = \Xi_i'(x) \circ
\Xi_i(x) = \psi(\Xi')_x \circ \psi(\Xi)_x = (\psi(\Xi') \odot \psi
(\Xi))_x$. Thus $\psi(\Xi' \diamond \Xi)=\psi(\Xi') \odot \psi(\Xi)$
and $\psi$ preserves compositions of modifications. Hence $\psi$ is
a functor.
\end{pf}

\begin{lem} \label{compositea}
The functor $\phi \circ \psi:PseudoCone(F,V) \rightarrow
PseudoCone(F,V)$ is the identity functor.
\end{lem}
\begin{pf}
First we verify this for objects, then for morphisms. Let $\pi':F
\Rightarrow \Delta_V$ be a pseudo natural transformation with
coherence isomorphisms $\tau'$\index{coherence isomorphism}.  Let
$b=\psi(\pi')$.  Then using the definitions of $b$ in Lemma
\ref{psifunctor} and the definition of $\pi$ above we evaluate
$\phi(\psi(\pi'))$ at each object $i$ of $\mathcal{J}$ and compare
the resulting functor $\phi(\psi(\pi'))_i$ to the functor $\pi'_i$.
Formally this is:
\begin{itemize}
\item
For all $x \in Obj \hspace{1mm}  A_i$, we have
$$\phi(\psi(\pi'))_ix=\phi(b)_ix=(b \circ \pi_i)x=bx=\pi_i'x.$$
\item
For all $x,y \in Obj \hspace{1mm}  A_i$ and all $g \in
Mor_{A_i}(x,y)$ we have $$\phi(\psi(\pi'))_ig=\phi(b)_ig=(b \circ
\pi_i)g=bg=\pi_i'g.$$
\end{itemize}
Thus $\phi(\psi(\pi'))=\pi'$ for all objects $\pi'$ of the category
$PseudoCone(F,V)$.  Hence $\phi \circ \psi$ is the identity on
objects.

Next we verify the lemma for morphisms. Let $\Xi:\sigma
\rightsquigarrow \sigma'$ be a morphism in the category
$PseudoCone(F,V)$, \ie $\Xi$ is a modification from the pseudo
natural transformation $\sigma:F \Rightarrow \Delta_V$ to the pseudo
natural transformation $\sigma':F \Rightarrow \Delta_V$. Let
$b=\psi(\sigma), b'=\psi(\sigma'):W \rightarrow V$ and
$\gamma=\psi(\Xi):b \Rightarrow b'$ for more convenient notation.
Then $\phi(\psi(\Xi))=\phi(\gamma):b \circ \pi \rightsquigarrow b'
\circ \pi$ is a modification from $\sigma$ to $\sigma'$ by the
result on objects.  For each $j \in Obj \hspace{1mm}  \mathcal{J}$
we have the natural transformation $\phi(\gamma)(j)=\gamma \ast
i_{\pi_j}:b \circ \pi_j \Rightarrow b' \circ \pi_j$.  But this
natural transformation is precisely $\Xi_j:\sigma_j \Rightarrow
\sigma_j'$ by the definition of $\gamma$ via $\psi$. Thus for all
morphisms $\Xi$ of the category $PseudoCone(F,V)$ we have
$\phi(\psi(\Xi))=\Xi$.  Hence $\phi \circ \psi$ is the identity on
morphisms.

\end{pf}

\begin{lem}
The composite functor $\psi \circ \phi:Mor_{\mathcal{C}}(W,V)
\rightarrow Mor_{\mathcal{C}}(W,V)$ is the identity functor.
\end{lem}
\begin{pf}
First we verify this for objects, then on generators for morphisms.
Let $b:W \rightarrow V$ be a functor and $x \in Obj \hspace{1mm} A_i
\subseteq Obj \hspace{1mm} W$. Then $\psi \circ \phi (b) x= \psi(b
\circ \pi)x = (b \circ \pi_i) x = bx$. Similarly for a morphism $g
\in Mor_{A_i}(x,y) \subseteq Mor_W (x,y)$ we have  $\psi \circ \phi
(b) g=$$ \psi(b \circ \pi)g =
$$(b \circ \pi_i) g  = bg$. For morphisms $h_{(x,f)}$, the
analogous calculation is $\psi \circ \phi (b) h_{(x,f)}= \psi(b
\circ \pi) h_{(x,f)} =(i_b*\tau_{i,j}(f))_x= b(\tau_{i,j}(f)_x) =
bh_{(x,f)}$. That follows because the coherence 2-cell up to which
$b \circ \pi$ is natural is
$(i_b*\tau_{i,j}(f))_x=b(\tau_{i,j}(f)_x)$, then we use the third
part of the definition of $\psi$ as well as the definition
$h_{(x,f)} = \tau_{i,j}(f)_x$. Thus $\psi \circ \phi (b) = b$ for
all objects $b$ of the category $Mor_{\mathcal{C}}(W,V)$. Hence
$\psi \circ \phi$ is the identity on the objects of the category
$Mor_{\mathcal{C}}(W,V)$.

Next we verify the lemma for morphisms.  Let $\gamma:b \Rightarrow
b'$ be a morphism in $Mor_{\mathcal{C}}(W,V)$, \ie a natural
transformation from some functor $b$ to some functor $b'$. Let
$\Xi=\phi(\gamma)$, $\sigma=\phi(b)$, and $\sigma'=\phi(b')$ for
more convenient notation. Then by definition $\Xi:\sigma = b \circ
\pi \rightsquigarrow b' \circ \pi = \sigma'$ is the modification
which takes $j \in \mathcal{J}$ to $\gamma*i_{\pi_j}$. Let $x \in
Obj \hspace{1mm}  A_j \subseteq Obj \hspace{1mm}  W$. Then
$\psi(\Xi)_x:\psi(\sigma)x=\sigma_jx \rightarrow
\sigma'_jx=\psi(\sigma')x$ is $\Xi_j(x)=(\gamma \ast i_{\pi_j})_x:(b
\circ \pi)_jx \rightarrow (b' \circ \pi)_jx$. This is described by
the following diagram.
$$\xymatrix@R=3pc@C=3pc{A_j \ar[rr]^{\pi_j} & \ar@{=>}[d]^{i_{\pi_j}} & W \ar[rr]^b  &
\ar@{=>}[d]^{\gamma} & V \\ A_j \ar[rr]_{\pi_j} & & W \ar[rr]_{b'} &
& V}$$

But by definition of $\phi$ and $(b \circ \pi)_j$, we see that
$\Xi_j(x)=\gamma_{\pi_jx}=\gamma_x$ is precisely $\gamma_x:bx
\rightarrow b'x$. Thus $\psi(\Xi)_x=\gamma_x$ and
$\psi(\phi(\gamma))=\psi(\Xi)=\gamma$. Hence $\psi \circ \phi$ is
the identity on the morphisms of the category
$Mor_{\mathcal{C}}(W,V)$.
\end{pf}
\begin{lem} \label{universalcolemma}
The small category $W$ and the pseudo natural transformation $\pi:F
\Rightarrow \Delta_W$ are universal in the sense that the functor
$\phi:Mor_{\mathcal{C}}(W,V) \rightarrow PseudoCone(F,V)$ defined by
$\phi b = b \circ \pi$ for objects $b$ is an isomorphism of
categories for all objects $V$ of $\mathcal{C}$.
\end{lem}
\begin{pf}
This follows immediately from the previous four lemmas because $V$
was an arbitrary object of the 2-category $\mathcal{C}$.
\end{pf}
\begin{lem}
The small category $W$ and the pseudo natural transformation
\newline $\pi:F \Rightarrow \Delta_W$ are a pseudo colimit of the pseudo functor $F:\mathcal{J} \rightarrow
\mathcal{C}$.
\end{lem}
\begin{pf}
This follows from Lemmas \ref{pilaxnatural} and
\ref{universalcolemma}.
\end{pf}

Thus every pseudo functor $F:\mathcal{J} \rightarrow \mathcal{C}$
from a small 1-category $\mathcal{J}$ to the 2-category
$\mathcal{C}$ of small categories admits a pseudo colimit. In other
words, the 2-category $\mathcal{C}$ of small categories admits
pseudo colimits. This completes the proof of Theorem
\ref{catlaxcolimits}.
\end{pf}

\begin{lem} \label{cattensorproducts}
The 2-category of $\mathcal{C}$ of small categories admits tensor
products.
\end{lem}
\begin{pf}
Let $J$ and $F$ be small categories. Then $J*F:=J \times F$ is a
tensor product\index{tensor product} of $J$ and $F$ with unit $\pi:J
\rightarrow Cat(F,J \times F)$ defined by
$$\pi(j)(x):=(j,x)$$ $$\pi(j)(f):=(1_j,f)$$ $$\pi(g)_x:=(g,1_x)$$
for $j \in Obj \hspace{1mm}  J, x \in Obj \hspace{1mm}  F, f \in Mor
\hspace{1mm}  F, g \in Mor \hspace{1mm}  J$. Alternatively one can
see that $Cat(J \times F, C)$ is isomorphic to $Cat(J,Cat(F,C))$ by
the usual adjunction.
\end{pf}

\begin{lem}
The 2-category $\mathcal{C}$ of small categories admits weighted
pseudo colimits.\index{colimit!pseudo
colimit|textbf}\index{colimit!weighted pseudo
colimit|textbf}\index{weighted|textbf}
\end{lem}
\begin{pf}
This 2-category admits pseudo coequalizers\index{coequalizer!pseudo
coequalizer}\index{coequalizer} by Theorem \ref{catlaxcolimits}. It
also admits tensor products\index{tensor product} by Lemma
\ref{cattensorproducts}. It is not difficult to construct
2-coproducts in this 2-category by using disjoint
union\index{disjoint union}. Hence, by the dual version of Theorem
\ref{streetpseudo}, the 2-category $\mathcal{C}$ admits weighted
pseudo limits.
\end{pf}

\begin{rem}
The 2-category of small groupoids admits weighted pseudo
colimits.\index{colimit!weighted pseudo
colimit|textbf}\index{colimit!pseudo
colimit|textbf}\index{groupoid}\index{weighted|textbf}
\end{rem}
\begin{pf}
The proof is the same as in the proof for the 2-category of small
categories except that we replace the free category by the free
groupoid.
\end{pf}

\begin{thm}
The 2-category of small categories and the 2-category of small
groupoids admit weighted bicolimits\index{bicolimit!weighted
bicolimit|textbf}\index{weighted|textbf}\index{groupoid}.
\end{thm}
\begin{pf}
These 2-categories admit weighted pseudo colimits. Every weighted
pseudo colimit is a weighted bicolimit.
\end{pf}

\chapter{Weighted Pseudo Limits in the 2-Category of Small Categories}
\label{limitsection} \label{sec:laxlimitsinCat} Not only does the
2-category $\mathcal{C}$ of small categories admit pseudo colimits,
but it also admits pseudo limits.  In fact we construct them
explicitly in the next proof. The notation remains the same as in
the previous chapter. This description is not new, since the
candidate $L$ in the proof below can be found in \cite{street3}.
Theorem \ref{streetpseudo} allows us to conclude that $\mathcal{C}$
admits weighted\index{weighted} pseudo limits.

\begin{thm} \label{catlaxlimits}
The 2-category $\mathcal{C}$ of small categories admits pseudo
limits.\index{limit!pseudo limit|textbf}
\end{thm}
\begin{pf}
Let $\mathcal{J}$ be a small 1-category and $F: \mathcal{J}
\rightarrow \mathcal{C}$ a pseudo functor. Recall that a pseudo cone
from $X$ to $F$ is a pseudo natural transformation  $\Delta_X
\Rightarrow F$ and that $PseudoCone(X,F)$ denotes the category with
objects the pseudo cones from $X$ to $F$ and morphisms the
modifications between them. A pseudo limit of $F$ is an object $L
\in Obj \hspace{1mm} \mathcal{C}$ with a pseudo cone $\pi:\Delta_L
\Rightarrow F$ which are universal in the sense that $(\pi \circ):
Mor_{\mathcal{C}}(V,L) \rightarrow PseudoCone(V,F)$ is an
isomorphism of categories for all small categories $V$.

First we define candidates $L \in Obj \hspace{1mm} \mathcal{C}$ and
$\pi:\Delta_L \Rightarrow F$. Then we show that they are universal.
For each $j \in Obj \hspace{1mm}  \mathcal{J}$ let $A_j$ denote the
small category $Fj$ as in the proof for the pseudo
colimit\index{colimit!pseudo colimit}. Then the candidate for the
pseudo limit is $L:=PseudoCone(\mathbf{1},F)$, also called the
category of pseudo cones to $F$ on a point. The pseudo natural
transformation candidate $\pi:\Delta_L \Rightarrow F$ is defined for
all objects $\eta: \Delta_{\mathbf{1}} \Rightarrow F$ of $L$ as
$\pi_i(\eta):=\eta_i(\ast)$  for all $i \in Obj \hspace{1mm}
\mathcal{J}$. For morphisms $\Theta:\eta \rightsquigarrow \eta'$ of
$L$ define $\pi_i(\Theta):=\Theta_i(\ast):\eta_i(\ast) \rightarrow
\eta_i'(\ast)$ for all $i \in Obj \hspace{1mm} \mathcal{J}$. Define
the coherence isos $\tau_{i,j}$
$$\xymatrix@R=3pc@C=3pc{Mor_{\mathcal{J}}(i,j) \ar[r]^{\Delta_L} \ar[d]_F &
   Mor_{\mathcal{C}}(L,L) \ar[d]^{\pi_j \circ}
\\ Mor_{\mathcal{C}}(A_i,A_j)  \ar[r]_{\circ \pi_i} \ar@{=>}[ur]^{\tau_{i,j}} &
Mor_{\mathcal{C}}(L,A_j)}$$ belonging to $\pi:\Delta_L \Rightarrow
F$ by $\tau_{i,j}(f)_{\eta}:= \tau_{i,j}^{\eta}(f)_{\ast}$ for all
$f \in Mor_{\mathcal{J}}(i,j)$ and all $\eta \in Obj \hspace{1mm} L$
where $\tau_{i,j}^{\eta}$ is the coherence natural isomorphism
belonging to $\eta:\Delta_{\mathbf{1}} \Rightarrow F$.
$$\xymatrix@R=3pc@C=3pc{Mor_{\mathcal{J}}(i,j) \ar[r]^{\Delta_{\mathbf{1}}} \ar[d]_F &
   Mor_{\mathcal{C}}(\mathbf{1},\mathbf{1}) \ar[d]^{\eta_j \circ}
\\ Mor_{\mathcal{C}}(A_i,A_j)  \ar[r]_{\circ \eta_i}
\ar@{=>}[ur]^{\tau_{i,j}^{\eta}} &
Mor_{\mathcal{C}}(\mathbf{1},A_j)}$$

\begin{lem}
The map $\pi:\Delta_L \Rightarrow F$ is a pseudo natural
transformation with coherence 2-cells given by $\tau$.
\end{lem}
\begin{pf}
First we show that for each $j \in Obj \hspace{1mm}  \mathcal{J}$ we
have a morphism $\pi_j:L=\Delta_L(j) \rightarrow Fj=A_j$ in the
2-category $\mathcal{C}$. We claim that $\pi_j$ is a morphism, \ie a
functor. Let $1_{\eta}=\Theta:\eta \rightsquigarrow \eta$ be the
identity modification of the pseudo cone $\eta: \Delta_{\mathbf{1}}
\Rightarrow F$. This means $\Theta_j = i_{\eta_j}:\eta_j \Rightarrow
\eta_j$ is the identity natural transformation for all $j \in Obj
\hspace{1mm}  \mathcal{J}$. Then
$\pi_j(1_{\eta})=\pi_j(\Theta)=\Theta_j(\ast)=
1_{\eta_j(\ast)}=1_{\pi_j(\eta)}$ and $\pi_j$ preserves identities.
Now let $\Theta, \Xi$ denote modifications in $L$ such that $\Xi
\diamond \Theta$ exists. Then $\pi_j(\Xi \diamond \Theta)=(\Xi
\diamond \Theta)_j(\ast) =\Xi_j \odot
\Theta_j(\ast)=\Xi_j(\ast)\circ \Theta_j(\ast)=\pi_j(\Xi) \circ
\pi_j(\Theta)$. Thus $\pi_j:L \rightarrow A_j$ is a functor.

Next we show that $\tau_{i,j}$ as defined above is a natural
transformation for all $i,j \in Obj \hspace{1mm}  \mathcal{J}$. By
inspecting the definition diagram for $\tau_{i,j}$ above we see that
for all $f \in Mor_{\mathcal{J}}(i,j)$ we should have an element
$\tau_{i,j}(f)$ of $Mor \hspace{1mm} Mor_{\mathcal{C}}(L,A_j)$. To
this end, we claim that $\tau_{i,j}(f): Ff \circ \pi_i \Rightarrow
\pi_j$ is a natural transformation. To see this, let $\Theta:\eta
\rightsquigarrow \eta'$ be a modification, \ie a morphism in the
category $L$. Then by taking $\gamma=i_f$ in the definition of
modification and evaluating the modification diagrams
(\ref{defnmodification1}) and (\ref{defnmodification2}) at $\ast \in
Obj \hspace{1mm} \mathbf{1}$ with
$\alpha=\eta,\beta=\eta',A=i,B=j,\Xi=\Theta$ we obtain the
commutivity of the diagram in the category $A_j$
$$\xymatrix@R=3pc@C=4pc{Ff(\eta_i(\ast)) \ar[r]^-{\tau_{i,j}^{\eta}(f)_{\ast}}
\ar[d]_{Ff(\Theta_i(\ast))} & \eta_j(\ast) \ar[d]^{\Theta_j(\ast)}
 \\ Ff(\eta_i'(\ast)) \ar[r]_-{\tau_{i,j}^{\eta'}(f)_{\ast}} & \eta_j'(\ast)}$$
where $\tau^{\eta}$ and  $\tau^{\eta'}$ denote the coherence natural
transformations belonging to the pseudo cones $\eta$ and $\eta'$
respectively. Using the definitions
$\tau_{i,j}(f)_{\eta}:=\tau_{i,j}^{\eta}(f)_{\ast},
\pi_i(\eta):=\eta_i(\ast)$, and $\pi_i(\Theta):=\Theta_i(\ast)$ we
see that this diagram is
$$\xymatrix@R=3pc@C=4pc{Ff \circ \pi_i(\eta) \ar[r]^-{\tau_{i,j}(f)_{\eta}}
\ar[d]_{Ff \circ \pi_i(\Theta)} & \pi_j(\eta) \ar[d]^{\pi_j(\Theta)}
 \\ Ff \circ \pi_i(\eta') \ar[r]_-{\tau_{i,j}(f)_{\eta'}} & \pi_j(\eta')}$$
which says precisely that $\eta \mapsto \tau_{i,j}(f)_{\eta}$ is
natural for fixed morphisms $f:i \rightarrow j$ of $\mathcal{J}$.
Thus $\tau_{i,j}(f):Ff \circ \pi_i \Rightarrow \pi_j$ is a natural
transformation. On the other hand, the assignment
$Mor_{\mathcal{J}}(i,j) \ni f \mapsto \tau_{i,j}(f)$ is vacuously
natural because the category $Mor_{\mathcal{J}}(i,j)$ is discrete.
Thus $\tau_{i,j}$ is a natural transformation for all $i,j \in Obj
\hspace{1mm} \mathcal{J}$.

The natural isomorphisms $\tau$ satisfy the unit axiom and
composition axiom involving $\delta$ and $\gamma$ because the
individual $\tau^{\eta}$ do.
\end{pf}

Now we must show that the small category $L$ and the pseudo natural
transformation $\pi:\Delta_L \Rightarrow F$ are universal in the
sense that the functor $\phi:Mor_{\mathcal{C}}(V,L) \rightarrow
PseudoCone(V,F)$ defined by $\phi b = \pi \circ b$ for objects $b$
is an isomorphism of categories for all objects $V$ of
$\mathcal{C}$. More precisely, $\phi$ is defined for $b \in Obj
\hspace{1mm} Mor_{\mathcal{C}}(V,L)$ and $j \in Obj \hspace{1mm}
\mathcal{J}$ as $\phi(b)(j)=\pi_j \circ b$. The natural
transformations for the pseudo cone $\phi b$ are $\tau_{i,j}(f) \ast
i_b$ for all $f:i \rightarrow j$ in $\mathcal{J}$. For morphisms
$\gamma:b \Rightarrow b'$ in $Mor \hspace{1mm} Mor_{C}(V,L)$ we
define $\phi(\gamma):\pi \circ b \rightsquigarrow \pi \circ b'$ to
be the modification which takes $j \in Obj \hspace{1mm}
\mathcal{J}$ to $\phi(\gamma)(j)=i_{\pi_j} \ast \gamma$. In the
following, $V$ is a fixed object of the 2-category $\mathcal{C}$ of
small categories.

\begin{lem}
The map $\phi:Mor_{\mathcal{C}}(V,L) \rightarrow PseudoCone(V,F)$ is
a functor.
\end{lem}
\begin{pf}
The proof is analogous to the proof for the $\phi$ of the pseudo
colimit.
\end{pf}

Now we construct a functor $\psi:PseudoCone(V,F) \rightarrow
Mor_{\mathcal{C}}(V,L)$ that is inverse to $\phi$. First we define
$\psi$ for objects, then for morphism. Finally we verify that it is
a functor and inverse to $\phi$. The key observation in the
construction is that we can get a pseudo cone on a point by
evaluating a pseudo cone on an object. This is the essence of the
identification we make below.

\begin{rem} \label{PL1}
Let $Obj \hspace{1mm}  P$ be the subset of the set \newline
$\{(a_i)_i \times (\varepsilon_f)_f \in \prod_{i \in Obj
\hspace{1mm} \mathcal{J}} Obj \hspace{1mm}  A_i \times \prod_{f \in
Mor \hspace{1mm}  \mathcal{J}} Mor \hspace{1mm}  A_{Tf} \vert$ $
\varepsilon_f:Ff(a_{Sf})\rightarrow a_{Tf}$ is iso
\newline for all $f \in Mor \hspace{1mm}  \mathcal{J} \}$ consisting of all
$(a_i)_i \times (\varepsilon_f)_f$ such that:
\begin{itemize}
\item
$\varepsilon_{1_j} \circ \delta^F_{j\ast}(a_j)=1_{a_j}$ for all $j
\in Obj \hspace{1mm}  \mathcal{J}$.
\item
$\varepsilon_g \circ (Fg(\varepsilon_f))=\varepsilon_{g \circ f}
\circ \gamma^F_{f,g}(a_{Sf})$ for all $f,g \in Mor \hspace{1mm}
\mathcal{J}$ such that $g \circ f$ exists.
\end{itemize}
Then $Obj \hspace{1mm}  L$ and $Obj \hspace{1mm}  P$ are in
bijective correspondence via the map $Obj \hspace{1mm}  L
\rightarrow Obj \hspace{1mm}  P$, $\eta \mapsto (\eta_i(\ast))_i
\times (\tau^{\eta}_{Sf,Tf}(f)_{\ast})_f$.
\end{rem}
\begin{pf}
The two conditions express exactly the required coherences for a
pseudo cone $\eta:\Delta_{\mathbf{1}} \Rightarrow F$. Any pseudo
cone $\eta:\Delta_{\mathbf{1}} \Rightarrow F$ is completely
determined by the data listed in the image sequence.
\end{pf}

\begin{rem} \label{PL2}
Let $\eta =(a_i)_i \times (\varepsilon_f)_f$ and $\eta'=(a_i')_i
\times (\varepsilon_f')_f$ be elements of $Obj \hspace{1mm}  P$. Let
$Mor_P(\eta,\eta')$ denote the set of $(\xi_i)_i \in \prod_{i \in
Obj \hspace{1mm}  \mathcal{J}} Mor_{A_i}(a_i,a_i')$ such that
$$\xymatrix@R=3pc@C=3pc{Ff(a_i) \ar[r]^-{\varepsilon_f} \ar[d]_{Ff(\xi_i)} & a_j
\ar[d]^{\xi_j} \\ Ff(a_i') \ar[r]_-{\varepsilon_f'} & a_j'}$$
commutes for all $f:i \rightarrow j$ in $\mathcal{J}$. Then
$Mor_L(\eta,\eta')$ and $Mor_P(\eta, \eta')$ are in bijective
correspondence via the map $Mor_L(\eta,\eta') \rightarrow
Mor_P(\eta, \eta')$, $\Theta \mapsto (\Theta_i(\ast))_i$. Moreover,
the composition $\Theta \diamond \Xi$ in $Mor_L(\eta,\eta')$
corresponds to the componentwise composition in $Mor_P(\eta,
\eta')$.
\end{rem}
\begin{pf}
The diagram is the result of evaluating the coherence stated in
diagrams (\ref{defnmodification1}) and (\ref{defnmodification2}) in
the definition of modification at $\ast$. The claim about
composition follows immediately from the definition of vertical
composition $\diamond$ of modifications.
\end{pf}

\begin{rem} \label{PL3a}
Under the identification above, $P$ is a category and $\pi_j$ is the
projection onto the $j$-th coordinate.
\end{rem}
\begin{pf}
This follows directly from the definition of $\pi$ and the
identification.
\end{pf}

We will use the identification without explanation. Now we define a
functor $\psi(\pi')=b:V \rightarrow L$ for any object $\pi'$ of
$PseudoCone(V,F)$. This will substantiate the comment that
evaluating a pseudo cone on an object gives a pseudo cone on a
point.

\begin{lem} \label{PL3}
Let $\pi':\Delta_V \Rightarrow F$ be a pseudo natural transformation
with coherence natural isomorphisms $\tau'$. For any fixed $x \in
Obj \hspace{1mm} V$ we have $\psi(\pi')(x):=b(x):=(\pi_i'(x))_i
\times (\tau'_{Sf,Tf}(f)_x)_f$ is an element of $Obj \hspace{1mm} P
= Obj \hspace{1mm} L$.
\end{lem}
\begin{pf}
Evaluating the coherences for $\tau$ involving $\delta$ and $\gamma$
at the object $x$ gives the coherences in the definition of $P$.
Thus $b(x) \in Obj \hspace{1mm}  P$ and $b(x)$ is a pseudo cone
$\Delta_{\mathbf{1}} \Rightarrow F$, in other words $b(x)$ is a
pseudo cone on a point.
\end{pf}

\begin{lem} \label{PL4}
Let $\pi':\Delta_V \Rightarrow F$ be a pseudo natural transformation
with coherence natural isomorphisms $\tau'$. Then for any fixed $h
\in Mor_V(x,y)$ we have a modification
$\psi(\pi')(h):=b(h):=(\pi_i'(h))_i:b(x) \rightsquigarrow b(y)$.
This notation means $b(h)_i(\ast):=\pi_i'(h)$.
\end{lem}
\begin{pf}
For notational convenience let $\eta:=b(x):\Delta_{\mathbf{1}}
\Rightarrow F$ and $\eta':= b(y):\Delta_{\mathbf{1}} \Rightarrow F$.
Let $\Theta=b(h)$. Then $\tau_{i,j}'(f)_x =
\tau^{\eta}_{i,j}(f)_{\ast}$ and $ \tau_{i,j}'(f)_y =
\tau^{\eta'}_{i,j}(f)_{\ast}$ and $\Theta_i(\ast)=\pi_i'(h)$ for all
$f:i \rightarrow j$ in $\mathcal{J}$ by the identification. The
naturality of $\tau_{i,j}'(f)$ says $\tau_{i,j}'(f)_y \circ
Ff(\pi_i'(h))=\pi_j'(h) \circ \tau_{i,j}'(f)_x$ for all $f:i
\rightarrow j$ in $\mathcal{J}$. Rewriting this identity using
$\eta, \eta'$, and $\Theta$ gives $\tau^{\eta'}_{i,j}(f)_{\ast}
\circ
 Ff(\Theta_i(\ast)) =
\Theta_j(\ast) \circ \tau^{\eta}_{i,j}(f)_{\ast}$. This last
identity says that the composition of natural transformations
(2-cells)
$$\xymatrix@R=3pc@C=3pc{\mathbf{1} \ar[rr]^{\eta_i} & \ar@{=>}[d]^{\Theta_i} &
 Fi \ar[rr]^{Ff} &
 \ar@{=>}[d]^{Fi_f}  & Fj \\
\mathbf{1} \ar[rr]_{\eta_i'} & & Fi \ar[rr]_{Ff}
\ar@{=>}[d]^{\tau^{\eta'}_{i,j}(f)} & & Fj \\ \mathbf{1}
\ar[rr]_{\Delta_{\mathbf{1}}(f)} & & \mathbf{1} \ar[rr]_{\eta_j'} &
& Fj}$$ is the same as the composition
$$\xymatrix@R=3pc@C=3pc{\mathbf{1} \ar[rr]^{\eta_i} & & Fi \ar[rr]^{Ff}
\ar@{=>}[d]^{\tau_{i,j}^{\eta}(f)} & & Fj \\
\mathbf{1} \ar[rr]^{\Delta_{\mathbf{1}}(f)} &
\ar@{=>}[d]^{\Delta_{\mathbf{1}}(i_f)} & \mathbf{1}
 \ar[rr]^{\eta_j} & \ar@{=>}[d]^{\Theta_j} & Fj \\
\mathbf{1} \ar[rr]_{\Delta_{\mathbf{1}}(f)} & & \mathbf{1}
\ar[rr]_{\eta_j'} & & Fj}$$ of natural transformations for all $f:i
\rightarrow j$ in $\mathcal{J}$. The only 2-cells in the category
$\mathcal{J}$ are of the form $i_f$. Therefore we have verified
diagrams (\ref{defnmodification1}) and (\ref{defnmodification2}) for
$\Theta$ to be a modification. Thus $\psi(\pi')(h)=b(h)=\Theta:\eta
\rightsquigarrow \eta'$ is a modification.
\end{pf}

\begin{lem} \label{PL5}
For any pseudo natural transformation $\pi':\Delta_V \Rightarrow F$
the map $\psi(\pi')=b:V \rightarrow L$ is a functor.
\end{lem}
\begin{pf}
For each $x \in Obj \hspace{1mm} V$ and all $j \in Obj \hspace{1mm}
\mathcal{J}$ we have $b(1_x)_j(\ast)=\pi_j'(1_x)=1_{\pi_j'x}$ since
$\pi_j':V \rightarrow A_j$ is a functor. Hence
$b(1_x)_j=i_{b(x)_j}$. Hence $b(1_x):b(x) \rightsquigarrow b(x)$ is
the identity modification. If $h$ and $\ell$ are morphisms in $V$
such that $\ell \circ h$ exists, then $b(\ell \circ
h)_j(\ast)=\pi_j'(\ell \circ h) = \pi_j'(\ell) \circ \pi_j'(h)
=b(\ell)_j(\ast) \circ b(h)_j(\ast)=(b(\ell)_j \odot b(h)_j)
(\ast)=(b(\ell) \diamond b(h))_j(\ast)$. Hence $b(\ell \circ h) =
b(\ell) \diamond b(h)$ and $b$ preserves compositions.
\end{pf}

\begin{lem} \label{PL6}
Let $\Xi:\alpha \rightsquigarrow \beta$ be a morphism in the
category $PseudoCone(V, F)$. Then $\psi(\Xi):\psi(\alpha)
\Rightarrow \psi(\beta)$ defined by $V \ni x \mapsto (\Xi_i(x))_i
\in Mor_L(\psi(\alpha)x, \psi(\beta)x)$ is a natural transformation.
As in Remark \ref{PL2} above, this definition
means $\psi(\Xi)(x)_i(\ast):=\Xi_i(x)$. 
\end{lem}
\begin{pf}
Since $\Xi:\alpha \rightsquigarrow \beta$ is a modification, for
each object $i$ of $\mathcal{J}$ there is a 2-cell of $\mathcal{C}$
(a natural transformation) $\Xi_i:\alpha_i \Rightarrow \beta_i$ and
these satisfy the condition listed in the definition of
modification. Evaluating this condition in diagrams
(\ref{defnmodification1}) and (\ref{defnmodification2}) at $x \in V$
we see that $(\Xi_i(x))_i:\psi(\alpha)x \rightsquigarrow
\psi(\beta)x$ is a modification. Hence $(\Xi_i(x))_i \in
Mor_L(\psi(\alpha)x, \psi(\beta)x)$.

We claim that $\psi(\Xi)$ is natural, \ie that the diagram
$$\xymatrix@R=3pc@C=4pc{\psi(\alpha)x \ar[r]^{(\Xi_i(x))_i}
\ar[d]_{\psi(\alpha)g=(\alpha_i(g))_i}
& \psi(\beta)x \ar[d]^{(\beta_i(g))_i=\psi(\beta)g}\\
\psi(\alpha)y \ar[r]_{(\Xi_i(y))_i}  & \psi(\beta)y}$$ in $L$
commutes. We only need to verify that the diagram commutes
componentwise, since the vertical composition of modifications
corresponds to the componentwise composition of these sequences
under the identification. But the diagram obviously commutes
componentwise because $\Xi_i:\alpha_i \Rightarrow \beta_i$ is a
natural transformation.
\end{pf}

\begin{thm} \label{calc1}
The map $\psi:PseudoCone(V,F) \rightarrow Mor_{\mathcal{C}}(V,L)$ as
defined in the previous lemmas is a functor.
\end{thm}
\begin{pf}
Suppose $\Xi:\alpha \rightsquigarrow \alpha$ is the identity
modification for a pseudo cone $\alpha:\Delta_V \Rightarrow F$. Then
$\Xi_j=i_{\alpha_j}:\alpha_j \Rightarrow \alpha_j$ for all $j \in
Obj \hspace{1mm}  \mathcal{J}$, so that
$\Xi_j(x)=(i_{\alpha_j})_x=1_{\alpha_j(x)}$. Then $x \mapsto
(1_{\alpha_j(x)})_j$ is the identity morphism $\psi(\alpha)
\rightarrow \psi(\alpha)$ in $Mor_{\mathcal{C}}(V,L)$.

If $\Xi,\Theta$ are modifications in $PseudoCone(V,F)$ such that
$\Theta \diamond \Xi$ exists, then for all $x \in V$ we have
$$\aligned \psi(\Theta \diamond \Xi)(x)
 &=((\Theta \diamond \Xi)_i(x))_i \\
 &=((\Theta_i \odot \Xi_i)(x))_i \\
 &=(\Theta_i(x) \circ \Xi_i(x))_i \\
 &=(\Theta_i(x))_i \diamond (\Xi_i(x))_i \\
 &=\psi(\Theta)(x) \diamond \psi(\Xi)(x) \\
 &=(\psi(\Theta) \odot \psi(\Xi))(x). \endaligned$$
Hence $\psi(\Theta \diamond \Xi)=\psi(\Theta) \odot \psi(\Xi)$ and
$\psi$ preserves compositions.

\end{pf}

Now that we have constructed the functor $\psi$, we prove that it is
inverse to $\phi$.

\begin{lem} \label{calc2}
The functor $\psi$ is a left inverse for $\phi$, \ie $\psi \circ
\phi=1_{Mor_{\mathcal{C}}(V,L)}$.
\end{lem}
\begin{pf}
First we verify the identity on objects. Let $b:V \rightarrow L$ be
an object of $Mor_{\mathcal{C}}(V,L)$. Recall that $\phi(b)$ is the
pseudo natural transformation $\pi \circ b$ with the coherence
natural transformations $\tau_{i,j}'(f)=\tau_{i,j}(f) \ast i_b$ for
all $f:i \rightarrow j$. For
$x \in V$ we have $$\aligned \psi \circ \phi (b)(x) &=\psi(\pi \circ b)(x) \\
&=(\pi_i \circ b(x))_i \times (\tau_{Sf,Tf}'(f)_x)_f \\
&=(\pi_i \circ b(x))_i \times ((\tau_{Sf,Tf}(f) \ast i_b)_x)_f \\
&=(b(x)_i(\ast))_i \times (\tau_{Sf,Tf}(f)_{b(x)})_f \\
&=(b(x)_i(\ast))_i \times (\tau_{Sf,Tf}^{b(x)}(f)_{\ast})_f \text{  by definition}\\
&=b(x) \text{ by the identification.} \endaligned $$  For $g:x
\rightarrow y$ in $V$ we have
$$\aligned \psi \circ \phi (b)(g) &=\psi(\pi \circ b)(g) \\
&=(\pi_i \circ b(g))_i \\
&=(b(g)_i(\ast))_i  \\
&=b(g) \text{ by the identification.}
\endaligned $$ Thus $\psi \circ \phi (b)$ and $b$ agree as
functors.

Next we verify the identity on morphisms. Let $\gamma:b \Rightarrow
b'$ be a natural transformation. Then for $x \in V$ we have
$$\aligned \psi \circ \phi (\gamma)_x
&=\psi(i_{\pi} \ast \gamma)_x \\
&=((i_{\pi_j} \ast \gamma)_x)_j \\
&=(\pi_j(\gamma_x))_j \\
&=(\gamma_{xj}(\ast))_j \\
&=\gamma_x \text{ by the identification.} \endaligned $$ Thus $\psi
\circ \phi (\gamma)$ and $\gamma$ agree as natural transformations
and $\psi \circ \phi=1_{Mor_{\mathcal{C}}(V,L)}$.

Another way to see this is to notice that $\pi_i$ is the projection
onto the $i$-th coordinate.
\end{pf}

\begin{lem} \label{calc3}
The functor $\psi$ is a right inverse for $\phi$, \ie $\phi \circ
\psi=1_{PseudoCone(V,F)}$.
\end{lem}
\begin{pf}
First we verify the identity on objects. Let $\pi':\Delta_V
\Rightarrow F$ be a pseudo cone. For $j \in Obj \hspace{1mm}
\mathcal{J}$ and $x \in V$ we have
$$\aligned (\phi \circ \psi(\pi'))_j(x) &=(\pi \circ \psi(\pi'))_j(x) \\
&=\pi_j \circ \psi(\pi')(x) \\
&=\pi_j((\pi_i'(x))_i \times (\tau_{Sf,Tf}'(f)_x)_f) \\
&=\pi_j'(x).
\endaligned$$
The last equality follows because $\pi_j$ is basically projection
onto the $j$-th coordinate under the identification.

Next we verify the identity on morphisms. Let $\Xi:\alpha
\rightsquigarrow \beta$ be a modification in $PseudoCone(V,F)$. For
$j \in Obj \hspace{1mm} \mathcal{J}$ and
$x \in V$ we have $$\aligned (\phi \circ \psi(\Xi))_j(x) &=(i_{\pi_j} \ast \psi(\Xi))_x\\
&=\pi_j(\psi(\Xi)_x) \\
&=\pi_j((\Xi_i(x))_i) \\
&=\Xi_j(x). \endaligned $$ Thus $\phi \circ \psi(\Xi) = \Xi$ and
$\phi \circ \psi=1_{PseudoCone(V,F)}$.
\end{pf}

\begin{lem}
The small category $L$ with the pseudo cone $\pi:\Delta_L
\Rightarrow F$ is a pseudo limit of the pseudo functor
$F:\mathcal{J} \rightarrow \mathcal{C}$.
\end{lem}
\begin{pf}
The functor $\phi:Mor_{\mathcal{C}}(V,L) \rightarrow
PseudoCone(V,F)$ is an isomorphism of categories by the previous
lemmas. Since $V$ was arbitrary we conclude that $L$ and $\pi$ are
universal.
\end{pf}

Thus every pseudo functor $F:\mathcal{J} \rightarrow \mathcal{C}$
from a small 1-category $\mathcal{J}$ to the 2-category
$\mathcal{C}$ of small categories admits a pseudo limit. In other
words, the 2-category $\mathcal{C}$ of small categories admits
pseudo limits. This completes the proof of Theorem
\ref{catlaxlimits}.
\end{pf}

\begin{lem} \label{catcotensor}
The 2-category $\mathcal{C}$ of small categories admits cotensor
products\index{cotensor product|textbf}.
\end{lem}
\begin{pf}
Let $J$ and $F$ be small categories. Then
$\{J,F\}:=\mathcal{C}(J,F)$ is a cotensor product of $J$ and $F$
with unit $\pi:J \rightarrow \mathcal{C}(\mathcal{C}(J,F),F)$
defined by evaluation.
\end{pf}

\begin{thm} \label{catweightedlaxlimits}
The 2-category $\mathcal{C}$ of small categories admits weighted
pseudo limits.\index{limit!weighted pseudo
limit|textbf}\index{limit!pseudo
limit|textbf}\index{weighted|textbf}
\end{thm}
\begin{pf}
This 2-category admits 2-products\index{2-product}. It also admits
cotensor products and pseudo equalizers\index{equalizer!pseudo
equalizer} by Lemma \ref{catcotensor} and Theorem
\ref{catlaxlimits}. Theorem \ref{streetpseudo} then implies that it
admits weighted pseudo limits.
\end{pf}

\begin{rem}
The 2-category of small groupoids admits weighted pseudo
limits.\index{limit!weighted pseudo limit|textbf}\index{limit!pseudo
limit|textbf}\index{weighted|textbf}\index{groupoid}
\end{rem}
\begin{pf}
The proof is exactly the same as the proof for small categories,
since $L=PseudoCone(\mathbf{1},F)$ is obviously a groupoid when the
target of $F$ is the 2-category of small groupoids.
\end{pf}

\begin{thm}
The 2-category of small categories and the 2-category of small
groupoids admit weighted
bilimits.\index{bilimit|textbf}\index{bilimit!weighted
bilimit|textbf}\index{weighted|textbf}\index{groupoid}
\end{thm}
\begin{pf}
They admit weighted pseudo limits, hence they also admit weighted
bilimits.
\end{pf}

\chapter{Theories and Algebras}
\label{sec:theories} The axioms for a group provide an example for
the concept of a {\it theory} and an example of a group is an {\it
algebra}\index{algebra} over the theory of
groups\index{theory!theory of groups}\index{group!theory of
groups}\index{group}. In this chapter we describe what this means.
Hu and Kriz point out in \cite{hu} that Lawvere's notion of a theory
\cite{lawvere} is equivalent to another notion of theory. We prove
this equivalence. It is well known that the category of algebras
over a theory $T$ is equivalent to the category of algebras for some
monad $C$. We present a version of this. Next we generalize theories
in two ways: theories on a set of objects and theories enriched in
groupoids\index{groupoid}. Theories on a set of objects allow us to
describe algebraic structures on more than one set, such as modules
or theories themselves. They also allow us to describe the free
theory\index{theory!free theory}\index{free theory} on a sequence of
sets. Theories enriched in groupoids will be used in Chapter
\ref{sec:laxTalgebras} to describe pseudo
algebras\index{algebra!pseudo algebra} over a theory $T$ as strict
algebras\index{algebra!algebra over a theory enriched in groupoids}
over a theory $\mathcal{T}$ enriched in groupoids\index{groupoid}.

A theory can also be described as a finitary monad\index{finitary
monad}\index{monad!finitary monad} on the category $Sets$ of small
sets as put forth in \cite{benabou}. Theories on more than one
object are called {\it many-sorted}\index{theory!many sorted
theory|textbf} in the monad\index{monad} description. Free finitary
monads\index{free finitary monad}\index{monad!free finitary monad}
in the enriched and many-sorted contexts can be found in
\cite{kelly3a} and \cite{lack2}. See \cite{street5} for monads in a
general 2-category.

\begin{defn}
A {\it theory}\index{theory|textbf}\index{Lawvere theory|textbf} is
a category $T$ with objects $0,1,2, \dots$ such that $n$ is the
product of $1$ with itself $n$ times in the category $T$ and each
$n$ is equipped with a limiting cone.
\end{defn}

This definition means for each $n \in Obj \hspace{1mm}  T$ we have
chosen morphisms $pr_i:n \rightarrow 1$ for $i=1, \dots,n$ with the
universal property: for any object $m \in Obj \hspace{1mm}  T$ and
morphisms $w_i:m \rightarrow 1$ for $i=1, \dots,n$ there exists a
unique morphism \newline $\prod_{j=1}^n w_j:m \rightarrow
n$ such that the diagram $$\xymatrix@C=3pc@R=3pc{n \ar[r]^{pr_i} & 1 \\
m \ar@{.>}[u]^{\prod_{j=1}^n w_j} \ar[ur]_{w_i} &}$$ commutes for
all $i=1, \dots,n$. In particular $0$ is the terminal
object\index{terminal object} of the category $T$. Note that we do
not require the projection $pr_1:1 \rightarrow 1$ to be the
identity, although it will automatically be an isomorphism. A useful
notation is $T(n):=Mor_T(n,1)$ for $n \in Obj \hspace{1mm}  T$.
Elements of $T(n)$ are called {\it words of arity
n}\index{arity|textbf}\index{word|textbf}.

Another relevant morphism is the following. Let $\iota_i:\{1, \dots,
n_i\} \rightarrow \{1, \dots, n_1 + n_2 + \dots + n_k\}$ be the
injective map which takes the domain to the $i$-th block and suppose
that $w_i:n_i \rightarrow 1$ is a morphism for all $i=1, \dots k$.
Then there exists a unique map denoted $(w_1, \dots ,w_k)$ such that
$$\xymatrix@R=3pc@C=3pc{k \ar[r]^{pr_i} & 1
\\ n_1 + n_2 + \cdots + n_k \ar@{.>}[u]^{(w_1, \dots ,w_k)}
 \ar[ur]_{w_i \circ \iota_i'} &}$$
commutes for all $i=1, \dots, k$ where $\iota_i':n_1 + n_2 + \cdots
+ n_k \rightarrow n_i$ is the unique morphism such that
$$\xymatrix@R=3pc@C=3pc{n_i \ar[r]^{pr_j}  & 1
\\ n_1 + n_2 + \cdots +n_k \ar@{.>}[u]^{\iota_i'}
\ar[ur]_{pr_{\iota_i j}} & } $$ commutes. One should keep in mind
that $n_1 + n_2 + \cdots + n_k$ is the product of $n_1, \dots, n_k$.
Note that $(w_1, \dots ,w_k)$ is not the same thing as the tuple
$w_1, \dots ,w_k$. The arrow $(w_1, \dots, w_k)$ is not the product
of $w_1, \dots, w_k$.

\begin{lem}
Let $T$ be a theory\index{theory}. Then $Mor_T(m,n)$ can be
identified with the set-theoretic product\index{product}
$\prod_{j=1}^n Mor_T(m,1)$ via the map which takes $w$ to the
tuple\index{tuple} with entries $pr_1 \circ w, \dots, pr_n \circ w$.
We identify $w$ with that tuple. In particular a theory is
determined up to isomorphism by the sets $T(0), T(1), T(2), \dots$.
\end{lem}
\begin{pf}
This follows directly from the definition of product in a category.
\end{pf}

\begin{examp} \label{Endexample1}
Let $X$ be a set. Then the {\it endomorphism theory
End(X)}\index{theory!endomorphism theory|textbf}\index{endomorphism
theory|textbf}\index{$End(X)$|(textbf} has objects $0,1,2, \dots$
and hom sets $Mor_{End(X)}(m,n)=$$Map(X^m, X^n)$. Composition is the
usual function composition. Here we readily see that $\{*\}$ is the
terminal object\index{terminal object} and that
$End(X)(0)=Mor_{End(X)}(0,1)$ can be identified with $X$.

Let $w \in End(X)(k)$ and $w_i \in End(X)(n_i)$ for $i=1, \dots,k$.
Then the composite function $\gamma(w,w_1, \dots, w_k):=w \circ
(w_1,\dots, w_k)$ is an element of $End(X)(n_1+ \dots + n_k)$. This
{\it composition}\index{composition} is associative. Let $1:=1_X \in
End(X)(1)$. Then apparently $w \circ (1, \dots, 1) =w$ and $1 \circ
w=w$, \ie the composition is also unital.

Let $\xymatrix@1{\{1, \dots,k\} \ar[r]^f & \{1, \dots,\ell \}
\ar[r]^g & \{1, \dots,m\} }$ be maps of sets. For a word $w \in
End(X)(k)$ we define a new word $w_f \in End(X)(\ell)$ by $w_f(x_1,
\dots, x_{\ell}) := w(x_{f1}, \dots, x_{fk})$ called the {\it
substituted word}\index{word!substituted word}\index{substituted
word}. Thus we have maps $$\xymatrix{End(X)(k) \ar[r]^{()_f} &
End(X)(\ell) \ar[r]^{()_g} & End(X)(m).}$$ If $e: \emptyset
\rightarrow \{1, \dots, k\}$ is the empty function and $x \in X =
End(X)(0)$, then the substituted word $x_e:X^k \rightarrow X$ is the
constant function $(x_1, \dots, x_k) \mapsto x$. There are no other
functions $\emptyset \rightarrow \{1, \dots, k\}$. We easily see
that $(w_f)_g=w_{g \circ f}$ and $w_{id_k}=w$ for the identity map
$id_k:\{1, \dots, k\} \rightarrow \{1, \dots, k\}$, \ie these {\it
substitution maps}\index{substitution maps} are functorial.

These substitution maps relate to the composition in two ways, which
we now describe. Let $f:\{1, \dots,k\} \rightarrow \{1, \dots,\ell
\}$, $w \in End(X)(k)$, and $w_i \in End(X)(n_i)$ for $i=1, \dots,
\ell$. Then $w_f \circ (w_1, \dots, w_{\ell}) = (w \circ (w_{f1},
\dots, w_{fk}))_{\bar{f}}$ where $$\xymatrix{\bar{f}:\{1,2, \dots,
n_{f1} + n_{f2} + \dots + n_{fk} \} \ar[r] & \{1,2, \dots, n_1 + n_2
+ \dots + n_{\ell} \}}$$ is the function obtained by parsing the
sequence $1,2, \dots, n_1+n_2 + \dots + n_{\ell}$ into consecutive
blocks $B_1, \dots, B_{\ell}$ of lengths $n_1, \dots, n_{\ell}$
respectively and then writing them in the order $B_{f1}, \dots,
B_{fk}$. For example, let $n_1=1,n_2=2,n_3=3,n_4=1, w \in T(3)$, and
$w_i \in T(n_i)$ for $i=1, \dots, 4$ and let $f:\{1,2,3\}
\rightarrow \{1, 2,3, 4\}$ be given by
$$\begin{pmatrix}
1 & 2 & 3 \\ 3 & 2 & 4
\end{pmatrix}.$$
Then $\bar{f}:\{1,2, \dots, 6\} \rightarrow \{1,2, \dots, 7\}$ is
given by
$$\begin{pmatrix}
1       & 2 & 3 & 4      & 5 & 6 \\
B_{f1}  &   &   & B_{f2} &   & B_{f3}
\end{pmatrix}
=
\begin{pmatrix}
1  & 2 & 3 & 4 & 5 & 6 \\
4  & 5 & 6 & 2 & 3 & 7
\end{pmatrix}.$$
We see that $$\aligned w_f \circ (w_1, w_2, w_3, w_4)(x_1, \dots,
x_7) &= w_f(w_1(x_1),w_2(x_2,x_3),w_3(x_4,x_5,x_6),w_4(x_7)) \\
&=w(w_3(x_4,x_5,x_6),w_2(x_2,x_3),w_4(x_7)) \\
&=w \circ (w_{f1},w_{f2},w_{f3})(x_4,x_5,x_6,x_2,x_3,x_7) \\
&=w \circ (w_{f1},w_{f2},w_{f3})(x_{\bar{f}1},x_{\bar{f}2}, \dots,
x_{\bar{f}6}). \endaligned$$ In other words we have $w_f \circ (w_1,
w_2, w_3) = (w \circ (w_{f1},w_{f2}, w_{f3}))_{\bar{f}}$. Note that
$\bar{f}$ depends not only on $f$, but also on the arity of the
words we are composing. The equality $w_f \circ (w_1, \dots,
w_{\ell})=(w \circ (w_{f1}, \dots, w_{fk}))_{\bar{f}}$ is the first
relationship between composition and the substitution maps $()_f$.

The second way the composition and the substitution maps relate
occurs in the following situation. If $w \in End(X)(k)$, $w_i \in
End(X)(n_i)$, and $g_i:\{1, \dots,
 n_i \} \rightarrow \{1, \dots, n_i'\}$ are functions for $i=1, \dots, k$,
then $w \circ ((w_1)_{g_1}, \dots, (w_k)_{g_k})= (w \circ (w_1,
\dots, w_k))_{g_1 + \dots + g_k}$ where $g_1 + \dots + g_k:\{1,2,
\dots, n_1+ \dots + n_k \} \rightarrow \{1,2, \dots, n_1 '+ \dots +
n_k' \}$ is the function obtained by placing $g_1, \dots, g_k$ next
to each other from left to right.
\end{examp}

\begin{examp} \label{Endexample2}
Let $X$ be a category. Then the {\it endomorphism theory
End(X)}\index{theory!endomorphism theory|textbf}\index{endomorphism
theory|textbf} has objects $0,1,2, \dots$ and it has hom sets
$Mor_{End(X)}(m,n)=$$Functors(X^m, X^n)$. We can proceed as in the
previous example and define substituted functors (substituted
words). Note that $End(X)$ can be made into a 2-category by taking
the 2-cells to be natural transformations, although we leave out the
2-cells for now. In most applications we will only be concerned with
the 1-category $End(X)$.
\end{examp}

\begin{examp}
Let $X$ be an object of a category with finite products. Then we
obtain a theory $End(X)$ with hom sets
$Mor_{End(X)}(m,n):=$$Mor(X^m, X^n)$\index{theory!endomorphism
theory|textbf}\index{endomorphism theory|textbf}.
\end{examp}
\index{$End(X)$|)}

We can abstract the essential properties of $End(X)$ in the previous
examples to get the following lemma for arbitrary theories.

\begin{lem} \label{tofunctortheory}
Let $T$ be a theory\index{theory|textbf}. Then for all $k,n_1,
\dots,n_k \in \{0,1, \dots \}$ there is a map $\gamma: T(k) \times
T(n_1) \times \dots \times T(n_k) \rightarrow T(n_1 + \cdots + n_k)$
called {\it composition}\index{composition|textbf} and for every
function $f:\{1, \dots, k\} \rightarrow \{1, \dots, \ell \}$ there
is a map $\xymatrix@1{T(k) \ar[r]^{()_f} & T(\ell)}$ called {\it
substitution}\index{substitution|textbf}. These maps have the
following properties.
\begin{enumerate}
\item
The $\gamma$'s are associative, \ie
$$\gamma(w,\gamma(w^1,w^1_1,\dots,w^1_{n_1}),\gamma(w^2,w^2_1,\dots,w^2_{n_2}),
\dots, \gamma(w^k,w^k_1,\dots,w^k_{n_k}))=$$ $$\gamma(\gamma(w, w^1,
\dots, w^k), w^1_1, \dots,w_{n_1}^1,w^2_1, \dots, w^2_{n_2}, \dots,
w^k_1, \dots,w_{n_k}^k).$$
\item
The $\gamma$'s are unital, \ie there exists an element $1 \in T(1)$
called the {\it unit}\index{unit} such that $$\gamma(w,1,
\dots,1)=w=\gamma(1,w)$$ for all $w \in T(k)$. Moreover, such an
element is unique.
\item
The $\gamma$'s are equivariant\index{equivariant} in the sense that
$$\gamma(w_f,w_1, \dots, w_{\ell})=\gamma(w,w_{f1}, \dots,
w_{fk})_{\bar{f}}$$ for all $f:\{1, \dots, k\} \rightarrow \{1,
\dots, \ell \}$ where $\bar{f}:\{1,2, \dots, n_{f1} + n_{f2} + \dots
+ n_{fk} \} \rightarrow \{1,2, \dots, n_1 + n_2 + \dots + n_{\ell}
\}$ is the function that moves entire blocks according to $f$ as
mentioned in the example above. Here $\bar{f}$ depends also on the
particular $\gamma$.
\item
The $\gamma$'s are equivariant\index{equivariant} in the sense that
$$\gamma(w,(w_1)_{g_1}, \dots, (w_k)_{g_k})= \gamma(w,w_1, \dots,
w_k)_{g_1 + \dots + g_k}$$ for all functions $g_i:\{1, \dots, n_i \}
\rightarrow \{1, \dots, n_i'\}$ where $g_1 + \dots + g_k:\{1,2,
\dots, n_1+ \dots + n_k \} \rightarrow \{1,2, \dots, n_1 '+ \dots +
n_k' \}$ is the function obtained by placing $g_1, \dots, g_k$ next
to each other from left to right.
\item
The substitution is functorial\index{functorial}, \ie for functions
\newline $\xymatrix@1{\{1, \dots,k\} \ar[r]^f & \{1, \dots,\ell \}
\ar[r]^g & \{1, \dots,m\} }$ the composition
$$\xymatrix@1@C=3pc{T(k)\ar[r]^{()_f} & T(\ell) \ar[r]^{()_g} & T(m)}$$
is the same as $$\xymatrix@1@C=3pc{T(k) \ar[r]^{()_{g \circ f}} &
T(m)}$$ and for the identity function $id_k:\{1, \dots,k \}
\rightarrow \{1, \dots,k \}$ the map
$$\xymatrix@1@C=3pc{T(k) \ar[r]^{()_{id_k}} & T(k)}$$ is equal to the
identity for all $k \geq 0$.
\end{enumerate}
\end{lem}
\begin{pf}
First we define the substitution. Let $f:\{1, \dots, k\} \rightarrow
\{1, \dots, \ell \}$ be a function. Then there exists a unique
morphism $f'$ such that the diagram
$$\xymatrix@R=3pc@C=3pc{k \ar[r]^{pr_i} & 1 \\ \ell \ar@{.>}[u]^{f'} \ar[ur]_{pr_{fi}}
&}$$ commutes for all $i=1, \dots,k$. For $w \in T(k)$ define
$w_f:=w \circ f'$. Thus the map $\xymatrix@1{T(k) \ar[r]^{()_f} &
T(\ell)}$ is defined by precomposition with $f'$.

Next we define the composition $\gamma:T(k) \times T(n_1) \times
\dots \times T(n_k) \rightarrow T(n_1+n_2+ \dots + n_k)$. Let $w \in
T(k),w_i \in T(n_i)$ for $i=1, \dots, k$.  Define $\gamma(w, w_1,
\dots, w_k):= w \circ (w_1, \dots, w_k)$ where the composition
$\circ$ is the composition of the category $T$ and $(w_1, \dots,
w_k)$ is the unique morphism such that
$$\xymatrix@R=3pc@C=3pc{k \ar[r]^{pr_i} & 1
\\ n_1 + n_2 + \cdots + n_k \ar@{.>}[u]^{(w_1, \dots ,w_k)}
 \ar[ur]_{(w_i)_{\iota_i}} &}$$ commutes as defined above.

\begin{enumerate}
\item
We claim that $\gamma$ is associative. \newline
$\gamma(w,\gamma(w^1,w^1_1,\dots,w^1_{n_1}),\gamma(w^2,w^2_1,\dots,w^2_{n_2}),
\dots, \gamma(w^k,w^k_1,\dots,w^k_{n_k}))=$ \newline $=w \circ (w^1
\circ (w^1_1,\dots,w^1_{n_1}), \dots, w^k \circ
(w^k_1,\dots,w^k_{n_k}))$ \newline $=w \circ ((w^1, \dots, w^k)
\circ ((w^1_1,\dots,w^1_{n_1}), \dots, (w^k_1,\dots,w^k_{n_k})))$
\newline $=(w \circ (w^1, \dots, w^k)) \circ
(w^1_1,\dots,w^1_{n_1}, \dots, w^k_1,\dots,w^k_{n_k})$ \newline
$=\gamma(\gamma(w, w^1, \dots, w^k), w^1_1, \dots,w_{n_1}^1,w^2_1,
\dots, w^2_{n_2}, \dots, w^k_1, \dots,w_{n_k}^k)$ \newline The
second to last equality follows by associativity of composition in
the category $T$ and by properties of products.
\item
We claim that $\gamma$ is unital. Let $1:1 \rightarrow 1$ be the
projection morphism of the object 1 in the category $T$, which is
not necessarily the identity morphism of the object 1. Then $(1,
\dots, 1):k \rightarrow k$ is the identity morphism of the object
$k$ because $1_{\iota_i}=1 \circ \iota_i' =1 \circ (pr_1^{-1} \circ
pr_i)=pr_1 \circ (pr_1^{-1} \circ pr_i) = pr_i$ in the diagram
$$\xymatrix@R=3pc@C=3pc{k \ar[r]^{pr_i} & 1 \\ \ar@{.>}[u]^{(1, \dots, 1)} k
\ar[ur]_{1_{\iota_i}} & }$$ for all $i=1, \dots, k$. Here
$\iota_i:\{1\} \rightarrow \{1, \dots, k\}$ is defined by
$\iota_i(1)=i$. Thus $\gamma(w, 1, \dots, 1)= w \circ (1, \dots,
1)=w \circ 1_k = w$.

To show $\gamma(1,w)=w$ we consider the diagram
$$\xymatrix@R=3pc@C=3pc{1 \ar[r]^{pr_1}  & 1 \\ n
\ar@{.>}[u]^{(w)} \ar[ur]_{w_{\iota_1}} & }$$ where $\iota_1:\{1,
\dots, n\} \rightarrow \{1, \dots, n\}$ is the identity. Then
$w_{\iota_1} =w$ and $(w)=pr_1^{-1} \circ w$. Thus $\gamma(1, w) = 1
\circ (w) = pr_1 \circ (pr_1^{-1} \circ w)=w$.

The uniqueness follows from $1=\gamma(1, 1')=1'$.

\item
Let $f:\{1, \dots, k\} \rightarrow \{1, \dots, \ell \}$ be a
function and $w_i \in T(i)$ for $i=1, \dots, \ell$. Using the
definitions of $\bar{f}:\{1, \dots, n_{f1} + \dots + n_{fk} \}
\rightarrow \{1, \dots, n_1 + \dots + n_{\ell} \}$ and $\iota_i$
from above we see that the following two diagrams
$$\xymatrix@R=3pc@C=3pc{k \ar[r]^{pr_i} & 1
\\ n_{f1} + \dots + n_{fk} \ar[u]^{(w_{f1}, \dots, w_{fk})}
\ar[ur]_{(w_{fi})_{\iota_i}} &
\\ n_1 + \dots + n_{\ell}  \ar@{=}[r] \ar[u]^{\bar{f}'} &
n_1 + \dots + n_{\ell} \ar[uu]_{(w_{fi})_{\iota_{fi}}}}$$
$$\xymatrix@R=3pc@C=3pc{k \ar[r]^{pr_i} & 1
\\ \ell \ar[u]^{f'} \ar[ur]_{pr_{fi}} &
\\ n_1 + \dots + n_{\ell}  \ar@{=}[r] \ar[u]^{(w_1, \dots, w_{\ell})} &
n_1 + \dots + n_{\ell} \ar[uu]_{(w_{fi})_{\iota_{fi}}}}$$ commute
for all $i=1, \dots, k$. Hence by the universal property of the
product $k$ we have $f' \circ (w_1, \dots, w_{\ell})=(w_{f1}, \dots,
w_{fk}) \circ \bar{f}'$. Using this we see that
$$\aligned \gamma(w_f, w_1, \dots, w_{\ell}) &=w \circ f' \circ (w_1, \dots,
w_{\ell}) \\
&=w \circ (w_{f1}, \dots, w_{fk}) \circ \bar{f}' \\
&=\gamma(w,w_{f1}, \dots, w_{fk})_{\bar{f}}. \endaligned$$
\item
Let $g_i:\{1, \dots, n_i \} \rightarrow \{1, \dots, n_i'\}$ be
functions for $i=1, \dots, k$. Then $$\aligned \gamma(w,(w_1)_{g_1},
\dots, (w_k)_{g_k})
&=w \circ (w_1 \circ g_1', \dots, w_k \circ g_k') \\
&=w \circ (w_1, \dots, w_k) \circ (g_1', \dots, g_k') \\
&=w \circ (w_1, \dots, w_k) \circ (g_1+ \dots + g_k)' \\
&=\gamma(w,w_1, \dots, w_k)_{g_1 + \dots + g_k}.
\endaligned$$
\item
Let $\xymatrix@1{\{1, \dots,k\} \ar[r]^f & \{1, \dots,\ell \}
\ar[r]^g & \{1, \dots,m\} }$ be functions. Then $f'$ and $g'$ make
the two small subdiagrams in
$$\xymatrix@R=3pc@C=3pc{k \ar[r]^{pr_i} & 1
\\ \ell \ar[u]^{f'} \ar[ur]_{pr_{fi}} &
\\ m \ar[u]^{g'} \ar@{=}[r] & m \ar[uu]_{pr_{gfi}}}$$
commute for all $i=1, \dots, k$. Thus the outer diagram commutes and
$(g \circ f)'=f' \circ g'$ by the universal property of the product.
We conclude $(w_f)_g=w \circ f' \circ g' = w \circ (g \circ f)' =
w_{g \circ f}$. The identity $1_k:k \rightarrow k$ makes
$$\xymatrix@R=3pc@C=3pc{k \ar[r]^{pr_i} & 1
\\ k \ar[u]^{1_k} \ar[ur]_{pr_{id_k(i)}} & }$$
commute for all $i=1, \dots, k$ where $id_k:\{1, \dots, k\}
\rightarrow \{1, \dots, k\}$ is the identity function. Hence
$(id_k)'=1_k$ and $w_{id_k}=w \circ (id_k)' = w \circ 1_k = w$ for
all $w \in T(k)$.
\end{enumerate}
We have verified all of the axioms.
\end{pf}

There is another description of a theory which can be formulated by
using the category $\Gamma$.

\begin{defn} \label{definitionofGamma}
Let $\Gamma$\index{$\Gamma$} be the category with objects
$\emptyset=0, 1,2, \dots$ where $k=\{1, \dots, k\}$. The morphisms
$k \rightarrow \ell$ are just maps of sets. In particular $0$ is the
initial object\index{initial object} since the only map $\emptyset
\rightarrow k$ is the empty function. There are no maps $k
\rightarrow \emptyset$ for $k \geq 1$. The object $1$ is the
terminal object\index{terminal object}. Let $+:\Gamma \times \Gamma
\rightarrow \Gamma$ denote the usual functor obtained by adding the
sets and placing maps side by side.
\end{defn}

\begin{rmk}
Let $T$ be a theory\index{theory|textbf}. Then by the previous lemma
$T$ defines a functor from $\Gamma$ to $Sets$ by $k \mapsto T(k)$
and $f \mapsto ()_f$. Moreover, this functor comes with maps
$\gamma: T(k) \times T(n_1) \times \dots \times T(n_k) \rightarrow
T(n_1 + \dots + n_k)$ which satisfy 1. through 5. The compositions
$\gamma$, unit $1$, and substitution are sometimes called the {\it
operations of theories}\index{operations of theories}. The relations
in 1. through 5. are sometimes called the {\it relations of
theories}\index{relations of theories}.
\end{rmk}

\begin{lem} \label{functortheory}
Let $T$ be a functor from $\Gamma$ to $Sets$ equipped with maps
$\gamma: T(k) \times T(n_1) \times \dots \times T(n_k) \rightarrow
T(n_1 + \dots + n_k)$ and an element $1 \in T(1)$ which satisfy (1)
through (5) where $T(f)=:()_f$ for functions $f:k \rightarrow \ell$.
Then $T$ determines a theory\index{theory|textbf} with
$Mor(n,1)=T(n)$ for all $n \geq 0$.
\end{lem}
\begin{pf}
Define the underlying category of the theory to formally have
objects $0,1,2 \dots$ and morphisms $Mor(m,n):=\prod_{i=1}^n
Mor(m,1)$. In particular $Mor(m,0)$ only has one element. We denote
a tuple of words $w_1, \dots, w_n \in Mor(m,1)$ by $\prod_{i=1}^n
w_i$. For $k,\ell \geq 0$ let $\iota_{\ell,k}:\{1, \dots, \ell k \}
\rightarrow \{1, \dots, k \}$ be the function such that
$\iota_{\ell,k}(i + jk)=i$ for $i=1, \dots, k$, in other words
$\iota_{\ell, k}$ wraps the domain around the codomain $\ell$ times.
Now define the composition of $\prod_{i=1}^{\ell} w_i \in
Mor(k,\ell)$ with $\prod_{i=1}^m v_i \in Mor(\ell,m)$ to be
$\prod_{i=1}^m v_i \circ \prod_{i=1}^{\ell} w_i := \prod_{i=1}^m
\gamma(v_i,w_1, \dots, w_{\ell})_{\iota_{\ell,k}}$. This composition
is associative because $\gamma$ is associative and equivariant.

Let $f_i:\{1\} \rightarrow \{1, \dots, n \}$ be the map $f_i(1)=i$.
Define $pr_i :=1_{f_i} \in T(n)$ where $1 \in T(1)$ is the
distinguished element whose existence we assumed. This notation is
slightly imprecise because we have different sequences $pr_1, \dots,
pr_n$ for different $n \geq 0$. From the context it will always be
clear which sequence of morphisms is meant. We claim that
$\prod_{i=1}^n pr_i \in Mor(n,n)$ is the identity on the object $n$.
Let $\prod_{i=1}^m w_i \in Mor(n,m)$. Then
$$\aligned \prod_{i=1}^m w_i \circ \prod_{i=1}^n pr_i
&= \prod_{i=1}^m \gamma(w_i, pr_1, \dots, pr_n)_{\iota_{n,n}} \\
&=\prod_{i=1}^m \gamma(w_i, 1_{f_1}, \dots, 1_{f_n})_{\iota_{n,n}} \\
&=\prod_{i=1}^m (\gamma(w_i, 1, \dots,1)_{f_1 + \dots + f_n})_{\iota_{n,n}} \text{ by equivariance} \\
&=\prod_{i=1}^m \gamma(w_i, 1, \dots, 1)_{\iota_{n,n}
\circ (f_1 + \dots + f_n)} \text{ by functoriality of $T$} \\
&=\prod_{i=1}^m w_i \text{ since $\gamma$ is unital, $\iota_{n,n}
\circ (f_1 + \dots + f_n)=id_n$,} \\
& \text{\hspace{.63in} and functoriality of $T$.}
\endaligned$$

Now for the other side let $\prod_{i=1}^n w_i \in Mor(m,n)$. Then
$$\aligned \prod_{i=1}^n pr_i \circ \prod_{i=1}^n w_i
&=\prod_{i=1}^n \gamma(pr_i,w_1, \dots, w_n)_{\iota_{n,m}} \text{
by definition} \\
&=\prod_{i=1}^n \gamma(1_{f_i},w_1, \dots, w_n)_{\iota_{n,m}}
\text{ by definition} \\
&=\prod_{i=1}^n (\gamma(1,w_i)_{\bar{f_i}})_{\iota_{n,m}} \text{
by equivariance} \\
&=\prod_{i=1}^n (w_i)_{\iota_{n,m} \circ \bar{f_i}} \text{ by
unitality of $\gamma$ and functoriality of $T$} \\
&=\prod_{i=1}^n w_i \text{ since $\iota_{n,m} \circ \bar{f_i} =
id_m$.} \endaligned $$ This can be seen by observing that
$\bar{f_i}:\{1, \dots,m\} \rightarrow \{1, \dots, nm \}$ has the
form
$$\begin{pmatrix}
1 & 2 & \dots & m
\\ (i-1)m+1 & (i-1)m+2 & \dots & (i-1)m + m
\end{pmatrix}$$
and by using the definition of $\iota_{n,m}$. Thus $\prod_{i=1}^n
pr_i \in Mor(n,n)$ is the identity on the object $n$.

Thus far we have shown that we have a category with objects $0,1,2,
\dots$ and morphisms $Mor(m,n)$. We claim that $n$ is the product of
$n$ copies of $1$ in this category with projections $pr_1, \dots,
pr_n:n \rightarrow 1$ introduced above. First note for
$\prod_{i=1}^n w_i \in Mor(m,n)$ we have $$\aligned pr_i \circ
\prod_{i=1}^n w_i
&= \gamma(pr_i, w_1, \dots, w_n)_{\iota_{n,m}} \\
&=\gamma(1_{f_i}, w_1, \dots, w_n)_{\iota_{n,m}} \text{ by definition} \\
&=\gamma(1,w_i)_{\iota_{n,m} \circ \bar{f_i}} \text{ by
equivariance and functoriality} \\
&=w_i \text{ since $\iota_{n,m} \circ \bar{f_i} = id_m$ and by
functoriality.} \endaligned $$ Now suppose we are given morphisms
$w_1, \dots, w_n \in Mor(m,1)$. Then
$$\xymatrix@R=3pc@C=3pc{n \ar[r]^{pr_i} & 1
\\ m \ar[u]^{\prod_{j=1}^n w_j} \ar[ur]_{w_i} & }$$
commutes for all $i=1, \dots, n$ by the remark just made. If
$\prod_{i=1}^n v_i \in Mor(m,n)$ is another morphism such that
$$\xymatrix@R=3pc@C=3pc{n \ar[r]^{pr_i} & 1
\\ m \ar[u]^{\prod_{j=1}^n v_j} \ar[ur]_{w_i} & }$$
commutes for all $i=1, \dots, n$, then by the remark $v_i=pr_i \circ
\prod_{j=1}^n v_j = w_i$ and hence $\prod_{j=1}^n v_j=\prod_{j=1}^n
w_j$ and the factorizing map is unique. Hence $n$ is the product of
$n$ copies of $1$.

We conclude that the functor $T$ with the maps $\gamma$ satisfying
the axioms (1) through (5) determines a theory with the indicated
hom sets.
\end{pf}

\begin{thm} \label{equivalenceoftheories}
A theory\index{theory|textbf} $T$ is determined by either of the
following equivalent collections of data:
\begin{enumerate}
\item
A category $T$ with objects $0,1,2, \dots$ such that $n$ is the
categorical product of $1$ with itself $n$ times and each $n$ is
equipped with a choice of projections.
\item
A functor $T:\Gamma \rightarrow Sets$ equipped with maps
$\gamma:T(k) \times T(n_1) \times \cdots \times T(n_k) \rightarrow
T(n_1 + \cdots + n_k)$ and a unit $1 \in T(1)$ which satisfy (1)
through (5) of Lemma \ref{tofunctortheory}.
\end{enumerate}
\end{thm}
\begin{pf}
In each description $Mor_T(n,1)$ is the same. By the universality of
products this determines the rest of the theory. The two processes
of Lemmas \ref{tofunctortheory} and \ref{functortheory} are
``inverse'' to one another by further inspection, provided we
identify $Mor_T(m,n)$ with $\prod_{i=1}^n T(m)$.
\end{pf}

\begin{defn}
Let $S$ and $T$ be theories. In the categorical description of $S$
and $T$ a {\it morphism of theories}\index{morphism!morphism of
theories|textbf} $\Phi:S \rightarrow T$ is a functor from the
category $S$ to the category $T$ such that $\Phi(n_S)=n_T$ and
$\Phi(pr_i)=pr_i$ for all projections.
\end{defn}

One easily sees that the theories form a category and we have a
suitable forgetful functor.

\begin{thm}
The forgetful functor\index{forgetful functor} from the category of
theories to $\prod_{n \geq 0}Sets$ given by $T \mapsto
(T(0),T(1),\dots)$ admits a left adjoint called the {\it free theory
functor}\index{free theory|textbf}\index{free theory!free theory
functor|textbf}\index{theory|free theory
functor|textbf}\index{functor!free theory functor|textbf}.
\end{thm}
\begin{pf}
On page \pageref{freetheorymonad} we will construct the free theory
on the sequence of sets
\newline $(T(0),T(1),\dots)$.
\end{pf}

To make later proofs easier, we need the following lemma.

\begin{lem} \label{theorymorphisms}
Let $\Phi:S \rightarrow T$ be a morphism\index{morphism!morphism of
theories|textbf} of theories.
\begin{enumerate}
\item
Let $f: \{1, \dots, k\} \rightarrow \{1, \dots, \ell\}$ be a
function. As usual, $f': \ell \rightarrow k$ denotes the unique
morphism in any theory such that
$$\xymatrix@R=3pc@C=3pc{k \ar[r]^{pr_i} & 1 \\ \ell \ar@{.>}[u]^{f'} \ar[ur]_{pr_{fi}}
}$$ commutes. Then $\Phi(f')=f'$.
\item
Let $f: \{1, \dots, k\} \rightarrow \{1, \dots, \ell\}$ be a
function and $w \in Mor_S(k,1)$. Then $\Phi(w_f)=\Phi(w)_f$.
\item
Let $w_1, \dots, w_n \in Mor_S(m,1)$. Then $\Phi(\prod_{j=1}^n w_j)=
\prod_{j=1}^n \Phi(w_j)$.
\item
Let $w_i \in Mor_S(n_i,1)$ for $i=1, \dots, k$. Then we also have
$\Phi(w_1, \dots, w_k)=(\Phi(w_1), \dots, \Phi(w_k))$.
\end{enumerate}
\end{lem}
\begin{pf}
\begin{enumerate}
\item
The diagram
$$\xymatrix@R=3pc@C=3pc{k \ar[r]^{pr_i} & 1
\\ \ell \ar[u]^{\Phi(f')} \ar[ur]_{pr_{fi}} & }$$
commutes for all $i=1, \dots, k$ by the properties of $\Phi$. Then
$\Phi(f')=f'$ by the universal property of the product.
\item
This follows from (1) and the definition $w_f=w \circ f'$.
\item
The properties of $\Phi$ imply that the diagram
$$\xymatrix@R=3pc@C=3pc{n \ar[r]^{pr_i} & 1
\\ m \ar[u]^{\Phi(\prod_{j=1}^n w_j)} \ar[ur]_{\Phi(w_i)}
}$$ commutes for all $i=1, \dots, n$. Then $\Phi(\prod_{j=1}^n w_j)=
\prod_{j=1}^n \Phi(w_j)$ by the universal property of the product.
\item
By (2) we have $\Phi((w_i)_{\iota_i})=\Phi(w_i)_{\iota_i}$. Hence,
the properties of $\Phi$ imply that the diagram
$$\xymatrix@R=3pc@C=3pc{k \ar[r]^{pr_i} & 1
\\ n_1 + n_2 + \cdots + n_k \ar[u]^{\Phi(w_1, \dots, w_k)}
\ar[ur]_{\Phi(w_i)_{\iota_i}} & }$$ commutes for all $i=1, \dots,
k$. Then $\Phi(w_1, \dots, w_k)=(\Phi(w_1), \dots, \Phi(w_k))$ by
the universal property of the product.
\end{enumerate}
\end{pf}

Just as a theory has a categorical description and a functorial
description, a morphism of theories\index{morphism!morphism of
theories|textbf} also has a second description. We work towards the
second description in the following two lemmas.

\begin{lem} \label{naturalmorphismoftheories}
Let $\Phi:S \rightarrow T$ be a morphism of theories, \ie a functor
such that $\Phi(n_S)=n_T$ and $\Phi(pr_i)=pr_i$ for all projections.
Then $\Phi$ determines a natural transformation $S \Rightarrow T$
also denoted by $\Phi$ such that
$$\xymatrix@C=9pc@R=4pc{S(k) \times S(n_1) \times \cdots \times S(n_k)
\ar[r]^{\Phi_k \times \Phi_{n_1} \times \cdots \times \Phi_{n_k}}
\ar[d]_{\gamma^S} & T(k) \times T(n_1) \times \cdots \times T(n_k)
\ar[d]^{\gamma^T}
\\ S(n_1 + \cdots + n_k) \ar[r]_{\Phi_{n_1 + \cdots + n_k}}
& T(n_1 + \cdots + n_k) }$$ commutes and $\Phi_1(1_S)=1_T$, where
$S,T: \Gamma \rightarrow Sets$ are the functors in the functorial
description of the theories $S$ and $T$.
\end{lem}
\begin{pf}
Let $\Phi_m:Mor_S(m,1) \rightarrow Mor_T(m,1)$ denote the map
obtained from the functor $\Phi$, \ie $\Phi_m(w):=\Phi(w)$ for $w
\in S(m)$. Then for $f:m \rightarrow n$ in $\Gamma$ and $w \in
S(m)$, we have $\Phi(w_f)=\Phi(w)_f$ by Lemma \ref{theorymorphisms}.
Hence $$\xymatrix@C=3pc@R=3pc{S(m) \ar[r]^{\Phi_m} \ar[d]_{S(f)} &
T(m) \ar[d]^{T(f)}
\\ S(n) \ar[r]_{\Phi_n} & T(n)}$$
commutes and $m \mapsto \Phi_m$ is natural.

Let $w \in S(k)$ and $w_i \in S(n_i)$ for $i=1, \dots, k$. Then
$$\aligned
\Phi_{n_1 + \cdots + n_k}(\gamma^S(w,w_1, \dots, w_k)) &=
\Phi(w \circ (w_1, \dots, w_k)) \\
&= \Phi(w) \circ (\Phi(w_1), \dots, \Phi(w_k))  \\
&= \gamma^T(\Phi(w), \Phi(w_1), \dots, \Phi(w_k)) \\
&= \gamma^T(\Phi_k(w), \Phi_{n_1}(w_1), \dots, \Phi_{n_k}(w_k) ).
\endaligned
$$
Hence the natural transformation $m \mapsto \Phi_m$ preserves the
$\gamma$'s.

Let $1_S \in S(1)$ and $1_T \in T(1)$ be the units in the respective
theories. Then $\Phi_1(1_S)=1_T$ because the functor $\Phi$
preserves projections.

Thus $\Phi:S \Rightarrow T$ is a natural transformation which
preserves the compositions and the units.
\end{pf}

\begin{lem}
Let $S, T: \Gamma \rightarrow Sets$ be theories. Let $\Phi:S
\Rightarrow T$ be a natural transformation preserving the $\gamma$'s
and their units as in Lemma
\ref{naturalmorphismoftheories}.
Then $\Phi$ determines a functor $S \rightarrow T$ also denoted
$\Phi$, where $S$ and $T$ are the categories in the categorical
description of the theories $S,T: \Gamma \rightarrow Sets$.
Moreover, the functor $\Phi:S \rightarrow T$ satisfies
$\Phi(n_S)=n_T$ and $\Phi(pr_i)=pr_i$ for all projections.
\end{lem}
\begin{pf}
We define $\Phi(n_S)=n_T$ for all $n_S \in Obj \hspace{1mm}  S$ and
$\Phi(\prod_{j=1}^{\ell} w_j):=\prod_{j=1}^{\ell} \Phi_k(w_j)$ for
all $\prod_{j=1}^{\ell}w_j \in Mor_S(k, \ell)$. Then for
$\prod_{i=1}^m v_i \in Mor_S(\ell,m)$ we have

$$\aligned \Phi(\prod_{i=1}^m v_i \circ
\prod_{j=1}^{\ell}w_j) &= \Phi(\prod_{i=1}^m \gamma(v_i, w_1, \dots,
w_{\ell})_{\iota_{\ell,k}}) \text{ from Lemma }
\ref{functortheory} \\
 &=\prod_{i=1}^m \gamma(\Phi_{\ell}(v_i),\Phi_k(w_1), \dots,
\Phi_k(w_{\ell}))_{\iota_{\ell,k}} \\
&=\prod_{i=1}^m \Phi_{\ell}(v_j) \circ \prod_{j=1}^{\ell} \Phi_k(w_j) \\
&= \Phi(\prod_{i=1}^m v_i) \circ \Phi(\prod_{j=1}^{\ell} w_j).
\endaligned
$$


Hence $\Phi$ preserves compositions.

We claim that $\Phi$ preserves projections. Let $f_i:\{1\}
\rightarrow \{1, \dots, n\}$ be the map $f_i(1)=i$. Then
$(1_S)_{f_i}=pr_i$ and
$$\aligned
\Phi(pr_i) &= \Phi_n((1_S)_{f_1}) \\
&=\Phi_n(1_S)_{f_i} \text{ by naturality} \\
&=(1_T)_{f_i} \text{ since $\Phi$ preserves the unit} \\
&= pr_i .
\endaligned
$$
Hence $\Phi$ preserves projections.

We claim that $\Phi$ preserves identities. Recall that
$\prod_{j=1}^n pr_j: n \rightarrow n$ is the identity on the object
$n$ of the category $S$. Then
$$\aligned
 \Phi(\prod_{j=1}^n pr_j ) &= \prod_{j=1}^n \Phi(pr_j) \text{ by definition} \\
&= \prod_{j=1}^n pr_j \text{ because $\Phi$ preserves projections.}
\endaligned$$
Thus $\Phi$ preserves identities and is a functor $S \rightarrow T$.
\end{pf}

Combining these two lemmas gives us the two descriptions of a
morphism of theories in the following theorem.

\begin{thm} \label{equivalenceoftheorymorphisms}
Let $S$ and $T$ be theories. Then a morphism\index{morphism!morphism
of theories|textbf} $S \rightarrow T$ of theories is given by either
of the following equivalent collections of data:
\begin{enumerate}
\item
A functor $\Phi:S \rightarrow T$ such that $\Phi(n_S)=n_T$ for all
$n_S \in Obj \hspace{1mm}  S$ and $\Phi(pr_i)= pr_i$ for all
projections.
\item
A natural transformation $\Phi:S \Rightarrow T$ of the functors
$S,T: \Gamma \rightarrow Sets$ which preserves the $\gamma$'s and
the units.
\end{enumerate}
\end{thm}
\begin{pf}
The processes of the previous two lemmas are ``inverse'' to each
other by inspection.
\end{pf}

\begin{thm} \label{theorydescriptions}
The category of theories\index{theory|textbf} with objects and
morphisms as in (1) of Theorems \ref{equivalenceoftheories} and
\ref{equivalenceoftheorymorphisms} is equivalent to the category
with objects and morphisms as in (2) of Theorems
\ref{equivalenceoftheories} and \ref{equivalenceoftheorymorphisms}.
\end{thm}
\begin{pf}
This relies on the bijection $Mor_T(m,n) \cong \prod_{j=1}^n
Mor_T(m,1)$.
\end{pf}

The concept of an {\it algebra}\index{algebra} is closely related to
the concept of theories. Roughly speaking, an algebra over a theory
is a category together with a rule that assigns an $n$-ary operation
on $X$ to every word of the theory of arity $n$ in such way that
compositions, substitutions, and identity 1 are preserved.

\begin{defn}
Let $X$ be a category and $T$ a theory. Then $X$ is a {\it
$T$-algebra}\index{algebra!$T$-algebra|textbf} if it is equipped
with a morphism of theories $T \rightarrow End(X)$, where $End(X)$
is the theory in Example \ref{Endexample2}. We also say $X$ is an
{\it algebra over the theory $T$}.\index{algebra!algebra over a
theory|textbf}
\end{defn}

Notice that if $X$ is a set viewed as a discrete category, this is
the usual definition of an algebra over a theory. Note also that we
have two versions of $T$-algebra, one is given by the
categorical\index{algebra!categorical $T$-algebra} description of
theories and the other by the functorial\index{algebra!functorial
$T$-algebra} description. A familiar example of an algebra is a
group, since a group is an algebra over the theory of groups as
follows.

\begin{examp}
Let $T$ be the theory\index{group}\index{group!theory of
groups|textbf}\index{theory!theory of groups|textbf} of groups, \ie
there are morphisms $e \in T(0), \nu \in T(1)$, and $\mu \in T(2)$
which satisfy the usual group axioms. The theory $T$ is the smallest
theory containing such $e,\nu,\mu$. A set $X$ is a group if there is
a morphism of theories $T \rightarrow End(X)$. This means we have
realizations of  $e,\nu,$ and $\mu$ on $X$.
\end{examp}

\begin{defn}
Let $X$ and $Y$ be $T$-algebras. Then a functor $H:X \rightarrow Y$
is a {\it morphism of $T$-algebras}\index{morphism!morphism of
$T$-algebras|textbf} in the categorical description if
$$\xymatrix@C=5pc@R=3pc{Mor_T(m,n) \ar[r] \ar[d] &
Mor_{End(X)}(m,n) \ar[d]^{H^{\times n} \circ }
\\ Mor_{End(Y)}(m,n) \ar[r]_-{\circ H^{\times m} } & Functors(X^m,Y^n) }$$
commutes for all $m$ and $n$. A functor $H:X \rightarrow Y$ is a
{\it morphism of $T$-algebras} in the functorial description if
$$\xymatrix@C=5pc@R=3pc{T(m) \ar[r] \ar[d] & End(X)(m) \ar[d]^{H \circ}
\\ End(Y)(m) \ar[r]_-{H^{\times m}} &
Functors(X \times \dots \times X, Y)}$$ commutes for all $m$.
\end{defn}

\begin{examp}
Let $T$ be the theory of groups and let $X$ and $Y$ be groups. Then
a set map $H:X \rightarrow Y$ is a morphism of $T$-algebras if and
only if it is a group homomorphism\index{group homomorphism}.
\end{examp}

\begin{thm} \label{star}
The category of categorical $T$-algebras is equivalent to the
category of functorial
$T$-algebras\index{algebra!$T$-algebra|textbf}\index{algebra!categorical
$T$-algebra}\index{algebra!functorial $T$-algebra}.
\end{thm}
\begin{pf}
The proof is similar to Theorem \ref{theorydescriptions}.
\end{pf}

Let $T$ be any theory. It is well known that $T$-algebras are
algebras for a monad\index{monad} $C$, which depends on $T$. See for
example \cite{maclane1} or \cite{power3}. We now present a version
of this in preparation for the 2-monad\index{2-monad} whose strict
algebras\index{algebra!strict algebra} are pseudo $T$-algebras. Let
$Cat_0$\index{$Cat_0$} denote the 1-category of small categories. We
define a functor $C:Cat_0 \rightarrow Cat_0$ as follows. For a small
category $X$, set  \label{2monad}
$$Obj \hspace{1mm}  CX:=\frac{(\mathop{\bigcup}_{n \geq 0} (T(n) \times Obj \hspace{1mm}
X^n))}{\Gamma}$$ where the quotient by $\Gamma$ means to mod out by
the smallest congruence\index{congruence} satisfying $(w_f, x_1,
\dots, x_n) \sim (w, x_{f1}, \dots, x_{fm})$ for all $m \in
\mathbb{N}_0$, $w \in T(m)$, and maps $f:m=\{1, \dots,m\}
\rightarrow \{1, \dots, n\}=n$. To define the morphisms of $CX$ we
note that $\bigcup_{n \geq 0} (T(n) \times X^n)$ is a category if we
interpret $T(n)$ as a discrete category for each $n$. Consider the
directed graph\index{directed graph} with objects $Obj \hspace{1mm}
CX$ and arrows from $[a]$ to $[b]$ given by the union
$$\bigcup Mor_{\cup_{n \geq 0} (T(n) \times X^n)}(a',b')$$ over all
$a' \sim a$ and $b' \sim b$. Next we take the free category on this
directed graph and mod out by the relations of $\bigcup_{n \geq 0}
(T(n) \times X^n)$ and the relations
$$(i_{w_f}, g_1, \dots, g_n)=(i_w,g_{f1}, \dots, g_{fm}).$$
This quotient category is $CX$. We define $C$ on functors $X
\rightarrow Y$ analogously. Then $C:Cat_0 \rightarrow Cat_0$ is a
functor because each step in the construction is functorial.

Next we define a natural transformation $\eta:1_{Cat_0} \Rightarrow
C$ by $\eta_X(x):=[1,x]$ for $x \in Obj \hspace{1mm} X$ and
$\eta_X(g):=[i_1,g]$ for a morphism $g$ in $X$. We also define a
natural transformation $\mu:C^2 \Rightarrow C$ by
$$\mu_X([w,[v^1,x^1_1, \dots, x^1_{j_1}],[v^2,x^2_1, \dots,
x^2_{j_2}], \dots, [v^k,x^k_1, \dots, x^k_{j_k}]]):=$$ $$
[\gamma(w,v^1,v^2, \dots, v^k),x^1_1, \dots, x^k_{j_k}]$$ for $w \in
T(k), v^i \in T(j_i)$,  and $(x^i_1, \dots, x^i_{j_i}) \in X^{j_i}$
for $i=1, \dots, k$. On morphisms we define it to be
$$\mu_X([i_w,[i_{v^1},g^1_1, \dots, g^1_{j_1}],[i_{v^2},g^2_1, \dots,
g^2_{j_2}], \dots, [i_{v^k},g^k_1, \dots, g^k_{j_k}]]):=$$ $$
[i_{\gamma(w,v^1,v^2, \dots, v^k)},g^1_1, \dots, g^k_{j_k}].$$ These
assignments make $\mu_X:C^2X \rightarrow CX$ into a well defined
functor because of the equivariances of $\gamma$. These natural
transformations commute appropriately to make $C$ into a {\it
monad}\index{monad} on the category $Cat_0$.

\begin{thm} \label{C=T}
The category of $C$-algebras\index{algebra!$C$-algebra} is
equivalent to the category of
$T$-algebras\index{algebra!$T$-algebra}.
\end{thm}
\begin{pf}
Let $\mathcal{C}_C$ and $\mathcal{C}_T$ denote the categories of
$C$-algebras and $T$-algebras respectively. We construct a functor
$\phi:\mathcal{C}_T \rightarrow \mathcal{C}_C$. Let $(X,\Phi)$ be a
$T$-algebra. Then $\Phi_n:T(n) \rightarrow Functors(X^n,X)$ is a
sequence of maps that is natural in $n$, preserves identity $1 \in
T(1)$, and preserves compositions $\gamma$. This sequence of maps
completely describes the algebraic structure. Let $h'$ denote the
element of $Functors(\bigcup_{n \geq 0}(T(n) \times X^n), X)$ that
corresponds to the sequence under the bijection
\begin{equation} \label{mapbijection}
Functors(\bigcup_{n \geq 0}(T(n) \times X^n), X) \leftrightarrow
 \prod_{n \geq 0} Functors(T(n), X^{X^n}).
 \end{equation} Then
$$h'(w_f, x_1, \dots, x_{n})=h'(w,x_{f1},\dots, x_{fm})$$
$$h'(i_{w_f}, g_1, \dots, g_{n})=h'(i_w,g_{f1},\dots, g_{fm})$$
because
$$\xymatrix@C=4pc@R=4pc{T(m) \ar[r]^-{\Phi_m} \ar[d]_{()_f} & Functors(X^m,X)
\ar[d]^{()_f} \\ T(n) \ar[r]_-{\Phi_n} & Functors(X^n,X)}$$
commutes. Hence $h':\bigcup_{n \geq 0} (T(n) \times X^n) \rightarrow
X$ induces a functor $h:CX \rightarrow X$, namely
$$[w,x_1, \dots, x_m] \mapsto \Phi_m(w)(x_1, \dots x_m)$$
$$[i_w, g_1, \dots, g_m] \mapsto \Phi_m(i_w)(y_1, \dots, y_m)
\circ \Phi_m(w)(g_1, \dots, g_m)$$
$$=\Phi_m(w)(g_1, \dots, g_m)$$
for $g_i:x_i \rightarrow y_i$. Then $h:CX \rightarrow X$ makes $X$
into a $C$-algebra because the diagrams
$$\xymatrix@R=3pc@C=3pc{C^2X \ar[d]_{\mu_X} \ar[r]^{Ch} & CX \ar[d]^{h}
& X \ar[r]^{\eta_X} \ar[dr]_{1_X}  & CX \ar[d]^h \\
CX \ar[r]_{h} & X & & X}$$ commute.

We define $\phi((X,\Phi)):=(X,h)$. For a morphism $H:(X, \Phi)
\rightarrow (Y, \Psi)$ of $T$-algebras, let $\phi(H):X \rightarrow
Y$ be the same functor as $H$ on the underlying categories. Then
$$\xymatrix@R=3pc@C=3pc{CX \ar[r]^{h_X} \ar[d]_{C\phi(H)} & X \ar[d]^{\phi(H)}
\\ CY \ar[r]_{h_Y}  & Y}$$
commutes. Then $\phi:\mathcal{C}_T \rightarrow \mathcal{C}_C$ is
obviously a functor.

An ``inverse'' to $\phi$ can easily be constructed using the
bijection (\ref{mapbijection}). For example, let $(X,h)$ be a
$C$-algebra. Then $h:CX \rightarrow X$ corresponds uniquely to a
functor $h':\bigcup_{n \geq 0}(T(n) \times X^n) \rightarrow X$ which
satisfies $$h'(w_f, x_1, \dots, x_{n})=h'(w,x_{f1},\dots, x_{fm})$$
$$h'(i_{w_f}, g_1, \dots, g_{n})=h'(w,g_{f1},\dots, g_{fm})$$
and $h'$ corresponds uniquely to some sequence $\Phi_n$ natural in
$n$ which preserves 1 and $\gamma$.

The equivalence of Theorem \ref{star} yields the desired result.
\end{pf}

The concept of theory can be generalized to handle algebraic
structures on more than one set, such as modules.

\begin{defn}
A {\it theory on a set of objects $J$}\index{theory!theory on a set
of objects|textbf}, also called a {\it many-sorted
theory}\index{theory!many sorted theory|textbf}\index{many sorted
theory|textbf}, is a category $\mathbf{T}$ whose objects are finite
sequences $(j_1^{m_1}, \dots, j_p^{m_p})$ with $j_1, \dots, j_p \in
J,p \geq 1$, and $m_1, \dots, m_p \in \mathbb{N}_0$ such that
$(j_1^{m_1}, \dots, j_p^{m_p})$ is a product of copies of $j \in J$
where each $j$ appears $\sum_{r:j_r=j} m_r$ times. Each sequence is
equipped with a limiting cone. Objects are equal to their reduced
form, \eg $(j^{m_1},j^{m_2})=(j^{m_1+m_2})$. We also abbreviate
$(j^1)=j$.
\end{defn}

\begin{examp}
An ordinary theory\index{theory} is a theory on one object, \ie on
the set $\{1\}$. We previously used $n$ to denote $1^n$ in the new
notation.
\end{examp}

\begin{examp}
Let $X_1$ and $X_2$ be categories. Then the {\it endomorphism theory
\newline $ End(X_j:j \in J)$ on $X_1$ and
$X_2$}\index{theory!endomorphism theory|textbf}\index{endomorphism
theory|textbf} is an example of a theory on the set $J=\{1,2\}$. The
morphisms are
$$Mor_{ End(X_j:j \in J)}((j_1^{m_1}, \dots, j_p^{m_p}),(k_1^{n_1}, \dots,
k_q^{n_q})):=$$
$$Functors(X_{j_1}^{m_1} \times \cdots \times
X_{j_p}^{m^p},X_{k_1}^{n_1} \times \cdots \times X_{k_q}^{n^q})$$
for $j_r,k_s \in \{1,2\}$ and $m_r,n_s \in \mathbb{N}_0$. We easily
see that $1^0$ and $2^0$ as well as $(1^0,2^0)$ and $(2^0,1^0)$ are
terminal objects\index{terminal object} and that $(j_1^{m_1}, \dots,
j_p^{m_p})$ is a product of $\sum_{r:j_r=1} m_r$ copies of $1$ and
$\sum_{r:j_r=2} m_r$ copies of $2$ equipped with the usual
projections. Note also that there is a bijective correspondence.
$$Mor_{ End(X_j:j \in J)}((j_1^{m_1}, \dots, j_p^{m_p}),(k_1^{n_1}, \dots,
k_q^{n_q}))$$ $$\updownarrow$$

\begingroup
\vspace{-2\abovedisplayskip} \Small
$$ \prod_{r:k_r=1}Mor_{ End(X_j:j \in
J)}((j_1^{m_1}, \dots, j_p^{m_p}),1)^{\times n_r} \times
\prod_{s:k_s=2}Mor_{ End(X_j:j \in J)}((j_1^{m_1}, \dots,
j_p^{m_p}),2)^{\times n_s}$$
\endgroup
\noindent In other words, the theory is determined by the sets
$$Mor_{ End(X_j:j \in J)}((j_1^{m_1}, \dots,
j_p^{m_p}),1)=: End(X_j:j \in J)_1(j_1^{m_1}, \dots, j_p^{m_p})$$
$$Mor_{ End(X_j:j \in J)}((j_1^{m_1}, \dots,
j_p^{m_p}),2)=: End(X_j:j \in J)_2(j_1^{m_1}, \dots, j_p^{m_p})$$
where $j_1, \dots, j_p \in \{1,2\}$ and $m_1, \dots, m_p \in
\mathbb{N}_0$ such that $j_r \neq j_{r+1}$ for all $1 \leq r \leq
p-1$.

Note also that for $n_1, \dots, n_q \in \mathbf{N}_0$ and $k_1,
\dots, k_q \in J$ and maps $$f:\sum_{r:j_r=1}m_r \rightarrow
\sum_{r:k_r=1}n_r$$
$$g:\sum_{s:j_s=2}m_s \rightarrow \sum_{s:k_s=2}n_s$$
in $\Gamma$ we have {\it substitution}\index{substitution} maps
$$\xymatrix@1@C=3pc{ End(X_j:j \in J)_1(j_1^{m_1}, \dots, j_p^{m_p})
\ar[r]^{()_{f,g}} &  End(X_j:j \in J)_1(k_1^{n_1}, \dots,
k_q^{n_q})}$$
$$\xymatrix@1@C=3pc{ End(X_j:j \in J)_2(j_1^{m_1}, \dots, j_p^{m_p})
\ar[r]^{()_{f,g}} &  End(X_j:j \in J)_2(k_1^{n_1}, \dots,
k_q^{n_q})}.$$ For example, let $w \in
 End(X_j:j \in J)_1(1^2,2^2,1^1,2^2)$ and $$f:=\begin{pmatrix} 1 & 2 &
3 \\ 1 & 1 & 1
\end{pmatrix}, \hspace{1pc} g:=\begin{pmatrix} 1 & 2 &
3 & 4 \\ 1 & 2 & 2 &1
\end{pmatrix},$$ where $(k_1^{n_1}, k_2^{n_2})=(1^1,2^2)$ so that
$$f:3 \rightarrow 1, \hspace{1pc} g:4 \rightarrow 2.$$ Then
$w_{f,g} \in  End(X_j:j \in J)_1(1^1,2^2)$ is defined by
$$\aligned
w_{f,g}(x^1_1, x^2_1, x^2_2) &:=
w(x^1_{f1},x^1_{f2},x^2_{g1},x^2_{g2},x^1_{f3},x^2_{g3},x^2_{g4})\\
&=w(x^1_1,x^1_1,x^2_1,x^2_2,x^1_1,x^2_2,x^2_1).
\endaligned$$
The notation $()_{f,g}$ suppresses the dependence of the map
$()_{f,g}$ on $(j_1^{m_1}, \dots, j_p^{m_p})$ and $(k_1^{n_1},
\dots, k_q^{n_q})$.

There are also two {\it compositions}\index{composition} $\gamma_1$
and $\gamma_2$. For example $$\gamma_1: End(X_j:j \in J)_1(1^2,2^2)
\times
 End(X_j:j \in J)_1(\bar{n}_1) \times  End(X_j:j \in J)_1(\bar{n}_2)
\times$$
$$\times  End(X_j:j \in J)_2(\bar{n}_3) \times  End(X_j:j \in J)_2(\bar{n}_4)
\rightarrow  End(X_j:j \in J)_1(\bar{n}_1 \cdot \bar{n}_2 \cdot
\bar{n}_3 \cdot \bar{n}_4)$$ and
$$\gamma_2: End(X_j:j \in J)_2(2^3,1^1) \times
 End(X_j:j \in J)_2(\bar{n}_1) \times  End(X_j:j \in J)_2(\bar{n}_2)
\times$$
$$\times  End(X_j:j \in J)_2(\bar{n}_3) \times  End(X_j:j \in J)_1(\bar{n}_4)
\rightarrow  End(X_j:j \in J)_2(\bar{n}_1 \cdot \bar{n}_2 \cdot
\bar{n}_3 \cdot \bar{n}_4)$$ where $\bar{n}_1 \cdot \bar{n}_2 \cdot
\bar{n}_3 \cdot \bar{n}_4$ means to concatenate the objects
$\bar{n}_1, \dots, \bar{n}_4$ and to reduce, \eg $(1^1,2^2) \cdot
(2^3,1^2)=(1^1,2^5,1^2)$.

There are also units $1_1 \in  End(X_j:j \in J)_1(1)$ and $1_2 \in
 End(X_j:j \in J)_2(2)$.

The compositions are associative, unital, and equivariant. The
substitution is also functorial. This example easily extends to
arbitrary $J$.
\end{examp}

\begin{defn}
Let $\Gamma_J$ denote the category whose objects are finite
sequences \newline $(j_1^{m_1}, \dots, j_p^{m_p})$ with $j_1, \dots,
j_p \in J,p \geq 1$, and $m_1, \dots, m_p \in \mathbb{N}_0$. Objects
are equal to their reduced form, \eg
$(j^{m_1},j^{m_2})=(j^{m_1+m_2})$. We also abbreviate $(j^1)=j$. The
morphisms are $$Mor_{\Gamma_J}((j_1^{m_1}, \dots,
j_p^{m_p}),(k_1^{n_1}, \dots, k_q^{n_q})):= \prod_{\ell \in J}
Mor_{\Gamma}(\sum_{r:j_r=\ell}m_r, \sum_{s:k_s=\ell}n_s)$$ where
$\Gamma$ denotes the category in Definition \ref{definitionofGamma}.

\end{defn}
In this definition the hom sets are assumed to be disjoint.

Several of the results on theories carry over to these generalized
theories on a set of objects.

\begin{thm}
A theory\index{theory|theory on a set of objects|textbf}
$\mathbf{T}$ on a set of objects $J$ is equivalent to a collection
of functors $\{\mathbf{T}_j:\Gamma_J \rightarrow Sets| j \in J \}$
equipped with compositions
$$\gamma_j:\mathbf{T}_j(j_1^{k_1}, \dots, j_p^{k_p}) \times
\mathbf{T}_{j_1}(\bar{n}_1^1) \times \cdots \times
\mathbf{T}_{j_1}(\bar{n}_{k_1}^1) \times$$
$$\times \mathbf{T}_{j_2}(\bar{n}_1^2) \times \cdots \times \mathbf{T}_{j_2}(\bar{n}_{k_2}^2)
\times$$
$$ \cdots $$
$$\times \mathbf{T}_{j_p}(\bar{n}_1^p) \times \cdots \times \mathbf{T}_{j_p}(\bar{n}_{k_p}^p)
\rightarrow \mathbf{T}_j(\bar{n}_1^1 \cdots \bar{n}_{k_1}^1 \cdot
\bar{n}^2_1 \cdots \bar{n}^2_{k_2} \cdots \bar{n}_1^{p} \cdots
\bar{n}_{k_p}^p)$$ for each $j \in J$ and $(j_1^{k_1}, \dots,
j_p^{k_p}), \bar{n}_1^1, \dots, \bar{n}_{k_p}^p \in Obj \hspace{1mm}
\Gamma_J$ and equipped with units $1_j \in \mathbf{T}_j(j)$ for each
$j \in J$ which satisfy analogues of (1) through (5) in Lemma
\ref{tofunctortheory}. Elements of $\mathbf{T}_j(\bar{n})$ are
called {\it words}.
\end{thm}
\begin{pf}
Set $\mathbf{T}_j(\bar{n}):=Mor_{\mathbf{T}}(\bar{n},j)$ and proceed
like in the case of a theory on the set $\{1\}$.
\end{pf}

\begin{examp} \label{theoryoftheoriesexample}
The {\it theory $\mathbf{R}$ of theories}\index{theory!theory of
theories|textbf} is a theory on the set $\mathbb{N}_0$. There are
three types of generating morphisms.
\begin{itemize}
\item
For each $k \geq 1$ and $n_1, \dots, n_k \geq 0$ there is a morphism
$\gamma:(k,n_1, \dots, n_k) \rightarrow (n_1+ \cdots + n_k)$ called
{\it composition}\index{composition|textbf}.
\item
For each $f:m \rightarrow n$ in $\Gamma$ there is a morphism
$()_f:(m) \rightarrow (n)$ called {\it
substitution}\index{substitution|textbf}.
\item
There is a morphism $1:(1^0) \rightarrow (1^1)$ called the {\it
unit}\index{unit|textbf}.
\end{itemize}
The substitution and unit are not to be confused with the
substitution and units with which every theory on a set of objects
is equipped. These morphisms must satisfy the relations of theories
in Lemma \ref{tofunctortheory}, namely associativity, equivariances,
unitality, and functoriality.
\end{examp}

Next we can speak of morphisms of theories on the set $J$ as well as
algebras for theories on the set $J$ just as in the case $J=\{1\}$.

\begin{defn}
A {\it morphism of theories on a set $J$}\index{morphism!morphism of
theories on a set of objects|textbf} is a functor $\Phi: \mathbf{S}
\rightarrow \mathbf{T}$ such that $\Phi(j_1^{m_1}, \dots,
j_p^{m_p})=(j_1^{m_1}, \dots, j_p^{m_p})$ and $\Phi(pr)=pr$ for
every projection.
\end{defn}

\begin{thm}
The analogue of Theorem \ref{theorydescriptions} holds for
theories\index{theory!theory on a set of objects} on a set of
objects $J$.
\end{thm}

\begin{defn}
Let $\mathbf{T}$ be a theory on the set $J$ and $\{X_j|j \in J\}$ a
collection of categories. Then $\{X_j\}_j$ form an {\it algebra over
$\mathbf{T}$} or a {\it
$\mathbf{T}$-algebra}\index{algebra!$\mathbf{T}$-algebra|textbf}\index{algebra!algebra
over a theory on a set of objects|textbf} if they are equipped with
a morphism $\Phi:\mathbf{T} \rightarrow End(X_j:j\in J)$ of theories
on $J$.
\end{defn}

\begin{examp}
Let $\mathbf{R}$ denote the theory of theories\index{theory!theory
of theories|textbf}. Let $T$ be a theory. Then $\{T(j)|j \in
\mathbb{N}_0\}$ form an $\mathbf{R}$-algebra. In other words, a
theory is an algebra over the theory of
theories\index{algebra!algebra over the theory of theories}. A
morphism of theories is nothing more than a morphism of algebras
over the theory of theories\index{morphism!morphism of theories}.
\end{examp}

\begin{thm}
The analogue of Theorem \ref{star} holds for a theory $\mathbf{T}$
on a set of objects.
\end{thm}

We can use the theory $\mathbf{R}$ of \label{freetheorymonad}
theories to construct a monad\index{monad} $C$ on the category
$\prod_{n\geq 0} Sets$ whose algebras are the usual theories. In
fact, $CT$ is the sequence of sets underlying the {\it free theory}
on $T$. This free theory is essential to several of the proofs in
this paper. Let $T=(T(n))_{n \geq 0}$ be an object of $\prod_{n\geq
0} Sets$ and $J:=\mathbb{N}_0$. Then the {\it free
theory}\index{theory!free theory|textbf}\index{free theory|textbf}
on $T$ is defined by
$$CT(n):=\frac{\bigcup_{\bar{m} \in Obj \hspace{1mm}  \Gamma_J}
\mathbf{R}_n(\bar{m})\times T(j_1)^{\times m_1} \times \cdots \times
T(j_p)^{\times m_p}}{\Gamma_J}$$ where $\bar{m}=(j_1^{m_1}, \dots,
j_p^{m_p})$.

We can generalize the notion of theory\index{theory} in yet another
direction. Instead of considering arbitrary sets $J$, we can
consider theories which are also 2-categories in which every 2-cell
is iso. We will use these to describe pseudo
algebras\index{algebra!pseudo algebra} in a compact way. See
\cite{power3} for a more general concept of enriched Lawvere
theory\index{Lawvere theory}\index{Lawvere theory!enriched Lawvere
theory|textbf}\index{theory!enriched theory|textbf}.

\begin{defn}
A {\it theory enriched in groupoids}\index{theory!theory enriched in
groupoids|textbf}\index{groupoid} is a 2-category $\mathcal{T}$ with
iso 2-cells and with objects $0,1,2, \dots$ such that $n$ is the
2-product\index{2-product} of $1$ with itself $n$ times in the
2-category $\mathcal{T}$ and each $n$ is equipped with a limiting
2-cone.
\end{defn}

This definition means for each $n \in Obj \hspace{1mm} \mathcal{T}$
we have chosen morphisms $\pi_i^n=pr_i:n \rightarrow 1$ for $i=1,
\dots, n$ with the universal property that
$$\xymatrix@1@C=3pc{Mor_{\mathcal{T}}(m,n) \ar[r]^-{\pi^n \circ} &
2-Cone(m,F)}$$ is an isomorphism for all $m \in Obj \hspace{1mm}
\mathcal{T}$, where $F:\{1, \dots, n\} \rightarrow \mathcal{T}$ is
the 2-functor which is constant 1. It is tempting to call such a
theory a 2-theory, but we reserve that name for something else. As
before, we use the notation $\mathcal{T}(n)$ for the category
$Mor_{\mathcal{T}}(n,1)$. Using the universal property, we can
construct $\prod$ and $(\dots)$ for the 2-cells. For any object $m
\in Obj \hspace{1mm}  T$, morphisms $w_i, v_i:m \rightarrow 1$, and
2-cells $\alpha_i:w_i \Rightarrow v_i$ for $i=1, \dots,n$, there
exists a unique 2-cell $\prod_{j=1}^n \alpha_j:\prod_{j=1}^n w_j
\Rightarrow \prod_{j=1}^n v_j$ such that
$$i_{pr_i} * \prod_{j=1}^n \alpha_j = \alpha_i$$ for all
$i=1, \dots n$. For any $k \in \mathbb{N}_0$, any morphisms
$w_i,v_i:n_i \rightarrow 1$, and any 2-cells $\alpha_i:w_i
\Rightarrow v_i$ for $i=1, \dots, k$, there is a unique 2-cell
$(\alpha_1, \dots, \alpha_k):(w_1, \dots, w_k) \Rightarrow (v_1,
\dots, v_k)$ such that
$$i_{pr_i} * (\alpha_1, \dots, \alpha_k) = (\alpha_i)_{\iota_i}$$
for all $i=1, \dots, k$.

\begin{examp}
Let $X$ be a category. Then the {\it endomorphism theory $End(X)$
enriched in groupoids}\index{theory!endomorphism theory enriched in
groupoids|textbf} has objects $0,1,2, \dots$, morphisms $Obj
\hspace{1mm} Mor_{End(X)}(m,n)=Functors(X^m,X^n)$ and 2-cells the
natural isomorphisms. \index{$End(X)$|textbf}
\end{examp}

Most of the work on theories carries over to the enriched context
with minor additions for the 2-cells. The statements of the relevant
theorems are as follows. The term {\it map} is simply replaced by
{\it functor}.

\begin{lem}
Let $\mathcal{T}$ be a theory enriched in groupoids. Then the
morphism category $Mor_{\mathcal{T}}(m,n)$ is isomorphic to the
product category $\prod_{j=1}^n Mor_{\mathcal{T}}(m,1)$.
\end{lem}

\begin{lem}
Let $\mathcal{T}$ be a theory enriched in
groupoids\index{theory!theory enriched in groupoids|textbf}. Then
for all $k,n_1, \dots,n_k \\ \in \{0,1, \dots \}$ there is a functor
$\gamma: \mathcal{T}(k) \times \mathcal{T}(n_1) \times \dots \times
\mathcal{T}(n_k) \rightarrow \mathcal{T}(n_1 + \cdots + n_k)$ called
{\it composition} and for every function $f:\{1, \dots, k\}
\rightarrow \{1, \dots, \ell \}$ there is a functor
$\xymatrix@1{\mathcal{T}(k) \ar[r]^{()_f} & \mathcal{T}(\ell)}$
called {\it substitution}. These functors satisfy the enriched
analogues of (1) through (5) in Lemma \ref{tofunctortheory}.
\end{lem}
\begin{pf}
Define $\gamma(w, w_1, \dots, w_k):=w \circ (w_1, \dots, w_k)$ as
before. Additionally, define $\gamma(\alpha, \alpha_1, \dots,
\alpha_k):=\alpha * (\alpha_1, \dots, \alpha_k)$ for 2-cells. Define
$w_f:=w \circ f'$ as before and $\alpha_f:= \alpha * i_{f'}$ where
$i_{f'}:f' \Rightarrow f'$ is the identity 2-cell of the morphism
$f'$ in $\mathcal{T}$ and $\alpha:w \Rightarrow v$ is a 2-cell. The
rest of proof is similar to Lemma \ref{tofunctortheory}.
\end{pf}

\begin{lem} \label{2functortheory}
Let $\mathcal{T}$ be a 2-functor from $\Gamma$ to the 2-category
$Cat$ of small categories equipped with functors $\gamma:
\mathcal{T}(k) \times \mathcal{T}(n_1) \times \dots \times
\mathcal{T}(n_k) \rightarrow \mathcal{T}(n_1 + \dots + n_k)$ and an
object $1 \in \mathcal{T}(1)$ which satisfy (1) through (5) of Lemma
\ref{tofunctortheory} where $()_f:=\mathcal{T}(f)$ for functions
$f:k \rightarrow \ell$. Then $\mathcal{T}$ determines a theory
enriched in groupoids\index{groupoid} with $Mor(n,1)=\mathcal{T}(n)$
for all $n \geq 0$.
\end{lem}

\begin{thm} \label{enrichedtheoryequivalence}
A theory $\mathcal{T}$ enriched in groupoids\index{theory!theory
enriched in groupoids|textbf} is determined by either of the
following equivalent collections of data:
\begin{enumerate}
\item
A 2-category $\mathcal{T}$ with objects $0,1,2, \dots$ such that $n$
is the 2-categorical product of $1$ with itself $n$ times and each
$n$ is equipped with a limiting 2-cone.
\item
A 2-functor $\mathcal{T}:\Gamma \rightarrow Cat$ equipped with
functors $\gamma:\mathcal{T}(k) \times \mathcal{T}(n_1) \times
\cdots \times \mathcal{T}(n_k) \rightarrow \mathcal{T}(n_1 + \cdots
+ n_k)$ and a unit $1 \in \mathcal{T}(1)$ which satisfy (1) through
(5) of Lemma \ref{tofunctortheory}.
\end{enumerate}
\end{thm}
\begin{pf}
In each description $Mor_{\mathcal{T}}(n,1)$ is the same. By the
universality of 2-products\index{2-product} this determines the rest
of the theory.
\end{pf}

\begin{defn}
Let $\mathcal{S}$ and $\mathcal{T}$ be theories enriched in
groupoids. In the 2-categorical description of $\mathcal{S}$ and
$\mathcal{T}$ a {\it morphism of theories enriched in
groupoids}\index{morphism!morphism of theories enriched in
groupoids|textbf} $\Phi:\mathcal{S} \rightarrow \mathcal{T}$ is a
2-functor from the 2-category $\mathcal{S}$ to the 2-category
$\mathcal{T}$ such that $\Phi(n_{\mathcal{S}})=n_{\mathcal{T}}$ and
$\Phi(pr_i)=pr_i$ for all projections.
\end{defn}

The analogue for Lemma \ref{theorymorphisms} incorporates the
2-cells below.

\begin{lem}
Let $\Phi:\mathcal{S} \rightarrow \mathcal{T}$ be a morphism of
theories enriched in groupoids\index{groupoid}.
\begin{enumerate}
\item
Let $f: \{1, \dots, k\} \rightarrow \{1, \dots, \ell\}$ be a
function. As usual, $f': \ell \rightarrow k$ denotes the unique
morphism in any theory such that
$$\xymatrix@R=3pc@C=3pc{k \ar[r]^{pr_i} & 1 \\ \ell \ar@{.>}[u]^{f'} \ar[ur]_{pr_{fi}}
}$$ commutes. Then $\Phi(f')=f'$.
\item
Let $f: \{1, \dots, k\} \rightarrow \{1, \dots, \ell\}$ be a
function and $w \in Mor_{\mathcal{S}}(k,1)$. Then
$\Phi(w_f)=\Phi(w)_f$.
\item
Let $w_j, v_j \in Mor_{\mathcal{S}}(m,1)$ and $\alpha_j:w_j
\Rightarrow v_j$ for $j=1, \dots, n$. Then $\Phi(\prod_{j=1}^n w_j)=
\prod_{j=1}^n \Phi(w_j)$ and $\Phi(\prod_{j=1}^n \alpha_j)=
\prod_{j=1}^n \Phi(\alpha_j)$.
\item
Let $w_j,v_j \in Mor_{\mathcal{S}}(n_j,1)$ for $j=1, \dots, k$. Then
we have $\Phi(w_1, \dots, w_k)=(\Phi(w_1), \dots, \Phi(w_k))$ and
$\Phi(\alpha_1, \dots, \alpha_k)=(\Phi(\alpha_1), \dots,
\Phi(\alpha_k))$.
\end{enumerate}
\end{lem}

\begin{thm} \label{enrichedmorphismequivalence}
Let $\mathcal{S}$ and $\mathcal{T}$ be theories enriched in
groupoids\index{groupoid}. Then a morphism $\mathcal{S} \rightarrow
\mathcal{T}$ of theories enriched in groupoids is given by either of
the following equivalent collections of data:
\begin{enumerate}
\item
A 2-functor $\Phi:\mathcal{S} \rightarrow \mathcal{T}$ such that
$\Phi(n_{\mathcal{S}})=n_{\mathcal{T}}$ for all $n_{\mathcal{S}} \in
Obj \hspace{1mm}  \mathcal{S}$ and $\Phi(pr_i)= pr_i$ for all
projections
\item
A 2-natural transformation $\Phi:\mathcal{S} \Rightarrow
\mathcal{T}$ of the 2-functors $\mathcal{S},\mathcal{T}: \Gamma
\rightarrow Cat$ which preserves the $\gamma$'s and the units.
\end{enumerate}
\end{thm}

\begin{thm} \label{enrichedtheorydescriptions}
The 2-category of theories enriched in groupoids\index{theory!theory
enriched in groupoids|textbf}\index{groupoid} with objects and
morphisms as in (1) of Theorems \ref{enrichedtheoryequivalence} and
\ref{enrichedmorphismequivalence} is 2-equivalent to the 2-category
with objects and morphisms as in (2) of Theorems
\ref{enrichedtheoryequivalence} and
\ref{enrichedmorphismequivalence}.
\end{thm}

We can now define algebras over theories enriched in groupoids in
analogy to algebras over theories.

\begin{defn}
Let $X$ be a category and $\mathcal{T}$ a theory enriched over
groupoids. Then $X$ is a {\it $\mathcal{T}$-algebra} if it is
equipped with a morphism of theories $\mathcal{T} \rightarrow
End(X)$ enriched in groupoids. We also say $X$ is an {\it algebra
over the theory
$\mathcal{T}$}\index{algebra!$\mathcal{T}$-algebra|textbf}\index{algebra!algebra
over a theory enriched in groupoids|textbf}\index{groupoid}.
\end{defn}

Our main example, pseudo $T$-algebras, will be given in the next
chapter as strict $\mathcal{T}$-algebras, where $\mathcal{T}$ is
obtained from the free theory\index{theory!free theory}\index{free
theory} on $T$.

\begin{thm} \label{star2}
The analogue of Theorem \ref{star} holds for theories enriched in
groupoids\index{groupoid}\index{theory!theory enriched in
groupoids}.
\end{thm}

\chapter{Pseudo $T$-Algebras}
\label{sec:laxTalgebras} In this chapter we introduce the 2-category
of pseudo $T$-algebras for a theory $T$. A {\it pseudo algebra} in
this paper is the same thing as a {\it lax
algebra}\index{algebra!lax algebra} in \cite{hu}, \cite{hu1}, and
\cite{hu2}. We construct from $T$ a theory $\mathcal{T}$ enriched in
groupoids\index{theory!theory enriched in groupoids}\index{groupoid}
and show that a pseudo algebra over $T$ is the same thing as an
algebra over $\mathcal{T}$. Theorem
\ref{laxTalgebras=strictCalgebras} says that the 2-category of
pseudo $T$-algebras and pseudo morphisms is 2-equivalent to the
2-category of strict $C$-algebras with pseudo morphisms for the
2-monad\index{2-monad} $C$ defined on page
\pageref{2monadprecisely}. This 2-category of strict $C$-algebras
and pseudo morphisms admits pseudo limits by a result of Blackwell,
Kelly, and Power in \cite{blackwell}. Hence the 2-category of pseudo
$T$-algebras admits pseudo limits. In the next chapter we give a
concrete construction of a pseudo limit. For more on pseudo algebras
over 2-monads\index{2-monad} see \cite{hermida}, \cite{lack2}, and
\cite{lack}.

\begin{defn} \label{laxalgebradefinition}
Let $T$ be a theory. A category $X$ is a {\it pseudo $T$-algebra} or
a {\it pseudo algebra over $T$}\index{algebra!pseudo algebra|textbf}
if it is equipped with {\it structure maps}  $\Phi_n : T(n)
\rightarrow Functors(X^n,X)$ for every $n \in \mathbb{N}$ as well as
the coherence isomorphisms\index{coherence isomorphism|(textbf}
below\index{algebra!pseudo $T$-algebra|textbf}\index{algebra!pseudo
algebra over $T$|textbf}. Moreover, the coherence isomorphisms are
required to satisfy the coherence diagrams\index{coherence
diagram|(textbf} below. We write simply $\Phi$ for all $\Phi_n$. The
coherence isomorphisms are indexed\index{indexed} by the operations
of theories\index{operations of theories} and are as follows:
\begin{enumerate}
\item
For every $k \in \mathbb{N}$, $w \in T(k)$, and all words $w_1,
\ldots, w_k$, there is a natural isomorphism\index{$c$} $c_{w, w_1,
\ldots ,w_k} : \Phi(\gamma(w,w_1, \dots, w_k)) \Rightarrow
\gamma(\Phi(w), \Phi(w_1), \ldots, \Phi(w_k))$. This means that
$\Phi$ preserves composition up to a natural isomorphism.

\item
There is a natural isomorphism\index{$I$} $I:\Phi(1) \Rightarrow
1_X$ where $1$ is the identity word and $1_X$ is the identity
functor $X \rightarrow X$. This means that $\Phi$ preserves the
identity up to a natural isomorphism.

\item
For every word $w \in T(m)$ and function $f:\{1, \ldots, m\}
\rightarrow \{1, \ldots, n\}$, there is a natural
isomorphism\index{$s$} $s_{w,f}: \Phi(w_f) \Rightarrow \Phi(w)_f$
where the substituted functor $\Phi(w)_f: X^n \rightarrow X$ is
defined in Examples and \ref{Endexample1} and \ref{Endexample2}.
This means that $\Phi$ preserves the substitution up to a natural
isomorphism.
\end{enumerate}
The coherence diagrams are indexed\index{indexed} by
relations\index{relations of theories} of theories and are as
follows. The commutivity of these diagrams means that they commute
when evaluated on every tuple of objects of $X$ of appropriate
length.

\begin{enumerate}
\item
The composition coherence isomorphisms are associative. For example,
for $u,v,w \in T(1)$ the diagram below must commute where $i_F$
means the identity natural transformation $F \rightarrow F$ for a
functor $F$.
\newpage
\begingroup
\vspace{-2\abovedisplayskip} \Small
$$
\xy (-38,0)*{ \xymatrix@R=3pc@C=.2pc{\hspace{2cm} \Phi(\gamma(w,
\gamma(v,u))) = \Phi(\gamma(\gamma(w, v),u))
\ar@{=>}[r]^-{c_{\gamma(w, v),u} } \ar@{=>}[d]_{c_{w, \gamma(v,u)}}
& \gamma(\Phi(\gamma(w, v)), \Phi(u))
\ar@{=>}[d]^{\gamma(c_{w,v},i_{\Phi(u)})} \\
\gamma(\Phi(w), \Phi(\gamma(v, u)))
\ar@{=>}[r]_-{\gamma(i_{\Phi(w)}, c_{v , u})} & \gamma(\Phi(w),
\gamma( \Phi (v), \Phi(u)))=\gamma(\gamma(\Phi(w),\Phi(v)), \Phi(u))
\hspace{3.2cm}} };\endxy
$$
\endgroup
\noindent \\ \\ \\ \\
\item
The natural isomorphism for the identity word commutes with the
natural isomorphism for the composition, \ie for every $n \in
\mathbb{N}$ and every word $w \in T(n)$ the diagram below must
commute where $1_X$ is the identity functor on $X$.

\begingroup
\vspace{-2\abovedisplayskip} \small
$$\xymatrix@C=6pc@R=3pc{\Phi(\gamma(w,1, \ldots, 1))
 \ar@{=}[r] \ar@{=>}[d]_{c_{w,1, \ldots ,1}} &
\Phi(w) \ar@{=}[d] \\ \gamma(\Phi(w),\Phi(1), \ldots, \Phi(1))
\ar@{=>}[r]_-{\gamma(i_{\Phi(w)}, I, \ldots, I)} &
\gamma(\Phi(w),1_X, \ldots, 1_X)}$$
\endgroup
\noindent

\item
The natural isomorphism for the identity word commutes with the
natural isomorphism for the composition also in the sense that for
every word $w \in T(n)$ the diagram below must commute.

\begingroup
\vspace{-2\abovedisplayskip} \small
$$\xymatrix@C=5pc@R=3pc{\Phi(\gamma(1,w)) \ar@{=}[r] \ar@{=>}[d]_{c_{1,w}}
& \Phi(w) \ar@{=}[d] \\ \gamma(\Phi(1), \Phi(w))
\ar@{=>}[r]_{\gamma(I,i_{\Phi(w)})} & \gamma(1_X, \Phi(w))}$$
\endgroup
\noindent

\item
Let $f:\{1, \dots, k\} \rightarrow \{1, \dots, \ell \}$ be a
function and let \newline $\bar{f}:\{1,2, \dots, n_{f1} + n_{f2} +
\dots + n_{fk} \} \rightarrow \{1,2, \dots, n_1 + n_2 + \dots +
n_{\ell} \}$ be the function that moves entire blocks according to
$f$ as in Example \ref{Endexample1}.  Then equivariance is preserved
in the sense that the diagram below must commute.

\begingroup
\vspace{-2\abovedisplayskip} \small
$$
\xymatrix@C=7pc@R=3pc{ \Phi(\gamma(w,w_{f1}, \dots,
w_{fk})_{\bar{f}}) \ar@{=}[d] \ar@{=>}[r]^{s_{\gamma(w,w_{f1},
\dots, w_{fk}), \bar{f}}}
 & \Phi(\gamma(w,
w_{f1}, \dots, w_{fk}))_{\bar{f}} \ar@{=>}[d]|{(c_{w,w_{f1},
\dots, w_{fk}})_{\bar{f}}} \\
\Phi(\gamma(w_f,w_1, \dots, w_{\ell})) \ar@{=>}[d]|{c_{w_f, w_1,
\dots, w_{\ell}}} & \gamma(\Phi(w), \Phi(w_{f1}),
\Phi(w_{fk}))_{\bar{f}} \ar@{=}[d]
\\ \gamma(\Phi(w_f),
\Phi(w_1), \dots, \Phi(w_{\ell}))
\ar@{=>}[r]_{\overset{\phantom{M}}{\gamma(s_{w,f},i_{\Phi(w_1)},
\dots, i_{\Phi(w_{\ell})})}} & \gamma(\Phi(w)_f, \Phi(w_1), \dots,
\Phi(w_{\ell})) }
$$
\endgroup
\noindent

\item
Let $g_i:\{1, \dots, n_i \} \rightarrow \{1, \dots, n_i'\}$ be
functions and let \newline $g_1 + \dots + g_k:\{1,2, \dots, n_1+
\dots + n_k \} \rightarrow \{1,2, \dots, n_1 '+ \dots + n_k' \}$ be
the function obtained by placing $g_1, \dots, g_k$ next to each
other from left to right. Then equivariance is preserved in the
sense that the diagram below must commute.

\begingroup
\vspace{-2\abovedisplayskip} \Small
$$\xymatrix@C=7pc@R=3pc{\Phi(\gamma(w,w_1, \dots,
w_k)_{g_1 + \dots + g_k}) \ar@{=}[d] \ar@{=>}[r]^{s_{\gamma(w,w_1,
\dots, w_k), g_1 + \cdots + g_k}} & \Phi(\gamma(w,w_1, \dots,
w_k))_{g_1+ \cdots +g_k} \ar@{=>}[d]|{(c_{w,w_1, \dots, w_k})_{g_1+
\cdots + g_k}}
\\ \Phi(\gamma(w,(w_1)_{g_1}, \dots, (w_k)_{g_k})) \ar@{=>}[d]|{c_{w,
(w_1)_{g_1}, \dots, (w_k)_{g_k}}} & \gamma(\Phi(w),\Phi(w_1),
\dots, \Phi(w_k))_{g_1 + \cdots + g_k} \ar@{=}[d] \\
\gamma(\Phi(w), \Phi((w_1)_{g_1}), \dots, \Phi((w_k)_{g_k}))
\ar@{=>}[r]_{\overset{\phantom{M}}{\gamma(i_{\Phi(w)},s_{w_1, g_1},
\dots, s_{w_k, g_k})}} & \gamma(\Phi(w), \Phi(w_1)_{g_1}, \dots,
\Phi(w_k)_{g_k})}$$
\endgroup
\noindent

\item
The substitution coherence isomorphisms are associative, \ie for
every word $w \in T(\ell)$ and functions $f:\{1, \ldots, \ell \}
\rightarrow \{1, \ldots, m\}$ and $g:\{1, \ldots,m\} \rightarrow
\{1, \ldots, n\}$ we mimic the equality $w_{g \circ f} = (w_f)_g $
by requiring the diagram below to commute. Here $(s_{w,f})_g$ is the
natural transformation which is defined for objects $A_1, \ldots,
A_n$ of $X$ by $(s_{w,f})_g(A_1, \ldots, A_n)= s_{w,f}(A_{g1},
\ldots, A_{gm})$.

\begingroup
\vspace{-2\abovedisplayskip} \small
$$\xymatrix@R=3pc@C=5pc{\Phi((w_f)_g)=\Phi(w_{g \circ f}) \ar@{=>}[r]^-{s_{w,g \circ
f}} \ar@{=>}[d]_{s_{(w_f),g}} & \Phi(w)_{g \circ f} \ar@{=}[d] \\
\Phi(w_f)_g \ar@{=>}[r]_{(s_{w,f})_g} & (\Phi(w)_f)_g}$$
\endgroup
\noindent

\item
For all $w \in T(k)$ and $id_k:\{1, \dots, k\} \rightarrow \{1,
\dots, k\}$ the natural transformation $s_{w,id_k}$ is the identity.
\end{enumerate}
\end{defn}

\begin{rmk}

One can compactly describe the concept of a pseudo
algebra\index{algebra!pseudo algebra} as follows. A category $X$ is
a pseudo $T$-algebra if it is equipped with a {\it pseudo morphism
of theories}\index{morphism!pseudo morphism of theories|textbf}
$\Phi:T \rightarrow End(X)$. The assignment $\Phi$ is pseudo in the
sense that the requirements of Lemma \ref{naturalmorphismoftheories}
are only satisfied up to coherence isos, namely the assignment
preserves $\gamma$ up to $c$, preserves the identity up $I$, and is
natural up to $s$ as in the diagrams below and these coherence isos
satisfy coherence diagrams.
$$\xymatrix@C=1.5pc@R=3pc{T(k) \times T(n_1)
\times \cdots \times T(n_k) \ar[r] \ar@{}@<1.5ex>[r]^-{\Phi(k)
\times \Phi(n_1) \times \cdots \times \Phi(n_k)} \ar[d]_{\gamma^T} &
End(X)(k) \times End(X)(n_1) \times \cdots \times End(X)(n_k)
\ar[d]^{\gamma^{End(X)}}
\\ T(n_1 + \cdots + n_k) \ar[r]_{\Phi(n_1 + \cdots + n_k)}
\ar@{=>}[ur]^c & End(X)(n_1 + \cdots + n_k) }$$
$$\xymatrix@C=3pc@R=3pc{X
\ar[d]_{\Phi(1)(1_T)} \ar@{=}[r] & X \ar[d]^{1_X} \ar@{}[d]^{\phantom{\Phi(1)(1_T)}} \\
X \ar@{=}[r] \ar@{=>}[ur]^I & X }$$
$$\xymatrix@C=4pc@R=3pc{T(m) \ar[r]^-{\Phi(m)} \ar[d]_{T(f)}
\ar@{}[d]_{\phantom{End(X)(f)}} & End(X)(m) \ar[d]^{End(X)(f)}
\\ T(n) \ar[r]_-{\Phi(n)} \ar@{=>}[ur]^{s_{-,f}} & End(X)(n)}$$
\end{rmk}

\begin{rmk} It is possible to describe the general form of these
coherence diagrams. In general, a relation\index{relation} $\alpha
\circ \beta = \alpha' \circ \beta'$ in the theory of
theories\index{theory!theory of theories} and a tuple $\bar{w}$ of
words gives rise to a coherence diagram
$$\xymatrix@C=4pc@R=3pc{\Phi(\alpha ' \circ \beta '(\bar{w}))
\ar@{=>}[r]^{\varepsilon_{\alpha'}(\beta'(\bar{w}))} \ar@{=}[d] &
\alpha'(\Phi(\beta'(\bar{w})))
\ar@{=>}[d]^{\alpha'(\varepsilon_{\beta'}(\bar{w}))}
\\ \Phi(\alpha \circ
\beta(\bar{w})) \ar@{=>}[d]_{\varepsilon_{\alpha}(\beta(\bar{w}))}
 & \alpha' \circ \beta'(\Phi(\bar{w})) \ar@{=}[d]
\\
\alpha(\Phi(\beta(\bar{w})))
\ar@{=>}[r]_{\alpha(\varepsilon_{\beta}(\bar{w}))} & \alpha \circ
\beta (\Phi (\bar{w})) }$$ where $\varepsilon_{\alpha},
\varepsilon_{\alpha'}, \varepsilon_{\beta}$, and
$\varepsilon_{\beta'}$ are the coherence isos associated to the
morphisms $\alpha$, $\alpha'$, $\beta$, and $\beta'$ respectively in
the theory of theories\index{theory!theory of theories} and
$\Phi(\bar{w})$ denotes the tuple of words obtained by applying
$\Phi$ to each of the constituents of $\bar{w}$. Note that
$\varepsilon_{\alpha}, \varepsilon_{\alpha'}, \varepsilon_{\beta}$,
and $\varepsilon_{\beta'}$ are tuples of the 2-cells $c,I,s$ and
identity 2-cells. In the definition of pseudo algebra above, the
morphisms $\beta, \beta'$ are tuples of generating morphisms in all
cases except in (4). In (4) the $\beta'$ is the result of applying a
substitution morphism\index{substitution morphism} in the theory of
theories to $\gamma$. This substitution morphism\index{substitution
morphism} can be written in terms of $f$ appropriately. In this case
we have $\varepsilon_{\beta'}(\bar{w})=c_{w,w_{f1}, \dots, w_{fk}}$.
\end{rmk}

\begin{defn} \label{defnlaxalgebramorphism}
Let $X$ and $Y$ be pseudo $T$-algebras and $H:X \rightarrow Y$ a
functor between the underlying categories. Denote the structure maps
of $X$ and $Y$ by $\Phi$ and $\Psi$ respectively. For all $n \in
\mathbb{N}$ and all $w \in T(n)$ let $\rho_w:H \circ \Phi(w)
\Rightarrow \Psi(w) \circ (H, \ldots, H)$ be a natural isomorphism.
Then $H$ is a {\it pseudo morphism of pseudo $T$-algebras with
coherence iso 2-cells\index{coherence isomorphism|)}\index{coherence
2-cell} $\rho_w$} (or just {\it morphism of pseudo
$T$-algebras}\index{morphism!morphism of pseudo
$T$-algebras|textbf}\index{morphism!pseudo morphism of pseudo
$T$-algebras|textbf} for short) if the following coherence
diagrams\index{coherence diagram|)} of natural isomorphisms are
satisfied.
\begin{enumerate}
\item
For all $k \in \mathbb{N}$, $w \in T(k)$, and all words $w_1, \dots,
w_k$ of $T$ the diagram below must commute.

\begingroup
\vspace{-2\abovedisplayskip} \normalsize
$$\xymatrix@C=2.6pc@R=3pc{H \circ \Phi(w \circ (w_1,\ldots,w_k)) \ar@{=>}[r]^-{i_H \ast
c_{w,w_1,\ldots,w_k}} \ar@{=>}[dd]_{\rho_{w \circ (w_1,\ldots,w_k)}}
& H \circ \Phi(w) \circ (\Phi(w_1), \ldots, \Phi(w_k))
\ar@{=>}[d]^{\rho_w \ast i_{(\Phi(w_1), \ldots, \Phi(w_k))}} \\   &
\Psi(w) \circ (H, \ldots, H) \circ (\Phi(w_1), \ldots, \Phi(w_k))
\ar@{=>}[d]^{i_{\Psi(w)} \ast ( \rho_{w_1}, \ldots, \rho_{w_k})}
 \\  \Psi(w \circ (w_1,\ldots,w_k)) \circ
(H, \dots, H)
\ar@{=>}[r]_{\overset{\phantom{l}}{c_{w,w_1,\ldots,w_k}\ast
 i_{(H, \ldots,H)}}} & \Psi(w) \circ
(\Psi(w_1), \ldots, \Psi(w_k)) \circ (H, \ldots, H) }$$
\endgroup
\noindent
\item
The diagram below must commute.

\begingroup
\vspace{-2\abovedisplayskip} \normalsize
$$\xymatrix@R=3pc@C=3pc{H \circ \Phi(1)
\ar@{=>}[r]^{i_H \ast I} \ar@{=>}[d]_{\rho_{1}} & H \circ 1_X
\ar@{=}[d] \\ \Psi(1) \circ H \ar@{=>}[r]_{I \ast i_H} & 1_Y \circ
H}$$
\endgroup
\noindent
\item
For every word $w \in T(m)$ and every function $f:\{1, \dots,m\}
\rightarrow \{1, \dots, n\}$ the diagram below must commute.

\begingroup
\vspace{-2\abovedisplayskip} \normalsize
$$\xymatrix@C=6pc@R=3pc{H \circ \Phi(w_f) \ar@{=>}[r]^{i_H \ast s_{w,f}}
\ar@{=>}[d]_{\rho_{w_f}} & H \circ \Phi(w)_f
\ar@{=>}[d]^{(\rho_w)_f} \\ \Psi(w_f) \circ (H, \dots, H)
\ar@{=>}[r]_-{s_{w,f}\ast i_{(H, \dots,H)}} & \Psi(w)_f \circ (H,
\dots, H)}$$
\endgroup
\noindent
\end{enumerate}
\end{defn}

\begin{examp}
Let $T$ be the theory of commutative monoids\index{commutative
monoid!theory of commutative monoids}\index{theory!theory of
commutative monoids} and let $FiniteSets$\index{$FiniteSets$} be the
category of finite sets and bijections\index{algebra!pseudo
algebra!examples of pseudo algebras}. Define $A \coprod B:=A \times
\{1\} \cup B \times \{2\}$ for finite sets $A$ and $B$. Define
coproduct similarly for morphisms of finite sets. Then
$\coprod:FiniteSets \times FiniteSets \rightarrow FiniteSets$ is a
functor which makes $FiniteSets$ into a pseudo $T$-algebra, \ie a
{\it pseudo commutative monoid}\index{commutative monoid!pseudo
commutative monoid|textbf}\index{commutative monoid!pseudo
commutative monoid!examples of pseudo commutative monoids}. More
generally, any symmetric monoidal category\index{symmetric monoidal
category} is a pseudo $T$-algebra.
\end{examp}

\begin{examp}
Let $T$ be the theory of commutative semi-rings\index{algebra!pseudo
algebra!examples of pseudo algebras}\index{commutative
semi-ring!theory of commutative semi-rings}\index{theory!theory of
commutative semi-rings}. Then the category of finite dimensional
complex vector spaces\index{vector space} is a pseudo $T$-algebra
whose structure is given by direct sum\index{direct sum} and tensor
product\index{tensor product}. We also say this category is a {\it
pseudo commutative semi-ring}\index{commutative semi-ring!pseudo
commutative semi-ring|textbf}\index{semi-ring|see{commutative
semi-ring}}\index{commutative semi-ring!pseudo commutative
semi-ring!example of pseudo commutative semi-ring}.
\end{examp}

\begin{defn}
Let $X,Y$, and $Z$ be pseudo $T$-algebras and $G:X \rightarrow Y,
H:Y \rightarrow Z$ morphisms of pseudo $T$-algebras with coherence
2-cells $\rho^G_w$ and $\rho^H_w$ respectively. Then the {\it
composition}\index{composition of morphisms of pseudo $T$-algebras}
$H \circ G$ is the composition of the underlying functors . It has
the coherence 2-cells\index{coherence 2-cell} $\rho_w^{H \circ
G}:=(\rho^H_w \ast i_{(G, \dots, G)}) \odot (i_H \ast \rho^G_w):H
\circ G \circ \Phi(w) \Rightarrow \Psi(w) \circ (H \circ G, \dots, H
\circ G)$ where $\Phi$ ad $\Psi$ denote the structure maps of $X$
and $Z$ respectively.
\end{defn}

\begin{lem}
The composition of morphisms of pseudo $T$-algebras is a morphism of
pseudo $T$-algebras.
\end{lem}
\begin{pf}
Immediate.
\end{pf}

\begin{defn}
Let $X$ and $Y$ be pseudo $T$-algebras with structure maps $\Phi$
and $\Psi$ respectively. Let $G,H:X \rightarrow Y$ be morphisms of
pseudo $T$-algebras. A natural transformation $\alpha:G \Rightarrow
H$ between the underlying functors is a {\it 2-cell in the
2-category of pseudo $T$-algebras}\index{2-cell!2-cell in the
2-category of pseudo $T$-algebras|textbf} if for all $n \in
\mathbb{N}$ and all $w \in T(n)$
$$\xymatrix@C=6pc@R=3pc{G \circ \Phi(w) \ar@{=>}[r]^{\alpha \ast i_{\Phi(w)}}
\ar@{=>}[d]_{\rho_w^G}
 & H \circ \Phi(w)
\ar@{=>}[d]^{\rho_w^H} \\ \Psi(w) \circ (G, \dots, G)
\ar@{=>}[r]_{i_{\Psi(w)} \ast (\alpha, \dots, \alpha)} & \Psi(w)
\circ(H, \dots,H) }$$ commutes. The vertical and horizontal
compositions of the 2-cells are just the vertical and horizontal
composition of the underlying natural transformations.
\end{defn}

\begin{lem}
The small pseudo $T$-algebras with morphisms and 2-cells defined
above form a 2-category.
\end{lem}
\begin{pf}
The axioms can be verified directly.
\end{pf}

Next we work towards a description of pseudo $T$-algebras as strict
algebras over a {\it 2-monad}\index{2-monad} $C$ by way of a theory
$\mathcal{T}$ enriched in groupoids\index{theory!theory enriched in
groupoids}\index{groupoid}. As mentioned in the last chapter, a
pseudo $T$-algebra is the same thing as a strict
$\mathcal{T}$-algebra. This was observed in \cite{hu2}. We can see
this as follows. Let $T'$ denote the free\index{theory!free
theory}\index{free theory} theory on the sequence of sets underlying
$T$. Recall that $T'$ was described in terms of the sets $T'(n)$ for
$n \geq 0$ and the compositions, substitutions, and identities. From
this description, the hom sets are $Mor_{T'}(m,n)=\prod_{j=1}^n
T'(m)$. There is a map of theories $T' \rightarrow T$ which gives
the theory structure on $T$. Let the underlying 1-category of the
2-category $\mathcal{T}$ be $T'$. For $v,w \in
\mathcal{T}(n)=Mor_{\mathcal{T}}(n,1)$ we define a unique iso 2-cell
between $v$ and $w$ if $v$ and $w$ map to the same element of $T(n)$
under the map of theories $T' \rightarrow T$. Otherwise there is no
2-cell between $v$ and $w$. With these definitions, the only 2-cell
between $w$ and $w$ is the identity and the vertical composition of
2-cells is uniquely defined. Thus $\mathcal{T}(n)$ is a category.
Next define $Mor_{\mathcal{T}}(m,n)$ to be the product category
$\prod_{j=1}^n \mathcal{T}(m)$ for all $m,n \in Obj \hspace{1mm}
\mathcal{T}$. From this it follows that there is a unique iso 2-cell
between $v,w \in Mor_{\mathcal{T}}(m,n)$ if they map to the same
element of $Mor_T(m,n)$ and otherwise there is no 2-cell. This
uniquely defines the horizontal composition of 2-cells and
$\mathcal{T}$ is a 2-category. From the definitions it also follows
easily that $n$ is the 2-product\index{2-product} of $n$ copies of
$1$ in $\mathcal{T}$. Hence $\mathcal{T}$ is a theory enriched in
groupoids\index{theory!theory enriched in
groupoids}\index{groupoid}. In \cite{hu2} $\mathcal{T}$ is denoted
$(Th(T),G(T))$.

We introduce the notation $c,I,s$\index{$c$}\index{$I$}\index{$s$}
for some of these 2-cells, which breaks the usual convention of
labelling 2-cells by lowercase Greek letters. Let
$$\xymatrix@1{c_{w,w_1, \dots, w_k}:(()_{id_{n_1+ \cdots + n_k}},
\gamma(w,w_1, \dots, w_k)) \ar@{=>}[r] & (\gamma,w,w_1, \dots,
w_k)}$$ denote the unique 2-cell for $w \in T(k), w_i \in T(n_i),
i=1, \dots, k$. The $\gamma$ on the right is a generator of the
theory of theories\index{theory!theory of theories} while the
$\gamma$ on the left is the composition in the theory $T$. The map
$id_{n_1 + \cdots + n_k}$ is the identity of the object $n_1 +
\cdots + n_k$ in the category $\Gamma$ of Definition
\ref{definitionofGamma}. Let
$$\xymatrix@1{I:(()_{id_1},1) \ar@{=>}[r] & (1,*)}$$ where
$(()_{id_1},1) \in \mathbf{R}_1(1) \times T(1)$ and $(1,*) \in
\mathbf{R}_1(1^0) \times T(1)^0$. Here $\mathbf{R}$ denotes the
theory of theories in Example \ref{theoryoftheoriesexample}. Let
$$\xymatrix@1{s_{w,f}:(()_{id_n},w_f) \ar@{=>}[r] &
(()_f,w)}$$ denote the unique 2-cell for $w \in T(m)$ and $f:m
\rightarrow n$ in $\Gamma$. We call these 2-cells as well as
identity 2-cells the {\it elementary 2-cells}\index{elementary
2-cell}. By the following inductive proof, every other 2-cell in
$\mathcal{T}$ can be obtained from these ones and their inverses.

\begin{lem} \label{basic2cell}
Let $\alpha$ be a word in the theory of theories, \ie $\alpha \in
\mathbf{R}_n(\bar{m})$ for some $n \in \mathbb{N}_0$,
$\bar{m}=(j_1^{m_1}, \dots, j_p^{m_p})$, and $m:=m_1+ \dots +m_p$.
Then the 2-cell
$$\xymatrix@1{(()_{id_n}, \alpha(v_1, \dots, v_m)) \ar@{=>}[r] &
(\alpha,v_1, \dots,v_m)}$$ in $\mathcal{T}$ can be expressed as a
vertical composition
$$ \sigma_s \odot \sigma_{s-1} \odot \cdots \odot \sigma_1$$
where each $\sigma_r$ is the result of applying a morphism in
$\mathbf{R}$ to a tuple of elementary 2-cells.
\end{lem}
\begin{pf}
Let $\alpha=\alpha_i \circ \cdots \circ \alpha_1$ where $\alpha_1,
\dots, \alpha_i$ are tuples of generating morphisms in the theory
$\mathbf{R}$ of theories such that $i$ is minimal. We induct on $i$.
If $i=1$, then $\alpha$ is a generating morphism for $\mathbf{R}$
and the 2-cell
$$\xymatrix@1{(()_{id_n}, \alpha(v_1, \dots, v_m)) \ar@{=>}[r] &
(\alpha,v_1, \dots,v_m)}$$ must be one of $c,I,$ or $s$. Now let $i
\geq 1$ and suppose the Lemma holds for all words that can be
expressed with $i$ terms or less. Suppose $\alpha \in
\mathbf{R}_n(\bar{m})$ has an expression with $i+1$ terms but not
does not have an expression with fewer terms. Then $\alpha=\beta
\circ (\beta_1, \dots, \beta_k)$ where $\beta$ is a generating
morphism for the theory of theories and $\beta_1, \dots, \beta_k$
are some words in the theory of theories, each with $i_1, \dots, i_k
\leq i$. Then the 2-cells
$$\aligned
& \xymatrix@1{\varepsilon_1:(()_{id},\beta_1(v_1, \dots))
\ar@{=>}[r] & (\beta_1, v_1, \dots)}  \\
& \xymatrix@1{\varepsilon_2:(()_{id},\beta_2(\dots)) \ar@{=>}[r] &
(\beta_2, \dots)}  \\  & \hspace{6pc} \vdots  \\
& \xymatrix@1{\varepsilon_k:(()_{id},\beta_k(\dots, v_m))
\ar@{=>}[r] & (\beta_k, \dots,v_m)}
\endaligned$$ can be obtained from elementary  2-cells
in the prescribed manner by the induction hypothesis. Here $id$ is
generically used to denote any identity morphism in $\Gamma$. Then
$$\xymatrix@R=3pc@C=1pc{(()_{id_n}, \alpha(v_1, \dots, v_m)) \ar@{=}[r] &
(()_{id_n}, \beta \circ (\beta_1, \dots, \beta_k)(v_1, \dots, v_m))
\ar@{=>}[d]^{e} \\ (\beta, \beta_1(w_1, \dots), \beta_2(\dots),
\dots, \beta_k(\dots,v_m)) \ar@{=}[r]
\ar@{=>}[d]_{\beta(\varepsilon_1, \dots, \varepsilon_k)} & (\beta,
(\beta_1, \dots, \beta_k)(v_1, \dots, v_m)) \\
(\beta,(\beta_1,w_1, \dots), (\beta_2, \dots), \dots, (\beta_k,
\dots, v_m)) \ar@{=}[r] & (\beta \circ (\beta_1, \dots, \beta_k),
v_1, \dots, v_m) \ar@{=}[d] \\ & (\alpha, v_1, \dots, v_m) }
$$ is also a composition of the prescribed type, where $e$ is an
elementary 2-cell.
\end{pf}

\begin{lem}
Let $\alpha$ and $\beta$ be words in the theory of
theories\index{theory!theory of theories}. Suppose that there is a
2-cell
$$\xymatrix@1{(\alpha,v_1, \dots, v_{m_1}) \ar@{=>}[r] & (\beta,w_1,
\dots, w_{m_2})}$$ in $\mathcal{T}$. Then this 2-cell is a vertical
composition of 2-cells obtained from elementary 2-cells and their
inverses by applying morphisms in the theory of
theories\index{elementary 2-cell}.
\end{lem}
\begin{pf}
From Lemma \ref{basic2cell} we have 2-cells
$$\xymatrix@R=3pc@C=3pc{(\alpha,v_1, \dots, v_{m_1}) \ar@{<=}[d] & (\beta,w_1, \dots,
w_{m_2})  \ar@{<=}[d] \\ (()_{id}, \alpha(v_1, \dots, v_{m_1}))
\ar@{=}[r] & (()_{id}, \beta(w_1, \dots, w_{m_2}))}$$ of the
prescribed type. We obtain the desired result by inverting the
2-cell on the left.
\end{pf}

\begin{thm}
There is a bijection between the set of small pseudo $T$-algebras
and the set of small
$\mathcal{T}$-algebras\index{algebra!$\mathcal{T}$-algebra|textbf}\index{algebra!pseudo
$T$-algebra|textbf}.
\end{thm}
\begin{pf}
Let $(X, \Phi)$ be a small pseudo $T$-algebra. Define a morphism
$\Psi: \mathcal{T} \rightarrow End(X)$ of theories enriched in
groupoids\index{theory!theory enriched in groupoids}\index{groupoid}
by the following sequence of functors $\Psi_n:\mathcal{T}(n)
\rightarrow End(X)(n)$. For notational convenience, the subscript
$n$ is usually left off below. For $(\alpha, w_1, \dots, w_{\ell})
\in \mathcal{T}(n)$ define
$$\Psi(\alpha,w_1, \dots, w_{\ell}):=\alpha(\Phi(w_1), \dots,
\Phi(w_{\ell})).$$ For elementary 2-cells, define
$$\Psi(c_{w,w_1, \dots, w_k}):=c_{w,w_1, \dots, w_k}$$
$$\Psi(I):=I$$
$$\Psi(s_{w,f}):=s_{w,f}$$
where the symbols on the right denote the coherence natural
isomorphisms from the pseudo $T$-algebra structure.

If $\alpha$ is a word in the theory of theories\index{theory!theory
of theories} and $\varepsilon_1, \dots, \varepsilon_k$ are
elementary 2-cells, then
$$\Psi(\alpha(\varepsilon_1, \dots,
\varepsilon_k)):=\alpha(\Psi(\varepsilon_1), \dots,
\Psi(\varepsilon_k)).$$ This is well defined, because if
$\alpha(\varepsilon_1, \dots, \varepsilon_k)=\beta(\varepsilon_1,
\dots, \varepsilon_k)$ with $\varepsilon_1, \dots, \varepsilon_k$
elementary, then $\alpha=\beta$.

Consider the 2-cell $$\xymatrix@1{(()_{id_n}, \alpha(v_1, \dots,
v_m)) \ar@{=>}[r] & (\alpha ,v_1, \dots, v_m)}$$ for some $\alpha
\in \mathbf{R}_n(\bar{m})$. By the above lemma, the word $\alpha$
can be expressed in the form $\sigma_s \odot \cdots \odot \sigma_1$
where each $\sigma_r$ is obtained from a tuple of elementary 2-cells
by applying a morphism in $\mathbf{R}$. Define
$$\Psi(\sigma_s \odot \cdots \odot \sigma_1):=\Psi(\sigma_s) \odot
\cdots \odot \Psi(\sigma_1)$$ where each $\Psi(\sigma_r)$ is defined
as in the previous paragraph. To see that this is well defined,
suppose $\sigma_s \odot \cdots \odot \sigma_1=\sigma_{s'}' \odot
\cdots \odot \sigma_1'$ where each $\sigma_{r'}'$ is obtained from a
tuple of elementary 2-cells by applying a morphism in $\mathbf{R}$.
Such a sequence gives rise to an expression $\alpha=\alpha_{s'}'
\circ \cdots \circ \alpha_1'$ where $\alpha_1', \dots, \alpha_{s'}'$
are tuples of generating morphisms. Let $\alpha=\alpha_s \circ
\cdots \circ \alpha_1$ be the expression that arose from $\sigma_s
\odot \cdots \odot \sigma_1$. It suffices to consider the case
$$\alpha=\alpha_4 \circ \alpha_3 \circ \alpha_2 \circ
\alpha_1=\alpha_4 \circ \alpha_3' \circ \alpha_2' \circ \alpha_1$$
with $\alpha_3 \circ \alpha_2=\alpha_3' \circ \alpha_2'$ because
$\alpha_{s'}' \circ \cdots \circ \alpha_1'$ can be obtained from
$\alpha_s \circ \cdots \circ \alpha_1$ by a finite number of
applications of the relations in the theory of theories. Then we
have the following diagram, whose vertical columns are
$\Psi(\sigma_4 \odot \sigma_3 \odot \sigma_2 \odot \sigma_1)$ and
$\Psi(\sigma_4 \odot \sigma_3' \odot \sigma_2' \odot \sigma_1)$
respectively.
$$\xymatrix@R=3pc@C=3pc{\Phi(\alpha_4 \circ \alpha_3 \circ \alpha_2 \circ \alpha_1(\bar{w}))
\ar@{=}[r] \ar@{=>}[d]_{\varepsilon_4(\alpha_3 \circ \alpha_2 \circ
\alpha_1(\bar{w}))} & \Phi(\alpha_4 \circ \alpha_3' \circ \alpha_2'
\circ \alpha_1(\bar{w})) \ar@{=>}[d]^{\varepsilon_4(\alpha_3' \circ
\alpha_2'
\circ \alpha_1(\bar{w}))} \\
\alpha_4 \Phi(\alpha_3 \circ \alpha_2 \circ \alpha_1(\bar{w}))
\ar@{=}[r] \ar@{=>}[d]_{\alpha_4(\varepsilon_3(\alpha_2 \circ
\alpha_1(\bar{w})))} & \alpha_4 \Phi(\alpha_3' \circ \alpha_2' \circ
\alpha_1(\bar{w})) \ar@{=>}[d]^{\alpha_4(\varepsilon_3'(\alpha_2'
\circ \alpha_1(\bar{w})))}
\\ \alpha_4 \circ \alpha_3 \Phi(\alpha_2 \circ \alpha_1(\bar{w}))
\ar@{=>}[d]_{\alpha_4 \circ
\alpha_3(\varepsilon_2(\alpha_1(\bar{w})))} & \alpha_4 \circ
\alpha_3' \Phi(\alpha_2' \circ \alpha_1(\bar{w}))
\ar@{=>}[d]^{\alpha_4 \circ
\alpha_3'(\varepsilon_2'(\alpha_1(\bar{w})))} \\
 \alpha_4
\circ \alpha_3 \circ \alpha_2 \Phi(\alpha_1(\bar{w})) \ar@{=}[r]
\ar@{=>}[d]_{\alpha_4 \circ \alpha_3 \circ
\alpha_2(\varepsilon_1(\bar{w}))} & \alpha_4 \circ \alpha_3' \circ
\alpha_2' \Phi(\alpha_1(\bar{w})) \ar@{=>}[d]^{{\alpha_4 \circ
\alpha_3' \circ \alpha_2'(\varepsilon_1(\bar{w}))}}
\\ \alpha_4
\circ \alpha_3 \circ \alpha_2 \circ \alpha_1\Phi(\bar{w}) \ar@{=}[r]
& \alpha_4 \circ \alpha_3' \circ \alpha_2' \circ
\alpha_1\Phi(\bar{w})}$$ Here $\varepsilon_i$ denotes the tuple of
elementary 2-cells needed to bring $\alpha_i$ past $\Phi$. The inner
square commutes because of the coherence diagrams. The top and
bottom squares commute because $\alpha_3 \circ \alpha_2 = \alpha_3'
\circ \alpha_2'$. Hence $$\Psi(\sigma_4 \odot \sigma_3 \odot
\sigma_2 \odot \sigma_1)=\Psi(\sigma_4 \odot \sigma_3' \odot
\sigma_2' \odot \sigma_1)$$ and $\Psi$ is well defined on any 2-cell
of the form
$$\xymatrix@1{(()_{id_n}, \alpha(v_1, \dots,
v_m)) \ar@{=>}[r] & (\alpha ,v_1, \dots, v_m)}.$$

Next we must define $\Psi$ on 2-cells of the form
$$\xymatrix@1{(\alpha ,v_1, \dots, v_{m_1}) \ar@{=>}[r] &
(\beta ,w_1, \dots, w_{m_2})}.$$ According to Lemma \ref{basic2cell}
we have 2-cells
$$\xymatrix@R=3pc@C=3pc{(\alpha,v_1, \dots, v_{m_1}) \ar@{<=}[d]_{\mu} &
(\beta,w_1, \dots, w_{m_2})  \ar@{<=}[d]_{\nu} \\ (()_{id},
\alpha(v_1, \dots, v_{m_1})) \ar@{=}[r] & (()_{id}, \beta(w_1,
\dots, w_{m_2}))}$$ on which $\Psi$ is already defined. Define
$$\Psi(\nu \odot \mu^{-1}):=\Psi(\nu) \odot \Psi(\mu)^{-1}.$$ To
see that this is well defined, suppose $$\xymatrix@1{\sigma_s \odot
\cdots \odot \sigma_1: (\alpha ,v_1, \dots, v_{m_1}) \ar@{=>}[r] &
(\beta ,w_1, \dots, w_{m_2}) }$$ is another expression where each
$\sigma_r$ is obtained by applying a morphism in $\mathbf{R}$ to a
tuple of elementary 2-cells or their inverses. Then $$\aligned
\Psi(\nu)= & \Psi(\sigma_s \odot \cdots \odot \sigma_1 \odot \mu) \\
\Psi(\nu)= & \Psi(\sigma_s \odot \cdots \odot \sigma_1) \odot
\Psi(\mu) \\
\Psi(\nu) \odot \Psi(\mu^{-1}) = & \Psi(\sigma_s \odot \cdots
\odot \sigma_1) \\
\Psi(\nu \odot \mu^{-1})= & \Psi(\sigma_s \odot \cdots \odot
\sigma_1)
\endaligned$$ and $\Psi$ is well defined on 2-cells.

By construction $\Psi_n:\mathcal{T}(n) \rightarrow End(X)(n)$ is a
functor and it preserves $\gamma, ()_g,$ and $(1,*)=1$. Hence $X$ is
a $\mathcal{T}$-algebra with structure maps given by $\Psi$. This
procedure $\Phi \mapsto \Psi$ defines a map $$\text{Pseudo
$T$-Algebras} \rightarrow \text{$\mathcal{T}$-Algebras}.$$

Now we define a map
$$ \text{$\mathcal{T}$-Algebras}\rightarrow \text{Pseudo $T$-Algebras} .$$ Let $(X, \Psi)$ be a $\mathcal{T}$-algebra.
Then define natural isomorphisms
$$c_{w,w_1, \dots, w_k}:=\Psi(c_{w,w_1, \dots, w_k})$$
$$I:=\Psi(I)$$
$$s_{w,f}:=\Psi(s_{w,f})$$
where the symbols $c,I,s$ on the right are 2-cells in $\mathcal{T}$.
Also define $$\Phi_n(w):=\Psi_n(()_{id_n},w)$$ for $w \in T(n)$.
Then the coherence diagrams are satisfied because
$\Psi_n:\mathcal{T}(n) \rightarrow End(X)(n)$ is a functor for every
$n$ and $\Psi$ preserves $\gamma,()_g,$ and $1$.

We can easily check that the two procedures are inverse to one
another and that they define a bijection.
\end{pf}

Next we can define a 2-monad\index{2-monad|textbf} $C:Cat
\rightarrow Cat$ like on page \pageref{2monad}. Define a 2-functor
$C$ by \label{2monadprecisely}
$$ CX:=\frac{(\mathop{\bigcup}_{n \geq 0} (\mathcal{T}(n) \times
X^n))}{\Gamma}$$ for any small category $X$. We can similarly define
2-natural transformations $\eta:1_{Cat} \Rightarrow C$ and $\mu:C^2
\Rightarrow C$.

\begin{thm} \label{laxTalgebras=strictCalgebras}
Let $\mathcal{C}_C$ denote the 2-category of small strict
$C$-algebras\index{algebra!$C$-algebra|textbf}, pseudo morphisms,
and 2-cells. Let $\mathcal{C}_T$ denote the 2-category of small
pseudo $T$-algebras. Then $\mathcal{C}_C$ and $\mathcal{C}_T$ are
2-equivalent.
\end{thm}
\begin{pf}
The small $C$-algebras are precisely the small
$\mathcal{T}$-algebras by a proof similar to Theorem \ref{C=T}. But
by the previous theorem, the small $\mathcal{T}$-algebras are
precisely the pseudo $T$-algebras. To see that the morphisms of the
2-categories $\mathcal{C}_C$ and $\mathcal{C}_T$ are the same, one
must compare the coherence isos of the morphisms. They are related
by $$\rho_{(\alpha,w_1, \dots, w_k) \times (\bar{x})}^C =
\alpha(\rho_{w_1}^T, \dots, \rho_{w_k}^T)(\bar{x}).$$ In diagram (1)
of Definition \ref{defnlaxalgebramorphism} the right vertical
composition can be replaced by the appropriate component of $\rho^C$
by the composition coherence diagram for coherence isos of pseudo
morphisms of $C$-algebras. Then (1) commutes by naturality of
$\rho^C$. In (2) of Definition \ref{defnlaxalgebramorphism}, the
right vertical equality can be replaced by the appropriate component
of $\rho^C$ by the unit coherence diagram for coherence isos of
pseudo morphisms of $C$-algebras. Then (2) commutes by the
naturality of $\rho^C$. Diagram (3) commutes by the naturality of
$\rho^C$. The 2-cells of the 2-categories $\mathcal{C}_C$ and
$\mathcal{C}_T$ are also the same.

Finally, the 2-equivalence of Theorem \ref{star2} yields the desired
2-equivalence.
\end{pf}

Power's Theorem 5.3 in \cite{power3} states that the 2-category of
strict $C$-algebras, pseudo morphisms and 2-cells is biequivalent to
the 2-category of strict $\mathcal{T}$-algebras, pseudo morphisms,
and 2-cells where $\mathcal{T}$ is a theory\index{theory!theory
enriched in categories} enriched in categories and $C$ is the
corresponding 2-monad in his construction. Power's theorem differs
from the above Theorem \ref{laxTalgebras=strictCalgebras} in several
regards. Theorem \ref{laxTalgebras=strictCalgebras} above uses
strict $C$-algebras to describe pseudo $T$-algebras, where $T$ is a
usual theory. Theorem \ref{laxTalgebras=strictCalgebras} also has a
2-equivalence rather than a biequivalence.

Theorem \ref{blackwellkellypower} states part of Theorem 2.6 from
\cite{blackwell}.

\begin{thm} \label{blackwellkellypower}
(Blackwell, Kelly, Power) Let $C$ be a
2-monad\index{2-monad|textbf}. Then the 2-category of small strict
$C$-algebras, pseudo morphisms, and 2-cells of pseudo morphisms
admits strictly weighted pseudo limits of strict
2-functors.\index{limit!weighted pseudo
limit|textbf}\index{weighted|textbf}
\end{thm}

We conclude the following completeness theorem from
\ref{blackwellkellypower}.

\begin{thm} \label{indexedlimitsstrict2functors}
Let $T$ be a theory. Then the 2-category of pseudo $T$-algebras
admits strictly weighted pseudo limits of strict
2-functors.\index{limit!weighted pseudo
limit|textbf}\index{weighted|textbf}
\end{thm}
\begin{pf}
A 2-equivalence of 2-categories preserves weighted pseudo limits
because it admits a left 2-adjoint. Then the result follows from the
previous two theorems.
\end{pf}

\chapter{Weighted Pseudo Limits in the 2-Category of Pseudo $T$-Algebras}
\label{sec:laxlimitsinlaxalgebras} In this chapter we show that the
2-category of pseudo $T$-algebras introduced in Chapter
\ref{sec:laxTalgebras} admits weighted\index{weighted} pseudo
limits\index{limit!weighted pseudo limit}. In Chapter
\ref{sec:laxlimitsinCat} we proved that the 2-category of small
categories admits weighted pseudo limits in Theorem
\ref{catlaxlimits}, Lemma \ref{catcotensor}, and Theorem
\ref{catweightedlaxlimits}. We modify the proofs in Chapter
\ref{sec:laxlimitsinCat} to obtain Theorem \ref{alglaxlimits}, Lemma
\ref{cotensorproductsT}, and Theorem \ref{weightedlaxlimitsT}.  Let
$\mathcal{C}$ denote the 2-category of small pseudo $T$-algebras in
this chapter. The existence of cotensor products\index{cotensor
product} in $\mathcal{C}$ allows us to conclude in Theorem
\ref{weightedlaxlimitsT} that $\mathcal{C}$ admits weighted pseudo
limits from a theorem of Street. This result is more general than
Theorem \ref{indexedlimitsstrict2functors} because it allows the
functors to be pseudo. The proof in this chapter for pseudo limits
is also constructive, whereas Theorem
\ref{indexedlimitsstrict2functors} is not.

\begin{thm} \label{alglaxlimits}
The 2-category $\mathcal{C}$ of small pseudo $T$-algebras admits
pseudo limits.\index{algebra!pseudo algebra!pseudo limits of pseudo
algebras|textbf}\index{limit!pseudo limit|textbf}
\end{thm}
\begin{pf}
Let $\mathcal{J}$ be a small 1-category and $F:\mathcal{J}
\rightarrow \mathcal{C}$ a pseudo functor. Let $\mathbf{1}$ denote
the terminal object\index{terminal object} of the 2-category of
small categories as in Theorem \ref{catlaxlimits}. Let $U$ denote
the forgetful 2-functor from the 2-category $\mathcal{C}$ of pseudo
$T$-algebras to the 2-category of small categories. The candidate
for the pseudo limit of $F$ is $L:=PseudoCone(\mathbf{1}, U \circ
F)$ as before. Note that these are pseudo cones into the 2-category
of small categories, not into the 2-category of pseudo $T$-algebras.
We define $\pi:\Delta_L \Rightarrow F$ as in Theorem
\ref{catlaxlimits}. We must show that $L$ has the structure of a
pseudo $T$-algebra, that $\pi$ is a pseudo natural transformation to
$F$, and that $L$ and $\pi$ are universal. These proofs will draw on
the analogous results for the pseudo limit of $U \circ F$.

\begin{lem} \label{limitisalgebra}
The small category $L$ admits a pseudo $T$-algebra structure.
\end{lem}
\begin{pf}
We first make the identification of the categories $P$ and $L$ as in
Remarks \ref{PL1} and \ref{PL2}. Let $\eta^{\ell}=(a_i^{\ell})_i
\times (\varepsilon_f^{\ell})_f \in Obj \hspace{1mm}  L$ and
$(\xi_i^{\ell})_i \in Mor \hspace{1mm} L$ for $1 \leq \ell \leq n$
and $w \in T(n)$. We denote the structure maps of the pseudo
$T$-algebra $Fi=A_i$ by $\Phi_i$ for all $i \in Obj \hspace{1mm}
\mathcal{J}$. Let $a_i :=\Phi_i(w)(a_i^1, \dots, a_i^n)$ and
$\varepsilon_f := \Phi_{Tf}(w)(\varepsilon_f^1, \dots,
\varepsilon_f^n) \circ \rho^{Ff}_w(a_{Sf}^1, \dots, a_{Sf}^n):
Ff(a_{Sf}) \rightarrow a_{Tf}$ as well as $\xi_i:=\Phi_i(w)(\xi_i^1,
\dots, \xi_i^n)$. Then the structure maps of the pseudo $T$-algebra
$L$ are defined by
 $\Phi(w)(\eta^1, \dots, \eta^n):=(a_i)_i \times (\varepsilon_f)_f$ and
$\Phi(w)((\xi_i^1)_i, \dots, (\xi_i^n)_i):=(\xi_i)_i$. We must
verify that these outputs belong to $L$.

We claim that $(a_i)_i \times (\varepsilon_f)_f \in Obj \hspace{1mm}
L$. We prove this by verifying the coherences in Remarks \ref{PL1}
and \ref{PL2} for a fixed word $w \in T(2)$. To avoid cumbersome
notation, we write $+$ for $\Psi(w)$ for any structure map $\Psi$.
The verification for a general word is the same. We abbreviate
$\rho^H_w$ as $\rho^H$ for any morphism $H$ of pseudo $T$-algebras.
The only word appearing in the following diagrams is $w$, so there
is no ambiguity. Let $\gamma_{f,g}:=\gamma_{f,g}^F$ and
$\delta_j:=\delta_j^F$. First we show that for all $j \in Obj
\hspace{1mm} \mathcal{J}$ the diagram.
\begin{equation} \label{1}
\xymatrix@R=3pc@C=4pc{a_j \ar[r]^-{\delta_{j\ast}(a_j)}
\ar[rd]_{1_{a_j}} & F1_j(a_j)
\ar[d]^{\varepsilon_{1_j}} \\
& a_j}
\end{equation}
commutes where $a_j=a_j^1+a_j^2$ and
$\varepsilon_{1_j}=(\varepsilon_{1_j}^1+\varepsilon_{1_j}^2) \circ
\rho^{F1_j}(a_j^1,a_j^2)$ as defined above. After writing this
diagram out we get
$$\xymatrix@C=9pc@R=3pc{a_j^1+a_j^2
\ar[d]_{\rho^{1_{Fj}}(a_j^1,a_j^2)=1_{a_j^1+a_j^2}}
\ar[r]^{\delta_{j\ast}(a_j^1+a_j^2)} &
 F1_j(a_j^1+a_j^2) \ar[d]^{\rho^{F1_j}(a_j^1,a_j^2)} \\
a_j^1+a_j^2 \ar[d]_{1_{a_j^1+a_j^2}} \ar[r]^{\delta_{j\ast}(a_j^1) +
\delta_{j\ast}(a_j^2)}
& F1_j(a_j^1)+F1_j(a_j^2) \ar[d]^{\varepsilon^1_{1_j} + \varepsilon^2_{1_j}}\\
a_j^1+a_j^2 \ar[r]_{1_{a_j^1+a_j^2}} & a_j^1+a_j^2 }$$ where the top
horizontal arrow is $\delta_{j\ast}(a_j)$ and the right vertical
composition is $\varepsilon_{1_j}$ by definition. The top square
commutes because $\delta_{j\ast}:1_{Fj} \Rightarrow F1_j$ is a
2-cell in the 2-category $\mathcal{C}$. The bottom square commutes
because $+$ is a functor and $\varepsilon_{1_j}^{\ell} \circ
\delta_{j\ast}(a_j^{\ell})=1_{a_j^{\ell}}$ for $\ell=1,2$. Hence
(\ref{1}) commutes. Next we show that for all $\xymatrix@1{i
\ar[r]^f & j \ar[r]^g & k}$ in $\mathcal{J}$ the diagram
\begin{equation} \label{2}
\xymatrix@R=3pc@C=9pc{Fg \circ Ff (a_i) \ar[r]^{\gamma_{f,g}(a_i)}
\ar[d]_{Fg(\varepsilon_f)} & F(g \circ f)(a_i)
\ar[d]^{\varepsilon_{g \circ f}}\\ Fg(a_j) \ar[r]_{\varepsilon_g} &
a_k}
\end{equation}
commutes where $\varepsilon_f=\varepsilon_f^1 + \varepsilon_f^2$
\etc After writing out this diagram we get the diagram below whose
outermost square is (\ref{2}). The upper left triangle commutes by
the definition of composition for morphisms of pseudo $T$-algebras.
The upper right quadrilateral commutes because $\gamma_{f,g}:Fg
\circ Ff \Rightarrow F(g \circ f)$ is a 2-cell in the 2-category of
pseudo $T$-algebras. The lower left square commutes because
$\rho^{Fg}:Fg(+) \Rightarrow Fg+Fg$ is a natural transformation. The
bottom right square commutes because $+$ is a functor and
$\varepsilon_g^{\ell} \circ (Fg(\varepsilon_f^{\ell}))
=\varepsilon^{\ell}_{g \circ f} \circ \gamma_{f,g}(a_{i}^{\ell})$
for $\ell=1,2$. Thus all four inner diagrams commute and (\ref{2})
commutes. Thus both coherences in Remark \ref{PL1} are satisfied and
$\eta^1 + \eta^2 = (a_i)_i \times (\varepsilon_f)_f$ is an object of
$L$.

\begingroup
\vspace{-2\abovedisplayskip} \small
$$\xymatrix@R=8pc@C=1pc{Fg \circ Ff(a_i^1+a_i^2)
\ar[rr]^{\gamma_{f,g}(a_i^1 + a_i^2)}
\ar[d]_{Fg(\rho^{Ff}(a_i^1,a_i^2))} \ar[rd]^{\hspace{.05in} \rho^{Fg
\circ Ff}(a_i^1, a_i^2)} & & F(g \circ f)(a_i^1+a_i^2)
\ar[d]^{\rho^{F(g \circ f)}(a_i^1,a_i^2)}
\\ Fg(Ff(a_i^1)+Ff(a_i^2))
\ar[r]_{\overset{\phantom{l}}{\rho^{Fg}(Ff(a_i^1),Ff(a_i^2))}}
\ar[d]_{Fg(\varepsilon_j^1 + \varepsilon_j^2)} &  Fg \circ
Ff(a_i^1)+ Fg \circ Ff(a_i^2)
\ar[r]_{\overset{\phantom{l}}{\gamma_{f,g}(a_i^1)+
\gamma_{f,g}(a_i^2)}}
\ar[d]^{Fg(\varepsilon_f^1)+Fg(\varepsilon_f^2)} & F(g \circ
f)(a_i^1)+F(g \circ f)(a_i^2) \ar[d]^{\varepsilon^1_{g \circ f} +
\varepsilon^2_{g \circ f} }
\\ Fg(a_j^1+a_j^2) \ar[r]_{\rho^{Fg}(a^1_j, a^2_j)}
& Fg(a_j^1)+Fg(a_j^2) \ar[r]_{\varepsilon^1_g + \varepsilon^2_g} &
a_k^1+a_k^2}$$
\endgroup
\noindent

 We claim that $(\xi_i^1)_i+(\xi_i^2)=(\xi_i)_i$ is a
morphism in $L$ where $(\xi_i^1)_i: (a_i^1)_i \times
(\varepsilon_f^1)_f \rightarrow (b_i^{1})_i \times (\zeta_f^1)_f$
and $(\xi_i^2)_i: (a_i^2)_i \times (\varepsilon_f^2)_f \rightarrow
(b_i^2)_i \times (\zeta_f^2)_f$ are morphisms in $L$. In other words
we must show that
\begin{equation} \label{3}
\xymatrix@R=3pc@C=3pc{Ff(a_i) \ar[r]^-{\varepsilon_f}
\ar[d]_{Ff(\xi_i)} & a_j \ar[d]^{\xi_j} \\ Ff(b_i) \ar[r]_-{\zeta_f}
& b_j}
\end{equation}
commutes for all morphisms $f:i \rightarrow j$ in $\mathcal{J}$,
where $a_i=a_i^1 +a_i^2$ \etc  If we write out the diagram we get
$$\xymatrix@R=5pc@C=5pc{Ff(a_i^1+a_i^2) \ar[r]^-{\rho^{Ff}(a_i^1,a_i^2)}
\ar[d]_{Ff(\xi_i^1+\xi_i^2)} & Ff(a_i^1)+Ff(a_i^2)
\ar[r]^-{\varepsilon_f^1 + \varepsilon_f^2} \ar[d]_{Ff(\xi^1_i) +
Ff(\xi_i^2)} & a_j^1+a_j^2 \ar[d]^{\xi_j^1 + \xi_j^2}
\\ Ff(b_i^1 +b_i^2) \ar[r]_-{\rho^{Ff}(b_i^1,b_i^2)}
& Ff(b_i^1) + Ff(b_i^2) \ar[r]_-{\zeta_f^1+\zeta_f^2} & b_j^1 +b_j^2
}$$ where the outermost square is (\ref{3}).  The square on the left
commutes because $\rho^{Ff}: Ff(+) \Rightarrow Ff + Ff$ is a natural
transformation. The right square commutes because the diagram
$$\xymatrix@R=3pc@C=3pc{Ff(a_i^{\ell}) \ar[r]^-{\varepsilon_f^{\ell}}
\ar[d]_{Ff(\xi_i^{\ell})} & a_j \ar[d]^{\xi_j^{\ell}} \\ Ff(b_i)
\ar[r]_-{\zeta_f^{\ell}} & b_j^{\ell}}$$ commutes for $\ell=1,2$ and
because $+$ is a functor. Hence $(\xi_i^1)_i+(\xi_i^2)=(\xi_i)_i$ is
a morphism in $L$. Thus $\Phi(w):L \times L \rightarrow L$.

The map $\Phi(w)$ preserves compositions and identities because the
individual components do. Thus $\Phi(w):L \times L \rightarrow L$ is
a functor. The same argument works for words in $T(n)$ for all $n
\in \mathbb{N}$. Thus $\Phi$ defines structure maps to make the
small category $L$ into a pseudo $T$-algebra.

We define the coherence isos\index{coherence isomorphism} for $\Phi$
to be those maps which have the coherence isos of $\Phi_i$ in the
$i$-th component. We can prove that they are morphisms of the
category $L$, \ie satisfy the diagram in Remark \ref{PL2}, by using
the coherence diagrams of $\rho$ with the respective coherence iso
as well as the naturality of the individual components.  The
coherence isos for $\Phi$ are natural because they are natural in
each component. The coherence isos for $\Phi$ satisfy the coherence
diagrams because the individual components do. Thus $L$ is a pseudo
$T$-algebra with structure maps $\Phi$.

\end{pf}

\begin{lem}
The map $\pi:\Delta_L \Rightarrow F$ is a pseudo natural
transformation with coherence iso 2-cells given by $\tau$.
\end{lem}
\begin{pf}
It is clear from the work on the small category case in Chapter
\ref{sec:laxlimitsinCat} that $\pi$ is a pseudo natural
transformation when we forget all the pseudo $T$-algebra structures.
Therefore it suffices to show that $\pi_j:L \rightarrow Fj$ is a
morphism of pseudo $T$-algebras for all $j \in Obj \hspace{1mm}
\mathcal{J}$ and that $\tau_{i,j}(f):Ff \circ \pi_i \Rightarrow
\pi_j$ is a 2-cell in the 2-category of pseudo $T$-algebras for all
morphisms $f:i \rightarrow j$ in $\mathcal{J}$.

Let $j \in Obj \hspace{1mm}  \mathcal{J}$. Then $\pi_j:L \rightarrow
Fj$ is a functor. We abbreviate $\Phi(w)$ for $w \in T(2)$ by $+$ as
above. Then for $\eta^{\ell}=(a_i^{\ell})_i \times
(\varepsilon_f^{\ell})_f \in Obj \hspace{1mm}  L$ for $\ell=1,2$ we
have $$\aligned \pi_j(\eta^1 + \eta^2) &= \pi_j((a_i^1+a_i^2)_i
\times ((\varepsilon_f^1 + \varepsilon_f^2) \circ
\rho^{Ff}_w(a_{Sf}^1,a_{Sf}^2))_f)
\\ &= a_j^1+a_j^2
\\ &= \pi_j(\eta^1)+\pi_j(\eta^2). \endaligned$$ The same calculation
works for words in $T(n)$ for all $n \in \mathbb{N}$. We conclude
that $\pi_j$ commutes with the structure maps for the pseudo
$T$-algebra structure. If we take $\rho_w^{\pi_j}=i_{\pi_j} \circ
i_{\Phi(w)}$ then $\pi_j$ is a morphism of pseudo $T$-algebras for
all $j \in \mathcal{J}$.

Let $f:i \rightarrow j$ be a morphism in $L$. To show that
$\tau_{i,j}(f)$ is a 2-cell, we must show that the diagram
\begin{equation} \label{4}
\xymatrix@C=9pc@R=4pc{Ff \circ \pi_i \circ \Phi(w)
\ar@{=>}[r]^{\tau_{i,j}(f) \ast
 i_{\Phi(w)}} \ar@{=>}[d]_{\rho_w^{Ff \circ \pi_i}}
 & \pi_j \circ \Phi(w)
\ar@{=>}[d]^{\rho_w^{\pi_j}}
\\ \Phi_j(w) \circ (Ff \circ \pi_i, \dots, Ff \circ \pi_i)
\ar@{=>}[r]_-{i_{\Phi_j(w)} \ast (\tau_{i,j}(f), \dots,
\tau_{i,j}(f))} & \Phi_j(w)
 \circ (\pi_j, \dots,\pi_j) }
\end{equation}
commutes for all words $w$. Recalling that
$\tau_{i,j}(f)_{\eta}:=\tau^{\eta}_{i,j}(f)$ and evaluating the
diagram on $(\eta^1,\eta^2)$ where $\eta^{\ell}=(a_i^{\ell})_i
\times (\varepsilon_f^{\ell})_f \in Obj \hspace{1mm}  L$ for $\ell
=1,2$ gives
$$\xymatrix@C=14pc@R=4pc{Ff(a_i^1+a_i^2)
\ar[d]_{\rho^{Ff}(a_i^1,a_i^2)}
 \ar[r]^-{(\varepsilon_f^1 + \varepsilon_f^2) \circ
\rho^{Ff}(a_i^1,a_i^2)} & a_j^1+a_j^2 \ar[d]^{1_{a_j^1+a_j^2}}
\\  Ff(a_i^1)+Ff(a_i^2) \ar[r]_-{\varepsilon_f^1+\varepsilon_f^2}
& a_j^1+a_j^2}$$ which obviously commutes. Hence $\tau_{i,j}(f)$ is
a 2-cell in the 2-category of pseudo $T$-algebras for all $f:i
\rightarrow j$ and $\pi$ is a pseudo natural transformation.

\end{pf}

Now we must show that the pseudo $T$-algebra $L$ and the pseudo
natural transformation $\pi:\Delta_L \Rightarrow F$ are universal in
the sense that the functor $\phi:Mor_{\mathcal{C}}(V,L) \rightarrow
PseudoCone(V,F)$ as defined in the small category case of Chapter
\ref{sec:laxlimitsinCat} is an isomorphism of categories for all
objects $V$ of $\mathcal{C}$. In the following, $V$ is a fixed
object of the 2-category $\mathcal{C}$ of pseudo $T$-algebras.

\begin{lem}
The map $\phi:Mor_{\mathcal{C}}(V,L) \rightarrow PseudoCone(V,F)$ is
a functor.
\end{lem}
\begin{pf}
The proof is analogous to the proof for the $\phi$ of the pseudo
colimit of small categories in Lemma \ref{phifunctor}. The only
difference is that here we have to verify that $\tau_{i,j}(f) \ast
i_b$ is a 2-cell of the 2-category $\mathcal{C}$ of pseudo
$T$-algebras for any morphism $b:V \rightarrow L$ as in the comments
just before Lemma \ref{phifunctor}. But that is immediate because
$i_b$ is obviously a 2-cell and the horizontal composition of
2-cells is again a 2-cell.
\end{pf}

Now we construct a functor $\psi:PseudoCone(V,F) \rightarrow
Mor_{\mathcal{C}}(V,L)$ that is inverse to $\phi$. First we define
$\psi$ for objects, then for morphisms. Finally we verify that it is
a functor and inverse to $\phi$. The next two lemmas define a
morphism $\psi(\pi'):V \rightarrow L$ in $\mathcal{C}$ for any
object $\pi'$ of $PseudoCone(V,F)$.

\begin{lem}
Let $\pi':\Delta_V \Rightarrow F$ be a pseudo natural transformation
with coherence 2-cells $\tau'$. For any fixed $x \in Obj
\hspace{1mm} V$ we have $\psi(\pi')(x):=b(x):=(\pi_i'(x))_i \times
(\tau'_{Sf,Tf}(f)_x)_f$ is an element of $Obj \hspace{1mm}  L$.
\end{lem}
\begin{pf}
This follows from Lemma \ref{PL3} by forgetting the pseudo
$T$-algebra structures. Thus $\psi(\pi')(x) \in Obj \hspace{1mm} L$.
\end{pf}

\begin{lem}
Let $\pi':\Delta_V \Rightarrow F$ be a pseudo natural transformation
with coherence 2-cells $\tau'$. Then for any fixed $h \in
Mor_V(x,y)$ we have a modification
$\psi(\pi')(h):=b(h):=(\pi_i'(h))_i:b(x) \rightsquigarrow b(y)$.
This notation means $b(h)_i(\ast):=\pi_i'(h)$.
\end{lem}
\begin{pf}
This is exactly the same as the proof of Lemma \ref{PL4} because the
pseudo $T$-algebra structure on $L$ makes no additional requirements
on the morphisms of the small category $L$.
\end{pf}

\begin{lem}
For any pseudo natural transformation $\pi':\Delta_V \Rightarrow F$
the map $\psi(\pi')=b:V \rightarrow L$ as defined above is a
morphism of pseudo $T$-algebras.
\end{lem}
\begin{pf}
By Lemma \ref{PL5} the map $b:V \rightarrow L$ is a functor between
the underlying small categories. We define a natural transformation
$\rho_w^b$ for $w \in T(2)$. We abbreviate the application of any
structure map to $w$ by $+$. Define
$\rho_w^b(x_1,x_2):=\rho^b(x_1,x_2):=
(\rho^{\pi_i'}(x_1,x_2))_i:b(x_1+x_2) \rightarrow b(x_1)+b(x_2)$ for
all $x_1,x_2 \in Obj \hspace{1mm} V$. We claim that
$\rho^b(x_1,x_2)$ is a morphism in $L$. Let $\tau_{i,j}'(f)$ denote
the coherence 2-cell of $\pi':\Delta_V \Rightarrow F$ for $f:i
\rightarrow j$ in $\mathcal{J}$. Since $\tau_{i,j}'(f) :Ff \circ
\pi'\Rightarrow \pi_j'$ is a 2-cell, we know that
$$\xymatrix@C=9pc@R=4pc{Ff \circ \pi_i'(x_1 +x_2) \ar[r]^{\tau_{i,j}'(f)_{x_1+x_2}}
\ar[d]_{\rho^{Ff \circ \pi_i'}(x_1,x_2)} & \pi_j'(x_1+x_2)
\ar[d]^{\rho^{\pi_j'}(x_1,x_2)}
\\ Ff \circ \pi_i'(x_1)+Ff \circ\pi_i'(x_2) \ar[r]_-{\tau_{i,j}'(f)_{x_1}+
\tau_{i,j}'(f)_{x_2}} & \pi_j'(x_1) + \pi_j'(x_2)}
$$
commutes. Rewriting the left vertical arrow and the bottom arrow
gives
$$\xymatrix@C=14pc@R=4pc{Ff(\pi_i'(x_1 +x_2)) \ar[r]^{\tau_{i,j}'(f)_{x_1+x_2}}
\ar[d]_{Ff \rho^{\pi_i'}(x_1,x_2)} & \pi_j'(x_1+x_2)
\ar[d]^{\rho^{\pi_j'}(x_1,x_2)}
\\ Ff(\pi_i'(x_1)+\pi_i'(x_2)) \ar[r]_-{(\tau_{i,j}'(f)_{x_1}+
\tau_{i,j}'(f)_{x_2}) \circ \rho^{Ff}(\pi_i'(x_1),\pi_i'(x_2))} &
\pi_j'(x_1) + \pi_j'(x_2)}
$$
which states precisely that
$\rho^b(x_1,x_2)=(\rho^{\pi_i'}(x_1,x_2))_i: b(x_1 +x_2) \rightarrow
b(x_1)+b(x_2)$ is a morphism in $L$ by Remark \ref{PL2}.  The map
$\rho^b$ is natural because each component is natural. Hence
$\rho^b$ is a natural transformation. If we define $\rho_w^b$
analogously for arbitrary words $w$ of the theory $T$, then the
coherences of Definition \ref{laxalgebradefinition} are satisfied
because they are satisfied componentwise. Hence $\psi(\pi')=b:V
\rightarrow L$ is a morphism of pseudo $T$-algebras.
\end{pf}

\begin{lem}
Let $\Xi:\alpha \rightsquigarrow \beta$ be a morphism in the
category $PseudoCone(V, F)$. Then $\psi(\Xi):\psi(\alpha)
\Rightarrow \psi(\beta)$ defined by $V \ni x \mapsto (\Xi_i(x))_i
\in Mor_L(\psi(\alpha)x, \psi(\beta)x)$ is a 2-cell in the
2-category of pseudo $T$-algebras. As in Lemma \ref{PL6}, this
definition means $\psi(\Xi)(x)_i(\ast):=\Xi_i(x)$.
\end{lem}
\begin{pf}
The map $\psi(\Xi)$ is a natural transformation by Lemma \ref{PL6}.
For all $i \in Obj \hspace{1mm} \mathcal{J}$ we have morphisms
$\alpha_i,\beta_i:V \rightarrow Fi$ and 2-cells $\Xi_i:\alpha_i
\Rightarrow \beta_i$. Hence
$$\xymatrix@C=9pc@R=4pc{\alpha_i(x_1+x_2) \ar[r]^{\Xi_i(x_1+x_2)}
\ar[d]_{\rho^{\alpha_i}(x_1,x_2)} & \beta_i(x_1+x_2)
\ar[d]^{\rho^{\beta_i}(x_1,x_2)}
\\ \alpha_i(x_1) + \alpha_i(x_2) \ar[r]_{\Xi_i(x_1)+\Xi_i(x_2)}
& \beta_i(x_1)+\beta_i(x_2) }$$ commutes. Since these are the
components for $\psi(\alpha)(x), \psi(\beta)(x)$, and
$\psi(\Xi)(x)$, we see that
$$\xymatrix@C=7pc@R=4pc{\psi(\alpha)(x_1+x_2) \ar[r]^{\psi(\Xi)(x_1+x_2)}
\ar[d]_{\rho^{\psi(\alpha)}(x_1,x_2)} & \psi(\beta)(x_1+x_2)
\ar[d]^{\rho^{\psi(\beta)}(x_1,x_2)}
\\ \psi(\alpha)(x_1)+\psi(\alpha)(x_2) \ar[r]_{\psi(\Xi)(x_1)+
\psi(\Xi)(x_2)} & \psi(\beta)(x_1)+\psi(\beta)(x_2) }$$ commutes.
Similar diagrams hold for arbitrary words $w$ in the theory $T$.
Thus $\psi(\Xi)$ is a 2-cell.
\end{pf}

\begin{thm} \label{isoofcategories}
The map $\psi:PseudoCone(V,F) \rightarrow Mor_{\mathcal{C}}(V,L)$ as
defined in the previous lemmas is an inverse functor to $\phi$.
\end{thm}
\begin{pf}
This follows from the calculations of Theorem \ref{calc1} and Lemmas
\ref{calc2} and \ref{calc3}.
\end{pf}

\begin{lem}
The pseudo $T$-algebra $L$ with the pseudo cone $\pi:\Delta_L
\Rightarrow F$ is a pseudo limit of the pseudo functor
$F:\mathcal{J} \rightarrow \mathcal{C}$.
\end{lem}
\begin{pf}
The functor $\phi:Mor_{\mathcal{C}}(V,L) \rightarrow
PseudoCone(V,F)$ is an isomorphism of categories by the previous
lemmas. Since $V$ was an arbitrary object of $\mathcal{C}$ we
conclude that $L$ and $\pi$ are universal.
\end{pf}

Thus every pseudo functor $F:\mathcal{J} \rightarrow \mathcal{C}$
from a small 1-category $\mathcal{J}$ to the 2-category
$\mathcal{C}$ of pseudo $T$-algebras admits a pseudo limit. Hence
$\mathcal{C}$ admits pseudo limits. This completes the proof of
Theorem \ref{alglaxlimits}.

\end{pf}

\begin{lem} \label{cotensorproductsT}
The 2-category $\mathcal{C}$ of small pseudo $T$-algebras admits
cotensor products\index{cotensor product|textbf}.
\end{lem}
\begin{pf}
Let $J \in Obj \hspace{1mm}  Cat$ and let $F$ be a pseudo
$T$-algebra. Let $U:\mathcal{C} \rightarrow Cat$ be the forgetful
functor. Define $P:=(UF)^{J}$, which is the 1-category of 1-functors
$J \rightarrow UF$. We claim that $P$ has the structure of a pseudo
$T$-algebra. Let $\Phi_n:T(n) \rightarrow Functors(F^n,F)$ denote
the structure maps for $F$. Define
\newline $\Phi^P_n:T(n) \rightarrow Functors(P^n,P)$ by
$$\Phi^P_n(w)(p_1, \dots, p_n)(j):=\Phi_n(w)(p_1(j), \dots,
p_n(j))$$ for $j\in Obj \hspace{1mm}  J$ and $p_1, \dots, p_n \in
Obj \hspace{1mm}  P$. Coherence isos are defined analogously. For
example, define $s^P_{w,f}:\Phi^P(w_f) \Rightarrow \Phi^P(w)_f$ for
$f:m \rightarrow n$ on $p_1, \dots, p_n \in  Obj \hspace{1mm} P$ as
the 1-natural transformation
$$\xymatrix@1{s_{w,f}^P(p_1, \dots, p_n):\Phi_n^P(w_f)(p_1, \dots,
p_n)  \ar@{=>}[r] & \Phi^P_n(w)_f(p_1, \dots, p_n)}$$ which is
$s_{w,f}^P(p_1, \dots, p_n)(j):=s_{w,f}(p_1(j), \dots, p_n(j))$ for
$j \in  Obj \hspace{1mm}  J$. Then all coherence diagrams are
satisfied because they are satisfied pointwise. Hence, $P$ has the
structure of a pseudo $T$-algebra.

We claim that $P$ is a cotensor product\index{cotensor product} of
$J$ and $F$. We use Remark \ref{cotensorunit}. Define a functor
$\pi:J \rightarrow \mathcal{C}(P,F)$ by $$\pi(j)(p):=p(j)$$
$$\pi(j)(\eta):=\eta(j)$$ $$\pi(g)(p):=p(g)$$ for $j$ an object of
$J$, $p$ a functor from $J$ to $UF$, $\eta$ a natural
transformation, and $g$ a morphism in $J$. Let $\sigma:J \rightarrow
\mathcal{C}(C,F)$ be a functor. Define a morphism $b:C \rightarrow
P$ of pseudo $T$-algebras by
$$b(c)(j):=\sigma(j)(c)$$
$$b(c)(f):=\sigma(f)(c)$$ $$b(m)(j):=\sigma(j)(m)$$  for $c \in
Obj \hspace{1mm} C$, $j \in Obj \hspace{1mm}  J$, $f \in Mor
\hspace{1mm}  J$, and $m \in Mor \hspace{1mm}  C$. Then $b$ is
strict and it is the unique morphism $C \rightarrow P$ such that
$\mathcal{C}(b,F) \circ \pi =\sigma$. A similar argument can be made
for 2-cells. Thus $P$ is a
 cotensor product of $J$ and $F$ with unit $\pi$.
\end{pf}

\begin{thm} \label{weightedlaxlimitsT}
The 2-category $\mathcal{C}$ of small pseudo $T$-algebras admits
weighted pseudo limits.\index{algebra!pseudo algebra!pseudo limits
of pseudo algebras|textbf}\index{limit!weighted pseudo
limit|textbf}\index{weighted|textbf}
\end{thm}
\begin{pf}
By Theorem \ref{alglaxlimits} it admits pseudo limits, and hence it
admits pseudo equalizers\index{equalizer!pseudo equalizer}. The
2-category $\mathcal{C}$ obviously admits
2-products\index{2-product}. By Lemma \ref{cotensorproductsT} it
admits cotensor products. Hence by Theorem \ref{streetpseudo} it
admits weighted pseudo limits.
\end{pf}

\begin{thm}
The 2-category $\mathcal{C}$ of small pseudo $T$-algebras admits
weighted bilimits.\index{bilimit|textbf}\index{bilimit!weighted
bilimit|textbf}\index{algebra!pseudo algebra!bilimits of pseudo
algebras|textbf}\index{weighted|textbf}
\end{thm}
\begin{pf}
It admits weighted pseudo limits and therefore admits weighted
bilimits.
\end{pf}

\chapter{Biuniversal Arrows and Biadjoints}
\label{sec:laxadjoints} After studying bilimits and bicolimits, we
turn our attention to another type of weakened structure called {\it
biadjoints}\index{biadjoint}. The concept of an
adjunction\index{adjunction} from 1-category theory consists of two
functors and a natural bijection between appropriate hom sets. Mac
Lane lists several equivalent ways of describing an adjunction in
\cite{maclane1} on pages 79-86. One of these ways involves a
universal arrow\index{universal arrow|textbf}\index{arrow!universal
arrow|textbf} for each object of the source category. To weaken
these concepts, we replace the functors by pseudo functors, the
natural bijection of hom sets by a pseudo natural equivalence of
categories, and the universal arrow by a biuniversal
arrow\index{biuniversal arrow}\index{arrow!biuniversal arrow}. The
main goal in this chapter is to prove that a
biadjunction\index{biadjunction} can be described via pseudo natural
equivalences or via biuniversal arrows. This is the meaning of
Theorem \ref{maintheorem1} and Theorem \ref{maintheorem2}.

A close result in the literature can be found in Gray's work
\cite{gray1}. His concept of {\it transcendental
quasiadjunction}\index{transcendental
quasiadjunction}\index{quasiadjunction}\index{quasiadjunction!transcendental
quasiadjunction} between two 2-functors on page 177 is similar to
the concept of biadjunction between two pseudo functors except that
the functors in a biadjoint are allowed to be pseudo. Gray remarks
on pages 180-181 that a transcendental
quasiadjunction{\index{quasiadjunction} gives rise to a certain
universal mapping property. The analogous concept for biadjoints is
a biuniversal arrow and the appropriate theorem is Theorem
\ref{maintheorem1}. On page 184 Gray remarks that under certain
hypothesis, the universal mapping property gives rise to a
quasiadjunction. The biadjoint version of this is Theorem
\ref{maintheorem2} in which the starting functor $G$ is allowed to
be a pseudo functor.

Kelly phrases a similar result in \cite{kelly2} on page 316 in terms
of homomorphisms\index{homomorphism of bicategories} of bicategories
and birepresentations\index{birepresentation}. His notion of
biadjoint is the same as in this paper, except that we are
considering only pseudo functors between 2-categories rather than
homomorphisms between bicategories. Kelly's statement is equivalent
to \ref{maintheorem2} after an application of Yoneda's
Lemma\index{Yoneda's Lemma for bicategories} for bicategories.
Yoneda's Lemma for bicategories can be found in \cite{street3}.

Street makes an observation on page 121 in \cite{street3} similar to
Theorem \ref{maintheorem2}: if each object admits a left
bilifting\index{bilifting} then a left biadjoint exists. The unit
for a left bilifting is the biuniversal arrow of Theorem
\ref{maintheorem2}.

MacDonald and Stone also have a weakened notion of adjunction in
\cite{macdonald1} called {\it soft adjunction}\index{soft
adjunction}\index{adjunction!soft adjunction}. In that article they
consider strict 2-functors and natural adjunctions between hom
categories. They prove theorems about the universality concepts that
arise in such a context.

We follow Mac Lane's presentation of adjoints except we account for
the 2-cells.  The notation in this study is analogous to the
notation in Mac Lane's book.  Recall the definition of a universal
arrow and its uniqueness.

\begin{defn} \label{defuniversalarrow}
Let $S:D \rightarrow C$ be a functor between 1-categories and $c \in
Obj \hspace{1mm}  C$. Then an object $r \in Obj \hspace{1mm} D$ and
a morphism $u \in Mor_C(c,Sr)$ are a {\it universal arrow from $c$
to $S$}\index{universal arrow|textbf}\index{arrow!universal
arrow|textbf} if for every $d \in Obj \hspace{1mm} D$ and every $f
\in Mor_C(c,Sd)$ there exists a unique morphism $f' \in Mor_D(r,d)$
such that $Sf' \circ u = f$. Pictorially this means for every $d$
and every $f$ as above, there exists a unique $f'$ making
$$\xymatrix@R=3pc@C=3pc{c \ar@{=}[d] \ar[r]^u & Sr  \ar@{.>}[d]^{Sf'}
 & r \ar@{.>}[d]^{f'} \\ c \ar[r]_{f} & Sd & d}$$
commute. This is equivalent to saying the assignment $f' \mapsto Sf'
\circ u$, \newline $Mor_D(r,d)\rightarrow Mor_C(c,Sd)$  is a
bijection of hom sets for every fixed $d \in Obj \hspace{1mm}  D$.
\end{defn}

\begin{lem}
Let $u:c \rightarrow Sr$ and $u':c \rightarrow Sr'$ be universal
arrows from the object $c$ to the functor $S$. Then there exists a
unique morphism $f':r \rightarrow r'$ such that $Sf' \circ u = u'$.
Moreover, the morphism $f':r \rightarrow r'$ is an isomorphism.
\end{lem}
\begin{pf}
There exist unique morphisms $f'$ and $g'$ such that the following
diagram commutes.
$$\xymatrix@R=3pc@C=3pc{c \ar@{=}[d] \ar[r]^u & Sr \ar@{.>}[d]^{Sf'} & r \ar@{.>}[d]^{f'}
\\ c \ar@{=}[d] \ar[r]^{u'} & Sr' \ar@{.>}[d]^{Sg'} & r' \ar@{.>}[d]^{g'}
\\ c \ar[r]_u & Sr & r}$$
The middle vertical column could be replaced by $S1_r$ to make the
outermost rectangle commutative. Hence by the uniqueness we have $g'
\circ f' = 1_r$. Similarly we can show that $f' \circ g' = 1_{r'}$.
Hence $f'$ is an isomorphism and $Sf' \circ u = u'$.
\end{pf}

Before weakening the concept of universal arrow, we prove a simple
lemma that will make it easier to visualize a biuniversal arrow.

\begin{lem} \label{visualization}
Let $\xymatrix@1{X  \ar@<.5ex>[r]^{\phi} & A \ar@<.5ex>[l]^{\psi}}$
be adjoint functors with unit $\theta:1_X \Rightarrow \psi \circ
\phi$ and counit $\mu: \phi \circ \psi \Rightarrow 1_A$. Suppose
that both the unit and the counit are natural isomorphisms. Let
$\nu:\phi (x) \rightarrow a$ be a morphism in $A$ and $x \in Obj
\hspace{1mm}  X, a \in Obj \hspace{1mm}  A$. Then there exists a
unique morphism $\nu':x \rightarrow \psi (a)$ such that
$$\xymatrix@C=4pc@R=3pc{x \ar@{.>}[d]_{\nu'}
& \phi (x) \ar[r]^-{\nu} \ar@{.>}[d]_{\phi(\nu')} & a \ar@{=}[d]
\\ \psi (a) & \phi (\psi (a)) \ar[r]_-{\mu(a)} & a}$$
commutes. Moreover, $\nu'$ is iso if and only if $\nu$ is
iso\index{iso}.
\end{lem}
\begin{pf}
The existence and uniqueness claims follow because $\mu(a)$ is a
universal arrow from $\phi$ to $a$. If $\nu'$ is iso, then
$\phi(\nu')$ is iso and so is $\nu=\mu(a) \circ \phi(\nu')$ because
$\mu(a)$ is iso by hypothesis. It only remains to show that $\nu'$
is iso if $\nu$ is iso. Suppose $\nu$ is iso. Then $\phi(\nu')$ is
iso from the commutivity of the diagram because $\mu(a)$ and $\nu$
are iso. By the naturality of $\theta$ we have
$$\xymatrix@C=4pc@R=3pc{x \ar[r]^-{\theta(x)} \ar[d]_{\nu'}
& \psi \circ \phi (x) \ar[d]^{\psi \circ \phi(\nu')}
\\ \psi(a) \ar[r]_-{\theta(\psi (a))}
& \psi \circ \phi (\psi (a))}$$ commutes. Then $\nu'$ is iso because
$\theta(x), \theta(\psi (a)),$ and $\psi(\phi(\nu'))$ are iso.
\end{pf}

To weaken the concept of universal arrow in the context of
2-categories, we replace the bijection of sets above by an
equivalence of the appropriate morphism categories.
\begin{defn} \label{deflaxuniversalarrow}
Let $S:\mathcal{D} \rightarrow \mathcal{C}$ be a pseudo functor
between 2-categories and $C \in Obj \hspace{1mm}  \mathcal{C}$. Then
an object $R \in Obj \hspace{1mm}  \mathcal{D}$ and a morphism $u
\in
 Mor_{\mathcal{C}}(C,SR)$ are a {\it biuniversal arrow}\index{biuniversal arrow|textbf}\index{arrow!biuniversal arrow|textbf}
from $C$ to $S$ if for every $D \in Obj \hspace{1mm}  \mathcal{D}$
the functor $\phi:Mor_{\mathcal{D}}(R,D) \rightarrow
Mor_{\mathcal{C}}(C,SD)$ defined by $f' \mapsto Sf' \circ u$ and
$\gamma \mapsto S \gamma \ast i_u$ is an equivalence of categories.
\end{defn}

We suppressed the dependence of $\phi$ on $D$ in the notation of the
definition. This definition implies that $\phi$ admits a right
adjoint $\psi$ such that the counit $\mu:\phi \circ \psi \Rightarrow
1_{Mor_{\mathcal{C}}(C,SD)}$ and unit are natural isomorphisms.
Pictorially the definition implies that for every object $D \in Obj
\hspace{1mm}  \mathcal{D}$ and every morphism $f: C \rightarrow SD$
in $\mathcal{C}$ there exists an $f'$ and a natural universal 2-cell
$\mu(f)$ which is iso (an arrow of the counit) as in the following
diagram.
$$\xymatrix@C=4pc@R=3pc{C \ar@{=}[d] \ar[r]^u & SR  \ar@{.>}[d]^{Sf'}
\ar@{:>}[dl]_{\mu(f)} & R \ar@{.>}[d]^{f'} \\ C \ar[r]_{f} & SD &
D}$$ The assignment $\psi:f \mapsto f'$ is functorial and $\mu:\phi
\circ \psi \Rightarrow 1_{Mor_{\mathcal{C}}(C,SD)}$ is a natural
transformation. This diagram is not equivalent to the definition
because it does not express the naturality of the 2-cells, nor does
it include the natural isomorphism (the unit) from the identity
functor on $Mor_{\mathcal{D}}(R,D)$ to $\psi \circ \phi$. The
universality of the 2-cell $\mu(f)$ from the functor $\phi$ to the
object $f$ means pictorially that the arrow $f'$ is unique up to
2-cell in the following way. If $\bar{f'}:R \rightarrow D$ is an
arrow in $\mathcal{D}$ and $\nu$ is a (not necessarily iso) 2-cell
as in
$$\xymatrix@C=4pc@R=3pc{C \ar@{=}[d] \ar[r]^u & SR  \ar[d]^{S\bar{f'}}
\ar@{=>}[dl]_{\nu} & R \ar[d]^{\bar{f}'} \\ C \ar[r]_{f} & SD & D}$$
then there exists a unique 2-cell $\nu':\bar{f'} \Rightarrow f'$
whose $\phi$ image factors $\nu$ via the universal arrow $\mu(f)$,
\ie $\nu'$ is such that \label{factorizing2cellnuprime}
$$\xymatrix@C=4pc@R=3pc{\bar{f}' \ar@{:>}[d]_{\nu'}
& S\bar{f}' \circ u \ar@{=>}[r]^-{\nu} \ar@{:>}[d]_{\phi(\nu')=S
\nu' * i_{u}} & f \ar@{=}[d]
\\ f' & Sf' \circ u \ar@{=>}[r]_-{\mu(f)} & f}$$
commutes. We also know that $\nu'$ is iso if and only if $\nu$ is
iso as in Lemma \ref{visualization}. Note that these diagrams are
dual to Definition \ref{defuniversalarrow}, although it is the same
concept of universal arrow.

One can ask if the equivalences of categories in the definition of
biuniversal arrow can be chosen in some natural way as in Remark
\ref{laxcolimitnatural}. They can in fact as the following theorem
shows.

\begin{thm}
Let $u:C \rightarrow SR$ be a biuniversal arrow from $C$ to the
pseudo functor $S$ as in Definition \ref{deflaxuniversalarrow}. Let
 $\phi_D:Mor_{\mathcal{D}}(R,D) \rightarrow Mor_{\mathcal{C}}(C,SD)$ be the
functor defined by $f' \mapsto Sf' \circ u$ and $\gamma \mapsto S
\gamma \ast i_u$. Then $D \mapsto \phi_D$ is a pseudo natural
transformation $Mor_{\mathcal{D}}(R,-) \Rightarrow
Mor_{\mathcal{C}}(C,S-)$. For $D \in Obj \hspace{1mm} \mathcal{D}$
let $\psi_D:Mor_{\mathcal{C}}(C,SD) \rightarrow
Mor_{\mathcal{D}}(R,D)$ be a right adjoint to $\phi_D$ such that the
unit $\eta_D:1_{Mor_{\mathcal{D}}(R,D)} \Rightarrow \psi_D \circ
\phi_D$ and the counit $\varepsilon_D:\phi_D \circ \psi_D
\Rightarrow 1_{Mor_{\mathcal{C}}(C,SD)}$ are natural isomorphisms.
Then $D \mapsto \psi_D$ is a pseudo natural transformation and $D
\mapsto \eta_D$ and $D \mapsto \varepsilon_D$ are iso modifications
$i_{Mor_{\mathcal{D}}(R,-)} \rightsquigarrow \psi \odot \phi$ and
$\phi \odot \psi \rightsquigarrow i_{Mor_{\mathcal{C}}(C,S-)}$ which
satisfy the triangle identities.
\end{thm}
\begin{pf}
Let $F,G:\mathcal{D} \rightarrow Cat$ be the pseudo functors defined
by $F(D)=Mor_{\mathcal{D}}(R,D)$ and $G(D)=Mor_{\mathcal{C}}(C,SD)$.
Then $F$ is a strict 2-functor. One can prove that $\phi:F
\Rightarrow G$ is a pseudo natural transformation by defining the
coherence 2-cell $\tau$ in terms of $\gamma^S$ and then using the
unit and composition axioms for $S$ to prove the unit and
composition axioms for $\phi$. After doing that, we are in the setup
of Lemma \ref{onevariable}, from which everything else follows.
\end{pf}

In analogy to the uniqueness statement for universal arrows, we have
a uniqueness statement for biuniversal arrows. It requires the
concept of pseudo isomorphism in a 2-category.

\begin{defn}
Let $\mathcal{D}$ be a 2-category and $f:R \rightarrow R'$ a
morphism in $\mathcal{D}$. Then $f$ is a {\it pseudo
isomorphism}\index{isomorphism}\index{isomorphism!pseudo
isomorphism|textbf} if there exists a morphism $g:R' \rightarrow R$
and iso 2-cells $g \circ f \Rightarrow 1_R$ and $g \circ f
\Rightarrow 1_{R'}$. A pseudo isomorphism is also called an {\it
equivalence}\index{equivalence|textbf}.
\end{defn}

\begin{lem} \label{laxuniqueness}
Let $S: \mathcal{D} \rightarrow \mathcal{C}$ be a pseudo functor.
Let $u_1:C \rightarrow SR_1$ and $u_2:C \rightarrow SR_2$ be
biuniversal arrows from $C$ to $S$. Then there exists a pseudo
isomorphism $g':R_1 \rightarrow R_2$ in $\mathcal{D}$ and an iso
2-cell as in (\ref{uniqueness}).
\begin{equation} \label{uniqueness}
\xymatrix@C=4pc@R=3pc{C \ar[r]^{u_1} \ar@{=}[d] & SR_1
\ar@{.>}[d]^{Sg'} \ar@{:>}[ld]_{\mu_1(u_2)} & R_1 \ar@{.>}[d]^{g'}
\\ C \ar[r]_{u_2} & SR_2 & R_2}
\end{equation}
Moreover, if $\bar{g}'$ and $\nu$ are a morphism and an iso 2-cell
that also fill in the diagram, then $\bar{g}'$ and $g'$ are
isomorphic via the unique 2-cell $\nu':\bar{g}' \rightarrow g'$ such
that $\mu_1(u_2) \circ (S \nu'*i_{u_1})=\nu$.
\end{lem}
\begin{pf}
The biuniversality of $u_1$ and $u_2$ guarantees the existence of
arrows $f',g'$, and $h'$ and iso 2-cells $\mu_1(u_2),\mu_2(u_1),$
and $\mu_1(u_1)$ to fill in the following diagrams.
\begin{equation} \label{two}
\xymatrix@C=4pc@R=3pc{C \ar@{=}[d] \ar[r]^{u_1} & SR_1
\ar@{.>}[d]^{Sf'} \ar@{:>}[ld]_{\mu_1(u_1)} & R_1 \ar@{.>}[d]^{f'}
\\ C \ar[r]_{u_1} & SR_1 & R_1}
\end{equation}
\begin{equation} \label{three}
\xymatrix@C=4pc@R=3pc{C \ar@{=}[d] \ar[r]^{u_1} & SR_1
\ar@{.>}[d]^{Sg'} \ar@{:>}[ld]_{\mu_1(u_2)} & R_1 \ar@{.>}[d]^{g'}
\\ C \ar@{=}[d] \ar[r]_{u_2}
& SR_2 \ar@{.>}[d]^{Sh'} \ar@{:>}[ld]^{\mu_2(u_1)} & R_2
\ar@{.>}[d]^{h'}
\\ C \ar[r]_{u_1} & SR_1 & R_1}
\end{equation}
The arrow $1_{R_1}$ also fills in the diagram
\begin{equation} \label{four}
\xymatrix@C=6pc@R=3pc{C \ar@{=}[d] \ar[r]^{u_1} & SR_1
\ar[d]^{S1_{R_1}} \ar@{=>}[ld]_{i_{u_1}*\delta^{-1}_{R_1 *}
\hspace{.1in}} & R_1 \ar[d]^{1_{R_1}}
\\ C \ar[r]_{u_1} & SR_1 & R_1}
\end{equation}
with an iso 2-cell. Diagram (\ref{three}) combined appropriately
with $(\gamma^S_{g',h'})^{-1}$ gives an iso 2-cell $h' \circ g'
\Rightarrow f'$ by the comments after the definition of biuniversal
arrow. Similarly, diagram (\ref{four}) gives an iso 2-cell $1_{R_1}
\Rightarrow f'$ for the same reason. Combining these two iso 2-cells
appropriately gives an iso 2-cell $h' \circ g' \Rightarrow 1_{R_1}$.
By a similar argument we obtain an iso 2-cell $g' \circ h'
\Rightarrow 1_{R_2}$. Thus $g':R_1 \rightarrow R_2$ is a pseudo
isomorphism. The iso 2-cell between $\bar{g}'$ and $g'$ is also
guaranteed by the comments after the definition of biuniversal arrow
in \ref{deflaxuniversalarrow}.
\end{pf}

After these preparations involving biuniversal arrows, we can now
introduce the main concept of this chapter.

\begin{defn} \label{defnlaxadjunction}
Let $\mathcal{X}$ and $\mathcal{A}$ be 2-categories. A {\it
biadjunction}\index{biadjunction|textbf} $\langle F,G,\phi \rangle:
\mathcal{X} \rightharpoonup \mathcal{A}$ consists of the following
data
\begin{itemize}
\item
Pseudo functors
$$\xymatrix@R=3pc@C=3pc{ \mathcal{X}  \ar@<.5ex>[r]^F & \mathcal{A} \ar@<.5ex>[l]^G}$$
between 2-categories
\item
For all $X \in Obj \hspace{1mm}  \mathcal{X}$ and all $A \in Obj
\hspace{1mm} \mathcal{A}$ an equivalence of categories
$\phi_{X,A}:Mor_{\mathcal{A}}(FX,A) \rightarrow
Mor_{\mathcal{X}}(X,GA)$ assigned in such a way to make $\phi$ into
a pseudo natural transformation in each variable between the
following pseudo functors of two variables.
$$\xymatrix@C=7pc@R=3pc{\mathcal{X}^{op} \times \mathcal{A} \ar[r]^{F^{op} \times
1_{\mathcal{A}}} & \mathcal{A}^{op} \times \mathcal{A} \ar[r]^-{Mor}
\ar@{=>}[d]^{\phi} & Cat
\\ \mathcal{X}^{op} \times \mathcal{A} \ar[r]^{1_{\mathcal{X}^{op}} \times G}
& \mathcal{X}^{op} \times \mathcal{X} \ar[r]^-{Mor} & Cat}$$
\end{itemize}
In this situation, $F$ is called a {\it left
biadjoint}\index{biadjoint|textbf}\index{biadjoint!left
biadjoint|textbf} for $G$ and $G$ is called a {\it right
biadjoint}\index{biadjoint!right biadjoint|textbf} for $F$.
\end{defn}

Recall again that a {\it biadjoint} is called a {\it lax
adjoint}\index{adjoint!lax adjoint|textbf} in \cite{hu}, \cite{hu1},
and \cite{hu2}. The degree of uniqueness of a left biadjoint (if a
left biadjoint exists), will be dealt with at the end of this
chapter. One can ask whether or not an adjoint functor
$\psi_{X,A}:Mor_{\mathcal{X}}(X,GA) \rightarrow
Mor_{\mathcal{A}}(FX,A)$ to $\phi_{X,A}$ can be chosen in a natural
way. This is similar to the question answered in Remark
\ref{laxcolimitnatural} for bicolimits. To show that right adjoints
can be chosen in a pseudo natural way, we need the following lemma.

\begin{lem} \label{onevariable}
Let $F,G:\mathcal{A} \rightarrow Cat$ be pseudo functors and $F$ a
strict 2-functor. Suppose we have a pseudo natural transformation
$\phi:F \Rightarrow G$ such that $\phi_A:FA \rightarrow GA$ is an
equivalence of categories for all $A \in Obj \hspace{1mm}
\mathcal{A}$. For each $A \in Obj \hspace{1mm}  \mathcal{A}$, let
$\psi_A:GA \rightarrow FA$ be a right adjoint to $\phi_A$ such that
the unit $\eta_A:1_{FA} \Rightarrow \psi_A \circ \phi_A$ and counit
$\varepsilon_A: \phi_A \circ \psi_A \Rightarrow 1_{GA}$ are natural
isomorphisms. Then $A \mapsto \psi_A$ is a pseudo natural
transformation $G \Rightarrow F$. The assignments $A \mapsto \eta_A$
and $A \mapsto \varepsilon_A$ define iso modifications $\eta:i_F
\rightsquigarrow \psi \odot \phi$ and $\varepsilon: \phi \odot \psi
\rightsquigarrow i_G$ respectively. Furthermore, $\eta$ and
$\varepsilon$ satisfy the triangle identities.
\end{lem}
\begin{pf}
For all $A \in Obj \hspace{1mm}  \mathcal{A}$ there exists such a
right adjoint $\psi_A$ because $\phi_A$ is an equivalence of
categories.

To show that $A \mapsto \psi_A$ is a pseudo natural transformation,
we need to define the coherence 2-cell $\tau_f'$ for each morphism
$f$ of $\mathcal{A}$, show that it is natural, it satisfies the unit
axiom, and that it satisfies the composition axiom.

For a morphism $f:A \rightarrow B$ in $\mathcal{A}$ let $\tau_f:Gf
\circ \phi_A \Rightarrow \phi_B \circ Ff$ denote the coherence
2-cell belonging to the pseudo natural transformation $\phi$. Define
$\tau_f':Ff \circ \psi_A \Rightarrow \psi_B \circ Gf$ to be the
composition of the 2-cells in diagram
(\ref{towardsdefiningtaufprime}).
\begin{equation} \label{towardsdefiningtaufprime}
\xymatrix@R=4pc@C=4pc{GA \ar[d]_{1_{GA}} \ar[r]^{\psi_A} & FA
\ar[d]^{1_{FA}} \ar@{=>}[ld]_{\varepsilon_A}
\\ GA \ar[d]_{Gf} & \ar[l]_{\phi_A} FA \ar[d]^{Ff}
\\ GB \ar[d]_{1_{GB}} &
\ar[l]^{\phi_B} FB \ar[d]^{1_{FB}} \ar@{=>}[ul]_{\tau_f^{-1}}
\ar@{=>}[ld]_{\eta_B}
\\ GB \ar[r]_{\psi_B} &  FB}
\end{equation}
We claim that the assignment $f \mapsto \tau_f'$ is natural in $f$.
To see this, let $f,g:A \rightarrow B$ be morphisms in $\mathcal{A}$
and $\mu:f \Rightarrow g$ a 2-cell in $\mathcal{A}$. Then $\tau_f'$
is the composition of the top row of 2-cells in diagram
(\ref{definingtaufprime}) and $\tau_g'$ is the bottom composition.

\begingroup
\vspace{-2\abovedisplayskip} \small
\begin{equation} \label{definingtaufprime}
\xymatrix@R=4pc@C=2pc{1_{FB} \circ Ff \circ \psi_A
\ar@{=>}[r]^{\underset{\phantom{\eta_B * i_{Ff}*i_{\psi_A}}} {\eta_B
* i_{Ff}*i_{\psi_A}}} \ar@{=>}[d]|{i_{1_{FB}}*F \mu
*i_{\psi_A}} & \psi_B \circ \phi_B \circ  Ff \circ \psi_A
\ar@{=>}[r]^{\underset{\phantom{\eta_B * i_{Ff}*i_{\psi_A}}}
{i_{\psi_B}*\tau_f^{-1}*i_{\psi_A}}} \ar@{=>}[d]|{i_{\psi_B \circ
\phi_B}*F \mu * i_{\psi_A}} & \psi_B \circ Gf \circ \phi_A \circ
\psi_A \ar@{=>}[r]^{\underset{\phantom{\eta_B * i_{Ff}*i_{\psi_A}}}
{i_{\psi_B}*i_{Gf}*\varepsilon_A}}
\ar@{=>}[d]|{i_{\psi_B}*G\mu*i_{\phi_A \circ \phi_A}} & \psi_B \circ
Gf \circ 1_{GA} \ar@{=>}[d]|{i_{\psi_B}*G \mu *i_{1_{GA}}}
\\ 1_{FB} \circ Fg \circ \psi_A
\ar@{=>}[r]_{\overset{\phantom{\eta_B * i_{Fg}*i_{\psi_A}}} {\eta_B
* i_{Fg}*i_{\psi_A}}} & \psi_B \circ \phi_B \circ  Fg \circ \psi_A
\ar@{=>}[r]_{\overset{\phantom{\eta_B * i_{Ff}*i_{\psi_A}}}
{i_{\psi_B}*\tau_g^{-1}*i_{\psi_A}}} & \psi_B \circ Gg \circ \phi_A
\circ \psi_A \ar@{=>}[r]_{\overset{\phantom{\eta_B *
i_{Ff}*i_{\psi_A}}} {i_{\psi_B}*i_{Gg}*\varepsilon_A}} & \psi_B
\circ Gg \circ 1_{GA}}
\end{equation}
\endgroup
\noindent
 The left square and the right square commute because of
the interchange law and the defining property of identity 2-cells.
The middle square commutes because $f \mapsto \tau_f$ is natural by
the definition of $\phi$ pseudo natural. Hence the outermost
rectangle commutes and $f \mapsto \tau_f'$ is natural.

We claim that $\tau'$ satisfies the unit axiom for pseudo natural
transformations. Since $F$ is strict, proving the coherence diagram
reduces to proving that $\tau_{1_A}'=i_{\psi_A}*\delta^G_{A*}$.
Using the definition of $\tau'$ above and the unit axiom for $\tau$
we see that $\tau_{1_A}'$ is the composition of 2-cells in diagram
(\ref{rightadjoint1}).
\begin{equation} \label{rightadjoint1}
\xymatrix@R=4pc@C=4pc{GA \ar[d]_{1_{GA}} \ar[r]^{\psi_A} & FA
\ar[d]^{1_{FA}} \ar@{=>}[ld]_{\varepsilon_A}
\\ GA \ar[d]_{G1_A} & \ar[l]_{\phi_A} FA \ar[d]^{F1_A}
\\ GA \ar[d]_{1_{GA}} &
\ar[l]^{\phi_A} FA \ar[d]^{1_{FA}}
\ar@{=>}[ul]_{\delta^G_{A*}*i_{\phi_A}} \ar@{=>}[ld]_{\eta_A}
\\ GA \ar[r]_{\psi_A} &  FA}
\end{equation}
But the composition of 2-cells in (\ref{rightadjoint1}) is the same
as the composition of 2-cells in (\ref{rightadjoint2}) by the
interchange law.
\begin{equation} \label{rightadjoint2}
\xymatrix@R=4pc@C=4pc{ & GA \ar[d]_{1_{GA}} \ar[r]^{\psi_A} & FA
\ar[d]^{1_{FA}} \ar@{=>}[ld]_{\varepsilon_A}
\\ GA \ar@{=}[r] \ar[d]_{G1_A} &  GA \ar[d]_{1_{GA}}
\ar@{=>}[ld]_{\delta^G_{A*}} & \ar[l]_{\phi_A} FA
\ar[d]^{F1_A=1_{FA}}
\\ GA \ar@{=}[r] & GA \ar[d]_{1_{GA}} &
\ar[l]^{\phi_A} FA \ar[d]^{1_{FA}} \ar@{=>}[ul]_{i_{\phi_A}}
\ar@{=>}[ld]_{\eta_A}
\\  &  GA \ar[r]_{\psi_A} &  FA}
\end{equation}
By one of the triangle identities we see that the right three
squares of (\ref{rightadjoint2}) collapse to $i_{\psi_A}$ and
therefore (\ref{rightadjoint1}) is the same as
$i_{\psi_A}*\delta^G_{A*}$. Hence
$\tau_{1_A}'=i_{\psi_A}*\delta^G_{A*}$ and the unit axiom is
satisfied.

We claim that $\tau'$ satisfies the composition axiom for pseudo
natural transformations. Let $\xymatrix@1{A \ar[r]^f & B \ar[r]^g &
C}$ be morphisms in $\mathcal{A}$. Since $F$ is a strict 2-functor,
proving the composition coherence reduces to proving that \newline
$\tau'_{g \circ f}= (i_{\psi_C} * \gamma^G_{f,g}) \odot (\tau_g' *
i_{Gf}) \odot (i_{Fg} * \tau_f')$. Following the same approach as
for the unit axiom, we write out $\tau_{g \circ f}'$ in
(\ref{rightadjoint3}).
\begin{equation} \label{rightadjoint3}
\xymatrix@R=4pc@C=4pc{GA \ar[d]_{1_{GA}} \ar[r]^{\psi_A} & FA
\ar[d]^{1_{FA}} \ar@{=>}[ld]_{\varepsilon_A}
\\ GA \ar[d]_{G(g \circ f)} & \ar[l]_{\phi_A} FA \ar[d]^{F(g \circ f)}
\\ GC \ar[d]_{1_{GC}} &
\ar[l]^{\phi_B} FC \ar[d]^{1_{FC}} \ar@{=>}[ul]_{\tau_{g \circ
f}^{-1}} \ar@{=>}[ld]_{\eta_C}
\\ GC \ar[r]_{\psi_C} &  FC}
\end{equation}
Using the composition axiom for $\tau$ and writing the 2-cells more
compactly we see that the composition of 2-cells in diagram
(\ref{rightadjoint3}) is the same as in diagram
(\ref{rightadjoint4}).
\begin{equation} \label{rightadjoint4}
\xymatrix@C=4pc@R=4pc{ & GA \ar[d]_{1_{GA}} \ar[r]^{\psi_A}
\ar@{}[dr]|(.25){\overset{\varepsilon_A}{\Leftarrow}} & FA
\ar[ld]^{\phi_A} \ar[d]^{Ff}
\\ GA \ar@{=}[r] \ar[ddd]_{G(g \circ f)}
& GA \ar@{}[rd]|(.75){\overset{\eta_B}{\Leftarrow}}
\ar@{}[r]|{\overset{\tau_f^{-1}}{\Leftarrow}} \ar[d]_{Gf} & FB
\ar[d]^{1_{FB}} \ar[ld]_{\phi_B}
\\ \ar@{}[rd]|{\overset{\gamma^G_{f,g}}{\Leftarrow}}
& GB \ar[r]_{\psi_B} \ar[d]_{1_{GB}} \ar@{}[rd]|(.25){
\overset{\varepsilon_B}{\Leftarrow}} & FB \ar[d]^{Fg}
\ar[ld]^{\phi_B}
\\ & GB \ar[d]_{Gg} \ar@{}[r]|{\overset{\tau_g^{-1}}{\Leftarrow}}
\ar@{}[rd]|(.75){ \overset{\eta_C}{\Leftarrow}} & FC \ar[d]^{1_{FC}}
\ar[ld]_{\phi_C}
\\ GC \ar@{=}[r] & GC \ar[r]_{\psi_C} & FC}
\end{equation}
The middle parallelogram involving $\eta_B$ and $\varepsilon_B$ is
the same as $i_{\phi_B}$ by the triangle identity. Hence
(\ref{rightadjoint4}) is
 $(i_{\psi_C} * \gamma^G_{f,g}) \odot (\tau_g' * i_{Gf})
\odot (i_{Fg} * \tau_f')$ and we conclude that $\tau'_{g \circ f}=
(i_{\psi_C} * \gamma^G_{f,g}) \odot (\tau_g' * i_{Gf}) \odot (i_{Fg}
* \tau_f')$ as required by the composition axiom.

Thus far we have shown that $A \mapsto \psi_A$ is a pseudo natural
transformation $G \Rightarrow F$. Next we show that $A \mapsto
\eta_A$ defines a modification $i_F \rightsquigarrow \psi \odot
\phi$.

Let $f,g:A \rightarrow B$ be morphisms in the 2-category
$\mathcal{A}$ and $\gamma:f \Rightarrow g$ a 2-cell. We claim that
the compositions in diagrams (\ref{defnmodification1}) and
(\ref{defnmodification2}) are the same, \ie that $\eta$ is a
modification. Our diagrams will of course have $F=G$, $\alpha=i_F$,
$\beta=\psi \odot \phi$, and the coherence iso belonging to $i_F$ is
trivial while the coherence iso for the composite pseudo natural
transformation $\psi \odot \phi$ is $(i_{\psi_B}*\tau_f) \odot
(\tau'_f * i_{\phi_A})$ by the remarks on page
\pageref{compositionoflaxnaturals} about coherence isos for a
vertical composition of pseudo natural transformations. Then we see
that the composition (\ref{defnmodification2}) is $\eta_B * F
\gamma$. We proceed by reducing (\ref{defnmodification1}) to $\eta_B
* F \gamma$. The composition in diagram (\ref{defnmodification1}) is
explicitly (\ref{etamodification1}), where we left off the vertical
equal signs.
\begin{equation} \label{etamodification1}
\xymatrix@R=3pc@C=3pc{FA \ar[rr]^{1_{FA}} & \ar@{=>}[d]^{\eta_A} &
FA \ar[r]^{Ff} \ar@{=>}@<5ex>[d]^{F\gamma}
 & FB
\\ FA \ar[r]^{\phi_A}
\ar@{=>}@<5ex>[d]^{i_{\phi_A}} & GA \ar[r]^{\psi_A} & FA \ar[r]_{Fg}
\ar@{=>}[d]^{\tau_g'} & FB
\\ FA \ar[r]_{\phi_A} & GA \ar[r]_{Gg} \ar@{=>}[d]^{\tau_g}
& GB \ar[r]^{\psi_B} \ar@{=>}@<5ex>[d]^{i_{\psi_B}} & FB
\\ FA \ar[r]_{Fg} & FB \ar[r]_{\phi_B} & GB \ar[r]_{\psi_B} & FB}
\end{equation}
Writing out the definition $\tau_g'$ in (\ref{etamodification1}) and
including some identities gives (\ref{etamodification2}).
\begin{equation} \label{etamodification2}
\xymatrix@R=3pc@C=3pc{FA \ar[rr]^{1_{FA}} \ar@{=}[d] &
\ar@{=>}[d]^{\eta_A} & FA \ar[r]^{Ff}  \ar@{=}[d]
\ar@{=>}@<5ex>[d]^{F\gamma} & FB \ar@{=}[d] \ar[r]^{1_{FB}}
\ar@{=>}@<5ex>[d]^{i_{1_{FB}}} & FB \ar@{=}[d]
\\ FA  \ar@{=}[d] \ar[r]^{\phi_A}
\ar@{=>}@<5ex>[d]^{i_{\phi_A}} & GA \ar[r]^{\psi_A}  \ar@{=}[d] & FA
\ar[r]_{Fg} \ar@{=>}[ld]_{\varepsilon_A} \ar[d]^{\phi_A} & FB
\ar@{=>}[ld]^{\tau_g^{-1}} \ar[r]_{1_{FB}} \ar[d]^{\phi_B} & FB
\ar@{=}[d] \ar@{=>}[ld]^{\eta_B}
\\ FA \ar[r]_{\phi_A} \ar@{=}[d] & GA \ar[r]_{1_{GA}} & GA \ar[r]_{Gg}
\ar@{=>}[d]^{\tau_g} & GB \ar[r]^{\psi_B}
\ar@{=>}@<5ex>[d]^{i_{\psi_B}} \ar@{=}[d] & FB \ar@{=}[d]
\\ FA \ar[rr]_{Fg} &  & FB \ar[r]_{\phi_B} & GB \ar[r]_{\psi_B} & FB}
\end{equation}
After cancelling $\tau_g$ with $\tau_g^{-1}$ and using one of the
triangle identities we see that (\ref{etamodification2}) is the same
as $\eta_B * F \gamma$. Thus we conclude that
(\ref{defnmodification1}) is the same as (\ref{defnmodification2})
and that $A \mapsto \eta_A$ is a modification.

One can similarly show that $A \mapsto \varepsilon_A$ is a
modification.

The modifications $\eta$ and $\varepsilon$ satisfy the triangle
identities because the individual 2-cells $\eta_A$ and
$\varepsilon_A$ do.
\end{pf}
Now we use this lemma to prove how the right adjoints
$\psi_{X,A}:Mor_{\mathcal{X}}(X,GA) \rightarrow
Mor_{\mathcal{A}}(FX,A)$ to $\phi_{X,A}$ can be chosen in a pseudo
natural way in the following theorem.

\begin{thm}
Let $\langle F,G,\phi \rangle: \mathcal{X} \rightharpoonup
\mathcal{A}$ be a biadjunction. For all $X \in Obj \hspace{1mm}
\mathcal{X}$ and all $A \in Obj \hspace{1mm}  \mathcal{A}$ let
$\psi_{X,A}:Mor_{\mathcal{X}}(X,GA) \rightarrow
Mor_{\mathcal{A}}(FX,A)$ be a right adjoint to $\phi_{X,A}$ such
that the unit $\eta_{X,A}:1_{Mor_{\mathcal{A}}(FX,A)} \Rightarrow
\psi_{X,A} \circ \phi_{X,A}$ and the counit
$\varepsilon_{X,A}:\phi_{X,A} \circ \psi_{X,A} \Rightarrow
1_{Mor_{\mathcal{X}}(X,GA)}$ are natural isomorphisms. Then the
assignment $(X,A) \mapsto \psi_{X,A}$ is pseudo natural in each
variable. Moreover, the assignments $(X,A) \mapsto \eta_{X,A}$ and
$(X,A) \mapsto \varepsilon_{X,A}$ comprise modifications in each
variable of the form $\eta:i_{Mor_{\mathcal{A}}(F-,-)}
\rightsquigarrow \psi \odot \phi$ and $\varepsilon:\phi \odot \psi
\rightsquigarrow i_{Mor_{\mathcal{X}}(-,G-)}$.
\end{thm}
\begin{pf}
We prove the pseudo naturality and modification in the second
variable. The first variable is similar. Let $\bar{F}$ respectively
$\bar{G}$ be the pseudo functor $\mathcal{A} \rightarrow Cat$
obtained by holding $X$ fixed in the first respectively second row
in Definition \ref{defnlaxadjunction}. See the proof of Lemma
\ref{fixedXlaxnatural} for a precise description of $\bar{F}$ and
$\bar{G}$. The pseudo functor $\bar{F}$ is actually a strict
2-functor because it is the composition of strict 2-functors. If we
drop the notation $X$ in all occurrences, we see that we are
precisely in the setup of Lemma \ref{onevariable}. This proves the
theorem for the second variable. To prove it for the first variable
we only need to prove an analogue of Lemma \ref{onevariable} for $F$
pseudo and $G$ strict.
\end{pf}

Next we prove a series of lemmas needed to prove Theorems
\ref{maintheorem1} and \ref{maintheorem2}.

\begin{lem} \label{adjunctionimpliesarrow}
Let $\mathcal{X}$ and $\mathcal{A}$ be 2-categories. Let $\langle
F,G,\phi \rangle: \mathcal{X} \rightharpoonup \mathcal{A}$ be a
biadjunction and let $\eta_X:=\phi_{X,FX}(1_{FX}):X \rightarrow
GFX$. Then $\eta_X:X \rightarrow G(FX)$ is a biuniversal
arrow\index{arrow!biuniversal arrow}\index{biuniversal arrow} from
$X$ to $G$.
\end{lem}
\begin{pf}
The assignment $(X,A) \mapsto \phi_{X,A}$  is pseudo natural in each
variable by assumption. Let $\tau$ denote the coherence 2-cells for
$\phi_{X,-}$. From the definition of pseudo natural transformation
$\phi_{X,-}$ we obtain for $f' \in Mor_{\mathcal{A}}(FX,D)$ the
following diagram in $Cat$.
$$\xymatrix@C=4pc@R=4pc{Mor_{\mathcal{A}}(FX,FX) \ar[r]^{\phi_{X,FX}}
\ar[d]_{f_{\ast}'} & Mor_{\mathcal{X}}(X,GFX) \ar[d]^{(Gf')_{\ast}}
\ar@{=>}[dl]_{\tau_{FX,D}(f') \hspace{.15in}}
\\ Mor_{\mathcal{A}}(FX,D) \ar[r]_{\phi_{X,D}} & Mor_{\mathcal{X}}(X,GD)}$$
Chasing $1_{FX}$ along this diagram gives a diagram in the
2-category $\mathcal{X}$.
$$\xymatrix@R=4pc@C=9pc{X \ar[r]^{\eta_X} \ar@{=}[d]  & G(FX)
\ar[d]^{Gf'} \ar@{=>}[dl]_{\tau_{FX,D}(f')(1_{FX}) \hspace{.25in}}
\\ X \ar[r]_{\phi_{X,D}(f')} & GD}$$
The map $Mor_{\mathcal{A}}(FX,D) \ni f' \mapsto
\tau_{FX,D}(f')(1_{FX})$ is natural.  This fact combined with the
diagram in $\mathcal{X}$ above says that we have a natural
isomorphism from the functor $Mor_{\mathcal{A}}(FX,D) \ni f' \mapsto
Gf' \circ \eta_X \in Mor_{\mathcal{X}} (X,GD)$ to the functor $f'
\mapsto \phi_{X,D}(f')$. From the definition of biadjunction,
$\phi_{X,D}$ is an equivalence of categories. Hence $f' \mapsto Gf'
\circ \eta_X$ is naturally isomorphic to an equivalence of
categories and is therefore itself an equivalence of categories
$Mor_{\mathcal{A}}(FX,D) \rightarrow Mor_{\mathcal{X}}(X,GD)$. We
conclude that $\eta_X$ is a biuniversal arrow.
\end{pf}

\begin{lem}
Let $\mathcal{X}$ and $\mathcal{A}$ be 2-categories. Let $\langle
F,G,\phi \rangle: \mathcal{X} \rightharpoonup \mathcal{A}$ be a
biadjunction and let $\eta_X:=\phi_{X,FX}(1_{FX}):X \rightarrow
GFX$. Then the assignment $X \mapsto \eta_X$ is a pseudo natural
transformation $1_{\mathcal{X}} \Rightarrow GF$.
\end{lem}
\begin{pf}
Let $f:X' \rightarrow X$ be a morphism in $\mathcal{X}$. Let $\tau$
respectively $\tau'$ denote the coherence 2-cells for the pseudo
natural transformation $\phi_{X',-}$ respectively $\phi_{-,FX}$. We
must show that we have a 2-cell
$$\xymatrix@R=3pc@C=3pc{X' \ar[r]^-{\eta_{X'}} \ar[d]_f & GFX' \ar[d]^{GFf}
 \ar@{=>}[dl]
\\ X \ar[r]_-{\eta_X}  & GFX }$$
in $\mathcal{X}$ which is natural in $f$ and satisfies the
coherences involving $\delta$ and $\gamma$. Since $\phi$ is pseudo
natural in each variable we have the diagram
$$\xymatrix@C=4pc@R=4pc{Mor_{\mathcal{A}}(FX',FX') \ar[r]^{(Ff)_{\ast}}
\ar[d]_{\phi_{X',FX'}} & Mor_{\mathcal{A}}(FX',FX)
\ar[d]^{\phi_{X',FX}} & Mor_{\mathcal{A}}(FX,FX)
\ar[d]^{\phi_{X,FX}} \ar[l]_{(Ff)^{\ast}}
\\ Mor_{\mathcal{X}}(X',GFX') \ar[r]_{(GFf)_{\ast}}
\ar@{=>}[ur]^{ \tau_{FX',FX}(Ff) \hspace{.2in}} &
Mor_{\mathcal{X}}(X',GFX) & Mor_{\mathcal{X}}(X,GFX)
\ar[l]^{f^{\ast}} \ar@{=>}[ul]_{\text{\hspace{.1in}
}\tau'_{X,X'}(f^{op})} }$$ in $Cat$. By chasing $1_{FX'}$ and
$1_{FX}$ from the upper corners of this diagram to the center and
then down we see that they both get mapped to $\phi_{X',FX}(Ff)$.
Chasing the identities in the opposite directions and evaluating the
natural transformations at the identities yields a diagram of
2-cells in $\mathcal{X}$.
$$\xymatrix@C=7pc@1{(GFf) \circ \eta_{X'}
\ar@{=>}[r]^-{\tau_{FX',FX}(Ff)(1_{FX'})} & \phi_{X',FX}(Ff) &
\eta_X \circ f \ar@{=>}[l]_-{\tau_{X,X'}'(f^{op})(1_{FX})}}$$ These
2-cells are invertible by hypothesis. Let $\tilde{\tau}_{X',X}(f)$
denote the composition from left to right obtained by inverting the
second 2-cell. $\tilde{\tau}_{X',X}$ is natural in $f$ because the
constituents are natural in $f$. The coherence 2-cells
$\tilde{\tau}$ satisfy the coherences with $\delta$ and $\gamma$
from $GF$ also because the individual constituents do. Hence
$$\xymatrix@C=6pc@R=4pc{X' \ar[r]^-{\eta_{X'}} \ar[d]_f & GFX'
\ar@{=>}[dl]_{ \tilde{\tau}_{X',X}(f) \hspace{2mm}} \ar[d]^{GFf}
\\ X \ar[r]_-{\eta_X}
 & GFX }$$
is natural in $f$ and satisfies the required coherences, so $X
\mapsto \eta_X$ is a pseudo natural transformation.
\end{pf}

Thus we have seen that given a biadjunction\index{biadjunction}
$\phi$ we get a pseudo natural transformation $\eta$ whose arrows
are biuniversal\index{biuniversal arrow}\index{arrow!biuniversal
arrow} arrows. Now we consider the converse of this statement.

\begin{lem} \label{converse}
Let $\mathcal{X}$ and $\mathcal{A}$ be 2-categories. Let
$\xymatrix@1{ \mathcal{X}  \ar@<.5ex>[r]^F & \mathcal{A}
\ar@<.5ex>[l]^G}$ be pseudo functors between 2-categories. Let
$\eta:1_{\mathcal{X}} \Rightarrow GF$ be a pseudo natural
transformation such that each arrow $\eta_X:X \rightarrow G(FX)$ is
a biuniversal\index{biuniversal arrow}\index{arrow!biuniversal
arrow} arrow from $X$ to $G$. Define $\phi_{X,A}(f):=Gf \circ
\eta_X$ for each $f:FX \rightarrow A$ and $\phi_{X,A}(\gamma):=G
\gamma \ast i_{\eta_X}$ for each $\gamma:f \Rightarrow f'$. Then
$\phi_{X,A}: Mor_{\mathcal{A}}(FX,A) \rightarrow
Mor_{\mathcal{X}}(X,GA)$ is an equivalence of categories for all $X
\in Obj \hspace{1mm}  \mathcal{X}$ and all $A \in Obj \hspace{1mm}
\mathcal{A}$.
\end{lem}
\begin{pf}
The functor $\phi_{X,A}$ is an equivalence since $\eta_X$ is a
biuniversal arrow.
\end{pf}

\begin{lem}
Let $\mathcal{X}$ and $\mathcal{A}$ be 2-categories. Let
$\xymatrix@1{ \mathcal{X}  \ar@<.5ex>[r]^F & \mathcal{A}
\ar@<.5ex>[l]^G}$ be pseudo functors between 2-categories. Let
$\eta:1_{\mathcal{X}} \Rightarrow GF$ be a pseudo natural
transformation such that each $\eta_X:X \rightarrow G(FX)$ is a
biuniversal\index{biuniversal arrow}\index{arrow!biuniversal arrow}
arrow from $X$ to $G$. Let $\phi_{X,A}$ be defined as in Lemma
\ref{converse} above. Then for fixed $A \in Obj \hspace{1mm}
\mathcal{A}$ the assignment $Obj \hspace{1mm} \mathcal{X}^{op} \ni X
\mapsto \phi_{X,A}$ denoted $\phi_{-,A}$ is pseudo natural.
\end{lem}
\begin{pf}
Let $A \in Obj \hspace{1mm}  \mathcal{A}$ be a fixed object
throughout this proof. Let $\bar{F}:\mathcal{X}^{op} \rightarrow
Cat$ denote the pseudo functor obtained by holding $A$ fixed in the
top row in the definition of biadjunction. This means
$\bar{F}(X)=Mor_{\mathcal{A}}(FX,A)$, $\bar{F}(f^{op})=(Ff)^{\ast}$,
and for $\alpha:f^{op} \Rightarrow (f')^{op}$ in $\mathcal{X}$ the
natural transformation $\bar{F}(\alpha):(Ff)^{\ast} \Rightarrow
(Ff')^{\ast}$ is $h \mapsto i_h \ast F \alpha$. Note that the
morphisms of $\mathcal{X}^{op}$ are formally the opposites of
morphisms of $\mathcal{X}$, but the 2-cells of $\mathcal{X}^{op}$
are precisely the same as the 2-cells in $\mathcal{X}$. The vertical
composition is the same in both $\mathcal{X}^{op}$ and
$\mathcal{X}$, although the horizontal compositions are switched.
The pseudo functor $\bar{F}$ is the composition of a pseudo functor
and a strict functor. For morphisms $\xymatrix@1{X \ar[r]^f & Y
\ar[r]^g & Z}$ in $\mathcal{X}$ we have
$\gamma_{g^{op},f^{op}}^{\bar{F}}: h \mapsto i_h \ast
\gamma^F_{f,g}$ and for $X \in Obj \hspace{1mm} \mathcal{X}^{op}$ we
have $\delta_{X \ast}^{\bar{F}}: h \mapsto i_h \ast \delta^F_{X
\ast}$ by the rules for composition of pseudo functors. Then
$\gamma_{g^{op},f^{op}}^{\bar{F}}: \bar{F}(f^{op}) \circ
\bar{F}(g^{op}) \Rightarrow \bar{F}(f^{op} \circ g^{op})$ and
$\delta_{X \ast}^{\bar{F}}:1_{\bar{F}X} \Rightarrow \bar{F}(1_X)$.
Let $\bar{G}$ denote the strict 2-functor obtained by holding $A$
fixed in the bottom row in the definition of biadjunction. This
means $\bar{G}(X)=Mor_{\mathcal{X}}(X,GA)$,
$\bar{G}(f^{op})=f^{\ast}$, and for $\alpha:f^{op} \Rightarrow
(f')^{op}$ in $\mathcal{X}$ the natural transformation
$\bar{G}(\alpha):\bar{G}(f^{op}) \Rightarrow \bar{G}((f')^{op})$ is
the natural transformation $h \mapsto i_h \ast \alpha$. The
2-functor $\bar{G}$ is the composition of two strict 2-functors and
is therefore strict.

In order to prove that $\phi_{-,A}$ is a pseudo natural
transformation from $\bar{F}$ to $\bar{G}$ we must display coherence
2-cells $\tau'$ up to which $\phi_{-,A}$ is natural and prove that
they satisfy the coherences involving $\delta$ and $\gamma$. Now we
describe this $\tau'$ and later prove the coherences. Let
$\tilde{\tau}$ denote the coherence 2-cells which make
$\eta:1_{\mathcal{X}} \Rightarrow GF$ pseudo natural, \ie for all
$f:X \rightarrow Y$ in $\mathcal{X}$ we have
$$\xymatrix@C=6pc@R=4pc{X \ar[r]^{\eta_X} \ar[d]_{f}
& GFX \ar[d]^{GFf} \ar@{=>}[ld]_{\hspace{.1in} \tilde{\tau}_{X,Y}(f)
\hspace{2mm}}
\\ Y \ar[r]_{\eta_Y} & GFY}$$
in $\mathcal{X}$.  Define a \label{firstcoherenceisoforphi} natural
isomorphism $\tau'_{f^{op}}=\tau'_{Y,X}(f^{op}):\bar{G}(f^{op})
\circ \phi_{Y,A} \Rightarrow \phi_{X,A} \circ \bar{F}(f^{op})$ by $h
\mapsto (\gamma^G_{Ff,h} \ast i_{\eta_X}) \odot (i_{Gh} \ast
(\tilde{\tau}_{X,Y}(f))^{-1})$ for $h \in Mor_{\mathcal{A}}(FY,A)$
as in the following diagram.
$$\xymatrix@R=4pc@C=6pc{X \ar[r]^{\eta_X} \ar[d]_f & GFX
\ar@{=>}[dl]_{\tilde{\tau}_{X,Y}(f) \hspace{2mm}} \ar[d]_{GFf} \ar@{=}[r] & GFX \ar[dd]^{G(h \circ Ff)} \\
Y \ar[r]_{\eta_Y} & GFY \ar[d]_{Gh}  &
\\ & GA \ar@{=}[r] \ar@{=>}[uru]|{\gamma_{Ff,h}^G} & GA }$$
The map $\tau'_{Y,X}(f^{op})$ is a natural transformation because
$\gamma^G_{Ff,h}$ is natural in $h$. The assignment $f^{op} \mapsto
\tau'_{Y,X}(f^{op})$ is also natural for a similar reason.

We claim that $\tau'$ satisfies the unit axiom for pseudo natural
transformations. We must show that the diagram of 2-cells in $Cat$
\begin{equation} \label{unit1}
\xymatrix@C=5pc@R=3pc{\phi_{X,A} \ar@{=}[r] \ar@{=}[d] &
1_{\bar{G}X} \circ \phi_{X,A} \ar@{=}[r] & \bar{G}(1_X) \circ
\phi_{X,A} \ar@{=>}[d]^{\tau_{1_X^{op}}'}
\\ \phi_{X,A} \circ 1_{\bar{F}X} \ar@{=>}[rr]_{i_{\phi_{X,A}} \ast
\delta^{\bar{F}}_{X \ast}} & & \phi_{X,A} \circ \bar{F}(1_X)}
\end{equation}
commutes for all $X \in Obj \hspace{1mm}  \mathcal{X}$. After we
evaluate this diagram on a morphism $h:FX \rightarrow A$ of
$\mathcal{A}$ we obtain the diagram of 2-cells
\begin{equation} \label{unit2}
\xymatrix@R=3pc@C=3pc{Gh \circ \eta_X \ar@{=}[r] \ar@{=}[dd] & Gh
\circ \eta_X \ar@{=}[r] & Gh \circ \eta_X \circ 1_X
\ar@{=>}[d]^{i_{Gh} \ast (\tilde{\tau}_{X,X}(1_X))^{-1}}
\\
& & Gh \circ GF1_X \circ \eta_X \ar@{=>}[d]^{\gamma^G_{F1_X,h} \ast
i_{\eta_X}}
\\G(h \circ 1_{FX}) \circ \eta_X \ar@{=>}[rr]_{G(i_h \ast \delta^F_{X \ast})
\ast i_{\eta_X}} & & G(h \circ F(1_X)) \circ \eta_X}
\end{equation}
in $\mathcal{X}$.  Since $\eta:1_{\mathcal{X}} \Rightarrow GF$ is a
pseudo natural transformation from the strict 2-functor to the
composition $G \circ F$ of pseudo functors, its unit axiom for
$\tilde{\tau}$ simplifies to the following commutative diagram.
$$\xymatrix@C=8pc@R=3pc{1_{GFX} \circ \eta_X \ar@{=>}[r]^{(G(\delta^F_{X\ast}) \odot
\delta^G_{FX \ast}) \ast i_{\eta_X}} \ar@{=}[rd]
 & GF1_X \circ \eta_X \ar@{=>}[d]^{\tilde{\tau}_{X,X}(1_X)}
\\ & \eta_X \circ 1_X}$$
Hence $(\tilde{\tau}_{X,X}(1_X))^{-1}=(G(\delta^F_{X\ast}) \odot
\delta^G_{FX \ast}) \ast i_{\eta_X}$ as 2-cells. Note also that
$\delta^{GF}_{X \ast}=(G(\delta^F_{X\ast}) \odot \delta^G_{FX
\ast})$ by the definition of composition of pseudo functors. Using
this, we see that diagram (\ref{unit2}) becomes the outermost
rectangle of the following diagram.
$$\xymatrix@C=5pc@R=3pc{Gh \circ 1_{GFX} \circ \eta_X
\ar@{=>}[rr]^{i_{Gh} \ast \delta^{GF}_{X \ast} \ast i_{\eta_X}}
\ar@{=>}[dd]_{i_{Gh} \ast i_{\eta_X}} \ar@{=>}[dr]^{\hspace{.2in}
i_{Gh} \ast \delta^G_{FX \ast} \ast i_{\eta_X}} & & Gh \circ GF1_X
\circ \eta_X \ar@{=>}[dd]^{\gamma^G_{F1_X, h} \ast i_{\eta_X}}
\\ &
Gh \circ G1_{FX} \circ \eta_X \ar@{=>}[ur]^{i_{Gh} \ast
G(\delta^F_{X \ast}) \ast i_{\eta_X} \hspace{.4in}}
\ar@{=>}[ld]^{\hspace{.2in} \gamma^G_{1_{FX},h} \ast i_{\eta_X}}
 &
\\ G(h \circ 1_{FX}) \circ \eta_X
\ar@{=>}[rr]_{G(i_h \ast \delta^F_{X \ast}) \ast i_{\eta_X}} & & G(h
\circ F1_X) \circ \eta_X }$$ The upper left vertex of this diagram
is the upper right vertex of diagram (\ref{unit2}) and the
composition of the top arrow and right vertical arrow of this
diagram is the right vertical arrow of diagram (\ref{unit2}). The
top triangle of this diagram commutes by definition. The left
triangle commutes by the unit axiom of the pseudo functor $G$
applied to the morphism $h:FX \rightarrow A$ of $\mathcal{A}$. The
right quadrilateral commutes by the naturality of $\gamma_{-,h}^G$
and because $G(i_h*\delta^F_{X*})=i_{Gh}*G(\delta^F_{X*})$. The
morphism $\eta_X$ and the 2-cell $i_{\eta_X}$ just tag along. Hence
the outermost rectangle commutes and diagram (\ref{unit2}) commutes.
This implies that diagram (\ref{unit1}) commutes. We conclude that
$\tau'$ satisfies the unit axiom required for $\phi_{-,A}$ to be a
pseudo natural transformation.

We claim that $\tau'$ satisfies the composition axiom required for
$\phi_{-,A}$ to be a pseudo natural transformation. We must prove
for all morphisms $\xymatrix@1{X \ar[r]^f & Y \ar[r]^g & Z}$ of
$\mathcal{X}$, \ie for all morphisms $\xymatrix@1{Z \ar[r]^{g^{op}}
& Y \ar[r]^{f^{op}} & X}$ of $\mathcal{X}^{op}$, the diagram of
2-cells in $Cat$
\begin{equation} \label{comp1}
\xymatrix@R=3pc@C=2pc{\bar{G}(f^{op}) \circ \bar{G}(g^{op}) \circ
\phi_{Z,A} \ar@{=>}[r] \ar@{=>}[d] & \bar{G}(f^{op}) \circ
\phi_{Y,A} \circ \bar{F}(g^{op}) \ar@{=>}[r] & \phi_{X,A} \circ
\bar{F}(f^{op}) \circ \bar{F}(g^{op}) \ar@{=>}[d]
\\ \bar{G}(f^{op} \circ g^{op}) \circ \phi_{Z,A}
\ar@{=>}[rr] & & \phi_{X,A} \circ \bar{F}(f^{op} \circ g^{op})}
\end{equation}
commutes. More precisely the diagram of 2-cells in $Cat$
\begin{equation} \label{comp2}
\xymatrix@C=4pc@R=3pc{f^* \circ g^* \circ \phi_{Z,A}
\ar@{=>}[r]^-{i_{f^*} \ast \tau'_{g^{op}}} \ar@{=}[d] & f^* \circ
\phi_{Y,A} \circ (Fg)^* \ar@{=>}[r]^-{\tau'_{f^{op}} \ast
i_{(Fg)^*}} & \phi_{X,A} \circ (Ff)^* \circ (Fg)^*
\ar@{=>}[d]^{i_{\phi_{X,A}} \ast \gamma^{\bar{F}}_{g^{op}, f^{op}}}
\\ (g \circ f)^* \circ \phi_{Z,A} \ar@{=>}[rr]_{\tau'_{f^{op} \circ g^{op}}}
& & \phi_{X,A} \circ (F(g \circ f))^*}
\end{equation}
must commute.  We evaluate this diagram on a morphism $h:FZ
\rightarrow A$ of $\mathcal{A}$, fill in the diagram with more
vertices, and cut the result down the middle column to get the left
respectively right half on page \pageref{biggestdiagram}. These are
diagrams of 2-cells in $\mathcal{X}$. Subdiagram (I) commutes by the
composition axiom applied to the morphisms $\xymatrix@1{X \ar[r]^f &
Y \ar[r]^g & Z}$ for the pseudo natural transformation
$\eta:1_{\mathcal{X}} \Rightarrow GF$ with its coherence 2-cells
$\tilde{\tau}$. Subdiagram (II) commutes by the composition axiom
applied to the morphisms $Ff,Fg,h$ for the pseudo functor $G$ with
its coherence 2-cells $\gamma^G$. The fifth arrow which is an
equality symbol was only drawn for convenience. Subdiagram (III)
commutes by the naturality of $\gamma^G$. All other subdiagrams
commute by definition or by the interchange law. Therefore the
outermost rectangle commutes when we put the two halves together.
This outermost rectangle is diagram (\ref{comp2}) evaluated on the
morphism $h: FZ \rightarrow A$ of $\mathcal{A}$. Hence (\ref{comp2})
and (\ref{comp1}) commute. We conclude that $\tau'$ satisfies the
composition axiom required for $\phi_{-,A}$ to be a pseudo natural
transformation.

Since $\phi_{-,A}$ with coherence 2-cells $\tau'$ satisfies the unit
axiom and composition axiom for pseudo natural transformations we
conclude that $\phi_{-,A}$ is a pseudo natural transformation for
fixed $A \in Obj \hspace{1mm}  \mathcal{A}$.

\fig{}{
\begingroup
\vspace{-2\abovedisplayskip} \normalsize
\label{biggestdiagram} $$\xymatrix@R=3pc@C=3pc{Gh \circ \eta_Z \circ
g \circ f \ar@{=>}[rr]^{\tau_{g^{op}}'(h)*i_f}
\ar@{=>}[dr]^{\hspace{.4in} i_{Gh}
* (\tilde{\tau}_{Y,Z}(g))^{-1} \ast i_f} \ar@{=}[dddd] & & G(h
\circ Fg) \circ \eta_Y \circ f
\\
& Gh \circ GFg \circ \eta_Y \circ f \ar@{=>}[ur]_{\hspace{.2in}
\gamma^G_{Fg,h} * i_{\eta_Y} * i_f} \ar@{=>}[dr]^{\hspace{.5in}
i_{Gh} *i_{GFg} * (\tilde{\tau}_{X,Y}(f))^{-1}} &
\\
& \text{(I)} & Gh \circ GFg \circ GFf \circ \eta_X
\ar@{=>}[d]^{i_{Gh} * \gamma^G_{Ff,Fg} * i_{\eta_X}}
\ar@{=>}@/_2pc/@<-10ex>[dd]_{i_{Gh} * \gamma^{GF}_{f,g}}
\\
& & Gh \circ G(Fg \circ Ff) \circ \eta_X \ar@{=>}[d]^{i_{Gh} *
G(\gamma^F_{f,g}) \ast i_{\eta_X}}
\\ Gh \circ \eta_Z \circ g \circ f
\ar@{=>}[rr]_{i_{Gh} * (\tilde{\tau}_{X,Z}(g \circ f))^{-1}}
\ar@{=>}[ruuu]_{i_{Gh} * (\tilde{\tau}_{Y,Z}(g))^{-1}*i_f} & & Gh
\circ GF(g \circ f) \circ \eta_X}
$$
$$
\xymatrix@R=3pc@C=1.8pc{G(h \circ Fg) \circ \eta_Y \circ f
\ar@{=>}[rr]^{\tau_{f^{op}}'(h \circ Fg)} \ar@{=>}[rd]_{i_{G(h \circ
Fg)}
* (\tilde{\tau}_{X,Y}(f))^{-1} \hspace{.5in}} & & G(h \circ Fg
\circ Ff) \circ \eta_X
\ar@{=>}[dddd]|{\overset{\phantom{e}}{\underset{\phantom{e}}{G(i_h
* \gamma^F_{f,g})
* i_{\eta_X}}}}
\\
& G( h \circ Fg) \circ GFf \circ \eta_X \ar@{=>}[ur]^{\gamma^G_{Ff,
h \circ Fg} * i_{\eta_X} \hspace{.2in}} \ar@{}[d]|{\text{(II)}} &
\\ Gh \circ GFg \circ GFf \circ \eta_X
\ar@{=>}[d]_{i_{Gh} * \gamma^G_{Ff,Fg} * i_{\eta_X}}
\ar@{=>}[ru]^{\gamma_{Fg,h}^G * i_{GFf} * i_{\eta_X} \hspace{.5in}}
& G(h \circ Fg \circ Ff) \circ \eta_X \ar@{=}[ruu] &
\\ Gh \circ G(Fg \circ Ff) \circ \eta_X
\ar@{=>}[ur]_{\hspace{.3in} \gamma^G_{Fg \circ Ff,h} * i_{\eta_X}}
\ar@{=>}[d]_{i_{Gh} * G(\gamma^F_{f,g}) \ast i_{\eta_X}} &
\text{(III)} &
\\ Gh \circ GF(g \circ f) \circ \eta_X
\ar@{=>}[rr]_{\gamma^G_{F(g \circ f),h} * i_{\eta_X}} & & G(h \circ
F(g \circ f)) \circ \eta_X }
$$
\endgroup
\noindent}
\end{pf}

\begin{lem} \label{fixedXlaxnatural}
Let $\mathcal{X}$ and $\mathcal{A}$ be 2-categories. Let
$\xymatrix@1{ \mathcal{X}  \ar@<.5ex>[r]^F & \mathcal{A}
\ar@<.5ex>[l]^G}$ be pseudo functors between 2-categories. Let
$\eta:1_{\mathcal{X}} \Rightarrow GF$ be a pseudo natural
transformation such that each $\eta_X:X \rightarrow G(FX)$ is a
biuniversal\index{biuniversal arrow}\index{arrow!biuniversal arrow}
arrow from $X$ to $G$. Let $\phi_{X,A}$ be defined as in Lemma
\ref{converse} above. Then for fixed $X \in Obj \hspace{1mm}
\mathcal{X}$ the assignment $Obj \hspace{1mm} \mathcal{A} \ni A
\mapsto \phi_{X,A}$ denoted $\phi_{X,-}$ is pseudo natural.
\end{lem}
\begin{pf}
Let $X$ be a fixed object of the 2-category $\mathcal{X}$ throughout
the proof. We introduce new pseudo functors $\bar{F}$ and $\bar{G}$
different from those in the previous proof.
 Let $\bar{F}:\mathcal{A} \rightarrow Cat$ be the strict
2-functor obtained by fixing $X$ in the top row in the definition of
biadjunction. This means $\bar{F}(A)=Mor_{\mathcal{A}}(FX,A)$,
$\bar{F}(f)=f_{\ast}$, and for $\alpha:f \Rightarrow f'$ we have
$\bar{F}(\alpha)$ is the natural transformation $e \mapsto \alpha
\ast i_e$. The 2-functor $\bar{F}$ is strict because it is the
composition of two strict 2-functors. Similarly let
$\bar{G}:\mathcal{A} \rightarrow Cat$ be the pseudo functor obtained
by fixing $X$ in the bottom row of the definition of biadjunction.
This means $\bar{G}(A)=Mor_{\mathcal{X}}(X,GA)$,
$\bar{G}(f)=(Gf)_{\ast}$, and for $\alpha:f \Rightarrow f'$ we have
$\bar{G}(\alpha)$ is the natural transformation $e \mapsto G(\alpha)
\ast i_e$.  The pseudo functor $\bar{G}$ is pseudo because it is the
composition of a pseudo functor and a strict functor. The definition
of composition of pseudo functors then says that the coherence
2-cells for $\bar{G}$ are $\gamma^{\bar{G}}_{f,g}:e \mapsto
\gamma^G_{f,g} \ast i_e$ for morphisms $f,g$ of $\mathcal{A}$ such
that $g \circ f$ exists and $\delta^{\bar{G}}_{A \ast}: e \mapsto
\delta^G_{A \ast} \ast i_e$ for $A \in Obj \hspace{1mm}
\mathcal{A}$. These are natural transformations, \ie 2-cells in
$Cat$, such that $\gamma_{f,g}^{\bar{G}}:\bar{G}(g) \circ \bar{G}(f)
\Rightarrow \bar{G}(g \circ f)$ and $\delta^{\bar{G}}_{A
\ast}:1_{\bar{G}(A)} \Rightarrow \bar{G}(1_A)$. They are natural in
$f$ and $g$ and they satisfy the required coherences for a pseudo
functor.

We must show that $\phi_{X,-}$ is a pseudo natural transformation
from $\bar{F}$ to $\bar{G}$. In other words we must display
coherence 2-cells $\tau$ up to which $\phi_{X,-}$ is natural and
satisfy the coherence diagrams involving $\gamma$ and $\delta$ from
$\bar{F}$ and $\bar{G}$. For morphisms $k:A \rightarrow A'$ of
$\mathcal{A}$ define $\tau_{A,A'}(k):e \mapsto \gamma^{G}_{e,k} \ast
i_{\eta_X}$ to fill in the \label{secondcoherenceisoforphi} diagram
$$\xymatrix@R=3pc@C=3pc{Mor_{\mathcal{A}}(FX,A) \ar[r]^{\phi_{X,A}}
\ar[d]_{k_{\ast}} & Mor_{\mathcal{X}}(X,GA) \ar[d]^{(Gk)_{\ast}}
\ar@{=>}[ld]_{\tau_{A,A'}(k) \hspace{.1in}}
\\ Mor_{\mathcal{A}}(FX,A') \ar[r]_{\phi_{X,A'}}
 & Mor_{\mathcal{A}}(X,GA')}$$
whose vertices are $\bar{F}(A),\bar{G}(A),\bar{G}(A')$, and
$\bar{F}(A')$ read clockwise. The map $\tau_{A,A'}(k)$ is a natural
transformation (2-cell in $Cat$) between the indicated functors
because $\gamma^G_{e,k}$ is natural in $e$. The assignment
$Mor_{\mathcal{A}}(A,A') \ni k \mapsto \tau_{A,A'}(k)$ is a natural
transformation $(\circ \phi_{X,A}) \circ \bar{G} \Rightarrow
(\phi_{X,A'} \circ) \circ \bar{F}$ because $\gamma^{G}_{e,k}$ is
natural in $k$. Hence this family $\tau$ of natural transformations
provides us with a candidate for the coherence 2-cells to make
$\phi_{X,-}$ into a pseudo natural transformation.

We claim that $\tau$ satisfies the unit axiom for pseudo natural
transformations. This requires a proof that the diagram of 2-cells
in $Cat$
$$\xymatrix@C=6pc@R=3pc{\phi_{X,A} \ar@{=>}[r]^-{i_{\phi_{X,A}}} \ar@{=}[d]
& 1_{\bar{G}A} \circ \phi_{X,A} \ar@{=>}[r]^{\delta^{\bar{G}}_{A
\ast} \ast i_{\phi_{X,A}}} & \bar{G}(1_A) \circ \phi_{X,A}
\ar@{=>}[d]^{\tau_{1_A}}
\\ \phi_{X,A} \circ 1_{\bar{F}A} \ar@{=}[rr]
& & \phi_{X,A} \circ \bar{F}(1_A)}$$ commutes for all $A \in Obj
\hspace{1mm} \mathcal{A}$. Evaluating this diagram on a morphism
$e:FX \rightarrow A$ of $\mathcal{A}$ results in the diagram of
2-cells
$$\xymatrix@C=6pc@R=3pc{Ge \circ \eta_X \ar@{=}[r] \ar@{=}[drr] & 1_{GA}
\circ Ge \circ \eta_X \ar@{=>}[r]^{\delta^G_{A \ast} \ast i_{Ge}
\ast i_{\eta_X}} & G(1_A) \circ Ge \circ \eta_X
\ar@{=>}[d]^{\gamma_{e,1_A}^G \ast i_{\eta_X}}
\\ & & G(1_A \circ e) \circ \eta_X}$$
in $\mathcal{X}$ which commutes because of the unit axiom for the
pseudo functor $G$. Hence $\tau$ satisfies the unit axiom for pseudo
natural transformations.

We claim that $\tau$ satisfies the composition axiom for pseudo
natural transformations. This requires us to prove for all morphisms
$\xymatrix@1{A \ar[r]^f & B \ar[r]^g & C}$ in $\mathcal{A}$ that the
diagram of 2-cells in $Cat$
$$\xymatrix@C=5pc@R=3pc{\bar{G}g \circ \bar{G}f \circ \phi_{X,A}
\ar@{=>}[r]^-{i_{\bar{G}g} \ast \tau_f}
\ar@{=>}[d]_{\gamma^{\bar{G}}_{f,g} \ast i_{\phi_{X,A}}} & \bar{G}g
\circ \phi_{X,B} \circ \bar{F}f \ar@{=>}[r]^{\tau_g \ast
i_{\bar{F}f}} & \phi_{X,C} \circ \bar{F}g \circ \bar{F}f \ar@{=}[d]
\\ \bar{G}(g \circ f) \circ \phi_{X,A} \ar@{=>}[rr]_{\tau_{g \circ f}}
& & \phi_{X,C} \circ \bar{F}(g \circ f)}$$ commutes. Evaluating this
diagram on a morphism $e:FX \rightarrow A$ of $\mathcal{A}$ results
in the diagram of 2-cells

\begingroup
\vspace{-2\abovedisplayskip} \normalsize
$$\xymatrix@C=3.6pc@R=3pc{Gg \circ Gf \circ Ge \circ \eta_X \ar@{=>}[r]^{\underset{\phantom{e}}{i_{Gg} \ast
\gamma^G_{e,f} \ast i_{\eta_X}}}  \ar@{=>}[d]_{\gamma^G_{f,g} \ast
i_{Ge} \ast i_{\eta_X}} & Gg \circ G(f \circ e) \circ \eta_X
\ar@{=>}[r]^-{\underset{\phantom{e}}{\gamma^G_{f \circ e,g} \ast
i_{\eta_X}}} & G(g \circ f \circ e) \circ \eta_X \ar@{=}[d]
\\ G(g \circ f) \circ G(e) \circ \eta_X \ar@{=>}[rr]_{\overset{\phantom{e}}{\gamma^G_{e,g \circ
f}*i_{\eta_x}}} & & G(g \circ f \circ e) \circ \eta_X}$$
\endgroup
\noindent in $\mathcal{X}$, which commutes by the composition axiom
for the pseudo functor $G$ applied to $\xymatrix@1{FX \ar[r]^e & A
\ar[r]^f & B \ar[r]^g & C}$. Hence $\tau$ satisfies the composition
axiom for pseudo natural transformations.

We conclude that $\phi_{X,-}$ is a pseudo natural transformation
from $\bar{F}$ to $\bar{G}$ with coherence 2-cells defined by
$\tau$.
\end{pf}

Now we can finally state and prove the two main theorems of this
chapter.

\begin{thm} \label{maintheorem1}
Let $\mathcal{X}$ and $\mathcal{A}$ be 2-categories. Let
$\xymatrix@R=3pc@C=3pc{ \mathcal{X}  \ar@<.5ex>[r]^F & \mathcal{A}
\ar@<.5ex>[l]^G}$ be pseudo functors. Then $F$ is a left
biadjoint\index{biadjoint!left biadjoint|textbf} for $G$ if and only
if there exists a pseudo natural transformation
$\eta:1_{\mathcal{X}} \Rightarrow GF$ such that $\eta_X:X
\rightarrow G(FX)$ is a biuniversal\index{biuniversal
arrow|textbf}\index{arrow!biuniversal arrow|textbf} arrow for all $X
\in Obj \hspace{1mm}  \mathcal{X}$.
\end{thm}
\begin{pf}
This follows immediately from the previous lemmas.
\end{pf}

\begin{thm} \label{maintheorem2}
Let $\mathcal{X}$ and $\mathcal{A}$ be 2-categories. Let
$\xymatrix@1{\mathcal{X}   & \mathcal{A} \ar[l]_G}$ be a pseudo
functor. Then there exists a left biadjoint\index{biadjoint!left
biadjoint|textbf} for $G$ if and only if for every object $X \in Obj
\hspace{1mm} \mathcal{X}$ there exists an object $R \in Obj
\hspace{1mm} \mathcal{A}$ and a biuniversal\index{biuniversal
arrow|textbf}\index{arrow!biuniversal arrow|textbf} arrow $\eta_X:X
\rightarrow G(R)$ from $X$ to $G$.
\end{thm}
\begin{pf}
By Lemma \ref{adjunctionimpliesarrow}, the existence of a left
biadjoint implies the existence of such a biuniversal arrow. Now we
prove the other direction. Suppose we have such a biuniversal arrow
for each $X \in Obj \hspace{1mm}  \mathcal{X}$. Define $FX:=R$. The
object $R \in Obj \hspace{1mm}  \mathcal{A}$ of course depends on
$X$. For $X \in Obj \hspace{1mm}  \mathcal{X}$ and $A \in Obj
\hspace{1mm}  \mathcal{A}$ let $\phi_{X,A}:Mor_{\mathcal{A}}(FX,A)
\rightarrow Mor_{\mathcal{X}}(X,GA)$ denote the functor $f' \mapsto
Gf' \circ \eta_X$ and $\alpha \mapsto G \alpha * i_{\eta_X}$. Let
$\psi_{X,A}:Mor_{\mathcal{X}}(X,GA) \rightarrow
Mor_{\mathcal{A}}(FX,A)$ denote a right adjoint equivalence, which
exists because $\eta_X$ is a biuniversal arrow. Let
$\mu_{X,A}:\phi_{X,A} \circ \psi_{X,A} \Rightarrow
1_{Mor_{\mathcal{X}}(X,GA)}$ denote a counit for these adjoint
functors.  All of this implies that for any morphism $f:X
\rightarrow GA$ there exists a morphism $f':=\psi_{X,A}(f)$ and a
2-cell $\mu_{X,A}(f)$ as in the diagram.
$$\xymatrix@C=5pc@R=3pc{X \ar@{=}[d] \ar[r]^{\eta_X} & G(FX)  \ar@{.>}[d]^{Gf'}
\ar@{:>}[dl]_{\mu_{X,A}(f) } & FX \ar@{.>}[d]^{f'} \\ X
 \ar[r]_{f}
 & GA & A}$$
Moreover, this 2-cell  $\mu_{X,A}(f)$ is a universal arrow from the
functor $\phi_{X,A} \circ \psi_{X,A}$ to the object $f$ because all
of the arrows of the counit of an adjunction are universal. This
means that for any other morphism $\bar{f}':FX \rightarrow A$ and
2-cell $\nu$ as in the diagram
$$\xymatrix@C=5pc@R=3pc{X \ar@{=}[d] \ar[r]^{\eta_X} & G(FX)  \ar[d]^{G
\bar{f}'} \ar@{=>}[dl]_{\nu } & FX \ar[d]^{\bar{f}'} \\ X
 \ar[r]_{f}
 & GA & A}$$
there exists a unique 2-cell $\nu':\bar{f}' \Rightarrow f'$ such
that the following diagram commutes.
$$\xymatrix@C=5pc@R=3pc{\bar{f}' \ar@{:>}[d]_{\nu'}
& G\bar{f}' \circ \eta_X \ar@{=>}[r]^-{\nu} \ar@{:>}[d]_{G \nu' *
i_{\eta_X}} & f \ar@{=}[d]
\\ f' & Gf' \circ \eta_X \ar@{=>}[r]_-{\mu_{X,A}(f)} & f}$$
If $\nu$ is iso, this 2-cell $\nu':\bar{f}' \Rightarrow f'$ is also
iso by the comments after Definition \ref{deflaxuniversalarrow}. The
uniqueness and iso property of $\nu'$ will be integral to defining
the coherence isomorphisms and proving the coherence diagrams below.

After setting up this notation, we define a left biadjoint candidate
$F$ for $G$. We already have $F$ defined for objects $X \in Obj
\hspace{1mm}  \mathcal{X}$ above.  For any morphism $h:X \rightarrow
Y$ in $\mathcal{X}$ define $Fh:=\psi_{X,FY}(\eta_Y \circ h)$. For
morphisms $h,h':X \rightarrow Y$ and any 2-cell $\alpha:h
\Rightarrow h'$ in $\mathcal{X}$ define
$F\alpha:=\psi_{X,FY}(i_{\eta_Y} \ast \alpha)$. Then the assignment
is obviously a functor on any fixed hom category because of the
interchange law and because $\psi_{X,FY}$ preserves identity 2-cells
and compositions of 2-cells. To define the coherence 2-cells
$\delta^F_X$ we now use the uniqueness described above. Note that
$F1_X=\psi_{X,FX}(\eta_X \circ 1_X)$ satisfies the diagram
$$\xymatrix@C=8pc@R=3pc{X \ar[r]^{\eta_X} \ar[d]_{1_X}
& GFX \ar@{:>}[dl]_{\mu_{X,FX}(\eta_X \circ 1_X) \hspace{.3in}}
\ar@{.>}[d]^{GF1_X} & FX \ar@{.>}[d]^{F1_X}
\\ X \ar[r]_{\eta_X}
& GFX & FX}$$ where $\mu_{X,FX}(\eta_X \circ 1_X)$ is universal. The
arrow $1_{FX}$ satisfies
$$\xymatrix@C=8pc@R=3pc{X \ar[r]^{\eta_X} \ar[d]_{1_X}
& GFX \ar@{=>}[dl]_{(\delta_{FX*}^G)^{-1} * i_{\eta_X}
\hspace{.3in}} \ar[d]^{G1_{FX}} & FX \ar[d]^{1_{FX}}
\\ X \ar[r]_{\eta_X}
& GFX & FX}$$ since $G$ is a pseudo functor. Let
$\delta_{X*}^F:1_{FX} \Rightarrow F1_X$ be the unique 2-cell whose
$\phi_{X,FX}$ image factors $(\delta_{FX*}^G)^{-1} * i_{\eta_X}$.
\begin{equation} \label{Fdelta}
\xymatrix@C=6pc@R=3pc{1_{FX} \ar@{:>}[d]_{\delta_{X*}^F} & G1_{FX}
\circ \eta_X \ar@{=>}[r]^-{(\delta_{FX*}^G)^{-1} * i_{\eta_X}}
\ar@{:>}[d]_{G(\delta^F_{X*}) \ast i_{\eta_X}} & 1_{GFX} \circ
\eta_X \ar@{=}[d]
\\ F1_X
& GF1_X \circ \eta_X \ar@{=>}[r]_-{\mu_{X,FX}(\eta_X \circ 1_X)} &
\eta_X \circ 1_X}
\end{equation}
It exists by the universality of $\mu_{X,FX}(\eta_X \circ 1_X)$. The
2-cell  $\delta_{X*}^F:1_{FX} \Rightarrow F1_X$ is iso because
$(\delta_{FX*}^G)^{-1} * i_{\eta_X}$ is iso. To define
$\gamma^F_{f,g}$ for $\xymatrix@1{X \ar[r]^f & Y \ar[r]^g & Z}$ in
$\mathcal{X}$ we similarly use the uniqueness. Note that $F(g \circ
f) =\psi_{X,FZ}(\eta_Z \circ g \circ f)$ satisfies the diagram
$$\xymatrix@C=4pc@R=3pc{X \ar[d]_f \ar[r]^{\eta_X}
& GFX \ar@{.>}[dd]^{GF(g \circ f)} \ar@{:>}[ldd]|-{\mu_{X,FZ}(\eta_Z
\circ g \circ f)} & FX \ar@{.>}[dd]^{F(g \circ f)}
\\ Y \ar[d]_g & &
\\ Z \ar[r]_{\eta_Z}
& GFZ & FZ}$$ where the 2-cell $\mu_{X,FZ}(\eta_Z \circ g \circ f)$
is universal. The arrow $Fg \circ Ff$ satisfies
\begin{equation} \label{gammaF1}
\xymatrix@C=4pc@R=4pc{X \ar[d]_f \ar[r]^{\eta_X} & GFX \ar@{=}[r]
\ar@{.>}[d]^{GFf} \ar@{:>}[dl]|-{\mu_{X,FY}(\eta_Y \circ f)} & GFX
\ar[dd]^{G(Fg \circ Ff)} \ar@{=>}[ddl]|-{(\gamma^G_{Ff,Fg})^{-1}} &
FX \ar@{.>}[d]_{Ff} \ar@/^1pc/[dd]^{Fg \circ Ff}
\\ Y \ar[d]_g \ar[r]^{\eta_Y}
& GFY \ar@{.>}[d]^(.4){GFg} \ar@{:>}[dl]|-{\mu_{Y,FZ}(\eta_Z \circ
g)} & & FY \ar@{.>}[d]_{Fg}
\\ Z \ar[r]_{\eta_Z}
& GFZ \ar@{=}[r] & GFZ & FZ}
\end{equation}
since $G$ is a pseudo functor. Let $\gamma^F_{f,g}: Fg \circ Ff
\Rightarrow F(g \circ f)$ be the unique 2-cell whose $\phi_{X,FZ}$
image factors the composition of the 2-cells in (\ref{gammaF1}) as
follows.
\begin{equation} \label{gammaF2}
\xymatrix@C=6pc@R=3pc{Fg \circ Ff  \ar@{:>}[d]_{\gamma_{f,g}^F} &
G(Fg \circ Ff) \circ \eta_X \ar@{=>}[r]
\ar@{:>}[d]_{G(\gamma^F_{f,g}) \ast i_{\eta_X}} &  \eta_Z \circ g
\circ f \ar@{=}[d]
\\ F(g \circ f)
& GF(g \circ f) \circ \eta_X \ar@{=>}[r]_-{\mu_{X,FZ}(\eta_Z \circ g
\circ f)} & \eta_Z \circ g \circ f}
\end{equation}
The top horizontal 2-cell in the previous diagram is the composition
of the 2-cells in (\ref{gammaF1}). The 2-cell $\gamma^F_{f,g}: Fg
\circ Ff \Rightarrow F(g \circ f)$ is iso because the composition of
2-cells in (\ref{gammaF1}) is iso. Thus we have completely defined a
left biadjoint candidate $F$ for $G$. Now we must show that the
2-cells do what they should in order for $F$ to be a pseudo functor.

We claim that $\gamma^F$ is natural in its two variables. We must
show for morphisms $\xymatrix@1{X \ar[r]^{f_i} & Y \ar[r]^{g_i} &
Z}$ in $\mathcal{X}$ and 2-cells $\alpha:f_1 \Rightarrow f_2$ and
$\beta:g_1 \Rightarrow g_2$ in $\mathcal{X}$ that
\begin{equation} \label{gammanatural}
\xymatrix@C=3pc@R=3pc{Fg_1 \circ Ff_1
\ar@{=>}[r]^{\gamma^F_{f_1,g_1}} \ar@{=>}[d]_{F \beta * F\alpha} &
F(g_1 \circ f_1) \ar@{=>}[d]^{F(\beta * \alpha)}
\\ Fg_2 \circ Ff_2 \ar@{=>}[r]_{\gamma^F_{f_2,g_2}}
& F(g_2 \circ f_2)}
\end{equation}
commutes.

Toward this end, consider diagrams (\ref{gammanatural1}) and
(\ref{gammanatural2}).
\begin{equation} \label{gammanatural1}
\xymatrix@C=7pc@R=3pc{G(Fg_1 \circ Ff_1) \circ \eta_X
\ar@{=>}[r]^-{\sigma_1} \ar@{=>}[d]_-{G(\gamma^F_{f_1,g_1})*
i_{\eta_X}} & \eta_Z \circ g_1 \circ f_1 \ar@{=}[d]
\\ G(F(g_1 \circ f_1)) \circ \eta_X
\ar@{=>}[r]^-{\mu_{X,FZ}(\eta_Z \circ g_1 \circ f_1)}
\ar@{=>}[d]_-{GF(\beta * \alpha ) * i_{\eta_X}} & \eta_Z \circ g_1
\circ f_1 \ar@{=>}[d]^-{i_{\eta_Z} * \beta * \alpha}
\\ GF(g_2 \circ f_2) \circ \eta_X
\ar@{=>}[r]_-{\mu_{X,FZ}(\eta_Z \circ g_2 \circ f_2)} & \eta_Z \circ
g_2 \circ f_2}
\end{equation}

\begin{equation} \label{gammanatural2}
\xymatrix@C=7pc@R=3pc{G(Fg_1 \circ Ff_1) \circ \eta_X
\ar@{=>}[r]^-{\sigma_1} \ar@{=>}[d]_-{G(F \beta * F \alpha )*
i_{\eta_X}} & \eta_Z \circ g_1 \circ f_1 \ar@{=>}[d]^-{i_{\eta_Z}
* \beta * \alpha}
\\ G(Fg_2 \circ Ff_2) \circ \eta_X \ar@{=>}[r]^-{\sigma_2}
\ar[d]_-{G(\gamma^F_{f_2, g_2}) * i_{\eta_x}} & \eta_Z \circ g_2
\circ f_2 \ar@{=}[d]
\\ G(F(g_2 \circ f_2)) \circ \eta_X \ar@{=>}[r]_-{\mu_{X,FZ}(\eta_Z \circ
g_2 \circ f_2)} & \eta_Z \circ g_2 \circ f_2}
\end{equation}
The top horizontal 2-cell $\sigma_1$ in both diagrams is the
composition of the 2-cells in diagram (\ref{gammaF1}) with $f,g$
replaced by $f_1, g_1$ respectively. The bottom horizontal 2-cell in
each diagram is $\mu_{X,FZ}(\eta_Z \circ g_2 \circ f_2)$. The center
horizontal 2-cell $\sigma_2$ in (\ref{gammanatural2}) is the
composition of the 2-cells in (\ref{gammaF1}) with $f,g$ replaced by
$f_2,g_2$ respectively. The top rectangle in (\ref{gammanatural1})
commutes because it is the analogue of (\ref{gammaF2}) for
$f_1,g_1$. The bottom rectangle in (\ref{gammanatural1}) commutes
because of the naturality of $\mu_{X,FZ}:\phi_{X,FZ} \circ
\psi_{X,FZ} \Rightarrow 1_{Mor_{\mathcal{X}}(X, GFZ)}$. Hence the
outer rectangle of (\ref{gammanatural1}) commutes. The top rectangle
of (\ref{gammanatural2}) commutes because of the naturality of
$(\gamma^G)^{-1}, \mu_{X,FY},$ and $\mu_{Y,FZ}$ by comparing with
the 2-cells of (\ref{gammaF1}). The bottom rectangle of
(\ref{gammanatural2}) commutes because it is the analogue of
(\ref{gammaF2}) for $f_2,g_2$. Hence the outer rectangle of
(\ref{gammanatural2}) commutes. From (\ref{gammanatural1}) and
(\ref{gammanatural2}) we conclude that both $F(\beta * \alpha) \odot
\gamma^F_{f_1, g_1}$ and $\gamma^F_{f_2,g_2} \odot (F \beta
* F \alpha)$ have $\phi_{X,FZ}$ images which fill in the right
diagram of (\ref{gammanatural3}).
\begin{equation} \label{gammanatural3}
\xymatrix@C=6pc@R=3pc{Fg_1 \circ Ff_1  \ar@{:>}[d] & G(Fg_1 \circ
Ff_1) \circ \eta_X \ar@{=>}[r]^-{(i_{\eta_Z} * \beta * \alpha) \odot
\sigma_1} \ar@{:>}[d] &  \eta_Z \circ g_2 \circ f_2 \ar@{=}[d]
\\ F(g_2 \circ f_2) & GF(g_2 \circ f_2) \circ \eta_X
\ar@{=>}[r]_-{\mu_{X,FZ}(\eta_Z \circ g_2 \circ f_2)} & \eta_Z \circ
g_2 \circ f_2}
\end{equation}
Since $\mu_{X,FZ}(\eta_Z \circ g_2 \circ f_2)$ is universal, we
conclude that $F(\beta * \alpha) \odot \gamma^F_{f_1,g_1}=
\gamma^F_{f_2,g_2} \odot (F\beta * F \alpha)$ and thus $\gamma^F$ is
natural in its two variables.

We claim that $\delta^F$ and $\gamma^F$ satisfy the unit axiom for
pseudo functors. Let $X \in Obj \hspace{1mm}  \mathcal{X}$ and let
$f:X \rightarrow Y$ be a morphism in $\mathcal{X}$. We must show
that $\gamma^F_{1_X,f}= (i_{Ff} * \delta^F_{X*})^{-1}$. By
definition, $\gamma^F_{1_X,f}$ is the unique 2-cell $Ff \circ F1_X
\Rightarrow F(f \circ 1_X)$ such that the composition of 2-cells
\begin{equation} \label{Funitaxiom1}
\xymatrix@R=3pc@C=3pc{X \ar[rr]^-{\eta_X} & \ar@{=>}[d]^{i_{\eta_X}}
& GFX \ar[rr]^{G(Ff \circ F1_X)} & \ar@{=>}[d]^{G(\gamma^F_{1_X,f})}
& GFY
\\ X \ar[rr]_-{\eta_X} &
& GFX \ar@{:>}[d]^{\mu_{X,FY}(\eta_Y \circ f \circ 1_X)}
\ar[rr]_{GF(f \circ 1_X)} & & GFY
\\ X \ar[rr]_-{f \circ 1_X} & & Y \ar[rr]_{\eta_Y} & & GFY}
\end{equation}
is the same as the composition of 2-cells
\begin{equation} \label{Funitaxiom2}
\xymatrix@C=6pc@R=3pc{X \ar[r]^{\eta_X}
\ar@{=>}@<9ex>[d]^{i_{\eta_X}}
 & GFX \ar[rr]^{G(Ff \circ F1_X)}
& \ar@{=>}[d]^{(\gamma^G_{F1_X,Ff})^{-1}} & GFY
\\ X \ar[r]_{\eta_X}
& GFX \ar[r]^{GF1_X} \ar@{:>}[d]^{\mu_{X,FX}(\eta_X \circ 1_X)} &
GFX \ar[r]^-{GFf} \ar@{=>}@<10ex>[d]^{i_{GFf}} & GFY
\\ X \ar[r]^{1_X} \ar@{=>}@<9ex>[d]^{i_{1_X}} & X \ar[r]_{\eta_X} &
GFX \ar@{:>}[d]^{\mu_{X,FY}(\eta_Y \circ f)} \ar[r]_-{GFf} & GFY
\\ X \ar[r]_{1_X} & X \ar[r]_f & Y \ar[r]_{\eta_Y} & GFY}
\end{equation}
where universal 2-cells are drawn with dotted double arrows for
clarity. We show that $(i_{Ff} * \delta^F_{X*})^{-1}$ is a 2-cell
with this defining property for $\gamma^F_{1_X,f}$.

Since $\gamma^G$ is natural we can rewrite the first horizontal
2-cell composition in (\ref{Funitaxiom2}) as the composition of the
first three 2-cells in the equal diagram (\ref{Funitaxiom3}).
\begin{equation} \label{Funitaxiom3}
\xymatrix@C=6pc@R=3pc{X \ar[r]^-{\eta_X}
\ar@{=>}@<8.5ex>[d]^{i_{\eta_X}}
 & GFX \ar[rr]^{G(Ff \circ F1_X)}
& \ar@{=>}[d]^{G((i_{Ff} * \delta^F_{X*})^{-1})} & GFY
\\ X \ar[r]^(.3){\eta_X} \ar@{=>}@<8.5ex>[d]^{i_{\eta_X}}
 & GFX \ar[rr]^(.3){G(Ff \circ 1_{FX})}
& \ar@{=>}[d]^{(\gamma^G_{1_{FX},Ff})^{-1}} & GFY
\\ X \ar[r]^(.3){\eta_X} \ar@{=>}@<8.5ex>[d]^{i_{\eta_X}}
 & GFX \ar[rr]^(.3){GFf \circ G1_{FX}}
& \ar@{=>}[d]^{G(i_{Ff})*G(\delta^F_{X*})} & GFY
\\ X \ar[r]_-{\eta_X}
& GFX \ar[r]^{GF1_X} \ar@{:>}[d]^{\mu_{X,FX}(\eta_X \circ 1_X)} &
GFX \ar[r]^{GFf} \ar@{=>}@<10ex>[d]^{i_{GFf}} & GFY
\\ X \ar[r]^-{1_X} \ar@{=>}@<8.8ex>[d]^{i_{1_X}} & X \ar[r]_{\eta_X} &
GFX \ar@{:>}[d]^{\mu_{X,FY}(\eta_Y \circ f)} \ar[r]_{GFf} & GFY
\\ X \ar[r]_-{1_X} & X \ar[r]_f & Y \ar[r]_{\eta_Y} & GFY}
\end{equation}
By the unit axiom for $G$, the definition of $\delta^F_{X*}$ in
(\ref{Fdelta}), and the interchange law we see that the second
horizontal composition in (\ref{Funitaxiom3}) is
$$\aligned (\gamma^G_{1_{FX},Ff})^{-1}*i_{\eta_X}
&=i_{GFf}*\delta^G_{FX*}*i_{\eta_X} \\ &=i_{GFf} *
(\mu_{X,FX}(\eta_X \circ 1_X) \odot
(G(\delta^F_{X*})*i_{\eta_X}))^{-1} \\ &=
(G(i_{Ff})*G(\delta^F_{X*})^{-1}*i_{\eta_X}) \odot (i_{GFf}
*\mu_{X,FX}(\eta_X \circ 1_X))^{-1}. \endaligned$$ Substituting
this in (\ref{Funitaxiom3}) for
$(\gamma^G_{1_{FX},Ff})^{-1}*i_{\eta_X}$ we see that the second
horizontal composition in (\ref{Funitaxiom3}) cancels with the third
and the fourth, leaving only
\begin{equation} \label{Funitaxiom4}
\xymatrix@R=3pc@C=3pc{X \ar[rr]^{\eta_X} & \ar@{=>}[d]^{i_{\eta_X}}
& GFX \ar[rr]^{G(Ff \circ F1_X)} &
\ar@{=>}[d]^{G((i_{Ff}*\delta^F_{X*})^{-1})} & GFY
\\ X \ar[rr]_{\eta_X} &
& GFX \ar@{:>}[d]^{\mu_{X,FY}(\eta_Y \circ f \circ 1_X)}
\ar[rr]_{GF(f \circ 1_X)} & & GFY
\\ X \ar[rr]_{f \circ 1_X} & & Y \ar[rr]_{\eta_Y} & & GFY
}.
\end{equation}
We see that the 2-cell compositions of
(\ref{Funitaxiom1}),(\ref{Funitaxiom2}), (\ref{Funitaxiom3}), and
(\ref{Funitaxiom4}) are all equal. Hence the 2-cell compositions
(\ref{Funitaxiom1}) and (\ref{Funitaxiom4}) are equal and by
universality of the 2-cell $\mu_{X,FY}(\eta_Y \circ f \circ 1_X)$ we
have $\gamma^F_{1_X,f}=(i_{Ff} * \delta^F_{X*})^{-1}$. The other
half of the unit axiom can be verified similarly. We conclude that
$\delta^F$ and $\gamma^F$ satisfy the unit axiom for pseudo
functors.

We claim that $\gamma^F$ satisfies the composition axiom for pseudo
functors. Let \newline $\xymatrix@1{W \ar[r]^{f} & X \ar[r]^{g} & Y
\ar[r]^{h} & Z }$ be morphisms of $\mathcal{X}$. We must show that
\newline $\gamma^F_{f,h \circ g}=\gamma^F_{g \circ f, h} \odot
(i_{Fh} * \gamma^F_{f,g}) \odot (\gamma^F_{g,h} * i_{Ff})^{-1}$. By
definition $(\gamma^F_{f,h \circ g})^{-1}$ is the unique 2-cell $F(h
\circ g \circ f) \Rightarrow F(h \circ g) \circ Ff$ such that the
composition of 2-cells
\begin{equation} \label{Fcompositionaxiom1}
\xymatrix@R=4pc@C=6pc{W \ar[r]^{\eta_W} \ar[d]_{f} & GFW \ar@{=}[r]
\ar[d]_{GFf} \ar@{:>}[ld] & GFW \ar@{=}[r] \ar[ddd]|-(.7){G(F(h
\circ g) \circ Ff)} \ar@{=>}[lddd]|-{(\gamma^G_{Ff,F(h \circ
g)})^{-1}} & GFW \ar[ddd]|-(.7){GF(h \circ g \circ f)}
\ar@{=>}[lddd]|-{G((\gamma^F_{f,h \circ g})^{-1})}
\\ X \ar[d]_g \ar[r]^{\eta_X} & GFX \ar[dd]|-(.53){GF(h \circ g)}
 \ar@{:>}[ldd] & &
\\ Y \ar[d]_h & & &
\\ Z \ar[r]_{\eta_Z}
& GFZ \ar@{=}[r] & GFZ \ar@{=}[r] & GFZ}
\end{equation}
is the same as the universal 2-cell $\mu_{X,FZ}(\eta_Z \circ h \circ
g \circ f)$. For clarity we continue to draw the universal 2-cells
as dotted double arrows. We prove that replacing $(\gamma^F_{f,h
\circ g})^{-1}$ in (\ref{Fcompositionaxiom1}) by $(\gamma^F_{g \circ
f, h} \odot (i_{Fh} * \gamma^F_{f,g}) \odot (\gamma^F_{g,h} *
i_{Ff})^{-1})^{-1}$ still gives $\mu_{X,FZ}(\eta_Z \circ h \circ g
\circ f)$. After that we conclude $\gamma^F_{f,h \circ
g}=\gamma^F_{g \circ f, h} \odot (i_{Fh} * \gamma^F_{f,g}) \odot
(\gamma^F_{g,h} * i_{Ff})^{-1}$ by the universality of the 2-cell
$\mu_{X,FZ}(\eta_Z \circ h \circ g \circ f)$. To this end, we claim
that the composition
\begin{equation} \label{Fcompositionaxiom2}
\xymatrix@R=4pc@C=4pc{W \ar[r]^{\eta_W} \ar[d]_{f} & GFW \ar@{=}[r]
\ar[d]_{GFf} \ar@{:>}[ld] & GFW \ar@{=}[r] \ar[ddd]|-(.6){G(F(h
\circ g) \circ Ff)} \ar@{=>}[lddd]|-{(\gamma^G_{Ff,F(h \circ
g)})^{-1}} & GFW \ar[ddd]|(.6){GF(h \circ g \circ f)}
\ar@{=>}[lddd]|-{}
\\ X \ar[d]_g \ar[r]^{\eta_X} & GFX \ar[dd]|-(.375){GF(h \circ g)}
 \ar@{:>}[ldd] & &
\\ Y \ar[d]_h & & &
\\ Z \ar[r]_{\eta_Z}
& GFZ \ar@{=}[r] & GFZ \ar@{=}[r] & GFZ}
\end{equation}
is the same as $\mu_{X,FZ}(\eta_Z \circ h \circ g \circ f)$ , where
the rightmost 2-cell is $G((\gamma^F_{g \circ f, h} \odot (i_{Fh} *
\gamma^F_{f,g}) \odot (\gamma^F_{g,h} * i_{Ff})^{-1})^{-1})$. We do
this by transforming (\ref{Fcompositionaxiom2}) to a diagram known
to be $\mu_{X,FZ}(\eta_Z \circ h \circ g \circ f)$. The naturality
of $\gamma^G$ guarantees that
$$\xymatrix@C=6pc@R=4pc{G(Fh \circ Fg) \circ GFf
\ar@{=>}[r]^{\gamma^G_{Ff,Fh \circ Fg}}
\ar@{=>}[d]_{G(\gamma^F_{g,h})* i_{GFf}} & G(Fh \circ Fg \circ Ff)
\ar@{=>}[d]^{G(\gamma^F_{g,h} * i_{Ff})}
\\ GF(h \circ g) \circ GFf \ar@{=>}[r]_{\gamma^G_{Ff,F(h \circ g)}}
& G(F(h \circ g) \circ Ff)}$$ commutes. Using this commutivity to
substitute for $(\gamma^G_{Ff,F(h \circ g)})^{-1}$ in
(\ref{Fcompositionaxiom2}) and cancelling $G(\gamma^F_{g,h}*
i_{Ff})^{-1} \odot G(\gamma^F_{g,h}* i_{Ff})$  gives
\begin{equation} \label{Fcompositionaxiom3}
\xymatrix@R=4pc@C=4pc{W \ar[r]^{\eta_W} \ar[d]_{f}
 & GFW \ar@{=}[r]
\ar[d]_{GFf} \ar@{:>}[ld] & GFW \ar@{=}[r] \ar[d]_{GFf}
\ar@{=>}[ld]|-{i_{GFf}}
 & GFW  \ar@{=}[r] \ar[ddd]|(.62){G(Fh \circ Fg \circ Ff)}
\ar@{=>}[lddd]|-{(\gamma^G_{Ff,Fh \circ Fg})^{-1}} & GFW
\ar[ddd]|(.62){GF(h \circ g \circ f)} \ar@{=>}[lddd]|-{}
\\ X \ar[d]_g \ar[r]^{\eta_X} & GFX \ar[dd]|-(.4){GF(h \circ g)}
 \ar@{:>}[ldd] \ar@{=}[r] & GFX \ar[dd]|-(.4){G(Fh \circ Fg)}
\ar@{=>}[ldd]|-{\overset{\phantom{.}}{G(\gamma^F_{g,h})}} & &
\\ Y \ar[d]_h & & & &
\\ Z \ar[r]_{\eta_Z}
& GFZ \ar@{=}[r] & GFZ \ar@{=}[r] & GFZ \ar@{=}[r] & GFZ}
\end{equation}
where the right 2-cell is $G((i_{Fh}*\gamma^F_{f,g})^{-1} \odot
(\gamma^F_{g \circ f,h})^{-1})$. We have also implicitly used the
fact that $G$ preserves the vertical composition of 2-cells. By the
definition of $\gamma^F_{g,h}$ in (\ref{gammaF1}) and
(\ref{gammaF2}), the lower left two rectangles of
(\ref{Fcompositionaxiom3}) can be rewritten to give the equal
composition (\ref{Fcompositionaxiom4}).
\begin{equation} \label{Fcompositionaxiom4}
\xymatrix@R=4pc@C=4pc{W \ar[r]^{\eta_W} \ar[d]_{f}
 & GFW \ar@{=}[r]
\ar[d]_{GFf} \ar@{:>}[ld] & GFW \ar@{=}[r] \ar[d]_{GFf}
\ar@{=>}[ld]|-{i_{GFf}}
 &
GFW  \ar@{=}[r] \ar[ddd]|-(.62){G(Fh \circ Fg \circ Ff)}
\ar@{=>}[lddd]|-{(\gamma^G_{Ff,Fh \circ Fg})^{-1}} & GFW
\ar[ddd]|(.62){GF(h \circ g \circ f)} \ar@{=>}[lddd]|-{}
\\ X \ar[d]_g \ar[r]^{\eta_X} & GFX \ar[d]_{GFg}
 \ar@{:>}[ld] \ar@{=}[r] & GFX \ar[dd]|-(.4){G(Fh \circ Fg)}
\ar@{=>}[ldd]|-{(\gamma^G_{Fg,Fh})^{-1}} & &
\\ Y \ar[d]_h \ar[r]_{\eta_Y} & GFY \ar[d]_{GFh} \ar@{:>}[ld]& & &
\\ Z \ar[r]_{\eta_Z}
& GFZ \ar@{=}[r] & GFZ \ar@{=}[r] & GFZ \ar@{=}[r] & GFZ}
\end{equation}
Recall that the composition axiom for the pseudo functor $G$
guarantees the commutivity of the following diagram.
$$
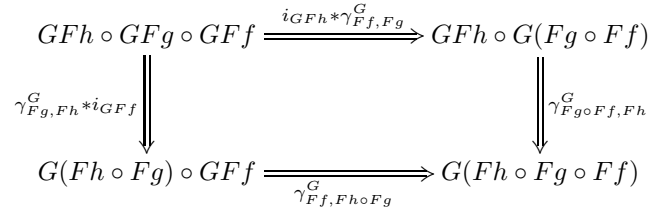
\xymatrix@C=5pc@R=3pc{GFh \circ GFg \circ GFf
\ar@{=>}[r]^{i_{GFh}*\gamma^G_{Ff,Fg}}
\ar@{=>}[d]_{\gamma^G_{Fg,Fh}*i_{GFf}} & GFh \circ G(Fg \circ Ff)
\ar@{=>}[d]^{\gamma^G_{Fg \circ Ff,Fh}}
\\  G(Fh \circ Fg) \circ GFf \ar@{=>}[r]_{\gamma^G_{Ff,Fh \circ Fg}}
& G(Fh \circ Fg \circ Ff)}$$ Using this composition axiom for the
pseudo functor $G$ we can replace the middle two columns of 2-cells
in (\ref{Fcompositionaxiom4}) to get the equal composition
(\ref{Fcompositionaxiom5}).
\begin{equation} \label{Fcompositionaxiom5}
\xymatrix@R=4pc@C=4pc{W \ar[r]^{\eta_W} \ar[d]_{f}
 & GFW \ar@{=}[r]
\ar[d]_{GFf} \ar@{:>}[ld] & GFW \ar@{=}[r] \ar[dd]|-(.6){G(Fg \circ
Ff)} \ar@{=>}[ldd]|-{(\gamma^G_{Ff,Fg})^{-1}}
 &
GFW  \ar@{=}[r] \ar[ddd]|(.62){G(Fh \circ Fg \circ Ff)}
\ar@{=>}[lddd]|-{(\gamma^G_{Fg \circ Ff,Fh})^{-1}} & GFW
\ar[ddd]|(.62){GF(h \circ g \circ f)} \ar@{=>}[lddd]|-{}
\\ X \ar[d]_g \ar[r]^{\eta_X} & GFX \ar[d]_{GFg}
 \ar@{:>}[ld]  &  & &
\\ Y \ar[d]_h \ar[r]_{\eta_Y} & GFY \ar[d]_{GFh} \ar@{:>}[ld]
\ar@{=}[r] & GFY \ar[d]_{GFh} \ar@{=>}[ld]|-{i_{GFh}} & &
\\ Z \ar[r]_{\eta_Z}
& GFZ \ar@{=}[r] & GFZ \ar@{=}[r] & GFZ \ar@{=}[r] & GFZ}
\end{equation}
In (\ref{Fcompositionaxiom5}) the right 2-cell is again
$G((i_{Fh}*\gamma^F_{f,g})^{-1} \odot (\gamma^F_{g \circ
f,h})^{-1})$ as in (\ref{Fcompositionaxiom3}) and
(\ref{Fcompositionaxiom4}). By the definition of $\gamma^F_{f,g}$ in
(\ref{gammaF1}) and (\ref{gammaF2}), we can rewrite the upper left
three rectangles of (\ref{Fcompositionaxiom5}) to obtain
(\ref{Fcompositionaxiom6}), which has
$G((i_{Fh}*\gamma^F_{f,g})^{-1} \odot (\gamma^F_{g \circ
f,h})^{-1})$ as its right 2-cell.
\begin{equation} \label{Fcompositionaxiom6}
\xymatrix@R=4pc@C=4pc{W \ar[r]^{\eta_W} \ar[d]_{f}
 & GFW \ar@{=}[r]
\ar[dd]|-(.4){GF(g \circ f)} \ar@{:>}[ldd] & GFW \ar@{=}[r]
\ar[dd]|-(.4){G(Fg \circ Ff)} \ar@{=>}[ldd]|-{G(\gamma^F_{f,g})}
 &
GFW  \ar@{=}[r] \ar[ddd]|-(.62){G(Fh \circ Fg \circ Ff)}
\ar@{=>}[lddd]|-{(\gamma^G_{Fg \circ Ff,Fh})^{-1}} & GFW
\ar[ddd]|(.62){GF(h \circ g \circ f)} \ar@{=>}[lddd]|-{}
\\ X \ar[d]_g &  &  & &
\\ Y \ar[d]_h \ar[r]_{\eta_Y} & GFY \ar[d]_{GFh} \ar@{:>}[ld]
\ar@{=}[r] & GFY \ar[d]_{GFh} \ar@{=>}[ld]|-{i_{GFh}} & &
\\ Z \ar[r]_{\eta_Z}
& GFZ \ar@{=}[r] & GFZ \ar@{=}[r] & GFZ \ar@{=}[r] & GFZ}
\end{equation}
The naturality of $\gamma^G$ implies that the diagram
$$\xymatrix@C=5pc@R=4pc{GFh \circ G(Fg \circ Ff)
\ar@{=>}[r]^{\gamma^G_{Fg \circ Ff, Fh}} \ar@{=>}[d]_{G(i_{Fh}) *
G(\gamma^F_{f,g})} & G(Fh \circ Fg \circ Ff) \ar@{=>}[d]^{G(i_{Fh}
* \gamma^F_{f,g})}
\\ GFh \circ G(F(g \circ f)) \ar@{=>}[r]_{\gamma^G_{F(g \circ f),Fh}}
& G(Fh \circ F(g\circ f))}$$ \vspace{1in} commutes. Using its
commutivity, we can rewrite (\ref{Fcompositionaxiom6}) by combining
its middle two columns of 2-cells with
$G((i_{Fh}*\gamma^F_{f,g})^{-1})$ from the last column to get
(\ref{Fcompositionaxiom7}).
\begin{equation} \label{Fcompositionaxiom7}
\xymatrix@R=4pc@C=4pc{W \ar[r]^{\eta_W} \ar[d]_{f}
 & GFW \ar@{=}[r]
\ar[dd]|-(.64){GF(g \circ f)} \ar@{:>}[ldd] & GFW
\ar[ddd]|-(.4){G(Fh \circ F(g \circ f))}
\ar@{=>}[lddd]|-{(\gamma^G_{F(g \circ f),Fh})^{-1}} \ar@{=}[r] & GFW
\ar[ddd]|(.4){GF(h \circ g \circ f)} \ar@{=>}[lddd]|-{G((\gamma^F_{g
\circ f,h})^{-1})}
\\ X \ar[d]_g &  & &
\\ Y \ar[d]_h \ar[r]_{\eta_Y} & GFY \ar[d]_{GFh} \ar@{:>}[ld]
 & &
\\ Z \ar[r]_{\eta_Z}
& GFZ \ar@{=}[r] & GFZ \ar@{=}[r] & GFZ}
\end{equation}
But by the definition of $\gamma^F_{g \circ f, h}$ in
(\ref{gammaF1}) and (\ref{gammaF2}), the composition of 2-cells in
(\ref{Fcompositionaxiom7}) is precisely $\mu_{X,FZ}(\eta_Z \circ h
\circ g \circ f)$. Since the compositions of 2-cells in the diagrams
(\ref{Fcompositionaxiom2}) through (\ref{Fcompositionaxiom7}) are
all equal, we conclude that the composition of 2-cells in
(\ref{Fcompositionaxiom2}) is $\mu_{X,FZ}(\eta_Z \circ h \circ g
\circ f)$. We conclude that $\gamma^F_{f,h \circ g}=\gamma^F_{g
\circ f, h} \odot (i_{Fh} * \gamma^F_{f,g}) \odot (\gamma^F_{g,h} *
i_{Ff})^{-1}$ by the universality of $\mu_{X,FZ}(\eta_Z \circ h
\circ g \circ f)$. Therefore $\gamma^F$ satisfies the composition
axiom for pseudo functors.

In summary, we have constructed a pseudo functor $F:\mathcal{X}
\rightarrow \mathcal{A}$ with natural coherence 2-cells $\delta^F$
and $\gamma^F$ and we have shown that they satisfy the unit axiom
and composition axiom for pseudo functors.

Next we have to show that $F$ is a left biadjoint using Theorem
\ref{maintheorem1}. By hypothesis we already have a morphism
$\eta_X:X \rightarrow G(FX)$ for all $X \in Obj \hspace{1mm}
\mathcal{X}$. We claim that the assignment $X \mapsto \eta_X$ is a
pseudo natural transformation from $1_{\mathcal{X}}$ to $GF$. We
need to define the 2-cells up to which $\eta$ is natural. For a
morphism $f:X \rightarrow Y$ of $\mathcal{X}$ define
\label{makeetalaxnatural} $\tau_f:=\mu_{X,FY}(\eta_Y \circ f)$. Then
the diagram
$$\xymatrix@C=3pc@R=3pc{X \ar[r]^-{\eta_X} \ar[d]_f & GFX \ar[d]^{GFf}
\ar@{=>}[dl]_{\tau_f} \\ Y \ar[r]_-{\eta_Y} & GFY}$$ illustrates the
source and target of the 2-cell. The map $f \mapsto \tau_f$ is
natural because $\mu_{X,FY}$ is a natural transformation. More
precisely let $\alpha:f_1 \Rightarrow f_2$ be a 2-cell in
$\mathcal{X}$ and let $f_1,f_2:X \rightarrow Y$ be morphisms in
$\mathcal{X}$. Then
\begin{equation} \label{taunatural1}
\xymatrix@C=7pc@R=3pc{\phi_{X,FY} (\psi_{X,FY}(\eta_Y \circ f_1))
\ar@{=>}[r]^-{\mu_{X,FY}(\eta_Y \circ f_1)} \ar@{=>}[d]|{\phi_{X,FY}
(\psi_{X,FY}(i_{\eta_Y}*\alpha))} & \eta_Y \circ f_1
\ar@{=>}[d]^{i_{\eta_Y}*\alpha}
\\ \phi_{X,FY} (\psi_{X,FY}(\eta_Y \circ f_2))
\ar@{=>}[r]_-{\mu_{X,FY}(\eta_Y \circ f_2)} & \eta_Y \circ f_2}
\end{equation}
commutes by the naturality of $\mu_{X,FY}$. By the definitions of
$F$, $\tau_{f_1}$, and $\tau_{f_2}$, diagram (\ref{taunatural1}) is
the same as the diagram
\begin{equation} \label{taunatural11}
\xymatrix@C=7pc@R=3pc{GFf_1 \circ \eta_X \ar@{=>}[r]^{\tau_{f_1}}
\ar@{=>}[d]_{GF\alpha * i_{\eta_X}} & \eta_Y \circ f_1
\ar@{=>}[d]^{i_{\eta_Y}*\alpha}
\\ GFf_2 \circ \eta_X
\ar@{=>}[r]_-{\tau_{f_2}} & \eta_Y \circ f_2}
\end{equation}
which says $f \mapsto \tau_f$ is natural. The map $f \mapsto \tau_f$
satisfies the unit axiom for pseudo natural transformations because
of (\ref{Fdelta}) and the definition of $\delta^{GF}$ for the
composite pseudo functor $GF$. The map $f \mapsto \tau_f$ satisfies
the composition axiom for pseudo natural transformations because of
(\ref{gammaF1}) and (\ref{gammaF2}), and the definition of
$\gamma^{GF}$ for the composite pseudo functor $GF$. Hence
$\eta:1_{\mathcal{X}} \Rightarrow GF$ is a pseudo natural
transformation with coherence 2-cells $\tau$.

By Theorem \ref{maintheorem1}, the constructed pseudo functor $F$ is
a left biadjoint because $\eta:1_{\mathcal{X}} \Rightarrow GF$ is a
pseudo natural transformation such that $\eta_X:X \rightarrow G(FX)$
is a biuniversal arrow for all $X \in Obj \hspace{1mm} \mathcal{X}$.
\end{pf}

We can summarize the previous two theorems in a way similar to Mac
Lane's theorem on page 83 of \cite{maclane1} as follows.

\begin{thm}
A biadjunction\index{biadjunction|textbf} $\langle F,G,\phi \rangle:
\mathcal{X} \rightharpoonup \mathcal{A}$ can be described up to
pseudo natural pseudo\index{isomorphism!pseudo isomorphism}
isomorphism (defined below) by either of the following data:
\begin{enumerate}
\item
Pseudo functors
$$\xymatrix@R=3pc@C=3pc{ \mathcal{X}  \ar@<.5ex>[r]^F & \mathcal{A} \ar@<.5ex>[l]^G}$$
and a pseudo natural transformation $\eta:1_{\mathcal{X}}
\Rightarrow GF$ such that each $\eta_X:X \rightarrow G(FX)$ is a
biuniversal\index{biuniversal arrow|textbf}\index{arrow!biuniversal
arrow|textbf} arrow from $X$ to $G$. Then $\phi_{X,A}$ is defined by
$\phi_{X,A}(f)=Gf \circ \eta_X$.
\item
A pseudo functor $G:\mathcal{A} \rightarrow \mathcal{X}$, for each
$X \in Obj \hspace{1mm}  \mathcal{X}$ an object $R \in \mathcal{A}$
depending on $X$, and for each $X \in Obj \hspace{1mm} \mathcal{X}$
a biuniversal\index{biuniversal
arrow|textbf}\index{arrow!biuniversal arrow|textbf} arrow $\eta_X:X
\rightarrow GR$ from $X$ to $G$. Then the pseudo functor $F$
satisfies $FX=R$ on objects and there is a natural iso 2-cell $GFh
\circ \eta_X \Rightarrow \eta_X' \circ h$ for morphisms $h:X
\rightarrow X'$.
\end{enumerate}
\end{thm}
\begin{pf}
Uniqueness will be proven below.
\end{pf}

Similar things can be formulated for bicounits. From 1-category
theory we know that any two left adjoints to a functor are naturally
isomorphic. A similar statement can be made for left biadjoints,
although we need the concept of {\it pseudo natural pseudo
isomorphism}.

\begin{defn}
Let $F,F':\mathcal{X} \rightarrow \mathcal{A}$ be pseudo functors.
Then a pseudo natural transformation $\alpha:F \Rightarrow F'$ is
called a {\it pseudo natural pseudo
isomorphism}\index{isomorphism!pseudo natural pseudo
isomorphism|textbf} or {\it pseudo natural
equivalence}\index{equivalence!pseudo natural equivalence|textbf} if
there exists a pseudo natural transformation $\alpha':F' \Rightarrow
F$ and there exist iso modifications $\alpha \odot \alpha'
\rightsquigarrow 1_{F'}$ and $\alpha' \odot \alpha \rightsquigarrow
1_F$.
\end{defn}

\begin{thm} \label{laxadjointuniqueness}
Let $F,F':\mathcal{X} \rightarrow \mathcal{A}$ be left biadjoints
for a pseudo functor $G:\mathcal{A} \rightarrow \mathcal{X}$. Then
there exists a pseudo natural pseudo isomorphism $\alpha:F
\Rightarrow F'$
\end{thm}
\begin{pf}
For $X \in Obj \hspace{1mm}  \mathcal{X}$, let $\eta_X:X \rightarrow
G(FX)$ and $\eta_X':X \rightarrow G(F'X)$ be the biuniversal arrows
obtained from the biadjunctions as in the theorems above. Then by
Lemma \ref{laxuniqueness} there exists a pseudo isomorphism
$\alpha_X:FX \rightarrow F'X$ and a pseudo inverse $\alpha_X':F'X
\rightarrow FX$ as well as 2-cells $\alpha_X' \circ \alpha_X
\Rightarrow 1_{FX}$ and $\alpha_X \circ \alpha_X' \Rightarrow
1_{F'X}$. It can be shown that the assignments $X \rightarrow
\alpha_X$ and $X \rightarrow \alpha_X'$ are pseudo natural and the
2-cells determine modifications $\alpha' \odot \alpha
\rightsquigarrow 1_F$ and $\alpha \odot \alpha' \rightsquigarrow
1_{F'}$.

For example, we construct the coherence 2-cell $\tau^{\alpha}$ up to
which $\alpha$ is natural. For $f \in Mor_{\mathcal{X}}(X,Y)$ we
have the following two diagrams.
\begin{equation} \label{alphanaturality1}
\xymatrix@C=4pc@R=3pc{X \ar[d]_f \ar[r]^{\eta_X} & GFX \ar[d]^{GFf}
\ar@{=>}[ld]_{\tau^{\eta}_f} & FX \ar[d]^{Ff}
\\ Y \ar@{=}[d] \ar[r]_{\eta_Y} & GFY \ar@{:>}[ld]^{\mu(\eta_Y')}
\ar@{.>}[d]^{G\alpha_Y} & FY \ar@{.>}[d]^{\alpha_Y}
\\ Y \ar[r]_{\eta_Y'} & GF'Y & F'Y
}
\end{equation}
\begin{equation} \label{alphanaturality2}
\xymatrix@C=4pc@R=3pc{X \ar@{=}[d] \ar[r]^{\eta_X} & GFX
\ar@{:>}[ld]_{\mu(\eta_X')} \ar@{.>}[d]^{G \alpha_X} & FX
\ar@{.>}[d]^{\alpha_X}
\\ X \ar[d]_{f} \ar[r]_{\eta_X'} & GF'X \ar@{=>}[ld]^{\tau^{\eta'}_f}
\ar[d]^{GF'f} & F'X \ar[d]^{F'f}
\\ Y \ar[r]_{\eta_Y'} & GF'Y & F'Y }
\end{equation}
But they can also be filled in as
\begin{equation} \label{alphanaturality3}
\xymatrix@C=4pc@R=3pc{X \ar@{=}[d] \ar[r]^{\eta_X} & GFX
\ar@{.>}[d]^{G \psi(\eta_Y' \circ f)} \ar@{:>}[ld] & FX
\ar@{.>}[d]^{\psi(\eta_Y' \circ f)}
\\ X \ar[r]_{\eta_Y' \circ f} & GF'Y & F'Y}
\end{equation}
where the dashed 2-cell is universal. The universality gives us iso
2-cells $\nu_f$ and $\nu_f'$ as in
$$\xymatrix@1@C=4pc@R=3pc{F'f \circ \alpha_X \ar@{=>}[r]^-{\nu_f'}
& \psi(\eta_Y' \circ f) & \alpha_Y \circ Ff \ar@{=>}[l]_-{\nu_f}}$$
whose $\phi$ images factor (via the universal 2-cell) the 2-cells in
(\ref{alphanaturality1}) and (\ref{alphanaturality2}) precomposed
with the appropriate $(\gamma^G)^{-1}$'s. Define
$\tau^{\alpha}_f:=\tau_f:=(\nu_f)^{-1} \odot \nu_f'$. This is the
coherence 2-cell up to which $\alpha$ will be natural.

A sketch of the naturality of $f \mapsto \tau_f$ goes as follows.
Let $\beta:f_1 \Rightarrow f_2$ be a 2-cell between $f_1,f_2:X
\rightarrow Y$. Then we must show that the outer rectangle of
$$\xymatrix@C=4pc@R=3pc{F'f_1 \circ \alpha_X \ar@{=>}[r]^-{\nu_{f_1}'}
 \ar@{=>}[d]_{F' \beta *
i_{\alpha_X}} & \psi(\eta_Y' \circ f_1)
\ar@{=>}[d]^{\psi(i_{\eta_Y'}* \beta)} & \alpha_Y \circ Ff_1
\ar@{=>}[l]_-{\nu_{f_1}} \ar@{=>}[d]^{i_{\alpha_Y} * F
 \beta}
\\ F'f_2 \circ \alpha_X \ar@{=>}[r]_-{\nu_{f_2}'} & \psi(\eta_Y' \circ f_2)
& \alpha_Y \circ Ff_2 \ar@{=>}[l]^-{\nu_{f_2}}}$$ commutes. We do
this by showing that the individual inner squares commute by
applying $\phi$ and using the universality and the fact that $\mu$
is a natural isomorphism. It also involves the naturality of the
$\gamma^G$'s.

We can also show that $\tau$ satisfies the composition and unit
axiom, although it is lengthy. Lastly we must verify that the 2-cell
assignments at the start actually give modifications $\alpha' \odot
\alpha \rightsquigarrow 1_F$ and $\alpha \odot \alpha'
\rightsquigarrow 1_{F'}$.

Thus, any two left biadjoints are pseudo naturally pseudo
isomorphic.
\end{pf}

There is a relationship between  bi(co)limits and biadjoints, just
like for (co)limits and adjoints.

\begin{rem}
Let $\mathcal{C}$ be a 2-category which admits
bicolimits\index{bicolimit} and bilimits\index{bilimit} and let
$\mathcal{J}$ be a 1-category. Let $\mathcal{C}^{\mathcal{J}}$ be
the 2-category with objects pseudo functors $\mathcal{J} \rightarrow
\mathcal{C}$, morphisms pseudo natural transformations, and 2-cells
the modifications. Let $\Delta:\mathcal{C} \rightarrow
\mathcal{C}^{\mathcal{J}}$ be the diagonal 2-functor\index{diagonal
2-functor}\index{2-functor!diagonal 2-functor}. Then
$bicolim:\mathcal{C}^{\mathcal{J}} \rightarrow \mathcal{C}$ is a
left biadjoint\index{biadjoint} for $\Delta$ and the arrows of the
biunit constructed in Theorem \ref{maintheorem1} are the universal
pseudo cones. Similarly, $bilim:\mathcal{C}^{\mathcal{J}}
\rightarrow \mathcal{C}$ is a right biadjoint for $\Delta$ and the
arrows of the bicounit are the universal pseudo cones.
\end{rem}

\chapter{Forgetful 2-Functors for Pseudo Algebras}
\label{sec:forgetfulfunctor} Next we show that forgetful 2-functors
for pseudo algebras\index{algebra!pseudo algebra} admit left
biadjoints. Let us consider the strict case as an example of what we
do below. Let $S$ be the theory\index{theory!theory of abelian
groups}\index{abelian groups}\index{abelian groups!theory of abelian
groups} of abelian groups and let $T$ be the theory of
rings\index{theory!theory of rings}\index{ring}\index{ring!theory of
rings}. Then we have an inclusion $S \hookrightarrow T$. Let $X$ be
a discrete $T$-algebra, \ie $X$ is a set and we have a morphism of
theories $T \rightarrow End(X)$. Then $X$ can be made into an
$S$-algebra by the composite map of theories $S \hookrightarrow T
\rightarrow End(X)$. This precomposition with the inclusion arrow
forgets the ring structure on the set $X$ and results in the
underlying abelian group. This precomposition with the inclusion
defines the {\it forgetful functor}\index{forgetful
functor}\index{functor!forgetful functor} from the category of rings
to the category of abelian groups. It admits a left adjoint which is
the appropriate {\it free functor}\index{free
functor}\index{functor!free functor}. Similarly, for any morphism of
theories $S \rightarrow T$ we have a forgetful 2-functor from pseudo
$T$-algebras to pseudo $S$-algebras and this 2-functor admits a left
biadjoint\index{biadjoint!left biadjoint}. Blackwell, Kelly, and
Power have shown that left biadjoints exist for the analogous
2-functors on 2-categories of strict algebras\index{algebra!strict
algebra} over 2-monads\index{2-monad} with pseudo morphisms in
\cite{blackwell}. Lack has given sufficient conditions in
\cite{lack} under which the inclusion of strict
algebras\index{algebra!strict algebra} over a 2-monad\index{2-monad}
into pseudo algebras\index{algebra!pseudo algebra} over the same
2-monad admits a left adjoint whose unit has components that are
equivalences. In such cases, every pseudo
algebra\index{algebra!pseudo algebra} over the 2-monad is equivalent
to a strict algebra\index{algebra!strict algebra} over the 2-monad.
Yanofsky has also studied quasiadjoints\index{quasiadjoint} to
forgetful 2-functors induced by morphisms of 2-theories in
\cite{yanofsky2000}, although his 2-theories are different from
those of \cite{hu}, \cite{hu1}, \cite{hu2}, and Chapter
\ref{sec:2-theories}.

\begin{defn} \label{defnoftheforgetful2functor}
Let $\phi:S \rightarrow T$ be a morphism of theories and let $X$ be
a pseudo $T$-algebra with structure maps $\Psi_n:T(n) \rightarrow
End(X)(n)$. Let $UX$ be the pseudo $S$-algebra which has $X$ as its
underlying category and $S$ structure maps defined by
$\Psi_n(\phi(w)):X^n \rightarrow X$ for $w \in S(n)$. Defining $U$
analogously for morphisms and 2-cells of the 2-category of pseudo
$T$-algebras yields a strict 2-functor $U$ from the 2-category of
pseudo $T$-algebras to the 2-category of pseudo $S$-algebras called
the {\it forgetful 2-functor associated to $\phi$}\index{forgetful
2-functor|textbf}\index{2-functor!forgetful 2-functor|textbf}.
\end{defn}

To show that the forgetful 2-functor associated to $\phi$ admits a
left biadjoint, we need to find a biuniversal arrow of the following
type: given a pseudo $S$-algebra $X$ there should exist a pseudo
$T$-algebra $R$ and a biuniversal arrow $\eta_X:X \rightarrow UR$ in
the category of pseudo $S$-algebras. We define this $R$ now.

\begin{notation} \label{primenotation}
Let $T$ be a theory. Let $T'$ denote the free theory\index{free
theory}\index{theory!free theory} on the sequence of sets
$T(0),T(1), \dots$ underlying the theory $T$. The category
$Alg'$\index{$Alg'$} is the category whose objects are small
$T'$-algebras and whose morphisms are morphisms of strict
$T'$-algebras. Let $Obj \hspace{1mm} Graph'$ be the collection of
small directed graphs\index{directed graph} whose object sets are
discrete $T'$-algebras. Let $Mor \hspace{1mm} Graph'$ be the
collection of morphisms of directed graphs whose object components
are morphisms of discrete $T'$-algebras. Then
$Graph'$\index{$Graph'$} is a category. We denote by $V'$ the left
adjoint to the forgetful\index{forgetful
functor}\index{functor!forgetful functor} functor $V:Alg'
\rightarrow Graph'$.
\end{notation}

The forgetful functor $V:Alg' \rightarrow Graph'$ admits a left
adjoint $V':Graph' \rightarrow Alg'$  by Freyd's Adjoint Functor
Theorem\index{Freyd's Adjoint Functor Theorem}. The functor $V'$ is
similar to taking the free\index{free category} category on a
directed graph\index{directed graph}, except the resulting category
is also a $T'$-algebra. The objects of the underlying directed graph
of $V'Y$ and the objects of the directed graph $Y$ are the same.

\begin{defn} \label{defnfreelaxTalgebra}
Let $\phi:S \rightarrow T$ be a morphism of theories. Let $X$ be a
pseudo $S$-algebra with structure maps $\Psi_n:S(n) \rightarrow
End(X)(n)$. We define the {\it free pseudo $T$-algebra $R$ on the
pseudo $S$-algebra $X$ associated to $\phi$}\index{algebra!free
pseudo $T$-algebra on a pseudo $S$-algebra|(textbf}\index{free
pseudo $T$-algebra on a pseudo $S$-algebra|(textbf} via intermediate
steps $R_{G'}$ and $R'$ as follows. Let $Obj \hspace{1mm}
R_{G'}$\index{$R_{G'}$|(} be the (discrete) free $T'$-algebra on the
discrete category $Obj \hspace{1mm} X$ and let $Mor \hspace{1mm}
R_{G'}$ be the collection of the following arrows:

\begin{enumerate}
\item
For every $n \in \mathbf{N}$, for all words $w \in T(n)$, $w_1 \in
T(m_1), \dots, w_n \in T(m_n)$, and for all objects $A^1_1,\dots,
A^1_{m_1}$,$A^2_1, \dots, A^2_{m_2}, \dots, A^n_1, \dots,
A^n_{m_n} \in Obj \hspace{1mm} R_{G'}$ there are arrows 
$$c_{w, w_1, \ldots, w_n}(A^1_1,\dots,A^n_{m_n}):$$
$$\xymatrix@1{ w \circ (w_1, \dots, w_n)(A^1_1,\dots,A^n_{m_n}) \ar[r] &
w(w_1(A^1_1,\dots, A^1_{m_1}), \dots, w_n(A^n_1, \dots,
A^n_{m_n}))}$$ $$c_{w, w_1,
\ldots,w_n}^{-1}(A^1_1,\dots,A^n_{m_n}):$$
$$\xymatrix@1{w(w_1(A^1_1,\dots, A^1_{m_1}), \dots, w_n(A^n_1, \dots,
A^n_{m_n})) \ar[r] & w \circ (w_1, \dots,
w_n)(A^1_1,\dots,A^n_{m_n})}.$$ Here $w \circ (w_1, \ldots, w_n)$ is
the composition in the original theory $T$. The target
$w(w_1(A^1_1,\dots, A^1_{m_1}), \dots, w_n(A^n_1, \dots,
A^n_{m_n}))$ is the result of composing in the free
theory\index{free theory} and applying it to the $A$'s in the free
algebra.
\item
For every $A \in Obj \hspace{1mm} R_{G'}$ there are arrows
$$\xymatrix@1{I_A:1(A) \ar[r] & A}$$ $$\xymatrix@1{I_A^{-1}:A \ar[r] &
1(A)}.$$ Here $1$ is the unit of the original theory $T$.
\item
For every word $w \in T(m)$, for every function $f:\{1, \dots,m\}
\rightarrow \{1,\dots,n\}$, and for all objects $A_1, \dots, A_n \in
Obj \hspace{1mm} R_{G'}$ there are arrows
$$\xymatrix@1{s_{w,f}(A_1, \ldots, A_n):w_f(A_1, \dots, A_n)
\ar[r] & w(A_{f1}, \dots ,A_{fm})}$$
$$\xymatrix@1{s_{w,f}^{-1}(A_1, \ldots, A_n):w(A_{f1}, \dots ,A_{fm})
\ar[r] & w_f(A_1, \dots, A_n)}.$$ The substituted word $w_f$ is the
substituted word in the original theory $T$. The target $w(A_{f1},
\dots ,A_{fm})$ is the result of substituting in $w$ in the free
theory and then evaluating on the $A$'s.
\item
For every word $w \in S(n)$ and objects $A_1, \ldots, A_n$ of $X$
there are arrows $$\xymatrix@1{\rho^{\eta}_w(A_1, \dots
A_n):\Psi(w)(A_1, \ldots, A_n) \ar[r] & \phi(w)(A_1, \ldots, A_n)}$$
$$\xymatrix@1{\rho^{\eta-1}_w(A_1, \dots A_n):\phi(w)(A_1, \ldots,
A_n) \ar[r] & \Psi(w)(A_1, \ldots, A_n)}.$$
\item
Include also all elements of $Mor \hspace{1mm}  \hspace{1mm} X$.
\end{enumerate}
Then $R_{G'}$ is an object of $Graph'$\index{$Graph'$}. Now we apply
$V'$ to $R_{G'}$ and we get a category $R'$ which is a $T'$-algebra.
The objects of $R_{G'}$\index{$R_{G'}$|)} and $R'$ are the same.

Let $K$ be the smallest congruence\index{congruence} on $R'$ with
the following properties:

\begin{enumerate}
\item
All of the relations necessary to make the coherence arrows
(including $\rho^{\eta}_w$) into natural transformations belong to
$K$. For example, if $A,B \in Obj\hspace{1mm} R'$ and $f:A
\rightarrow B$ is a morphism of $R'$, then the relation $I_A \circ f
= 1(f) \circ I_B$ belongs to $K$.
\item
All of the relations necessary to make the coherence arrows
(including $\rho^{\eta}_w$) into isos are in $K$. For example, for
every $A \in Obj \hspace{1mm} R'$ the relations $I_A \circ I_A^{-1}
= 1_A$ and $I_A^{-1} \circ I_A = 1_A$ are in $K$.
\item
All of the relations for pseudo algebras\index{algebra!pseudo
algebra} listed in Definition \ref{laxalgebradefinition} belong to
$K$, where the objects range over the objects of $R'$.
\item
The original composition relations in the category $X$ belong to
$K$.
\item
The coherence diagrams necessary to make the inclusion $\eta_X:X
\rightarrow UR$ into a morphism of pseudo $S$-algebras are in $K$.
These diagrams are listed in Definition
\ref{defnlaxalgebramorphism}. Note that these coherence diagrams
will involve the arrows $\rho^{\eta}_w(A_1,\dots,A_n):\Psi(w)(A_1,
\ldots, A_n) \rightarrow \phi(w)(A_1,\ldots, A_n)$ for $w \in S(n)$
and objects $A_1,\dots,A_n \in Obj \hspace{1mm} X$.
\item
If the relations $f_1=g_1, \dots ,f_n = g_n$  are in $K$ and $w \in
T'(n)$, then the relation $w(f_1, \dots, f_n) = w(g_1, \dots g_n)$
is also in $K$.
\end{enumerate}
Next mod out by the congruence\index{congruence} $K$ in $R'$ to
obtain the quotient category $R$ called the {\it free pseudo
$T$-algebra on the pseudo $S$-algebra $X$ associated to
$\phi$}\index{algebra!free pseudo $T$-algebra on a pseudo
$S$-algebra|)}\index{free pseudo $T$-algebra on a pseudo
$S$-algebra|)}. We do not use a capital Greek letter to denote the
structure maps of the pseudo $T$-algebra $R$. Instead we write the
words directly.
\end{defn}

In all of the following lemmas in this chapter we use the notation
just introduced in Definition \ref{defnoftheforgetful2functor},
Notation \ref{primenotation}, and Definition
\ref{defnfreelaxTalgebra}.
\begin{lem}
In the notation of the previous definition, the free pseudo
$T$-algebra $R$ on the pseudo $S$-algebra $X$ associated to $\phi$
is a pseudo $T$-algebra.
\end{lem}
\begin{pf}
First we note that $R$ is a (strict) $T'$-algebra. The functor from
the word $w \in T'(n)$ induces a functor on the quotient by relation
6 and the composition and identities in $T'$ are preserved.  The
structure maps have the coherence isos required of a pseudo
$T$-algebra because of the arrows we threw in. The coherence isos
satisfy the required coherence diagrams because of relations 1 and
2. Hence $R$ is a pseudo $T$-algebra.
\end{pf}

\begin{lem}
The inclusion functor denoted $\eta_X:X \rightarrow UR$ is a
morphism of pseudo $S$-algebras.
\end{lem}
\begin{pf}
The inclusion is a functor because of relation 4. It is a morphism
of pseudo $S$-algebras because for all $w \in S(n)$ the natural
transformation $\rho^{\eta}_w:\eta_X \circ \Psi(w) \Rightarrow
\phi(w)(\eta_X, \dots, \eta_X)$ satisfies the required coherences by
the relations in 1. and 5.
\end{pf}

\begin{lem}
For every pseudo $T$-algebra $D$ and every morphism $H:X \rightarrow
UD$ of pseudo $S$-algebras, there exists a morphism $H':R
\rightarrow D$ of pseudo $T$-algebras such that
$$\xymatrix@R=3pc@C=3pc{X \ar@{=}[d] \ar[r]^{\eta_X}
& UR \ar@{.>}[d]^{UH'} & R \ar@{.>}[d]^{H'}
\\ X \ar[r]_{H} & UD & D}$$
commutes.
\end{lem}
\begin{pf}
Let $\Phi$ denote the structure maps of the pseudo $T$-algebra $D$.
As above, $\Psi$ denotes the structure maps of the pseudo
$S$-algebra $X$ and we suppress the capital Greek letter when
denoting the structure maps of the pseudo $T$-algebra $R$. Note that
$D$ is a strict $T'$-algebra and we can therefore apply the
forgetful 2-functor $V:Alg' \rightarrow Graph'$ to it. We also use
$\Phi$ to denote the structure maps of the strict $T'$ algebra $D$.
To construct the morphism $H'$, we define a morphism $H_0':R_{G'}
\rightarrow VD$ in $Graph'$, which induces a morphism $H_1':R'
\rightarrow D$ in $Alg'$ by the definition of the left adjoint to
$V$. Then we show that $H_1'$ preserves the
congruence\index{congruence} $K$ and therefore induces a functor
$H':R \rightarrow D$. Lastly we show that $H'$ is a morphism of
pseudo $T$-algebras such that the desired diagram commutes.

We now define a morphism $H_0':R_{G'} \rightarrow VD$ in $Graph'$.
Defining $H_0'A:=HA$ for $A \in Obj \hspace{1mm} X$ induces a map
$H_0':Obj \hspace{1mm} R_{G'} \rightarrow Obj \hspace{1mm} D$ of
discrete $T'$ algebras. For $f \in Mor \hspace{1mm} X$ define
$H_0'f:=Hf$. For every $w \in S(n)$ and objects $A_1, \dots, A_n \in
Obj\hspace{1mm} X$ let $H_0'$ map the arrows
 $\rho^{\eta}_w(A_1, \dots A_n):\Psi(w)(A_1, \ldots, A_n)
\rightarrow \phi(w)(A_1, \ldots, A_n)$ to the coherence isos
 $\rho^H_w(A_1, \dots A_n):H(\Psi(w)(A_1, \ldots, A_n))
\rightarrow \Phi(\phi(w))(HA_1, \ldots, HA_n)$. Note that the source
and target of $\rho^H_w(A_1, \dots A_n)$ are equal to
$H_0'(\Psi(w)(A_1, \ldots, A_n))$ and $H_0'(\phi(w)(A_1, \ldots,
A_n))$ respectively.  Let $H_0'$ map the other coherence arrows 1
through 3. to the analogous ones in $Mor \hspace{1mm}  D$ with $H_0$
applied to sources and targets. Thus we have defined a morphism
$H_0':R_{G'} \rightarrow VD$ in $Graph'$.

The morphism  $H_0':R_{G'} \rightarrow VD$ in $Graph'$ induces a
morphism $H_1':R' \rightarrow D$ of $Alg'$ by the definition of the
left adjoint to $V$. We claim that $H_1'$ preserves the
congruence\index{congruence} $K$. It suffices to check the relations
1 through 6. We verify them in order of the list above.
\begin{enumerate}
\item
These are satisfied because the analogous arrows for $D$ and $H$ are
natural transformations and $H_1'$ maps coherence arrows to
coherence arrows.
\item
These are satisfied because the analogous arrows for $D$ and $H$ are
isos and $H_1'$ maps coherence arrows to coherence arrows.
\item
The target category $D$ is a pseudo $T$-algebra so these are
satisfied.
\item
The functor $H$ preserves the relations of the category $X$ and
$H_1'$ is defined in terms of $H$, which implies that these are
satisfied.
\item
These are satisfied because $\rho^H_w$ satisfies the coherences and
$H_1'(\rho^{\eta}_w)=\rho^H_w$.
\item
This is by induction. The base case is showing 1 through 5. as was
just done. Suppose the relations $f_1=g_1, \dots, f_n = g_n$ are in
$K$ and $H_1'f_i=H_1'g_i$ for all $i=1,\dots,n$. That is our
induction hypothesis.  Then $$\aligned H_1'(w(f_1, \dots,f_n))&
=\Phi(w)(H_1'(f_1),\dots,H_1'(f_n)) \text{ since
$H_1'$} \\ & \text{\hspace{.2in} is a morphism of $T'$-algebras} \\
&=\Phi(w)(H_1'g_1, \dots,H_1'g_n) \text{ by induction hypothesis} \\
&=H_1'(w(g_1,\dots,g_n)) \text{ since $H_1'$} \\ &
\text{\hspace{.2in} is a morphism of $T'$-algebras.} \endaligned
$$ Thus $H_1'(w(f_1, \dots,f_n))=H_1'(w(g_1,\dots,g_n))$ and
$H_1'$ satisfies this relation.
\end{enumerate}

Since $H_1'$ satisfies the relations, we conclude that $H_1': R'
\rightarrow D$ induces a functor $H':R \rightarrow D$ such that
$H_1'=H' \circ Q$ where $Q:R' \rightarrow R$ is the projection
functor onto the quotient category. The functor $H':R \rightarrow D$
is a morphism of strict $T'$-algebras because for $w \in T'(n)$,
$A_1, \dots,A_n \in Obj \hspace{1mm} R$, and for morphisms $f_1,
\dots, f_n \in Mor \hspace{1mm}  R$ we have
$$\aligned H'(w(A_1, \dots, A_n))&=H_1'(w(A_1, \dots, A_n)) \\
&=\Phi(w)(H_1'A_1, \dots,H_1'A_n) \\
&=\Phi(w)(H'A_1, \dots,H'A_n) \endaligned$$ since $H_1'$ and $H'$
agree on objects. We also have $$\aligned H'(w(f_1, \dots,
f_n))&=H_1'(w(f_1, \dots, f_n)) \\
&=\Phi(w)(H_1'f_1, \dots,H_1'f_n) \\
&=\Phi(w)(H'f_1, \dots,H'f_n) \endaligned $$ where $H_1'$ is
actually applied to representatives of $w(f_1, \dots, f_n), f_1,
\dots, f_n$. Hence $H'$ is a morphism of strict $T'$-algebras and
also a morphism of pseudo $T$-algebras, since $T(n) \subseteq T'(n)$
although this inclusion is not necessarily a map of theories.
According to these two demonstrations, the coherence 2-cells for the
morphism $H'$ of pseudo $T$-algebras are just identities.
\label{coherencetrivial}

We claim that
$$\xymatrix@R=3pc@C=3pc{X \ar@{=}[d] \ar[r]^{\eta_X}
& UR \ar@{.>}[d]^{UH'} & R \ar@{.>}[d]^{H'}
\\ X \ar[r]_{H} & UD & D}$$
commutes. It is sufficient to check this for the underlying functors
and the coherence 2-cells. The underlying functor of $H'$ is the
same as the underlying functor of $UH'$. Let $A \in Obj \hspace{1mm}
X$. Then $UH' \circ \eta_X(A)=UH'(A)=H'A=HA$. Similarly, for $f \in
Mor \hspace{1mm}  X$ we have $UH' \circ \eta_X(f)=UH'(f)=H'f=Hf$.
Hence the diagram commutes. The coherence 2-cells also commute
because $H'(\rho^{\eta}_w)=\rho^H_w$ and because the coherence
2-cells of $H'$ are identities.
\end{pf}

\begin{lem} \label{forgetfuluniversal} \label{inclusionuniversal}
The inclusion morphism $\eta_X:X \rightarrow UR$ is a
biuniversal\index{biuniversal arrow}\index{arrow!biuniversal arrow}
arrow from $X$ to the forgetful 2-functor\index{forgetful
2-functor}\index{2-functor!forgetful 2-functor}.
\end{lem}
\begin{pf}
Let $D$ be a pseudo $T$-algebra. Let $Mor_S(X,UD)$ denote the
category of morphisms of pseudo $S$-algebras from $X$ to $UD$. Let
$Mor_T(R,D)$ denote the category of morphisms of pseudo $T$-algebras
from $R$ to $D$.  Let $\phi:Mor_T(R,D) \rightarrow Mor_S(X,UD)$ be
the functor defined by $H' \mapsto UH' \circ \eta_X$ and $\gamma
\mapsto U \gamma * i_{\eta_X}$. Define a functor $\psi:Mor_S(X,UD)
\rightarrow Mor_T(R,D)$ as follows. For $H \in Obj \hspace{1mm}
Mor_S(X,UD)$ let $\psi H := H'$ where $H':R \rightarrow D$ is the
morphism of pseudo $T$ algebras constructed in the previous lemma.

If $H,J \in Obj \hspace{1mm} Mor_S(X,UD)$ and $\beta:H \Rightarrow
J$ is a 2-cell in the 2-category of pseudo $S$-algebras, define
$\psi(\beta)=\beta':H' \Rightarrow J'$ inductively as follows. If $A
\in Obj\hspace{1mm} X$ then define $\beta'A$ to make
$$\xymatrix@R=3pc@C=3pc{H'A \ar[r]^{\beta'A} \ar@{=}[d] & J'A \ar@{=}[d] \\ HA \ar[r]_{\beta A} &
JA}$$ commute. If $w \in T'(n)$ and $\beta'$ is already defined for
$A_1, \dots, A_n \in Obj \hspace{1mm}  R$, then $\beta'(w(A_1,
\dots, A_n)) := \Phi(w)(\beta'A_1, \dots, \beta' A_n)$. The
following inductive proof shows that $\beta':H' \Rightarrow J'$ is a
natural transformation. For $f \in Mor \hspace{1mm}  X$ the
naturality of $\beta'$ is guaranteed by the naturality of $\beta:H
\Rightarrow J$. The naturality of $\beta'$ for the coherence isos
thrown into the category $R$ during its construction follows because
$H'$ and $J'$ take coherence isos of $R$ to analogous ones in $D$
and the coherences isos in $D$ are natural. That concludes the base
case for the induction. Now suppose $\beta'$ is natural for
morphisms $f_i \in Mor_R(A_i,B_i)$ for $i=1, \dots, n$ and $w \in
T'(n)$. Then
$$\xymatrix@C=6pc@R=4pc{H'w(A_1, \dots, A_n) \ar[r]^{\beta'w(A_1, \dots, A_n)}
\ar[d]_{H'w(f_1, \dots, f_n)} & J'w(A_1, \dots, A_n)
\ar[d]^{J'w(f_1, \dots, f_n)}
\\ H'w(B_1, \dots, B_n) \ar[r]_{\beta'w(B_1, \dots, B_n)}
& J'w(B_1, \dots, B_n)}$$ commutes because $w$ commutes with
everything in the diagram by definition and because we apply the
functor $\Phi(w)$ to each of the individual naturality diagrams for
$f_i:A_i \rightarrow B_i$ and $i=1, \dots, n$. Hence $\beta'$ is
natural for any morphism in $R$ by this inductive proof. Moreover,
the natural transformation commutes appropriately with $\rho^{H'}$
and $\rho^{J'}$ because they are trivial and $\beta'(w(A_1, \dots,
A_n)$$=\Phi (w)(\beta'A_1, \dots, \beta'A_n)$. Hence
$\psi(\beta)=\beta'$ is a 2-cell in the 2-category of pseudo
$T$-algebras.

It is routine to check inductively that the assignment
$\psi:Mor_S(X,UD) \rightarrow Mor_T(R,D)$ preserves identities and
compositions and is thus a functor.

We claim that $\psi$ is a right adjoint for $\phi$. By the previous
lemma $\phi \circ \psi (H)=H$ for all $H \in Obj \hspace{1mm}
Mor_S(X,UD)$. We easily see that $\phi \circ \psi (\beta) = \beta$
for all $\beta \in Mor \hspace{1mm} Mor_S(X,UD)$. Hence the counit
$\mu:\phi \circ \psi \Rightarrow 1_{Mor_S(X,UD)}$ is the identity
natural \label{counittrivial} transformation, which is of course a
natural isomorphism. Next we define a unit $\theta:1_{Mor_T(R,D)}
\Rightarrow \psi \circ \phi$. For $J' \in Mor_T(R,D)$ let $H':=\psi
\circ \phi(J')$. Recall that $H'$ is strict, \ie $\rho^{H'}$ is
trivial, while $J'$ may not be strict. We define a 2-cell
$\theta(J'):J' \Rightarrow H'=\psi \circ \phi(J')$ in the category
of pseudo $T$-algebras inductively. For $A \in Obj \hspace{1mm} X
\subseteq Obj \hspace{1mm} R$ set $\theta(J')(A) :=1_{J'A}$. Suppose
$w \in T'(n)$ and $\theta(J')$ is already defined for $A_1, \dots,
A_n \in Obj \hspace{1mm} R$. Then define $\theta(J')(w(A_1,
\dots,A_n):J'(w(A_1, \dots, A_n)) \rightarrow H'(w(A_1, \dots,
A_n))$ by $\Phi(w)(\theta(J')A_1, \dots, \theta(J')A_n) \circ
\rho^{J'}_w(A_1, \dots, A_n)$. An inductive proof, similar to the
one above but also using the naturality of $\rho^{J'}_w$, shows that
$\theta(J')$ is a natural transformation and commutes with
$\rho^{J'}$ and $\rho^{H'}$ appropriately, \ie $\theta(J'):J'
\Rightarrow H'$ is a 2-cell. It is also iso by induction. The
assignment $J' \mapsto \theta(J')$ is natural by an inductive
argument that uses the diagram in the definition of 2-cell in the
2-category of pseudo $T$-algebras. Hence $\theta:1_{Mor_T(R,D)}
\Rightarrow \psi \circ \phi$ is a natural isomorphism. If we can
show that $\theta$ and $\mu$ satisfy the triangular identities, then
we can conclude that $\psi$ is a right adjoint for $\phi$

We claim that the unit $\theta$ and the counit $\mu$ satisfy the
triangular identities. First we show that
\begin{equation} \label{triangle1}
\xymatrix@1@C=4pc{\psi \ar@{=>}[r]^-{\theta * i_{\psi}} & \psi \circ
\phi \circ \psi \ar@{=>}[r]^-{i_{\psi}*\mu} & \psi}
\end{equation}
is the identity natural transformation $i_{\psi}:\psi \Rightarrow
\psi$. Let $H \in Obj \hspace{1mm} Mor_S(X, UD)$. Then
$$\aligned (i_{\psi} * \mu) \odot (\theta * i_{\psi})(H) &= \psi(\mu_H) \circ \theta_{\psi
H} \text{ by definition} \\ &=\theta_{\psi H} \text{ since $\mu_H$
is trivial.} \endaligned$$ But $\theta_{\psi H}=\theta(\psi H)$ is
the trivial 2-cell $\psi H \Rightarrow \psi H$ because $\psi H$ is a
strict morphism of pseudo $T$-algebras, \ie $\rho^{\psi H}_w$ is
trivial. Hence (\ref{triangle1}) is $i_{\psi}:\psi \Rightarrow
\psi$. Next we show that
\begin{equation} \label{triangle2}
\xymatrix@1@C=4pc{\phi \ar@{=>}[r]^-{i_{\phi}*\theta} & \phi \circ
\psi \circ \phi \ar@{=>}[r]^-{\mu * i_{\phi}} & \phi}
\end{equation}
is the identity natural transformation $i_{\phi}:\phi \Rightarrow
\phi$. Let $J' \in Obj \hspace{1mm} Mor_T(R,D)$. Then
$$\aligned (\mu * i_{\phi}) \odot (i_{\phi} * \theta)(J')&=\mu_{\phi J'} \circ
\phi(\theta_{J'}) \text{ by definition} \\ &=\phi(\theta_{J'})
\text{ since $\mu_{\phi J'}$ is trivial} \\ &=\theta_{J'} *
i_{\eta_X} \text{ by definition.} \endaligned $$ But $\theta_{J'}
* i_{\eta_X}$ is the trivial 2-cell $\phi(J')=J' \circ \eta_X
\Rightarrow J' \circ \eta_X$ because
$\theta_{J'}(A)=\theta(J')(A)=1_{J'A}$ for all $A \in Obj
\hspace{1mm} X$ and $\eta_X:X \rightarrow R$ is the inclusion
functor. Hence (\ref{triangle2}) is the identity natural
transformation $i_{\phi}:\phi \Rightarrow \phi$. Thus the unit and
counit satisfy the triangular identities and $\psi$ is a right
adjoint for $\phi$. Moreover, $\phi$ is an equivalence because the
unit and counit are natural isomorphisms. We conclude that $\eta_X:
X \rightarrow UR$ is a biuniversal arrow from $X$ to the 2-functor
$U$.
\end{pf}

\begin{rmk}
Although it is not necessary, we can construct the
factorizing\index{factorizing 2-cell} 2-cell $\nu'$ on page
\pageref{factorizing2cellnuprime} as follows. Let $H:X \rightarrow
UD$ be a morphism of pseudo $S$-algebras. Then $\psi(H)=H'$
satisfies
$$\xymatrix@C=4pc@R=3pc{X \ar[r]^{\eta_X} \ar@{=}[d]
& UR \ar@{:>}[ld]_{\mu(H)} \ar@{.>}[d]^{UH'} & R \ar@{.>}[d]^{H'}
\\ X \ar[r]_H & UD & D}$$
and $\mu(H)$ is the identity 2-cell. Suppose $\bar{H}':R \rightarrow
D$ is another morphism of pseudo $T$-algebras and $\nu$ is a 2-cell
as follows.
$$\xymatrix@C=4pc@R=3pc{X \ar@{=}[d] \ar[r]^{\eta_X} & U(R)  \ar[d]^{U
\bar{H}'} \ar@{=>}[dl]_{\nu } & R \ar[d]^{\bar{H}'} \\ X
 \ar[r]_{H}
 & JD & D}$$
Define a 2-cell $\nu':\bar{H}' \Rightarrow H'$ as follows. For $A
\in Obj \hspace{1mm} X \subseteq Obj \hspace{1mm} R$, $\nu' A:=\nu
A$. If $w \in T'(n)$ and $\nu'$ is already defined for $A_1, \dots,
A_n \in Obj \hspace{1mm} R$, then $\nu'(w(A_1, \dots, A_n)):=
\Phi(w)(\nu'A_1, \dots, \nu'A_n) \circ \rho^{\bar{H}'}(A_1, \dots,
A_n)$. By induction $\nu'$ is a natural transformation. It also
commutes with $\rho^{H'}$ and $\rho^{\bar{H}'}$ appropriately by
construction. Hence $\nu'$ is a 2-cell in the 2-category of pseudo
$T$-algebras. By construction we see that
\begin{equation} \label{nu'unique}
\xymatrix@C=4pc@R=3pc{\bar{H}' \ar@{:>}[d]_{\nu'} & U\bar{H}' \circ
\eta_X \ar@{=>}[r]^-{\nu} \ar@{:>}[d]_{U \nu' * i_{\eta_X}} & H
\ar@{=}[d]
\\ H' & UH' \circ \eta_X \ar@{=>}[r]_-{\mu(H)} & H}
\end{equation}
commutes. Such a 2-cell $\nu'$ is unique by the requirement that
(\ref{nu'unique}) commutes and by the commutivity with
$\rho^{\bar{H}'}$ and $\rho^{H'}$ required of 2-cells $\bar{H'}
\Rightarrow H'$. More precisely, the commutivity of
(\ref{nu'unique}) says that $\nu'A = \nu A$ for all $A \in Obj
\hspace{1mm} X$ and the appropriate commutivity with
$\rho^{\bar{H}'}$ and $\rho^{H'}$ specifies what $\nu'$ does to
objects of the form $w(A_1, \dots, A_n)$ for $A_1, \dots, A_n \in
Obj \hspace{1mm} R$. If $\nu$ is iso, then so is $\nu'$ by the
construction and the fact that $\rho^{\bar{H}'}$ is iso.
\end{rmk}

\begin{thm} \label{forgetfuladjoint}
Let $S$ and $T$ be theories and $\phi:S \rightarrow T$ a morphism of
theories. Then the forgetful 2-functor\index{forgetful
2-functor|textbf}\index{2-functor!forgetful 2-functor|textbf} $U$
associated to $\phi$ from the 2-category of small pseudo
$T$-algebras to the 2-category of small pseudo $S$-algebras admits a
left biadjoint\index{biadjoint!left biadjoint|textbf} denoted $F$.
Moreover, this pseudo functor $F$ is actually a strict 2-functor.
\end{thm}
\begin{pf}
For every pseudo $S$-algebra $X$ there exists a pseudo $T$-algebra
$R$ and a biuniversal arrow $\eta_X:X \rightarrow UR$ by Lemma
\ref{forgetfuluniversal}. This guarantees the existence of a left
biadjoint by Theorem \ref{maintheorem2}.

We can prove that $F$ is strict by inspecting its coherence isos
constructed in the general theory of Theorem \ref{maintheorem2}. Let
$\mathcal{X}$ be the 2-category of pseudo $S$-algebras, let
$\mathcal{A}$ be the 2-category of pseudo $T$-algebras, and let
$G:=U:\mathcal{A} \rightarrow \mathcal{X}$ be the forgetful
2-functor. For any pseudo $S$-algebra $X \in Obj \hspace{1mm}
\mathcal{X}$, we define $FX$ to be the free\index{algebra!free
pseudo $T$-algebra on a pseudo $S$-algebra}\index{free pseudo
$T$-algebra on a pseudo $S$-algebra} pseudo $T$-algebra $R$ on the
pseudo $S$-algebra $X$ associated to the morphism of theories
$\phi:S \rightarrow T$. The co-unit $\mu$ for the biuniversal
$\eta_X:X \rightarrow UR$ is the identity as we observed in Lemma
\ref{inclusionuniversal}. The pseudo functor $U=G$ is actually a
strict 2-functor, so $\delta^G$ and $\gamma^G$ are identity natural
transformations. After inspecting diagram (\ref{Fdelta}) on page
\pageref{Fdelta}, we see that $\delta^F_*$ must be trivial because
$(\delta^G_{FX*})^{-1}*i_{\eta_X}$ and $\mu_{X,FX}(\eta_X \circ
1_X)= \mu(\eta_X \circ 1_X)$  are trivial. Hence $F$ preserves
identities.

Similarly, each of the 2-cells in diagram (\ref{gammaF1}) on page
\pageref{gammaF1} is trivial, and therefore their composition is
trivial. After inspecting diagram (\ref{gammaF2}) on page
\pageref{gammaF2}, we see that $\gamma^F_{f,g}$ must also be trivial
because both the horizontal top and bottom arrows are trivial.
Therefore $F$ preserves compositions.

Since $F$ preserves compositions and identities, it is a strict
2-functor.
\end{pf}

\begin{thm}
The biuniversal\index{biuniversal arrow}\index{arrow!biuniversal
arrow} arrows $\eta_X:X \rightarrow UFX$ define a strict 2-natural
transformation $\eta:1_{\mathcal{X}} \Rightarrow U \circ F$, where
$\mathcal{X}$ is the 2-category of pseudo $S$-algebras.
\end{thm}
\begin{pf}
Recall that the counits $\mu$ for the biuniversal arrows $\eta_X$
are all trivial as indicated on page \pageref{counittrivial} in
Lemma \ref{inclusionuniversal}. In the proof of Theorem
\ref{maintheorem2} on page \pageref{makeetalaxnatural} the
biuniversal arrows $\eta_X:X \rightarrow UFX$ are made into a pseudo
natural transformation by defining $\tau_f:=\mu_{X,FY}(\eta_Y \circ
f)$ for $f:X \rightarrow Y$. We see that $\tau_f$ is trivial because
$\mu_{X,FY}$ is trivial. Hence $\eta$ is strictly 2-natural.
\end{pf}

Theorem \ref{forgetfuladjoint} can be sharpened. Let $\mathcal{A}$
denote the 2-category of pseudo $T$-algebras and let $\mathcal{X}$
denote the 2-category of pseudo $S$-algebras. Then the equivalence
of categories $Mor_{\mathcal{A}}(FX,A) \rightarrow
Mor_{\mathcal{X}}(X, UA)$ implicit in Theorem \ref{forgetfuladjoint}
is strictly 2-natural in each variable. However, it can be shown
that a left 2-adjoint does not exist in specific cases. The
equivalence in the other direction $Mor_{\mathcal{X}}(X, UA)
\rightarrow Mor_{\mathcal{A}}(FX,A) $ in Theorem
\ref{forgetfuladjoint} is not strictly 2-natural in each variable.
In fact, there is an example where there does not exist an
equivalence $Mor_{\mathcal{X}}(X, UA) \rightarrow
Mor_{\mathcal{A}}(FX,A) $ which is strictly 2-natural in each
variable, even after replacing $F$ by another biadjoint $F'$.
Counterexamples will be given after presenting Theorem
\ref{sharperadjoint}, which is a sharper version of Theorem
\ref{forgetfuladjoint}.

\begin{thm} \label{sharperadjoint}
Let $S$ and $T$ be theories. Let $U:\mathcal{A} \rightarrow
\mathcal{X} $ be the forgetful 2-functor\index{2-functor!forgetful
2-functor}\index{forgetful 2-functor} associated to a morphism $S
\rightarrow T$ of theories. Let $F$ denote the left
biadjoint\index{biadjoint!left biadjoint} to $U$ introduced in
Theorem \ref{forgetfuladjoint}. Then the equivalence of categories
$\phi_{X,A}:Mor_{\mathcal{A}}(FX,A) \rightarrow Mor_{\mathcal{X}}(X,
UA)$ from Theorem \ref{forgetfuladjoint} defined by
$\phi_{X,A}(f):=Uf \circ \eta_X$ is strictly 2-natural in each
variable.
\end{thm}
\begin{pf}
The universal arrow $\eta_X:X \rightarrow UFX$ is the inclusion
morphism. The functor $\phi_{X,A}:Mor_{\mathcal{A}}(FX,A)
\rightarrow Mor_{\mathcal{X}}(X, UA)$ is defined by
$\phi_{X,A}(f):=Uf \circ \eta_X$ as in Lemma \ref{converse}. The
functor $\phi_{X,A}$ is an equivalence of categories for all $X \in
Obj \hspace{1mm} \mathcal{X}$ and all $A \in Obj \hspace{1mm}
\mathcal{A}$ because $\eta_X$ is a biuniversal arrow. The coherence
isos $\tau'$ for the pseudo naturality of $\phi_{-,A}$ are defined
on page \pageref{firstcoherenceisoforphi} in terms of some trivial
2-cells, $\gamma^G$, and $\tilde{\tau}$, where $\tilde{\tau}$ is the
coherence iso for $\eta$. But $\gamma^G$ is trivial for $G=U$
because $U$ is a strict 2-functor. The coherence iso $\tilde{\tau}$
is also trivial because $\eta$ is a strict 2-natural transformation.
Hence $\tau'$ is also trivial and $\phi_{-,A}$ is strictly
2-natural, \ie $\phi$ is 2-natural in the first variable.

The coherence isos $\tau$ for $\phi_{X,-}$ are defined on page
\pageref{secondcoherenceisoforphi} for morphisms $k:A \rightarrow
A'$ by $\tau_{A,A'}(k):e \mapsto \gamma^G_{e,k}*i_{\eta_X}$. But
$G=U$ is a strict functor and $\gamma^G$ is trivial, hence $\tau$ is
also trivial. Therefore $\phi_{X,-}$ is strictly 2-natural, \ie
$\phi$ is 2-natural in the second variable.
    We conclude that $X,A \mapsto \phi_{X,A}$ is strictly 2-natural in
each variable.
\end{pf}

Before proving that Theorem \ref{forgetfuladjoint} cannot be further
improved to a left 2-adjoint, we need a theorem which states that we
can change a morphism of pseudo $T$-algebras in a specific way and
still have a morphism of pseudo $T$-algebras.

\begin{thm} \label{replacement}
Let $X,Y$ be pseudo $T$-algebras and $H:X \rightarrow Y$ a
morphism\index{morphism!morphism of pseudo $T$-algebras} of pseudo
$T$-algebras. Suppose that $J_0(x) \in Obj \hspace{1mm} Y$ and
$\alpha_0(x):J_0(x) \rightarrow H(x)$ is an isomorphism for each $x
\in Obj \hspace{1mm} X$. Then there exists a morphism $J:X
\rightarrow Y$ of pseudo $T$-algebras whose object function is $J_0$
and there exists an iso 2-cell $\alpha:J \Rightarrow H$ of pseudo
$T$-algebras such that $\alpha(x)=\alpha_0(x)$ for all $x \in Obj
\hspace{1mm} X$. Moreover, such $J$ and $\alpha$ are unique.
\end{thm}
\begin{pf}
For $x \in Obj \hspace{1mm}  X$ define $J(x):=J_0(x)$ and
$\alpha(x):=\alpha_0(x)$. For a morphism $f:x_1 \rightarrow x_2$ of
$X$ define $J(f):=\alpha(x_2)^{-1} \circ H(f) \circ \alpha(x_1)$. We
easily see that $J$ is a functor and $\alpha$ is natural
transformation from $J$ to the functor underlying $H$.

For $w \in T(n)$ let $\rho^H_w:H \circ \Phi(w) \Rightarrow \Psi
\circ (H, \dots, H)$ denote the coherence isomorphism for $H$, where
$\Phi$ and $\Psi$ denote the structure maps of $X$ and $Y$
respectively.  Define a natural isomorphism $\rho^J_w:J \circ
\Phi(w) \Rightarrow \Psi \circ (J, \dots, J)$ by the following
diagram.
$$\xymatrix@C=6pc@R=4pc{J \circ \Phi(w) \ar@{=>}[r]^{\alpha \ast i_{\Phi(w)}}
\ar@{=>}[d]_{\rho_w^J}
 & H \circ \Phi(w)
\ar@{=>}[d]^{\rho_w^H} \\ \Psi(w) \circ (J, \dots, J)
\ar@{=>}[r]_-{i_{\Psi(w)} \ast (\alpha, \dots, \alpha)} & \Psi(w)
\circ(H, \dots,H) }$$ In other words $\rho^J_w:=(i_{\Psi(w)} \ast
(\alpha^{-1}, \dots, \alpha^{-1})) \odot \rho_w^H \odot (\alpha \ast
i_{\Phi(w)})$. This is a natural transformation because it consists
of horizontal and vertical compositions of natural transformations.

We claim that $\rho^J_w$ satisfies the coherence diagrams required
to make $J$ a morphism of pseudo $T$-algebras. We can prove the
commutivity of any $J$ coherence diagram from the commutivity of the
analogous $H$ coherence diagram by using the following procedure.
First we draw the commutative $H$ coherence diagram and then we
circumscribe it with the analogous $J$ coherence diagram. Next we
draw the obvious isomorphisms between respective $J$ and $H$
vertices. All of the resulting inner diagrams commute because of the
interchange law, because of the definition of $\rho^J_w$, or because
of the diagram for $H$.

We present the substitution diagram to clarify the process. Let
$f:\{1, \dots, m \} \rightarrow \{1, \dots, n\}$ be a function and
$w \in T(m)$.

\begingroup
\vspace{-2\abovedisplayskip} \small
$$\xymatrix@R=4pc@C=1pc{J \circ \Phi(w_f) \ar@{=>}[ddd]_{\rho^J_{w_f}}
\ar@{=>}[rrr]^{i_J*s_{w,f}}  \ar@{=>}[dr]^{\alpha*i_{\Phi(w_f)}} & &
& J \circ \Phi(w)_f \ar@{=>}[ddd]^{(\rho^J_w)_f}
\ar@{=>}[dl]_{\alpha*i_{\Phi(w)_f} \hspace{1em}}
\\ & H \circ \Phi(w_f) \ar@{=>}[r]^{i_H*s_{w,f}} \ar@{=>}[d]_{\rho^H_{w_f}}
 & H \circ \Phi(w)_f \ar@{=>}[d]^{(\rho^H_w)_f} &
\\ & \Psi(w_f) \circ (H, \dots, H) \ar@{=>}[r]_{\overset{\phantom{l}}{s_{w,f}*i_{(H,\dots,H)}}}
\ar@{=>}[dl]^{\hspace{2em} i_{\Psi(w_f)}*(\alpha^{-1},\dots,
\alpha^{-1})} & \Psi(w)_f \circ (H, \dots, H)
\ar@{=>}[dr]_{i_{\Psi(w)_f}*(\alpha^{-1},\dots, \alpha^{-1})
\hspace{2em}} &
\\ \Psi(w_f) \circ (J, \dots, J) \ar@{=>}[rrr]_{s_{w,f}*i_{(J,\dots,J)}}
& & & \Psi(w)_f \circ (J, \dots, J)}$$
\endgroup
\noindent
 The top and bottom squares commute because of the
interchange law. The left and right squares commute because of the
definitions of $\rho^J_{w_f}$ and $\rho^J_w$. The innermost square
commutes because $H$ is a morphism of pseudo $T$-algebras. Hence the
outer rectangle commutes and $J$ satisfies the substitution
coherence diagram.

The other diagrams can be verified using the same procedure. The
only subtlety in this procedure occurs in the right hand vertical
composition of the composition axiom. We reproduce the right hand
part of that diagram obtained by the procedure mentioned above.

\begingroup
\vspace{-2\abovedisplayskip} \small
$$\xymatrix@R=4pc@C=2pc{... \ar@{=>}[r]^-{i_J*c_{w, w_1, \dots, w_n}}
& J \circ \Phi(w) \circ (\Phi(w_1), \dots, \Phi(w_n))
\ar@{=>}[dd]^{\rho^J_w*i_{(\Phi(w_1), \dots, \Phi(w_n))}}
\ar@{=>}[dl]_{\alpha*i_{\Phi(w)}*i_{(\Phi(w_1), \dots, \Phi(w_n))}
 \hspace{4em} \hspace{5mm}}
\\  H \circ \Phi(w) \circ (\Phi(w_1), \dots \Phi(w_n))
\ar@{=>}[d]_{\rho^H_w * i_{(\Phi(w_1),\dots,\Phi(w_n))}} & &
\\  \Psi(w) \circ (H, \dots, H) \circ (\Phi(w_1), \dots, \Phi(w_n))
\ar@{=>}[r]^{\underset{\phantom{l}}{i_{\Psi(w)}*(\alpha^{-1}, \dots,
\alpha^{-1})* i_{(\Phi(w_1), \dots, \Phi(w_n))}}}
\ar@{=>}[d]_{i_{\Psi(w)}*(\rho^H_{w_1}, \dots, \rho^H_{w_n})} &
\Psi(w) \circ (J, \dots, J) \circ (\Phi(w_1), \dots, \Phi(w_n))
\ar@{=>}[dd]^{i_{\Psi(w)}*(\rho^J_{w_1}, \dots, \rho^J_{w_n})}
\\  \Psi(w) \circ (\Psi(w_1), \dots, \Psi(w_n)) \circ (H, \dots, H)
\ar@{=>}[dr]_{i_{\Psi(w)}*i_{(\Psi(w_1), \dots,
\Psi(w_n))}*(\alpha^{-1}, \dots, \alpha^{-1}) \hspace{6em}
\hspace{5mm}} & &
\\... \ar@{=>}[r]_-{c_{w, w_1, \dots, w_n}*i_{(J, \dots, J)}}
 & \Psi(w) \circ (\Psi(w_1), \dots, \Psi(w_n)) \circ (J, \dots,
 J)}$$
\endgroup
\noindent

The upper right quadrilateral results from the diagram defining
$\rho^J_w$ by horizontally composing with $i_{(\Phi(w_1), \dots,
\Phi(w_n))}$. Then the upper right square commutes by iterated use
of the interchange law.

The bottom right quadrilateral results from the defining diagrams of
$\rho^J_{w_1}, \dots, \rho^J_{w_n}$ by taking their product,
horizontally composing with the identity 2-cell $$i_{(\Psi(w_1),
\dots, \Psi(w_n))}=(i_{\Psi(w_1)}, \dots, i_{\Psi(w_n)}),$$ and
finally reversing one of the arrows. The commutivity then follows
from the interchange law.

The other parts of the diagram are easily seen to commute, and we
conclude that $J$ satisfies the composition coherence.

The commutivity of all of these coherence diagrams implies that $J$
is a morphism of pseudo $T$-algebras. We conclude that $\alpha$ is a
2-cell in the 2-category of pseudo $T$-algebras by looking at its
defining diagram.

Now we turn to the uniqueness. Suppose $J':X \rightarrow Y$ is a
morphism of pseudo $T$-algebras and $\alpha':J' \Rightarrow H$ is a
2-cell in the 2-category of pseudo $T$-algebras such that for all $x
\in Obj \hspace{1mm} X$ we have $J'(x)=J_0(x)$ and
$\alpha'(x)=\alpha_0(x)$. Then for a morphism $f:x_1 \rightarrow
x_2$ in $X$ the diagram
$$\xymatrix@R=3pc@C=3pc{J_0(x_1) \ar[r]^{\alpha(x_1)} \ar[d]_{J'(f)} & H(x_1)
\ar[d]^{H(f)}
\\ J_0(x_1) \ar[r]_{\alpha(x_2)} & H(x_2)}$$
commutes. Hence $J'(f)=\alpha(x_2)^{-1} \circ H(f) \circ
\alpha(x_1)=J(f)$. For a word $w \in T(n)$, the diagram
$$\xymatrix@C=6pc@R=4pc{J' \circ \Phi(w) \ar@{=>}[r]^{\alpha \ast i_{\Phi(w)}}
\ar@{=>}[d]_{\rho_w^{J'}}
 & H \circ \Phi(w)
\ar@{=>}[d]^{\rho_w^H} \\ \Psi(w) \circ (J', \dots, J')
\ar@{=>}[r]_{i_{\Psi(w)} \ast (\alpha, \dots, \alpha)} & \Psi(w)
\circ(H, \dots,H) }$$ commutes. Hence $\rho^{J'}_w=(i_{\Psi(w)} \ast
(\alpha^{-1}, \dots, \alpha^{-1})) \odot \rho_w^H \odot (\alpha \ast
i_{\Phi(w)}) =\rho^J_w$. We conclude $J'=J$ as morphisms of pseudo
$T$-algebras.
\end{pf}

\begin{lem} \label{not2natural}
The functor $\psi_{X,A}:Mor_{\mathcal{X}}(X, UA) \rightarrow
Mor_{\mathcal{A}}(FX,A) $ in Theorem \ref{forgetfuladjoint} is not
strictly 2-natural in each variable.
\end{lem}
\begin{pf}
Suppose $\psi$ is strictly 2-natural. Then for any morphism of
pseudo $T$-algebras $J:FX \rightarrow FX$ the following diagram must
commute.
\begin{equation} \label{strictoutput}
\xymatrix@C=5pc@R=3pc{Mor_{\mathcal{A}}(FX, FX) \ar[d]_{J_*} &
\ar[l]_{\psi_{X,FX}} \ar[d]^{(UJ)_*} Mor_{\mathcal{X}}(X,UFX)
\\ Mor_{\mathcal{A}}(FX, FX) & Mor_{\mathcal{X}}(X,UFX) \ar[l]^{\psi_{X,FX}}
}
\end{equation}
According to page \pageref{coherencetrivial}, the output
$\psi_{X,FX}(H)$ is always a strict morphism of pseudo $T$-algebras
for all morphisms $H:X \rightarrow UFX$ of pseudo $S$-algebras. Let
$a \in Obj \hspace{1mm} FX$. Let $w$ be the trivial word in the
theory $T$. Then $w(\psi_{X,FX}(\eta_X)(a))$ is isomorphic to (but
not equal to) $\psi_{X,FX}(\eta_X)(a)$ via a coherence isomorphism.
By Theorem \ref{replacement} we can construct from this data a
morphism $J:FX \rightarrow FX$ of pseudo $T$-algebras such that
$J(w(\psi_{X,FX}(\eta_X)(a)))=\psi_{X,FX}(\eta_X)(a)$ and $J$ is the
identity on all other objects. Chasing $\eta_X$ along diagram
(\ref{strictoutput}) from the top right corner, we see that
$\psi_{X,FX}(UJ \circ \eta_X)=J \circ \psi_{X,FX}(\eta_X)$ and $J
\circ \psi_{X,FX}(\eta_X)$ must be strict because $\psi_{X,FX}(UJ
\circ \eta_X)$ is. But $J \circ \psi_{X,FX}(\eta_X)$ is not strict
because it does not commute with the application of $w$ by the
construction of $J$.
\end{pf}

In fact, we present an example where there is no pseudo natural
transformation $\psi$ as in Lemma \ref{not2natural} that is strictly
2-natural in the second variable, even after replacing $F$ by
another left biadjoint to $U$. The reason is that our morphisms of
pseudo algebras\index{algebra!pseudo algebra} are not required to be
strict, \ie they are not required to commute with the structure
maps.

\begin{examp} \label{nostrictification}
Let $S$ be the trivial theory\index{theory!trivial
theory}\index{trivial theory} and let $T$ be the theory of
commutative monoids\index{theory!theory of commutative
monoids}\index{commutative monoid!theory of commutative monoids}.
Let $\mathcal{X}$ be the 2-category of pseudo $S$-algebras and let
$\mathcal{A}$ be the 2-category of pseudo $T$-algebras. Let
$U:\mathcal{A} \rightarrow \mathcal{X}$ be the forgetful
2-functor\index{forgetful 2-functor}\index{2-functor!forgetful
2-functor} associated to the trivial map of theories $S \rightarrow
T$. Then there does not exist a left biadjoint\index{biadjoint!left
biadjoint} $F': \mathcal{X} \rightarrow \mathcal{A}$ which admits
equivalences of categories $\psi'_{X,A}:Mor_{\mathcal{X}}(X,UA)
\rightarrow Mor_{\mathcal{A}} (F'X,A)$ that are strictly 2-natural
in the second variable.
\end{examp}
\begin{pf}
First we prove that our constructed left biadjoint $F:\mathcal{X}
\rightarrow \mathcal{A}$ does not admit equivalences $\psi_{X,A}'$
that are strictly 2-natural in the second variable. Suppose for each
$X \in Obj \hspace{1mm} \mathcal{X}$ there exist equivalences
$\psi'_{X,A}:Mor_{\mathcal{X}}(X,UA) \rightarrow
Mor_{\mathcal{A}}(FX,A)$ that are strictly natural in $A$, the
second variable. Let $\phi_{X,A}'$ be a functor such that
$\phi_{X,A}' \circ \psi_{X,A}'$ and $\psi_{X,A}' \circ \phi_{X,A}'$
are naturally isomorphic to the respective identities.

Let $X$ be the pseudo $S$-algebra with only one object $*$ and no
nontrivial morphisms. Let $A$ be the category of finite sets with a
choice of disjoint union\index{disjoint union}. This makes $A$ into
a pseudo $T$-algebra.

We claim that there exists a morphism $H:X \rightarrow UA$ of pseudo
$S$-algebras such that $\psi_{X,A}'(H)(*) \neq \emptyset$. Suppose
not. Then for every morphism $H:X \rightarrow UA$, we have
$\psi_{X,A}'(H)(w(*, \dots, *)) \cong w(\emptyset, \dots,
\emptyset)=\emptyset$ and thus $\psi_{X,A}'(H)$ is constant
$\emptyset$. By the equivalence, every morphism $K:FX \rightarrow A$
of pseudo $T$-algebras is isomorphic to $\psi_{X,A}' \circ
\phi_{X,A}' (K)$. This implies that $K$ must also be constant
$\emptyset$. But this is a contradiction, since there are nontrivial
morphisms $FX \rightarrow A$. Thus there exists a morphism $H:X
\rightarrow UA$ of pseudo $S$-algebras such that $\psi_{X,A}'(H)(*)
\neq \emptyset$.

We claim that there exists an object $x \in Obj \hspace{1mm} FX$
such that $\psi_{X,A}'(H)(x) \neq H(*)$. Let $n \in \mathbb{N}$ be
large enough that
$$n \cdot | \psi_{X,A}'(H)(*) | > |H(*) |. $$
This is possible because $|\psi'_{X,A}(H)(*)| \neq 0$ from above.
Let $x = *+(*+(*+\cdots))$ where there are $n$ copies of $*$. Then $
| \psi_{X,A}'(H)(x) | =  n \cdot | \psi_{X,A}'(H)(*) | $ because
$\psi_{X,A}'(H)$ is a morphism of pseudo $T$-algebras and
isomorphisms in $A$ are bijections of sets. Thus $\psi_{X,A}'(H)(x)
\neq H(*)$.

Let $J_0(\psi_{X,A}'(H)(x))$ be any set of the same cardinality as
$\psi_{X,A}'(H)(x)$ but not equal to $\psi_{X,A}'(H)(x)$. Let
$\alpha_0(\psi_{X,A}'(H)(x)):J_0(\psi_{X,A}'(H)(x)) \rightarrow
\psi_{X,A}'(H)(x)$ be a bijection. Let $J_0(a)=a$ for all $a \in Obj
\hspace{1mm} A$ such that $a \neq \psi_{X,A}'(H)(x)$. Then by
Theorem \ref{replacement} there exists a morphism $J:A \rightarrow
A$ of pseudo $T$-algebras which is the identity except on the object
$\psi_{X,A}'(H)(x)$. In particular $J(H(*))=H(*)$ because $H(*) \neq
\psi_{X,A}'(H)(x)$ from above.

The 2-naturality in the second variable implies that
\begin{equation} \label{strictoutput1}
\xymatrix@C=5pc@R=3pc{Mor_{\mathcal{A}}(FX, A) \ar[d]_{J_*} &
\ar[l]_{\psi_{X,A}'} \ar[d]^{(UJ)_*} Mor_{\mathcal{X}}(X,UA)
\\ Mor_{\mathcal{A}}(FX, A) & Mor_{\mathcal{X}}(X,UA) \ar[l]^{\psi_{X,A}'}
}
\end{equation}
commutes, \ie $J \circ \psi_{X,A}'(H)=\psi_{X,A}'(UJ \circ H)$. But
$UJ \circ H = H$ because $J(H(*))=H(*)$. Hence $J \circ
\psi_{X,A}'(H)= \psi_{X,A}'(H)$. Evaluating this on $x$ gives
$$J(\psi_{X,A}'(H)(x))=\psi_{X,A}'(H)(x)$$ which contradicts
$$J(\psi_{X,A}'(H)(x)) \neq \psi_{X,A}'(H)(x).$$

Thus there cannot exist such a $\psi'_{X,A}:Mor_{\mathcal{X}}(X,UA)
\rightarrow Mor_{\mathcal{A}}(FX,A)$ and the reason is that we allow
morphisms which are not strict.

Let $F':\mathcal{X} \rightarrow \mathcal{A}$ be any left biadjoint
for $U:\mathcal{A} \rightarrow \mathcal{X}$. Suppose it admits
equivalences of categories $\psi'_{X,A}:Mor_{\mathcal{X}}(X,UA)
\rightarrow Mor_{\mathcal{A}}(F'X,A)$ that are strictly 2-natural in
the second variable. Since $F$ and $F'$ are left biadjoints for $U$,
there exists for each $X$ a pseudo isomorphism $FX \rightarrow FX'$
by the biuniversal arrow argument in Lemma \ref{laxuniqueness} and
Theorem \ref{laxadjointuniqueness}. This pseudo isomorphism induces
an equivalence of categories $Mor_{\mathcal{A}}(F'X,A) \rightarrow
Mor_{\mathcal{A}}(FX,A)$ which is strictly 2-natural in $A$.
Composing this with $\psi_{X,A}'$ gives an equivalence of categories
$Mor_{\mathcal{X}}(X,UA) \rightarrow Mor_{\mathcal{A}}(FX,A)$ which
is strictly 2-natural in $A$, the second variable. But it was shown
above that such a 2-natural equivalence cannot exist. Hence we have
arrived at a contradiction and we conclude that $F'$ does not admit
equivalences $\psi'_{X,A}:Mor_{\mathcal{X}}(X,UA) \rightarrow
Mor_{\mathcal{A}}(F'X,A)$ that are strictly 2-natural in the second
variable.
\end{pf}

We can build on the previous example to show that there does not
exist a left 2-adjoint to the forgetful 2-functor in that situation.

\begin{examp}
Let $S$ be the trivial theory\index{theory!trivial theory} and let
$T$ be the theory of commutative monoids\index{theory!theory of
commutative monoids}\index{commutative monoid!theory of commutative
monoids}. Let $\mathcal{X}$ be the 2-category of pseudo $S$-algebras
and let $\mathcal{A}$ be the 2-category of pseudo $T$-algebras. Let
$U:\mathcal{A} \rightarrow \mathcal{X}$ be the forgetful
2-functor\index{forgetful 2-functor}\index{2-functor!forgetful
2-functor} associated to the trivial map of theories $S \rightarrow
T$. Then there does not exist a left
2-adjoint\index{2-adjoint}\index{2-adjoint!left 2-adjoint} to $U$,
\ie there does not exist a 2-functor $F':\mathcal{X} \rightarrow
\mathcal{A}$ which admits isomorphisms of categories
$\phi_{X,A}:Mor_{\mathcal{A}}(F'X,A) \rightarrow
Mor_{\mathcal{X}}(X,UA)$ that are strictly 2-natural in each
variable.
\end{examp}
\begin{pf}
Suppose such a $\phi$ existed. Let $\psi_{X,A}:=\phi_{X,A}^{-1}$.
Then $\psi_{X,A}$ is strictly 2-natural in the second variable $A$
and is an equivalence of categories. But this is impossible by the
previous example.
\end{pf}

\chapter{Weighted Bicolimits of Pseudo $T$-Algebras}
\label{sec:laxcolaxT} In this chapter we show that the 2-category of
pseudo $T$-algebras admits weighted\index{weighted} bicolimits. The
proof builds on the free pseudo\index{free pseudo $T$-algebra on a
pseudo $S$-algebra}\index{algebra!free pseudo $T$-algebra on a
pseudo $S$-algebra} $T$-algebra construction from Chapter
\ref{sec:forgetfulfunctor} as well as the construction of pseudo
colimits in the 2-category of small categories from Chapter
\ref{sec:laxcolimitsinCat}. The present construction of bicolimits
does not capture pseudo colimits\index{colimit!pseudo colimit}
because of the equivalence of morphism categories inherent to the
construction of the free pseudo $T$-algebra. This equivalence arises
because the morphisms of pseudo $T$-algebras are pseudo morphisms of
pseudo $T$-algebras rather than strict morphisms. After proving that
this 2-category admits bicolimits and bitensor
products\index{bitensor product}, we conclude that it admits
weighted bicolimits\index{bicolimit!weighted bicolimit}.

\begin{thm} \label{Talgebralaxcolimits}
The 2-category $\mathcal{C}$ of small pseudo $T$-algebras admits
bicolimits.\index{algebra!pseudo algebra!bicolimits of pseudo
algebras}\index{bicolimit|textbf}
\end{thm}
\begin{pf}
Let $\mathcal{J}$ be a small 1-category and $F:\mathcal{J}
\rightarrow \mathcal{C}$ a pseudo functor. In the following
construction we use notation similar to the construction of the
biuniversal\index{biuniversal arrow}\index{arrow!biuniversal arrow}
arrows for forgetful 2-functors in Chapter
\ref{sec:forgetfulfunctor}.

First we define candidates $W \in Obj \hspace{1mm}  \mathcal{C}$ and
$\pi:F \Rightarrow \Delta_W$. Let $T'$ denote the free
theory\index{free theory}\index{theory!free theory} on the sequence
of sets $T(0), T(1), \dots $ underlying the theory $T$. Let
$Alg'$\index{$Alg'$} be the category of small $T'$-algebras. Let
$Graph'$\index{$Graph'$} be the category of small directed
graphs\index{directed graph} whose object sets are discrete $T'$
algebras. Then there is a forgetful functor $Alg' \rightarrow
Graph'$ and it admits a left adjoint $V'$ by Freyd's Adjoint Functor
Theorem\index{Freyd's Adjoint Functor Theorem}.

Let $Obj \hspace{1mm}  R_{G'}$\index{$R_{G'}$|(} be the free
(discrete) $T'$ algebra on the set $\coprod_{j \in Obj \hspace{1mm}
\mathcal{J}}Obj \hspace{1mm}  Fj$. Let $Mor \hspace{1mm}
\hspace{1mm}  R_{G'}$ be the collection of the following arrows:

\begin{enumerate}
\item
For every $n \in \mathbf{N}$, for all words $w \in T(n)$, $w_1 \in
T(m_1), \dots, w_n \in T(m_n)$, and for all objects $A^1_1,\dots,
A^1_{m_1}$,$A^2_1, \dots, A^2_{m_2}, \dots, A^n_1, \dots, A^n_{m_n}
\in Obj \hspace{1mm}  R_{G'}$ there are arrows $$c_{w, w_1, \ldots,
w_n}(A^1_1,\dots,A^n_{m_n}):$$ $$\xymatrix@1{w \circ (w_1, \dots,
w_n)(A^1_1,\dots,A^n_{m_n}) \ar[r] & w(w_1(A^1_1,\dots, A^1_{m_1}),
\dots, w_n(A^n_1, \dots, A^n_{m_n}))}$$  $$c_{w, w_1, \ldots,
w_n}^{-1}(A^1_1,\dots,A^n_{m_n}):$$
$$\xymatrix@1{w(w_1(A^1_1,\dots, A^1_{m_1}), \dots, w_n(A^n_1,
\dots, A^n_{m_n})) \rightarrow w \circ (w_1, \dots,
w_n)(A^1_1,\dots,A^n_{m_n})}.$$ Here $w \circ (w_1, \ldots, w_n)$ is
the composition in the original theory $T$. The target
$w(w_1(A^1_1,\dots, A^1_{m_1}), \dots, w_n(A^n_1, \dots,
A^n_{m_n}))$ is the result of composing in the free
theory\index{theory!free theory}\index{free theory} and applying it
to the $A$'s in the free algebra.
\item
For every $A \in Obj \hspace{1mm}  R_{G'}$ there are arrows
$$\xymatrix@1{I_A:1(A) \ar[r] & A}$$ $$\xymatrix@1{I_A^{-1}:A \ar[r]
& 1(A)}.$$ Here $1$ is the unit of the original theory $T$.
\item
For every word $w \in T(m)$, for every function $f:\{1, \dots,m\}
\rightarrow \{1,\dots,n\}$, and for all objects $A_1, \dots, A_n \in
Obj \hspace{1mm}  R_{G'}$ there are arrows
$$\xymatrix@1{s_{w,f}(A_1, \ldots, A_n):w_f(A_1, \dots, A_n) \ar[r]
& w(A_{f1}, \dots ,A_{fm})}$$ $$\xymatrix@1{s_{w,f}^{-1}(A_1,
\ldots, A_n):w(A_{f1}, \dots ,A_{fm}) \ar[r] & w_f(A_1, \dots,
A_n)}.$$ The substituted word $w_f$ is the substituted word in the
original theory $T$. The target $w(A_{f1}, \dots ,A_{fm})$ is the
result of substituting in $w$ in the free theory\index{theory!free
theory}\index{free theory} and then evaluating on the $A$'s.
\item
For every word $w \in T(n)$, $j \in Obj \hspace{1mm} \mathcal{J}$,
and objects $A_1, \dots, A_n$ of $Fj$ there are arrows
$$\xymatrix@1{\rho^{\pi_j}_w(A_1, \dots, A_n):\Phi_j(w)(A_1, \dots, A_n)
\ar[r] & w(A_1, \dots, A_n)}$$
$$\xymatrix@1{(\rho^{\pi_j}_w)^{-1}(A_1, \dots, A_n):w(A_1, \dots,
A_n) \ar[r] & \Phi_j(w)(A_1, \dots, A_n),}$$ where $\Phi_j$ denotes
the structure maps of the pseudo $T$-algebra $Fj$.
\item
Include all elements of $\coprod_{j \in \mathcal{J}} Mor
\hspace{1mm}  Fj$ in $Mor \hspace{1mm} R_{G'}$.
\item
For every morphism $f:i \rightarrow j$ of $\mathcal{J}$ and every $x
\in Obj \hspace{1mm}  Fi$ we include arrows
$$\xymatrix@1{h_{(x,f)}:x \ar[r] & a_f(x)}$$ $$\xymatrix@1{h_{(x,f)}^{-1}:a_f(x) \ar[r] & x}$$
as in the proof of Theorem \ref{catlaxcolimits}, where $a_f=Ff:Fi
\rightarrow Fj$.
\end{enumerate}

With these arrows, $R_{G'}$ is an object of
$Graph'$\index{$R_{G'}$|)}\index{$Graph'$}. Now we apply the functor
$V'$ to the directed graph\index{directed graph} $R_{G'}$ to get a
category $R'$ which is a $T'$-algebra.

Let $K$ be the smallest congruence\index{congruence} on the category
$R'$ with the following properties:

\begin{enumerate}
\item
All of the relations necessary to make the coherence arrows
(including $\rho^{\pi_j}_w$) into natural transformations belong to
$K$. For example, if $A,B \in Obj \hspace{1mm}  R'$ and $f:A
\rightarrow B$ is a morphism of $R'$,then the relation $I_A \circ f
= 1(f) \circ I_B$ belongs to $K$.
\item
All of the relations necessary to make the coherence arrows
(including $\rho^{\pi_j}_w$) into isos are in $K$. For example, for
every $A \in Obj \hspace{1mm}  R'$ the relations $I_A \circ I_A^{-1}
= 1_A$ and $I_A^{-1} \circ I_A = 1_A$ are in $K$.
\item
All of the relations for pseudo algebras\index{algebra!pseudo
algebra} listed in Definition \ref{laxalgebradefinition} belong to
$K$, where the objects range over the objects of $R'$.
\item
The original composition relations in each of the categories $Fj$
belong to $K$ for all $j \in Obj \hspace{1mm}  \mathcal{J}$.
\item
The coherence diagrams necessary to make the inclusion $\pi_j:Fj
\rightarrow R'$ into a morphism of pseudo $T$-algebras belong to
$K$. These diagrams are listed in Definition
\ref{defnlaxalgebramorphism}. Note that these coherence diagrams
will involve the arrows $\rho_w^{\pi_j}(A_1, \dots, A_n)$ for $w \in
T(n)$.
\item
All of the relations in the proof of Theorem \ref{catlaxcolimits}
are in $K$.
\item
If the relations $f_1=g_1, \dots ,f_n = g_n$  are in $K$ and $w \in
T'(n)$, then the relation $w(f_1, \dots, f_n) = w(g_1, \dots g_n)$
is also in $K$.
\end{enumerate}
Next we mod out by the congruence\index{congruence} $K$ in $R'$ and
we get a pseudo $T$-algebra $R=:W \in Obj \hspace{1mm} \mathcal{C}$.

We define a pseudo natural transformation $\pi:F \Rightarrow
\Delta_W$ as follows. For $j \in Obj \hspace{1mm}  \mathcal{J}$,
define $\pi_j:Fj \rightarrow W$ to be the inclusion functor. The
functor $\pi_j$ is a morphism of pseudo $T$-algebras because of the
relations we modded out by. Define $\tau_{i,j}(f)_x:\pi_i(x)
\rightarrow \pi_j \circ a_f(x)$ by $\tau_{i,j}(f)_x := h_{(x,f)}$ as
in the proof of Theorem of \ref{catlaxcolimits}. Then $x \mapsto
\tau_{i,j}(f)_x$ is a 2-cell $\pi_i \Rightarrow \pi_j \circ a_f$ in
the 2-category of pseudo $T$-algebras because of the relations we
modded out by and because of the work in the proof of Theorem
\ref{catlaxcolimits}. By an argument similar to Lemma
\ref{pilaxnatural} we conclude that $\pi:F \Rightarrow \Delta_W$ is
a pseudo natural transformation. The candidate for the
bicolimit\index{bicolimit} of $F$ is $W \in Obj \hspace{1mm}
\mathcal{C}$ with the pseudo cone $\pi:F \Rightarrow \Delta_W$. This
concludes the definition of the candidate for the bicolimit of $F$.

Let $V \in Obj \hspace{1mm}  \mathcal{C}$. Define the functor
$\phi:Mor_{\mathcal{C}}(W,V) \rightarrow PseudoCone(F,V)$ by $b
\mapsto b \circ \pi$ as before. We need to see that $\phi$ is an
equivalence of categories.

\begin{lem}
There is a functor $\psi:PseudoCone(F,V) \rightarrow
Mor_{\mathcal{C}}(W,V)$.
\end{lem}
\begin{pf}
First we define $\psi$ on objects. Let $\pi':F \Rightarrow \Delta_V$
be a pseudo natural transformation which is natural up to the
coherence iso 2-cells $\tau'$. From $\pi'$ we get a map of sets
$$\coprod_{j \in Obj \hspace{1mm}  \mathcal{J}} Obj \hspace{1mm}  Fj
\rightarrow Obj \hspace{1mm}  V$$ which induces a map $$d:Obj
\hspace{1mm}  R_{G'} \rightarrow Obj \hspace{1mm}  V$$ of discrete
$T'$ algebras. Define $d$ on arrows of $R_{G'}$ as follows:
\begin{itemize}
\item
$dg:=\pi_j'g$ for all $g \in Mor \hspace{1mm}  Fj$ and all $j \in
Obj \hspace{1mm} \mathcal{J}$
\item
$dh_{(x,f)}:=\tau_{i,j}'(f)_x$ and
$dh_{(x,f)}^{-1}:=(\tau_{i,j}'(f)_x)^{-1}$ for $f:i \rightarrow j$
in $\mathcal{J}$ and $x \in Obj \hspace{1mm}  Fi$
\item
$d$ takes a coherence arrow in $R_{G'}$ to the analogous coherence
iso in $V$
\item
$d(\rho^{\pi_j}_w):=\rho^{\pi_j'}_w$ where $\rho_w^{\pi_j'}$ is the
coherence iso of the morphism $\pi_j':Fj \rightarrow V$ of pseudo
$T$-algebras, and similarly
$d((\rho^{\pi_j}_w)^{-1}):=(\rho^{\pi_j'}_w)^{-1}$.
\end{itemize}
This defines a morphism $d:R_{G'} \rightarrow V$ of the category
$Graph'$, where part of the structure of the $T'$-algebra $V$ is
forgotten. The adjoint $Graph' \rightarrow Alg'$ to the forgetful
functor $Alg' \rightarrow Graph'$ gives us a morphism $R'
\rightarrow V$, which we also denote by $d$. Furthermore,
$d:R'\rightarrow V$ preserves the relations in $K$. Hence $d$
induces a map $b:R \rightarrow V$ on the quotient and $d$ is a
morphism of pseudo $T$-algebras. Note that the coherence isos of $b$
are trivial. This is how we define $\psi$ on objects:
$\psi(\pi'):=b$.

Let $\sigma, \sigma' \in Obj \hspace{1mm}  PseudoCone(F,V)$ and let
$\Xi:\sigma \rightsquigarrow \sigma'$ be a morphism in the category
$PseudoCone(F,V)$. Then define a 2-cell $\psi(\Xi):\psi(\sigma)
\Rightarrow \psi(\sigma')$ by $\psi(\Xi)_x := \Xi_j(x)$ for $x \in
Obj \hspace{1mm}  Fj$ and continue the definition inductively by
$$\psi(\Xi)_{w(x_1, \dots, x_n)} := \Psi(w)(\psi(\Xi)_{x_1}, \dots,
\psi(\Xi)_{x_n}),$$ where $\Psi$ denotes the structure maps of the
pseudo $T$-algebra $V$. Another inductive argument shows that this
assignment preserves compositions and identities.
\end{pf}
\begin{lem}
The functor $\phi \circ \psi:PseudoCone(F,V) \rightarrow
PseudoCone(F,V)$ is the identity functor.
\end{lem}
\begin{pf}
This is similar to Lemma \ref{compositea}. The only difference here
is that we must prove that the coherence isos for the morphism
$\pi_j':Fj \rightarrow V$ of pseudo $T$-algebras are the same as the
coherence isos for $(\phi \circ \psi(\pi'))_j$. But this is true
because the coherence isos of $\psi(\pi')$ are trivial.
\end{pf}

\begin{lem}
The composite functor $\psi \circ \phi: Mor_{\mathcal{C}}(W,V)
\rightarrow Mor_{\mathcal{C}}(W,V)$ is naturally isomorphic to the
identity functor.
\end{lem}
\begin{pf}
We construct a natural isomorphism $\eta:1_{Mor_{\mathcal{C}}(W,V)}
\Rightarrow \psi \circ \phi$. Let $b \in Obj \hspace{1mm}
Mor_{\mathcal{C}}(W,V)$. We define $\eta_b=:\alpha$ inductively. For
all $j \in Obj \hspace{1mm} \mathcal{J}$ and all $x \in Obj
\hspace{1mm} Fj \subseteq Obj \hspace{1mm}  W$ we have $\psi \circ
\phi(b)(x)=b(x)$. Define
$$\alpha_x:b(x) \rightarrow \psi \circ \phi (b)(x)$$ to be the
identity for such $x$. For $w \in T(n)$ and $x_1, \dots, x_n \in
\coprod_{j \in Obj \hspace{1mm}  \mathcal{J}} Obj \hspace{1mm} Fj$
define
$$\alpha_{w(x_1, \dots, x_n)}:=\rho^b_w(x_1, \dots, x_n).$$
Now let $x_1, \dots, x_n \in Obj \hspace{1mm}  W$ and $w \in T(n)$.
Suppose $\alpha_{x_1}, \dots, \alpha_{x_n}$ are already defined.
Then define
$$\alpha_{w(x_1, \dots, x_n)}:b(w(x_1, \dots, x_n)) \rightarrow
\psi \circ \phi(b)(w(x_1, \dots, x_n))$$ to be the composition
$$\xymatrix@R=3pc@C=3pc{b(w(x_1, \dots, x_n)) \ar[d]^{\rho_w^b(x_1, \dots, x_n)}
\\ \Psi(w)(bx_1, \dots, bx_n) \ar[d]^{\Psi(w)(\alpha_{x_1}, \dots,
\alpha_{x_n})} \\ \Psi(w)(\psi \circ \phi(b)x_1, \dots, \psi \circ
\phi(b)x_n).}$$

Then the assignment $x \mapsto \alpha_x$ is a 2-cell in the category
of pseudo $T$-algebras because it is natural and commutes with the
coherence isos of $b$ and $\psi \circ \phi (b)$ by an inductive
argument (recall the coherence isos of $\psi \circ \phi (b)$ are
trivial). An inductive argument also shows that $b \mapsto \eta_b$
is natural.
\end{pf}

\begin{lem}
The functor $\phi:Mor_{\mathcal{C}}(W,V) \rightarrow
PseudoCone(F,V)$ defined by $b \mapsto b \circ \pi$ is an
equivalence of categories.
\end{lem}
\begin{pf}
This follows immediately from the previous two lemmas.
\end{pf}

\begin{lem}
The object $W \in Obj \hspace{1mm}  \mathcal{C}$ and the pseudo cone
$\pi:F \Rightarrow \Delta_W$ comprise a bicolimit of $F$.
\end{lem}
\begin{pf}
This follows immediately from the previous lemma.
\end{pf}

This completes the proof that the 2-category of small pseudo
$T$-algebras admits bicolimits.
\end{pf}

\begin{lem}
The 2-category $\mathcal{C}$ of pseudo $T$-algebras admits bitensor
products\index{bitensor product|textbf}.
\end{lem}
\begin{pf}
Let $J$ be a category and $F$ a pseudo $T$-algebra. First we define
an object $R_{G'}$ of $Graph'$. Let $Obj \hspace{1mm} R_{G'}$ be the
free discrete $T'$-algebra on the set $Obj \hspace{1mm} J \times Obj
\hspace{1mm}  F$, where $T'$ is the free theory on $T$. Let $Mor
\hspace{1mm}  R_{G'}$ be the collection of the following arrows.

\begin{enumerate}
\item
For every $n \in \mathbf{N}$, for all words $w \in T(n)$, $w_1 \in
T(m_1), \dots, w_n \in T(m_n)$, and for all objects $A^1_1,\dots,
A^1_{m_1}$,$A^2_1, \dots, A^2_{m_2}, \dots, A^n_1, \dots, A^n_{m_n}
\in Obj \hspace{1mm}  R_{G'}$ there are arrows $$c_{w, w_1, \ldots,
w_n}(A^1_1,\dots,A^n_{m_n}):$$ $$\xymatrix@1{w \circ (w_1, \dots,
w_n)(A^1_1,\dots,A^n_{m_n}) \ar[r] & w(w_1(A^1_1,\dots, A^1_{m_1}),
\dots, w_n(A^n_1, \dots, A^n_{m_n}))}$$  $$c_{w, w_1, \ldots,
w_n}^{-1}(A^1_1,\dots,A^n_{m_n}):$$
$$\xymatrix@1{w(w_1(A^1_1,\dots, A^1_{m_1}), \dots, w_n(A^n_1,
\dots, A^n_{m_n})) \ar[r] & w \circ (w_1, \dots,
w_n)(A^1_1,\dots,A^n_{m_n})}.$$ Here $w \circ (w_1, \ldots, w_n)$ is
the composition in the original theory $T$. The target
$w(w_1(A^1_1,\dots, A^1_{m_1}), \dots, w_n(A^n_1, \dots,
A^n_{m_n}))$ is the result of composing in the free
theory\index{theory!free theory}\index{free theory} and applying it
to the $A$'s in the free algebra.
\item
For every $A \in Obj \hspace{1mm}  R_{G'}$ there are arrows
$$\xymatrix@1{I_A:1(A) \ar[r] & A}$$ $$\xymatrix@1{I_A^{-1}:A
\ar[r] & 1(A)}.$$ Here $1$ is the unit of the original theory $T$.
\item
For every word $w \in T(m)$, for every function $f:\{1, \dots,m\}
\rightarrow \{1,\dots,n\}$, and for all objects $A_1, \dots, A_n \in
Obj \hspace{1mm}  R_{G'}$ there are arrows
$$\xymatrix@1{s_{w,f}(A_1, \ldots, A_n):w_f(A_1, \dots, A_n) \ar[r] &
w(A_{f1}, \dots ,A_{fm})}$$ $$\xymatrix@1{s_{w,f}^{-1}(A_1, \ldots,
A_n):w(A_{f1}, \dots ,A_{fm}) \ar[r] & w_f(A_1, \dots, A_n)}.$$ The
substituted word $w_f$ is the substituted word in the original
theory $T$. The target $w(A_{f1}, \dots ,A_{fm})$ is the result of
substituting in $w$ in the free theory\index{theory!free
theory}\index{free theory} and then evaluating on the $A$'s.
\item
For every word $w \in T(n)$, $j \in Obj \hspace{1mm}  J$, and
objects $x_1, \dots, x_n$ of $F$ there are arrows
$$\xymatrix@1{\rho^{\pi(j)}_w((j,x_1), \dots, (j,x_n)):(j,\Phi(w)(x_1, \dots,
x_n)) \ar[r] & w((j,x_1), \dots, (j,x_n))}$$
$$\xymatrix@1{(\rho^{\pi(j)}_w)^{-1}:w((j,x_1), \dots, (j,x_n))
\ar[r] & (j,\Phi(w)(x_1, \dots, x_n))},$$ where $\Phi$ denotes
structure maps of the pseudo $T$-algebra $F$.
\item
Include all elements of $Mor \hspace{1mm}  J \times Mor \hspace{1mm}
F$ in $Mor \hspace{1mm} R_{G'}$.
\end{enumerate}

With these arrows, $R_{G'}$ is an object of $Graph'$. Now we apply
the free $T'$-algebra functor to the directed graph\index{directed
graph} $R_{G'}$ to get a category $R'$ which is a $T'$ algebra. Let
$K$ be the smallest congruence\index{congruence} on the category
$R'$ with the following properties:

\begin{enumerate}
\item
All of the relations necessary to make the coherence arrows
(including $\rho^{\pi(j)}_w$) into natural transformations belong to
$K$. For example, if $A,B \in Obj \hspace{1mm}  R'$ and $f:A
\rightarrow B$ is a morphism in $R'$,then the relation $I_A \circ f
= 1(f) \circ I_B$ belongs to $K$.
\item
All of the relations necessary to make the coherence arrows
(including $\rho^{\pi(j)}_w$) into isos are in $K$. For example, for
every $A \in Obj \hspace{1mm}  R'$ the relations $I_A \circ I_A^{-1}
= 1_A$ and $I_A^{-1} \circ I_A = 1_A$ are in $K$.
\item
All of the relations for pseudo algebras\index{algebra!pseudo
algebra} listed in Definition \ref{laxalgebradefinition} belong to
$K$, where the objects range over the objects of $R'$.
\item
The original composition relations in the category $J \times F$
belong to $K$.
\item
For each $j \in J$, the coherence diagrams necessary to make the
inclusion $F \rightarrow R'$, $x \mapsto (j,x)$ into a morphism of
pseudo $T$-algebras belong to $K$. These diagrams are listed in
Definition \ref{defnlaxalgebramorphism}. Note that these coherences
will involve the arrows $$\rho_w^{\pi(j)}((j,x_1), \dots,
(j,x_n)):(j,\Phi(w)(x_1, \dots, x_n)) \rightarrow w((j,x_1), \dots,
(j,x_n)).$$
\item
For any $g:j_1 \rightarrow j_2$ in $J$ and $x_1, \dots, x_n$ in $F$
we include the relation
$$\xymatrix@C=8pc@R=4pc{(j_1,\Phi(w)(x_1, \dots, x_n))
\ar[r]^{(g,1_{\Phi(w)(x_1, \dots, x_n)})}
\ar[d]_{\rho^{\pi(j_1)}_w(x_1, \dots, x_n)} & (j_2, \Phi(w)(x_1,
\dots, x_n)) \ar[d]^{\rho^{\pi(j_2)}_w(x_1, \dots, x_n)}
\\ w((j_1,x_1), \dots, (j_1,x_n)) \ar[r]_{w((g,x_1), \dots, (g,x_n))} &
w((j_2,x_1), \dots, (j_2,x_n)).}$$
\item
If the relations $f_1=g_1, \dots ,f_n = g_n$  are in $K$ and $w \in
T'(n)$, then the relation $w(f_1, \dots, f_n) = w(g_1, \dots g_n)$
is also in $K$.
\end{enumerate}

Next we mod out by the congruence\index{congruence} $K$ in $R'$ and
we get a pseudo $T$-algebra $J*F \in Obj \hspace{1mm}  \mathcal{C}$.
We define a functor $\pi:J \rightarrow \mathcal{C}(F,J*F)$ by
$$\pi(j)(x):=(j,x)$$
$$\pi(j)(f):=(1_j,f)$$
$$(\pi(g))_x:=(g,1_x)$$ for $j\in Obj \hspace{1mm}  J,x \in Obj \hspace{1mm}  F, f \in Mor \hspace{1mm}  F,$
and $g \in Mor \hspace{1mm}  J$. Then $\pi(j):F \rightarrow J*F$ is
a morphism of pseudo $T$-algebras with coherence isos
$\rho^{\pi(j)}$ and $\pi(g):\pi(j_1) \Rightarrow \pi(j_2)$ is a
2-cell in the 2-category of pseudo $T$-algebras because of the
relations. The relations also imply that $\pi$ is a functor.

We claim that $\pi$ induces an equivalence
$$\xymatrix@C=3pc@1{\mathcal{C}(J*F,C) \ar[r]^-{\phi} &
Cat(J,\mathcal{C}(F,C))}$$ $$b \mapsto \mathcal{C}(F,b) \circ \pi$$
$$ \alpha \mapsto \mathcal{C}(F,\alpha) * i_{\pi}$$ of categories.
Define a functor $\psi:Cat(J,\mathcal{C}(F,C)) \rightarrow
\mathcal{C}(J*F,C)$ as follows. For a functor $\sigma:J \rightarrow
\mathcal{C}(F,C)$, we have a map of sets $$Obj \hspace{1mm}  J
\times Obj \hspace{1mm}  F \rightarrow Obj \hspace{1mm}  C$$
$$(j,x) \mapsto \sigma(j)(x)$$ which induces a map $\psi(\sigma): Obj \hspace{1mm}  R_{G'}
\rightarrow Obj \hspace{1mm}  C$ of discrete $T'$-algebras
satisfying
$$\psi(\sigma)(j,x):=\sigma(j)(x)$$
$$\psi(\sigma)(w((j_1,x_1), \dots,
(j_n,x_n))):=\Phi^C(w)(\sigma(j_1)(x_1), \dots, \sigma(j_n)(x_n))$$
for $(j,x),(j_1,x_1), \dots, (j_n,x_n) \in J \times F$. Define
$\psi(\sigma)$ on arrows of $R_{G'}$ by
$$\psi(\sigma)(c_{w,w_1, \dots, w_n}(A_1^1, \dots, A^n_{m_n}))
:=c_{w,w_1, \dots, w_n}(\psi(\sigma)(A_1^1), \dots,
\psi(\sigma)(A^n_{m_n}))$$
$$\psi(\sigma)(I_A):=I_{\psi(\sigma)(I_A)}$$
$$\psi(\sigma)(s_{w,f}(A_1, \dots,
A_n)):=s_{w,f}(\psi(\sigma)(A_1), \dots, \psi(\sigma)(A_n))$$
$$\psi(\sigma)(g,f):=\sigma(j_2)(f) \circ \sigma(g)_{x_1} =
\sigma(g)_{x_2} \circ \sigma(j_1)(f)$$ for $A^k_{\ell},A,A_i \in Obj
\hspace{1mm}  R_{G'}$, $f:m \rightarrow n$, $g:j_1 \rightarrow j_2$
in $J$, and $f:x_1 \rightarrow x_2$ in $F$. We define $\psi(\sigma)$
similarly for $c_{w,w_1, \dots, w_n}^{-1}, I_A^{-1}, s_{w,f}^{-1}$.
Then $\psi(\sigma):R_{G'} \rightarrow C$ is a morphism in $Graph'$,
which induces a morphism $R' \rightarrow C$ in $Alg'$. It preserves
the relations and therefore induces a morphism $\psi(\sigma):J*F
\rightarrow C$ of pseudo $T$-algebras on the quotient. This is
actually a strict morphism of pseudo $T$-algebras. For a natural
transformation $\Xi: \sigma \Rightarrow \sigma'$ define a 2-cell
$\psi(\Xi):\psi(\sigma) \Rightarrow \psi(\sigma')$ inductively by
$$\psi(\Xi)_{(j,x)}:=(\Xi_j)_x$$
for $(j,x) \in Obj \hspace{1mm}  J \times Obj \hspace{1mm}  F$ and
$$\psi(\Xi)_{w(A_1,\dots,A_n)}:=\Phi^C(w)(\psi(\Xi)_{A_1}, \dots,
\psi(\Xi)_{A_n})$$ whenever $\psi(\Xi)_{A_1}, \dots,
\psi(\Xi)_{A_n}$ are already defined. From these definitions we can
conclude that $\psi$ is a functor and $\phi \circ \psi = 1_{Cat(J,
\mathcal{C}(F,C))}$. For example,
$$\aligned (\phi \circ \psi (\sigma))(j)(x) &= (\psi(\sigma)
\circ \pi(j))(x) \\
&=\psi(\sigma)(j,x) \\ &= \sigma(j)(x)\endaligned$$ and also
$$\aligned ((\phi \circ \psi (\Xi))_{j})_x &= ((\psi(\Xi)
* i_{\pi})_j)_x \\ &= \psi(\Xi)_{\pi(j) (x)} \\
&=(\Xi_j)_x. \endaligned$$

We construct a natural isomorphism $\eta:1_{\mathcal{C}(J*F,C)}
\Rightarrow \psi \circ \phi$. Let $b:J*F \rightarrow C$ be a
morphism of pseudo $T$-algebras. We define $\eta_b=:\alpha$
inductively. For all $(j,x) \in Obj \hspace{1mm}  J \times Obj
\hspace{1mm}  F$ we have
$$\aligned \psi \circ \phi(b)(j,x)&= \psi(\mathcal{C}(F,b)
\circ \pi)(j,x) \\ &=(\mathcal{C}(F,b) \circ \pi)(j)(x) \\
&=(b \circ \pi(j))(x) \\ &= b(j,x).
\endaligned$$ Define
$$\alpha_{(j,x)}:b(j,x) \rightarrow \psi \circ \phi (b)(j,x)$$ to be the
identity for such $(j,x)$. For $w \in T(n)$ and $(j_1,x_1), \dots,
(j_n,x_n) \in Obj \hspace{1mm}  J \times Obj \hspace{1mm}  F$ define
$$\alpha_{w((j_1,x_1), \dots, (j_n,x_n))}:=\rho^b_w((j_1,x_1), \dots, (j_n,x_n)).$$
For $A_1, \dots, A_n \in Obj \hspace{1mm}  R_{G'}=Obj \hspace{1mm}
J*F$ and $w \in T(n)$, define
$$\alpha_{w(A_1, \dots, A_n)}:b(w(A_1, \dots, A_n)) \rightarrow
\psi \circ \phi(b)(w(A_1, \dots, A_n))$$ to be the composition
$$\xymatrix@R=3pc{b(w(A_1, \dots, A_n)) \ar[d]^{\rho_w^b(A_1, \dots, A_n)}
\\ \Psi(w)(bA_1, \dots, bA_n) \ar[d]^-{\Psi(w)(\alpha_{A_1},
\dots, \alpha_{A_n})} \\
\Psi(w)(\psi \circ \phi(b)A_1, \dots, \psi \circ \phi(b)A_n). }
$$

Then the assignment $x \mapsto \alpha_x$ is a 2-cell in the category
of pseudo $T$-algebras because it is natural and commutes with the
coherence isos of $b$ and $\psi \circ \phi (b)$ by an inductive
argument (recall the coherence isos of $\psi \circ \phi (b)$ are
trivial). An inductive argument also shows that $b \mapsto \eta_b$
is natural.

By Remark \ref{laxtensorunit}, this implies that $J*F$ is a bitensor
product of $J$ and $F$\index{bitensor product}.
\end{pf}

\begin{thm}
The 2-category $\mathcal{C}$ of pseudo $T$-algebras admits weighted
bicolimits.\index{algebra!pseudo algebra!bicolimits of pseudo
algebras}\index{bicolimit!weighted
bicolimit|textbf}\index{weighted|textbf}\index{bicolimit}
\end{thm}
\begin{pf}
The 2-category $\mathcal{C}$ admits bicoproducts\index{bicoproduct}
and bicoequalizers\index{bicoequalizer} by Theorem
\ref{Talgebralaxcolimits}. It admits bitensor
products\index{bitensor product} by the previous lemma. Hence by
Theorem \ref{streetlaxcolimit} it admits weighted bicolimits.
\end{pf}

\chapter{Stacks}
\label{sec:stacks} In this chapter we introduce the language of
stacks\index{stack} in analogy to sheaves, since stacks generalize
sheaves. A stack\index{stack} is a contravariant pseudo
functor\index{contravariant}\index{functor!pseudo functor} from a
Grothendieck topology to a 2-category which takes Grothendieck
covers to bilimits in the sense described below. The target
2-category is required to admit bilimits. We have shown that the
2-category of pseudo algebras over a theory admits
bilimits\index{algebra!pseudo algebra!bilimits of pseudo algebras},
so we can speak of stacks of pseudo algebras\index{algebra!pseudo
algebra!stacks of pseudo algebras}. Some references for stacks are
\cite{breen}, \cite{fantechi1}, \cite{giraud1}, \cite{moerdijk1},
and \cite{vistoli1}. We are interested in stacks because we want to
capture the algebraic structure of
holomorphic\index{holomorphic}\index{rigged surface!holomorphic
families of rigged surfaces}\index{holomorphic families of rigged
surfaces} families of rigged surfaces\index{rigged surface} as in
Section \ref{sec:rigged}.

\begin{defn}
A {\it basis for a Grothendieck topology}\index{Grothendieck
topology|textbf}\index{basis for a Grothendieck topology|textbf} on
a category $\mathcal{B}$ with pullbacks\index{pullback} is a
function $K$ which assigns to each object $B$ of $\mathcal{B}$
  a collection of families of morphisms with codomain $B$ such
  that:
\begin{enumerate}
\item
If $g:B' \rightarrow B$ is an isomorphism, then $\{g\} \in K(B)$.
\item
If $\{g_i:B_i \rightarrow B | i \in I \} \in K(B)$, then for any
morphism $g:D \rightarrow B$ the family $\{\pi_i^2: B_i \times_B D
\rightarrow D | i \in I \}$ of pullbacks of the $g_i$ along $g$ is
in $K(D)$.
\item
If $\{g_i:B_i \rightarrow B | i \in I \} \in K(B)$ and
$\{f_{ij}:D_{ij} \rightarrow B_i | j \in J_i \} \in K(B_i)$ for all
$i$, then the composite family $\{g_i \circ f_{ij}:B_{ij}
\rightarrow B | i \in I, j \in J_i\}$ is in $K(B)$.
\end{enumerate}
\end{defn}
The second axiom is called the {\it stability axiom}\index{stability
axiom|textbf} because it says that $K$ is stable under
pullbacks\index{pullback}. The third axiom is called the {\it
transitivity axiom}\index{transitivity axiom|textbf}. Often we refer
to the basis as well as the category $\mathcal{B}$ as a Grothendieck
topology\index{Grothendieck topology|textbf}. We follow this
convention. Some authors call a Grothendieck topology a Grothendieck
site\index{Grothendieck site|textbf}. The elements of $K(B)$ are
called {\it Grothendieck covers}\index{Grothendieck cover|textbf}.

\begin{defn}
Let $\mathcal{B}$ be a Grothendieck topology and $\mathcal{C}$ a
concrete category. Then a {\it
$\mathcal{C}$-sheaf}\index{sheaf|textbf} on $\mathcal{B}$ is a
contravariant functor $G: \mathcal{B} \rightarrow \mathcal{C}$ which
takes Grothendieck covers to limits, \ie for any object $B$ of
$\mathcal{B}$ and for any Grothendieck cover $\{g_i:B_i \rightarrow
B | i \in I \} \in K(B)$ the following diagram is an equalizer,
\begin{equation} \label{sheafequalizer}
\xymatrix@C=3pc{G(B) \ar[r]^-e & \prod_{i \in I} G(B_i)
 \ar@<.5ex>[r]^-{p_1} \ar@<-.5ex>[r]_-{p_2} & \prod_{i,j \in I}
 G(B_i \times_B B_j)}
\end{equation}
where $e(a)=\{G(g_i)a\}_{i \in I}$ and $p_1(\{a_k\}_{k
 \in I})_{ij} = G(\pi_{ij}^1)a_i$ and $p_2(\{a_k\}_{k \in I})_{ij} =
  G(\pi_{ij}^2)a_j$. Here $\pi_{ij}^1,\pi_{ij}^2$ are the morphisms in
 the pullback\index{pullback} diagrams for $B_{ij}:=B_i \times_B B_j$.
 $$\xymatrix@R=3pc@C=3pc{B_i \times_B B_j \ar[r]^-{\pi_{ij}^1} \ar[d]_{\pi_{ij}^2} &
 B_i \ar[d]
 \\ B_j \ar[r] & B}$$

\end{defn}

See \cite{maclane2} for a thorough discussion of Grothendieck
topologies and sheaves. Diagram (\ref{sheafequalizer}) is an
equalizer\index{equalizer} if and only if it is {\it
exact}\index{exact}. Usually we speak of a $\mathcal{C}$-sheaf as a
sheaf of objects of $\mathcal{C}$. For example, if $\mathcal{C}$ is
the category of sets, then we speak of a sheaf of sets.  Next we
speak of stacks\index{stack} of categories and then generalize to
stacks of objects with algebraic structure.

Let $Cat$\index{$Cat$} denote the 2-category of small categories.
Suppose $\mathcal{B}$ is a Grothendieck topology. Let $G:
\mathcal{B} \rightarrow Cat$ be a contravariant\index{contravariant}
pseudo functor. Let $B$ be an object of $\mathcal{B}$ and $\{g_i:B_i
\rightarrow B | i \in I \} \in K(B)$ a Grothendieck
cover\index{Grothendieck cover}. Consider the diagram
\begin{equation} \label{stackdefn}
\xymatrix@C=3pc{\prod_{i \in I} G(B_i) \ar@<.5ex>[r]^-{p_1}
\ar@<-.5ex>[r]_-{p_2} & \prod_{i,j \in I} G(B_i \times_B B_j)
\ar[r]|-{p_{13}}  \ar@<1.2ex>[r]^-{p_{12}} \ar@<-1.2ex>[r]_-{p_{23}}
& \prod_{i,j,k \in I} G(B_i \times_B B_j \times_B B_k) }
\end{equation}
 where the arrows are defined as
$$p_1(\{a_k\}_{k})_{ij} : = G(\pi_{ij}^1)a_i$$
$$p_2(\{a_k\}_{k})_{ij} : = G(\pi_{ij}^2)a_j$$
$$p_{12}(\{a_{\ell m}\}_{\ell m})_{ijk} : = G(\pi^{12}_{ijk})a_{ij}$$
$$p_{13}(\{a_{\ell m}\}_{\ell m})_{ijk}: = G(\pi^{13}_{ijk})a_{ik}$$
$$p_{23}(\{a_{\ell m}\}_{\ell m})_{ijk}: = G(\pi^{23}_{ijk})a_{jk}.$$
Here $\pi^{12}_{ijk},\pi^{13}_{ijk},\pi^{23}_{ijk}$ are the
morphisms for the triple fiber product $B_i \times_B B_j \times_B
B_k$ as in the following commutative diagram from \cite{vistoli1}.
The unlabelled arrows are $g_i, g_j,$ and $g_k$ from the
Grothendieck cover.
$$\xymatrix@R=3pc@C=3pc{& B_{ijk} \ar[dl]_{\pi_{ijk}^{12}}
\ar[dd]|{}^(.3){\pi^{13}_{ijk}} \ar[rr]^(.3){\pi_{ijk}^{23}} & &
B_{jk} \ar[dl]_{\pi^1_{jk}} \ar[dd]^{\pi^2_{jk}} \\ B_{ij}
\ar[dd]_{\pi_{ij}^1} \ar[rr]_(.3){\pi_{ij}^2} &  & B_j \ar[dd] &
\\ & B_{ik} \ar[rr]|{}^(.3){\pi_{ik}^2} \ar[ld]_(.3){\pi^1_{ik}} & & B_k
\ar[ld] \\ B_i \ar[rr] & & B & }$$ Every face in this diagram is a
pullback\index{pullback} square. The object $B_{ijk}$ is the limit
of the diagram obtained from this one by deleting $B_{ijk}$ and the
arrows emanating from it.

Diagram (\ref{stackdefn}) can be interpreted as the image of a
pseudo functor $F:\mathcal{J} \rightarrow Cat$ as follows. Let
$\mathcal{J}$ be the free 1-category on the directed
graph\index{directed graph}
\begin{equation} \label{sourcediagram}
\xymatrix@C=3pc@R=3pc{X \ar@<.5ex>[r]^-{f_1} \ar@<-.5ex>[r]_-{f_2} &
Y \ar[r]|-{f_{13}}
  \ar@<1.2ex>[r]^-{f_{12}}
\ar@<-1.2ex>[r]_-{f_{23}} & Z }
\end{equation}
modded out by the relations below.
$$
\xymatrix@R=3pc@C=3pc{
   X \ar[r]^{f_1} \ar[d]_{f_1} & Y \ar[d]^{f_{12}} &   &
 & X \ar[r]^{f_2} \ar[d]_{f_2} & Y \ar[d]^{f_{13}}
\\ Y \ar[r]_{f_{13}} & Z &   &   & Y \ar[r]_{f_{23}} & Z
\\   &   & X \ar[r]^{f_1} \ar[d]_{f_2} & Y \ar[d]^{f_{23}}&   &
\\   &   & Y \ar[r]_{f_{12}} & Z &   &}
$$
Define a covariant pseudo functor $F:\mathcal{J} \rightarrow Cat$
which takes diagram (\ref{sourcediagram}) to diagram
(\ref{stackdefn}) and takes identity morphisms to identity
morphisms. The pseudo functor $F$ is defined on all possible
composites of nontrivial morphisms as:
$$F(f_{12} \circ f_1)(\{a_{\ell}\}_{\ell})_{ijk}:=
G(\pi_{ij}^1 \circ \pi_{ijk}^{12})a_i $$
$$F(f_{13} \circ f_2)(\{a_{\ell}\}_{\ell})_{ijk}:=
G(\pi_{ik}^2 \circ \pi_{ijk}^{13})a_k $$
$$F(f_{23} \circ f_1)(\{a_{\ell}\}_{\ell})_{ijk}:=
G(\pi_{jk}^1 \circ \pi_{ijk}^{23})a_j. $$ The identity coherence
isos $\delta^F$ for $F$ are equalities because $F$ takes identity
morphisms to identity morphisms. The coherence isos $\gamma^F$ for
composites of non-identity morphisms are defined as tuples of the
composition coherence isos for $G$. For example, the coherence iso
$\gamma^F_{f_1, f_{12}}\{a_{\ell}\}_{\ell}:F(f_{12}) \circ
F(f_1)\{a_{\ell}\}_{\ell} \rightarrow F(f_{12} \circ
f_1)\{a_{\ell}\}_{\ell}$ is defined as
$$\{\gamma^G_{\pi^{12}_{ijk},\pi^1_{ij}}a_i \}_{ijk}:
\{G(\pi_{ijk}^{12}) \circ G(\pi_{ij}^1) a_i \}_{ijk} \rightarrow
\{G(\pi_{ij}^1 \circ \pi_{ijk}^{12})a_i \}_{ijk}.
$$
The coherence isos $\gamma^F$ for composites involving one or more
identity morphisms are defined to be equalities. For example, the
coherence iso $$\gamma^F_{1_X, f_1}\{a_{\ell}\}_{\ell}:F(f_{1})
\circ F(1_X)\{a_{\ell}\}_{\ell} \rightarrow F(f_{1} \circ
1_X)\{a_{\ell}\}_{\ell}$$ is equality. The coherence diagram in the
pseudo functor unit axiom for $\delta^F$ is satisfied because of
this definition. The coherence diagram in the pseudo functor
composition axiom for $\gamma^F$ is satisfied because of the
diagrams for $\gamma^G$ and also because of this definition. The
coherence isos are also natural because $\mathcal{J}$ has no
nontrivial 2-cells. Thus $F:\mathcal{J} \rightarrow Cat$ is a pseudo
functor whose image is diagram (\ref{stackdefn}). By a {\it bilimit
of diagram\index{bilimit!bilimit of a diagram} (\ref{stackdefn})} we
mean a bilimit\index{bilimit} of this functor $F$.

In the context of stacks there is a canonical candidate for the
bilimit of $F$, namely $G(B)$. The candidate for the universal
pseudo cone $\pi':\Delta_{G(B)} \Rightarrow F$ is defined on objects
as follows.
$$\pi_X':G(B) \rightarrow \prod_iG(B_i)$$
$$\pi_X'(a):=\{G(g_i)a \}_i $$
$$ \pi_Y':G(B) \rightarrow \prod_{i,j} G(B_i \times_B
B_j)$$
$$\pi_Y'(a):=\{G( g_i \circ \pi^1_{ij})a \}_{ij} $$
 $$ \pi_Z':G(B) \rightarrow \prod_{i,j,k} G(B_i
\times_B B_j \times_B B_k)$$
$$\pi_Z'(a):=\{G(g_i \circ \pi_{ij}^1 \circ \pi_{ijk}^{12} )a\}_{ijk}$$
The coherence isos $\tau'_f:Ff \circ \pi_{Sf}' \Rightarrow
 \pi_{Tf}'\circ \Delta_{G(B)}(f)$ for the pseudo cone $\pi'$ and
non-identity morphisms $f$ in $\mathcal{J}$  are defined in terms of
$\gamma^G$. For example, for $f_1:X \rightarrow Y$ we have
$$\xymatrix@R=4pc@C=4pc{G(B) \ar[r]^{\pi_X'} \ar[d]_{1_{G(B)}=\Delta_{G(B)}(f_1)}
& \prod_i G(B_i) \ar[d]^{F(f_1)=p_1} \ar@{=>}[dl]^{\tau_{f_1}'}
\\ G(B) \ar[r]_-{\pi_Y'} & \prod_{i,j} G(B_i \times_B B_j)} $$
defined by $\tau_{f_1}'a:=\{\gamma^G_{\pi^1_{ij},g_i} a \}_{ij}:
\{G(\pi_{ij}^1) \circ G(g_i)a \}_{ij} \rightarrow \{G(g_i \circ
\pi_{ij}^1)a \}_{ij}$ for all objects $a$ of $G(B)$. For the
identity morphisms $1_X, 1_Y,$ and $1_Z$ of $\mathcal{J}$ we define
$\tau_{1_X}', \tau_{1_Y}',$ and $\tau_{1_Z}'$ to be equalities. The
coherence diagram for the unit axiom of pseudo natural
transformations is satisfied because of this definition. The
composition axiom for $\tau'$ and nontrivial morphisms is satisfied
because of the composition axiom for $\gamma^G$ and because
$\gamma^{\Delta_{G(B)}}$ is an equality. The composition axiom for
$\tau'$ whenever one or more of the morphisms is trivial follows
trivially. Thus $\pi':\Delta_{G(B)} \Rightarrow F$ is a pseudo
natural transformation with coherence isos $\tau'$. After these
preliminary remarks, we can finally define stack of categories.
\begin{defn}
Let $Cat$ denote the 2-category of small categories. Suppose
$\mathcal{B}$ is a Grothendieck topology. A {\it stack of
categories}\index{stack|textbf}\index{stack!stack of
categories|textbf} is a contravariant\index{contravariant|textbf}
pseudo functor $G: \mathcal{B} \rightarrow Cat$ which takes
Grothendieck covers\index{Grothendieck cover} to
bilimits\index{bilimit}, \ie for any object $B$ of $\mathcal{B}$ and
any Grothendieck cover $\{g_i:B_i \rightarrow B | i \in I \} \in
K(B)$ the diagram
$$\xymatrix@C=3pc@R=3pc{\prod_{i \in I} G(B_i)
\ar@<.5ex>[r]^-{p_1} \ar@<-.5ex>[r]_-{p_2} & \prod_{i,j \in I} G(B_i
\times_B B_j)  \ar[r]|-{p_{13}}  \ar@<1.2ex>[r]^-{p_{12}}
\ar@<-1.2ex>[r]_-{p_{23}} & \prod_{i,j,k \in I} G(B_i \times_B B_j
\times_B B_k) }$$ has $G(B)$ as a bilimit with universal pseudo cone
 $\pi':\Delta_{G(B)} \Rightarrow F$ as defined above.
\end{defn}
One common way to define a stack is via descent
objects\index{descent object} as in \cite{fantechi1},
\cite{giraud1}, \cite{moerdijk1}, or \cite{vistoli1}.

\begin{defn}
Let $\mathcal{B}$ be a Grothendieck topology and $G: \mathcal{B}
\rightarrow Cat$ a contravariant pseudo
functor\index{contravariant}\index{functor!pseudo functor}. Suppose
that $\{B_i \rightarrow B\}_i$ is a Grothendieck cover. Then an {\it
object with descent data on $\{B_i \rightarrow B\}_i$}\index{object
with descent data|textbf}\index{descent data|textbf} consists of an
object $\{a_i\}_i \in \prod_{i\in I} G(B_i)$ and isomorphisms
$\phi_{ij}:G(\pi^2_{ij})a_j \rightarrow G(\pi_{ij}^1)a_i$ in $G(B_i
\times_B B_j)$ which satisfy the {\it cocycle
condition}\index{cocycle condition|textbf}
$$G(\pi_{ijk}^{13})\phi_{ik}=G(\pi_{ijk}^{12})\phi_{ij}
\circ G(\pi_{ijk}^{23})\phi_{jk}$$ in $G(B_i \times_B B_j \times_B
B_k)$ up to the coherence isos of the pseudo functor $G$. See below.
A {\it morphism of descent objects}\index{morphism!morphism of
descent objects|textbf} $\{\xi_i\}_i:\{a_i\}_i \rightarrow
\{a_i'\}_i$ is a morphism in $\prod_{i \in I} G(B_i)$ such that the
diagram
$$\xymatrix@R=3pc@C=3pc{G(\pi_{ij}^2)a_j \ar[r]^{\phi_{ij}} \ar[d]_{G(\pi_{ij}^2)\xi_j}
& G(\pi_{ij}^1)a_i \ar[d]^{G(\pi_{ij}^1)\xi_i}
\\ G(\pi_{ij}^2)a_j' \ar[r]_{\phi_{ij}'} & G(\pi_{ij}^1)a_i' }$$
commutes in $G(B_i \times_B B_j)$. These objects and morphisms form
the {\it category of descent data on the cover $\{B_i \rightarrow
B\}_i$ }\index{descent data|textbf}\index{category of descent
data|textbf}. This category is denoted $G(\{B_i \rightarrow B\}_i)$.
There is a functor $G(B) \rightarrow G(\{B_i \rightarrow B\}_i)$
defined by $a \mapsto \{G(g_i)a\}_i$ where $g_i:B_i \rightarrow B$
are the morphisms from the Grothendieck cover. The $\phi_{ij}$
belonging to the image of $a$ under this functor are
$\phi_{ij}:=(\gamma^G_{\pi_{ij}^1,g_i}a)^{-1} \circ
(\gamma^G_{\pi_{ij}^2,g_j}a)$.
\end{defn}

The {\it cocycle condition}\index{cocycle condition|textbf} can be
stated explicitly as the requirement that the following diagram
commutes.
$$\xymatrix@C=4pc@R=4pc{G(\pi_{ijk}^{23})G(\pi^2_{jk})a_k
\ar[r]^{G(\pi_{ijk}^{23})\phi_{jk}} \ar[d]_{\gamma_{\pi_{ijk}^{23},
\pi^2_{jk}}a_k} & G(\pi_{ijk}^{23})G(\pi_{jk}^1)a_j
\ar[r]^{\gamma_{\pi_{ijk}^{23}, \pi^1_{jk}}a_j}    & G(\pi^1_{jk}
\circ \pi^{23}_{ijk}) a_j \ar@{=}[d] \\ G(\pi_{jk}^2 \circ
\pi^{23}_{ijk})a_k \ar@{=}[d] & & G(\pi_{ij}^2 \circ
\pi_{ijk}^{12})a_j \ar[d]^{\gamma^{-1}_{\pi_{ijk}^{12},
\pi_{ij}^2}a_j}
\\ G(\pi_{ik}^2 \circ \pi^{13}_{ijk})a_k
\ar[d]_{\gamma^{-1}_{\pi_{ijk}^{13},\pi_{ik}^2}a_k} & &
G(\pi_{ijk}^{12}) \circ G(\pi^2_{ij})a_j
\ar[d]^{G(\pi_{ijk}^{12})\phi_{ij}}
\\ G(\pi^{13}_{ijk}) \circ G(\pi_{ik}^2) a_k \ar[d]_{G(\pi_{ijk}^{13})
\phi_{ik}} & & G(\pi_{ijk}^{12}) \circ G(\pi_{ij}^1)a_i
\ar[d]^{\gamma_{\pi_{ijk}^{12},\pi_{ij}^1}a_i }
\\ G(\pi_{ijk}^{13}) \circ G(\pi_{ik}^1)  a_i
& G(\pi_{ik}^1 \circ \pi_{ijk}^{13})a_i
\ar[l]^-{\gamma^{-1}_{\pi_{ijk}^{13}, \pi_{ik}^1}a_i} & \ar@{=}[l]
G(\pi_{ij}^1 \circ \pi^{12}_{ijk})a_i }$$ This diagram is another
reason why we require our pseudo functors to have coherence arrows
that are iso: if $\gamma$ were not invertible, the cocycle condition
cannot be stated.

\begin{defn}
If $\mathcal{B}$ is a Grothendieck topology\index{Grothendieck
topology}, then a {\it Giraud stack of categories on
$\mathcal{B}$}\index{Giraud stack|textbf}\index{stack!Giraud
stack|textbf} is a contravariant pseudo
functor\index{contravariant}\index{functor!pseudo functor} $G:
\mathcal{B} \rightarrow Cat$ such that for any object $B$ of
$\mathcal{B}$ and any Grothendieck cover\index{Grothendieck cover}
$\{B_i \rightarrow B \}_i$ of $B$, the functor $G(B) \rightarrow
G(\{B_i \rightarrow B\}_i)$ is an equivalence of
categories.\footnote{This is not standard terminology. We have only
introduced it to distinguish the two definitions in the proof of
their equivalence. }
\end{defn}

\begin{thm}
Let $G:\mathcal{B} \rightarrow Cat$ be a contravariant pseudo
functor from a Grothendieck topology\index{Grothendieck topology} to
the 2-category of small categories. Then $G$ is a
stack\index{stack|textbf} if and only if it is a Giraud
stack\index{Giraud stack|textbf}\index{stack!Giraud stack|textbf}.
\end{thm}
\begin{pf}
From Chapter \ref{sec:laxlimitsinCat} we know that the category
$L:=PseudoCone(\mathbf{1},F)$ is a pseudo limit of $F$. It is
described as a subcategory of an appropriate product in Remarks
\ref{PL1} and \ref{PL2} in such a way that the pseudo cone
$\pi:\Delta_L \Rightarrow F$ consists of projections as in Remark
\ref{PL3a}.

We claim that the category $L$ of pseudo cones on a point is
equivalent to the category $G(\{B_i \rightarrow B\}_i)$ of descent
data by a functor $H:L \rightarrow G(\{B_i \rightarrow B\}_i)$.
Recall from Remark \ref{PL1} that each object of $L$ corresponds to
a tuple
$$\{a_i \}_i \times \{ a_{ij} \}_{ij} \times \{a_{ijk} \}_{ijk}
\times \{\varepsilon_f\}_f$$ of objects $$\{a_i \}_i \in
\prod_iG(B_i),$$ $$\{ a_{ij} \}_{ij} \in \prod_{ij}G(B_i \times_B
B_j),$$ $$\{a_{ijk} \}_{ijk} \in \prod_{ijk} G(B_i \times_B B_j
\times_B B_k),$$ and morphisms $\varepsilon_f$ indexed by morphisms
$f$ of $\mathcal{J}$ appropriately. For example,
$\varepsilon_{f_1}:F(f_1)\{a_i\}_i \rightarrow \{a_{ij} \}_{ij} $.
These morphisms satisfy the two axioms listed in Remark \ref{PL1}.
Each morphism in $L$ corresponds to a tuple
$$\{\xi_i \}_i \times \{ \xi_{ij} \}_{ij} \times \{\xi_{ijk}
\}_{ijk}$$ of morphisms in the product categories above and this
tuple commutes with the morphisms $\varepsilon_f$ appropriately.
Define
$$H( \{a_i \}_i \times \{ a_{ij} \}_{ij} \times \{a_{ijk} \}_{ijk} \times
\{\varepsilon_f\}_f):=\{a_i\}_i$$
$$H(\{\xi_i \}_i \times \{\xi_{ij} \}_{ij} \times \{\xi_{ijk} \}_{ijk}):=
\{\xi_i\}_i.$$ The descent data for $\{a_i\}_i$ are defined as the
components of $\{\phi_{ij}\}_{ij}:=(\varepsilon_{f_1})^{-1} \circ
\varepsilon_{f_2}$. Morphisms of $L$ map to morphisms of $G(\{B_i
\rightarrow B\}_i)$ because the outer diagram of
\begin{equation} \label{descentdata}
\xymatrix@R=3pc@C=3pc{F(f_1)\{a_i \}_i \ar[r]_-{\varepsilon_{f_1}}
\ar[d]_{F(f_1) \{\xi_i\}_i } & \{a_{ij}\}_{ij}
\ar[d]^{\{\xi_{ij}\}_{ij} } & F(f_2)\{a_i \}_i \ar[d]^{F(f_2)
\{\xi_i\}_i } \ar[l]^{\varepsilon_{f_2}} \ar@/_1pc/[ll]_{\{\phi_{ij}\}_{ij}} \\
F(f_1)\{a_i' \}_i
 \ar[r]^-{\varepsilon_{f_1}'} &  \{a_{ij}'\}_{ij}& F(f_2)\{a_i' \}_i
 \ar[l]_{\varepsilon_{f_2}'} \ar@/^1pc/[ll]^{\{\phi_{ij}'\}_{ij}}}
\end{equation}
commutes by Remark \ref{PL2}. To see that the $\phi_{ij}$ satisfy
the cocycle condition\index{cocycle condition}, consider the diagram
below.
\begin{equation} \label{towardscocycle}
\xymatrix@C=4pc@R=4pc{G(\pi_{ijk}^{23})G(\pi^2_{jk})a_k
\ar[r]^{G(\pi_{ijk}^{23})\phi_{jk}} \ar[d]_{\gamma_{\pi_{ijk}^{23},
\pi^2_{jk}}a_k} \ar[dr]_{G(\pi^{23}_{ijk})\varepsilon^{f_2}_{jk}} &
G(\pi_{ijk}^{23})G(\pi_{jk}^1)a_j \ar[r]^{\gamma_{\pi_{ijk}^{23},
\pi^1_{jk}}a_j}  \ar[d]^{G(\pi_{ijk}^{23} )\varepsilon^{f_1}_{jk} }
& G(\pi^1_{jk} \circ \pi^{23}_{ijk}) a_j \ar@{=}[d] \\
G(\pi_{jk}^2 \circ \pi^{23}_{ijk})a_k \ar@{=}[d] &
G(\pi_{ijk}^{23})a_{jk}  & G(\pi_{ij}^2 \circ \pi_{ijk}^{12})a_j
\ar[d]^{\gamma^{-1}_{\pi_{ijk}^{12}, \pi_{ij}^2}a_j}
\\ G(\pi_{ik}^2 \circ \pi^{13}_{ijk})a_k
\ar[d]_{\gamma^{-1}_{\pi_{ijk}^{13},\pi_{ik}^2}a_k} &
G(\pi_{ijk}^{12})a_{ij} & G(\pi_{ijk}^{12}) \circ G(\pi^2_{ij})a_j
\ar[d]^{G(\pi_{ijk}^{12})\phi_{ij}}
\ar[l]_{G(\pi_{ijk}^{12})\varepsilon^{f_2}_{ij}}
\\ G(\pi^{13}_{ijk}) \circ G(\pi_{ik}^2) a_k \ar[d]_{G(\pi_{ijk}^{13})
\phi_{ik}} \ar[r]^{G(\pi_{ijk}^{13})\varepsilon_{ik}^{f_2}} &
G(\pi_{ijk}^{13})a_{ik} & G(\pi_{ijk}^{12}) \circ G(\pi_{ij}^1)a_i
\ar[d]^{\gamma_{\pi_{ijk}^{12},\pi_{ij}^1}a_i }
\ar[ul]^{G(\pi_{ijk}^{12})\varepsilon^{f_1}_{ij}}
\\ G(\pi_{ijk}^{13}) \circ G(\pi_{ik}^1)  a_i
\ar[ur]_{G(\pi_{ijk}^{13})\varepsilon_{ik}^{f_1}} & G(\pi_{ik}^1
\circ \pi_{ijk}^{13})a_i \ar[l]^-{\gamma^{-1}_{\pi_{ijk}^{13},
\pi_{ik}^1}a_i} & \ar@{=}[l] G(\pi_{ij}^1 \circ \pi^{12}_{ijk})a_i }
\end{equation}
We want to show that the outer rectangle\index{cocycle condition}
commutes. The small triangles commute by definition of $\phi_{ij}$.
Next we draw another vertex $a_{ijk}$ inside the rectangle but
outside the triangles. Then we draw the arrows $\varepsilon_{ijk}^f$
for all non-identity morphisms $F$ of the category $\mathcal{J}$
with target $Z$. All of these arrows terminate at $a_{ijk}$. Each of
the resulting subdiagrams commutes because of the relations in
$\mathcal{J}$ or because of the second axiom on the morphisms
$\varepsilon_f$ in Remark \ref{PL1} . Note that we are using the
notation $\varepsilon_f=\{\varepsilon_{ijk}^{f} \}_{ijk}$. The outer
rectangle commute because all of the subdiagrams commute and
everything is iso. Hence the $\phi_{ij}$'s satisfy the cocycle
condition and $H$ maps $L$ into $G(\{B_i \rightarrow B\}_i)$. These
assignments obviously define a functor $H$.

The functor $H$ is faithful. Suppose
$$H(\{\xi_i \}_i \times \{\xi_{ij} \}_{ij} \times
\{\xi_{ijk} \}_{ijk})=H(\{\xi_i' \}_i \times \{\xi_{ij}' \}_{ij}
\times \{\xi_{ijk}' \}_{ijk}).$$ Then $\{\xi_i \}_i=\{\xi_i' \}_i$.
From this we conclude $\{\xi_{ij} \}_{ij}=\{\xi_{ij}' \}_{ij}$ by
diagram (\ref{descentdata}). A similar diagram with objects
$\{a_{ijk}\}$ and $\{a_{ijk}'\}$ in the center and arrows
$\varepsilon_{f_{12}},\varepsilon_{f_{23}}$ and
$\varepsilon_{f_{12}}',\varepsilon_{f_{23}}'$ pointing inward shows
that $\{\xi_{ijk} \}_{ijk}=\{\xi_{ijk}' \}_{ijk}$.

The functor $H$ is also full. Let $\{\xi_i \}_i$ be a morphism in
the category of descent data. Suppose further that its source and
target lie in the image of $H$. Then the outer diagram of diagram
(\ref{descentdata}) commutes and we define $\{ \xi_{ij} \}_{ij}$ to
be the unique arrow that makes diagram (\ref{descentdata}) commute.
It exists because the horizontal arrows are iso. We can also define
$\{\xi_{ijk} \}_{ijk}$ similarly, although we need to use diagram
(\ref{descentdata}) several times and the naturality of $\gamma^G$
to show that the necessary diagrams in Remark \ref{PL2} commute.

The functor $H$ is also surjective on objects. Suppose $\{a_i\}$ is
an object with descent data $\phi_{ij}$. Define
$a_{ij}:=G(\pi_{ij}^1)a_i$ and $a_{ijk}:=G(\pi_{ik}^1 \circ
\pi_{ijk}^{13})a_i$. Define $\varepsilon^{f_1}_{ij}:G(\pi_{ij}^1)a_i
\rightarrow a_{ij}$ to be the identity and
$\varepsilon_{ij}^{f_2}:=\phi_{ij}$. Let $\varepsilon_{ijk}^{f_{13}
\circ f_1}:G(\pi_{ik}^1 \circ \pi_{ijk}^{13})a_i \rightarrow
a_{ijk}$ also be the identity. Any $\varepsilon$ indexed by an
identity morphism is also trivial. Consider diagram
(\ref{towardscocycle}) with the additional vertex $a_{ijk}$ and the
additional $\varepsilon$'s mentioned just after diagram
(\ref{towardscocycle}). Requiring the inner diagrams to commute
uniquely defines the other $\varepsilon$'s which we did not define
yet. The commutivity of these smaller diagrams guarantees that the
tuple $$\{a_i \}_i \times \{ a_{ij} \}_{ij} \times \{a_{ijk}
\}_{ijk} \times \{\varepsilon_f\}_f$$ we have just defined is an
object of $L$. This object obviously maps under $H$ to $\{a_i\}_i$
with the correct descent data.

We conclude $H$ is an equivalence because it is faithfully full and
essentially surjective. Hence the category $L$ of pseudo cones is
equivalent to the category $G(\{B_i \rightarrow B\}_i)$ of descent
data.

There is also a functor $G(B) \rightarrow L$ defined like the
functor $G(B) \rightarrow G(\{B_i \rightarrow B\}_i)$ that makes the
diagrams
$$\xymatrix@R=3pc@C=3pc{G(\{B_i \rightarrow B\}_i) & L \ar[l]_-H &  \Delta_L
\ar@{=>}[d]^{\pi} \\
G(B) \ar[ur] \ar[u] & \Delta_{G(B)} \ar@{=>}[ur] \ar@{=>}[r]_{\pi'}
& F }$$ commute. Suppose $G$ is a Giraud stack. Then the left
vertical arrow is an equivalence. Hence the functor $G(B)
\rightarrow L$ is an equivalence and $\pi'$ makes $G(B)$ into a
bilimit of $F$ because $L$ is a bilimit of $F$ with pseudo limiting
cone $\pi$. Hence $G$ is a stack.

Suppose $G$ is a stack. Then $\pi'$ makes $G(B)$ into a bilimit of
$F$. Then the functor $G(B) \rightarrow L$ is an equivalence because
$L$ is also a bilimit and the right diagram commutes. Hence the
functor $G(B) \rightarrow G(\{B_i \rightarrow B\}_i)$ is also an
equivalence and $G$ is a Giraud stack.

This completes the proof that the two definitions of stack are
equivalent.

\end{pf}

Lastly, we define stacks of objects in a 2-category which admits
bilimits, such as the 2-category of pseudo
algebras\index{algebra!pseudo algebra} over a theory.

\begin{defn}
Let $\mathcal{C}$ be a 2-category whose objects have underlying
categories. Suppose $\mathcal{B}$ is a Grothendieck
topology\index{Grothendieck topology} and $\mathcal{C}$ admits
bilimits. A {\it stack of objects of
$\mathcal{C}$}\index{stack|textbf} is a
contravariant\index{contravariant} pseudo
functor\index{functor!pseudo functor} $G: \mathcal{B} \rightarrow
\mathcal{C}$ which takes Grothendieck covers\index{Grothendieck
cover} to bilimits\index{bilimit|textbf}, \ie for any object $B$ of
$\mathcal{B}$ and any Grothendieck cover $\{g_i:B_i \rightarrow B |
i \in I \} \in K(B)$ the diagram
$$\xymatrix@C=3pc{\prod_{i \in I} G(B_i)
\ar@<.5ex>[r]^-{p_1} \ar@<-.5ex>[r]_-{p_2} & \prod_{i,j \in I} G(B_i
\times_B B_j)  \ar[r]|-{p_{13}}  \ar@<1.2ex>[r]^-{p_{12}}
\ar@<-1.2ex>[r]_-{p_{23}} & \prod_{i,j,k \in I} G(B_i \times_B B_j
\times_B B_k) }$$ has $G(B)$ as a bilimit with universal pseudo cone
 $\pi':\Delta_{G(B)} \Rightarrow F$ as defined above.
\end{defn}

For example, a stack of pseudo algebras\index{algebra!pseudo
algebra!stacks of pseudo algebras}\index{stack|stack of pseudo
algebras|textbf} over a theory\index{theory} $T$ is a contravariant
pseudo functor\index{contravariant}\index{functor!pseudo functor}
from a Grothendieck topology\index{Grothendieck topology} into the
2-category of pseudo $T$-algebras which takes Grothendieck covers to
bilimits\index{bilimit} in the above sense.

\chapter{2-Theories, Algebras, and Weighted Pseudo Limits}
\label{sec:2-theories} The algebraic structure of the category of
rigged surfaces\index{rigged surface} can be described as a pseudo
algebra\index{algebra!pseudo algebra} over a certain
2-theory\index{2-theory} as in \cite{hu}, \cite{hu1}, and
\cite{hu2}. A {\it pseudo algebra over a
2-theory}\index{algebra!pseudo algebra over a 2-theory} in this
paper is the same as a {\it lax\index{lax} algebra over a
2-theory}\index{algebra!lax algebra over a 2-theory} in \cite{hu},
\cite{hu1}, and \cite{hu2}. However, the 2-theories of
\cite{yanofskyproposal}, \cite{yanofsky2000}, and
\cite{yanofsky2001} are different from the 2-theories in this paper.
In this chapter we review the relevant terminology and prove results
about limits. Before giving the definition of a 2-theory, we
motivate it with an example in the first section.

\section{The 2-Theory $End(X)$ Fibered over the Theory $End(I)$}
Let $I$ be a category and $k$ a positive
integer\index{2-theory!endomorphism
2-theory|(textbf}\index{theory!endomorphism
2-theory|(textbf}\index{2-theory!$(End(X),End(I))$|(textbf}\index{$(End(X),End(I))$|(textbf}\index{$End(X)$|(textbf}\index{endomorphism
2-theory|(textbf}. Suppose $X:I^k \rightarrow Cat$ is a strict
2-functor from the category $I^k$ to the 2-category $Cat$ of small
categories. Here $I^k$ is interpreted as a 2-category where the hom
sets are discrete categories. We will now describe the {\it 2-theory
End(X) fibered over the theory End(I)}, which is a contravariant
functor $End(I) \rightarrow Cat$ satisfying certain properties.

Recall that the theory $End(I)$\index{$End(I)$} is the category with
objects $0=\{*\}, 1=I,2=I^2, 3=\dots$ and morphisms
$Mor_{End(I)}(m,n)=Functors(I^m,I^n)$. Here $\{*\}$ denotes the
terminal object\index{terminal object} in the category of small
categories. As with any theory\index{theory}, the theory $End(I)$
can be completely described by the sets
$End(I)(n):=Mor_{End(I)}(n,1)$, a composition\index{composition},
substitution\index{substitution}, and a unit\index{unit} which
satisfy a list of axioms. See Theorem \ref{equivalenceoftheories} or
\cite{hu} for details.

From the theory $End(I)$ we can obtain another category denoted
$End(I)^k$, which also turns out to be a theory. It has objects
$0=\{*\} \times \cdots \times \{*\}, 1=I \times \cdots \times I,
2=I^2 \times \cdots \times I^2,
 3=\dots$ ($k$ copies in each product)
and it has morphisms $Mor_{End(I)^k}(m,n):=Mor_{End(I)}(m,n)^{\times
k}$. For example, $v \in Mor_{End(I)^k}(m,1)$ is a functor
$v:(I^m)^k \rightarrow I^k$ that is a $k$-tuple of functors $I^m
\rightarrow I$. For $n \in \mathbb{N}$ and $1\leq i \leq n$, let
$pr_i^{\times k}:(I^n)^k \rightarrow I^k$ be the morphism
$pr_i^{\times k} \in Mor_{End(I)^k}(n,1)$ whose $k$ components are
each the projection functor $pr_i:I^n \rightarrow I$ onto the $i$-th
coordinate. We can easily check that $n \in Obj \hspace{1mm}
End(I)^k$ is the product in $End(I)^k$ of $n$ copies of $1$ with
projection morphisms $pr_1^{\times k}, \dots, pr_n^{\times k}$.
Hence $End(I)^k$ is itself a theory and $Mor_{End(I)^k}(m,n)$ is in
bijective correspondence with $\prod_{i=1}^n Mor_{End(I)^k}(m,1)$.
We identify these two sets via the usual bijection. In other words,
for $k$-tuples $w_1, \dots, w_n \in Mor_{End(I)^k}(m,1)$ we let
$\prod_{j=1}^n w_j$ denote the unique morphism $m \rightarrow n$ of
$End(I)^k$ such that
$$\xymatrix@R=3pc@C=3pc{n \ar[r]^{pr_i^{\times k}} & 1
\\ m \ar@{.>}[u]^{\prod_{j=1}^n w_j} \ar[ur]_{w_i} &}$$
commutes for all $i=1, \dots, n$. This notation differs from
\cite{hu}, in which the notation $(w_1, \dots, w_n)$ is used instead
of the product. We reserve $(w_1, \dots, w_n)$ for a different
morphism. The reason for our choice will become clear later. Using
our convention, we have $w=\prod^n_{j=1} pr_j^{\times k} \circ w$
for $w \in Mor_{End(I)^k}(m,n)$.

Since $End(I)^k$ is a theory, it has a
substitution\index{substitution} and a
composition\index{composition} with unit\index{unit} which satisfy
certain axioms described in Chapter \ref{sec:theories} and
\cite{hu}. If $f:\{1, \dots, p\} \rightarrow \{1, \dots, q\}$ is a
function and $w \in$ {\small $ End(I)^k(p) = Mor_{End(I)^k}(p,1) =
Mor_{End(I)}(p,1)^ {\times k}$}, then the substituted
word\index{word!substituted word} $w_f$ is obtained by substituting
by $f$ in each of the words in the $k$-components of $w$. The
composition is also done componentwise. The unit\index{unit}
$1^{\times k}: I \times \cdots \times I \rightarrow  I \times \cdots
\times I$ is $k$ copies of the unit $1:I \rightarrow I$ in the
theory $End(I)$. These explicit descriptions of substitution,
composition, and unit follow from the definitions of the projections
in the theory $End(I)^k$ by the work in Chapter \ref{sec:theories}.

We follow the conventions of Chapter \ref{sec:theories} to define a
morphism $(w_1, \dots, w_n)$. Let $w_i \in End(I)^k(m_i)$ for $i=1,
\dots, n$. Let $\iota_i:\{1, \dots, m_i \} \rightarrow \{1, \dots,
m_1 + m_2 + \cdots + m_n\}$ be the injective function which takes
the domain to the $i$-th block. Then there exists a unique morphism
$(w_1, \dots,w_n)$ such that
$$\xymatrix@R=3pc@C=3pc{n \ar[r]^{pr_i^{\times k}} & 1
\\ m_1 + m_2 + \cdots + m_n \ar@{.>}[u]^{(w_1, \dots,w_n)}
\ar[ur]_{\hspace{2mm}(w_i)_{\iota_i}} &}$$ commutes for all $i=1,
\dots, n$. Explicitly, the morphism $(w_1, \dots, w_n)$ is obtained
by doing an analogous process in each of the $k$ components.

The strict 2-functor $X:I^k \rightarrow Cat$ gives rise to a
contravariant\index{contravariant} functor $End(X): End(I)
\rightarrow Cat$ as follows. For $m \in Obj \hspace{1mm}  End(I)$
the category $End(X)(m)$ has objects $Obj \hspace{1mm}  End(X)(m) =
\coprod_{n \geq 0} Mor_{End(I)^k}(m,n)$, in other words, the objects
of $End(X)(m)$ are the arrows of $End(I)^k$ with domain $m$. For
$\prod_{i=1}^p v_i, \prod_{i=1}^q w_i \in Obj \hspace{1mm}
End(X)(m)$ where $v_1, \dots, v_p,w_1, \dots, w_q \in
Mor_{End(I)^k}(m,1)$ we define the set of morphisms
$Mor_{End(X)(m)}(\prod_{i=1}^p v_i, \prod_{i=1}^q w_i)$ to be the
collection of natural transformations
\begin{equation} \label{morphism}
\alpha: X \circ v_1 \circ d^m \times \cdots \times X \circ v_p \circ
d^m \Rightarrow X \circ w_1 \circ d^m \times \cdots \times X \circ
w_q \circ d^m
\end{equation}
where $d^m: I^m \rightarrow (I^m)^k$ is the diagonal functor. Note
that $X \circ v_1 \circ d^m \times \cdots \times X \circ v_p \circ
d^m$ and $X \circ w_1 \circ d^m \times \cdots \times X \circ w_q
\circ d^m$ are functors $I^m \rightarrow Cat$. The composition of
morphisms in $End(X)(m)$ is the vertical composition of natural
transformations. With these definitions, $End(X)(m)$ is a category.
We must still define the contravariant functor $End(X)$ on morphisms
and verify that it preserves identities and compositions. For any
morphism $u:I^{\ell} \rightarrow I^m$ of the theory $End(I)$, define
$u^{\times k}:(I^{\ell})^k \rightarrow (I^{m})^k$ to be the functor
which is $u$ in each of the $k$ components. Note that $u^{\times k}
\circ d^{\ell} = d^m \circ u:I^{\ell} \rightarrow (I^m)^k$. The
functor $End(X)(u): End(X)(m) \rightarrow End(X)(\ell)$ is defined
on objects by $End(X)(u)(\prod_{i=1}^p v_i) := \prod_{i=1}^p v_i
\circ u^{\times k}$ and on morphisms $\alpha$ in (\ref{morphism}) by
$End(X)(u)(\alpha ) := \alpha * i_u$ where $*$ denotes the
horizontal composition of natural transformations and $i_u:u
\Rightarrow u$ is the trivial natural transformation. This makes
sense because
$$\aligned (X \circ v_1 \circ d^m \times \cdots \times X \circ v_p \circ d^m)
\circ u &= X \circ v_1 \circ d^m \circ u \times \cdots \times X
\circ v_p
\circ d^m \circ u \\
&=X \circ v_1 \circ u^{\times k} \circ  d^{\ell} \times \cdots
\times X \circ v_p \circ u^{\times k} \circ d^{\ell}
\endaligned
$$
and
$$\alpha * i_u: X \circ v_1 \circ u^{\times k} \circ  d^{\ell} \times \cdots
\times X \circ v_p \circ u^{\times k} \circ d^{\ell} \Rightarrow X
\circ w_1 \circ u^{\times k} \circ  d^{\ell} \times \cdots \times X
\circ w_q \circ u^{\times k} \circ d^{\ell}$$ really is a morphism
$$End(X)(u)(\prod_{i=1}^p v_i) = \prod_{i=1}^p v_i \circ u^{\times k}
\rightarrow \prod_{i=1}^p w_i \circ u^{\times
k}=End(X)(u)(\prod_{i=1}^q w_i).$$ If $u:I^{\ell} \rightarrow I^m$
is the identity functor $I^m \rightarrow I^m$, then
$End(X)(u):End(X)(m) \rightarrow End(X)(m)$ is also the identity
functor because $v_i \circ u^{\times k} = v_i$ for $i=1, \dots, p$
and $w_i \circ u^{\times k} =w_i$ for $i=1, \dots, q$ and also
$\alpha * i_u = \alpha$. If $\xymatrix@1{I^j \ar[r]^{u_1} & I^m
\ar[r]^{u_2} & I^{\ell}}$ are morphisms in $End(I)$, then
$u_2^{\times k} \circ u_1^{\times k} = (u_2 \circ u_1)^{\times k}$
and
$$(\alpha * i_{u_2})*i_{u_1}=\alpha * (i_{u_2}*i_{u_1})
=\alpha*i_{u_2 \circ u_1},$$ which together imply that
$$End(X)(u_2 \circ u_1)=End(X)(u_1) \circ End(X)(u_2).$$ Thus
$End(X):End(I) \rightarrow Cat$ preserves identities and
compositions and is a contravariant functor.

The category $End(X)(m)$ also admits certain products, which will be
a feature of a general 2-theory. For $v_1, \dots, v_p \in
Mor_{End(I)^k}(m,1)$ and $\prod_{i=1}^p v_i \in Mor_{End(I)^k}(m,p)
\subseteq Obj \hspace{1mm}  End(X)(m)$ define projections
$pr_j:\prod_{i=1}^p v_i \rightarrow v_j$ for $j=1, \dots, p$ to be
the projection natural transformations
$$X \circ v_1 \circ d^m \times \cdots \times X \circ v_p \circ d^m
\Rightarrow X \circ v_j \circ d^m.$$ Then $\prod_{i=1}^p v_i$ is
obviously the product of $v_1, \dots, v_p$ in the category
$End(X)(m)$ with these projections. This explains the choice of
notation $\prod_{i=1}^p v_i$. This product property will also be
required of a general 2-theory. We record for later use how these
products allow us to define morphisms $\iota'$ for every function
$\iota:\{1, \dots, p\} \rightarrow \{1, \dots, q\}$. Let $w_1,
\dots, w_q \in Mor_{End(I)^k}(m,1) \subseteq Obj \hspace{1mm}
End(X)(m)$. Then for a function $\iota:\{1, \dots, p\} \rightarrow
\{1, \dots, q\}$ there exists a unique morphism $\iota'$ such that
\begin{equation} \label{iotaprime}
\xymatrix@R=3pc@C=3pc{\prod_{i=1}^p w_{\iota(i)} \ar[r]^-{pr_{\ell}}
& w_{\ell}
\\ \prod_{i=1}^q w_i \ar@{.>}[u]^{\iota'} \ar[ur]_{pr_{\iota(\ell)}}
}
\end{equation}
commutes for all $\ell=1, \dots, p$. The arrows of the natural
transformation $\iota':X \circ w_1 \circ d^m \times \cdots \times X
\circ w_q \circ d^m \Rightarrow X \circ w_{\iota(1)} \circ d^m
\times \cdots \times X
 \circ w_{\iota(p)} \circ d^m$ have the appropriate projections as
their components.

The 2-theory $End(X)$ has several operations\index{operations of
2-theories} on it which any general 2-theory will also have, once we
define the notion of 2-theory. To make the description of these
operations easier, we follow the notation introduced by P. Hu and I.
Kriz in \cite{hu}. For objects $w, w_1, \dots, w_q \in
Mor_{End(I)^k}(m,1) \subseteq Obj \hspace{1mm}  End(X)(m)$ we set
$$End(X)(w;w_1, \dots, w_q):=Mor_{End(X)(m)}(\prod_{i=1}^q w_i,w).$$
The operations of P. Hu and I. Kriz are collated in the following
theorem.

\begin{thm} \label{2-theoryoperations}
The contravariant\index{contravariant} functor $End(X):End(I)
\rightarrow Cat$ has the following operations.\index{operations of
2-theories}
\begin{enumerate}
\item
For each $w \in T^k(m)$ there exists a unit\index{unit} $1_w \in
End(X)(w;w)$.
\item
For all $w, w_i, w_{ij} \in Mor_{End(I)^k}(m,1)$ there is a function
called $End(X)$-composition\index{$End(X)$-composition}.

\begingroup
\vspace{-2\abovedisplayskip} \small
$$\gamma:End(X)(w;w_1, \dots, w_q) \times End(X)(w_1;w_{11}, \dots, w_{1p_1})
\times \cdots \times End(X)(w_q;w_{q1}, \dots, w_{qp_q})$$
$$\rightarrow End(X)(w;w_{11}, \dots, w_{qp_q})$$
\endgroup
\noindent
\item
Let $w, w_1, \dots, w_q \in Mor_{End(I)^k}(m,1)$. For any function
$\iota:\{1, \dots, p\} \rightarrow \{1, \dots, q\}$ there is a
function
$$()^{\iota}:End(X)(w;w_{\iota(1)}, \dots,w_{\iota(p)}) \rightarrow
End(X)(w;w_1, \dots, w_q)$$ called
$End(X)$-functoriality\index{$End(X)$-functoriality}.
\item
Let $w, w_1, \dots, w_q \in Mor_{End(I)^k}(m,1)$. For any function
$f:\{1, \dots,m\} \rightarrow \{1, \dots, \ell \}$ there is a
function $$()_f:End(X)(w;w_1, \dots, w_q) \rightarrow End(X)(w_{f};
(w_1)_{f}, \dots, (w_q)_{f})$$ where $w_{f}$ means to substitute $f$
in each of the words in the $k$-tuple $w$. This function is called
$End(I)$-functoriality\index{$End(I)$-functoriality}. Note that
$End(X)(w;w_1, \dots, w_q)$ is a
hom set in the category $End(X)(m)$ while on the other hand \\
$End(X)(w_{f}; (w_1)_{f}, \dots, (w_q)_{f})$ is a hom set in the
category $End(X)(\ell)$.
\item
For $u_i \in End(I)(k_i), i=1, \dots, m$ and $w, w_1, \dots, w_q \in
Mor_{End(I)^k}(m,1)$ let $v_j:=\gamma^{\times k}(w_j;u_1^{\times k},
\dots,u_m^{\times k})$ for $j=1, \dots, q$ and furthermore let
$v:=\gamma^{\times k}(w;u_1^{\times k}, \dots,u_m^{\times k})$. Then
there is a function
$$(u_1,\dots,u_m)^*:End(X)(w;w_1, \dots, w_q) \rightarrow
End(X)(v;v_1, \dots, v_q)$$ called
$End(I)$-substitution\index{$End(I)$-substitution}. Here
$\gamma^{\times k}$ means to use the composition of the theory
End(I) in each of the $k$ components, which coincides with
composition in the theory $End(I)^k$. Note that $End(X)(w;w_1,
\dots, w_q)$ is a hom set in the category $End(X)(m)$ while
$End(X)(v;v_1, \dots, v_q)$ is a hom set in the category $End(X)(k_1
+ \cdots + k_m)$.
\end{enumerate}
\end{thm}
\begin{pf}
\begin{enumerate}
\item
The unit $1_w: X \circ w \circ d^m \Rightarrow X \circ w \circ d^m$
is the identity natural transformation $i_{X \circ w \circ d^m}:X
\circ w \circ d^m \Rightarrow X \circ w \circ d^m$.
\item
Let $\alpha:\prod_{i=1}^q w_i \rightarrow w$ and
$\alpha_i:\prod_{j=1}^{p_i} w_{ij} \rightarrow w_i$ for
$i=1,\dots,q$ be morphisms of $End(X)(m)$. Let $\iota_{\ell}: \{1,
\dots, p_{\ell}\} \rightarrow \{1, \dots, p_1+p_2 + \cdots + p_{q}
\}$ be the injective function which takes the domain to the
$\ell$-th block. We take the product $\prod_{i=1}^q
\prod_{j=1}^{p_i} w_{ij}$ to be
$$\prod_{i=1}^q \prod_{j=1}^{p_i} w_{ij}=w_{11} \times w_{12} \times \cdots
w_{1p_1} \times w_{21} \times \cdots \times w_{2p_2} \times w_{31}
\times \cdots \times w_{qp_q}.$$ Then there exists a unique morphism
$(\alpha_1, \dots, \alpha_q)$ such that
$$\xymatrix@R=3pc@C=3pc{\prod_{i=1}^q w_i \ar[r]^-{pr_{\ell}^{\times k}} & w_{\ell}
\\ \prod_{i=1}^q \prod_{j=1}^{p_i} w_{ij} \ar@{.>}[u]^{(\alpha_1, \dots,
\alpha_q)} \ar[ur]_{(\alpha_{\ell})_{\iota_{\ell}}} & }$$ commutes
for all $\ell=1, \dots, q$. This means that $$(\alpha_1, \dots,
\alpha_q): X \circ w_{11} \circ d^m \times \cdots \times X \circ
w_{qp_q} \circ d^m \Rightarrow$$
$$ X \circ w_1 \circ
d^m \times X \circ w_2 \circ d^m \times \cdots \times X \circ w_q
\circ d^m$$ is the natural transformation which is $\alpha_{\ell}$
on $X \circ w_{\ell 1} \circ d^m \times \cdots \times X \circ
w_{\ell p_{\ell}} \circ d^m$. Define
$$\gamma(\alpha;\alpha_1, \dots, \alpha_q):= \alpha \circ (\alpha_1, \dots
\alpha_q)$$ where the composition is in the category $End(X)(m)$.
\item
Let $w_1, \dots, w_q \in Mor_{End(I)^k}(m,1)$ and $\iota:\{1, \dots,
p\} \rightarrow \{1, \dots, q \}$ be a function. Let $\iota':
\prod_{i=1}^q w_i \rightarrow \prod_{i=1}^q w_{\iota(i)}$ be the
morphism defined in diagram (\ref{iotaprime}). Then we define
$End(X)$-functoriality $$End(X)(w;w_{\iota (1)}, \dots, w_{\iota
(p)}) \rightarrow End(X)(w;w_1, \dots, w_q)$$ by $\alpha \mapsto
\alpha \circ \iota'$.
\item
A function $f:\{1, \dots, m\} \rightarrow \{1, \dots, \ell \}$
induces a morphism $f':\ell \rightarrow m$ in $End(I)$ which in turn
gives rise to a morphism $(f')^{\times k}:(I^{\ell})^k \rightarrow
(I^m)^k$ in $End(I)^k$. Then $w_f=w \circ (f')^{\times k}$ by
definition and the functor $End(X)(f'):End(X)(m) \rightarrow
End(X)(\ell)$ gives us a map of hom sets
$$()_f:End(X)(w;w_1, \dots, w_q) \rightarrow End(X)(w_{f};
(w_1)_{f}, \dots, (w_q)_{f}).$$
\item
Let $\iota_i:\{1, \dots, k_i\} \rightarrow \{1, \dots, k_1 + k_2+
\cdots + k_m\}$ be the injective map which takes the domain to the
$i$-th block. Let $(u_1^{\times k}, \dots, u_m^{\times k})$ denote
the unique morphism in $End(I)^k$ such that
$$\xymatrix@R=3pc@C=3pc{m \ar[r]^{pr_i^{\times k}} & 1
\\  k_1 + k_2+ \cdots + k_m \ar@{.>}[u]^{(u_1^{\times k}, \dots,
 u_m^{\times k})} \ar[ur]_{(u_i)^{\times k}_{\iota_i}}}$$
commutes. Then we know from the general theory of
theories\index{theory!theory of theories} that
\newline $\gamma^{\times k}(w;u_1^{\times k}, \dots, u_m^{\times
k})= w \circ (u_1^{\times k}, \dots, u_m^{\times k})$ where the
composition ``$\circ$'' is the composition in the category
$End(I)^k$. Then $End(X)(u_1^{\times k}, \dots, u_m^{\times
k})(w)\\=v$ and the functor $End(X)(u_1^{\times k}, \dots,
u_m^{\times k})$ gives us the desired map of hom sets.
\end{enumerate}
\end{pf}

These operations on $End(X)$ satisfy certain relations.

\begin{thm} \label{2-theoryrelations}
The operations on the contravariant\index{contravariant} functor
$End(X):End(I) \rightarrow Cat$ satisfy the following
relations.\index{relations of 2-theories}
\begin{enumerate}
\item
$End(X)$-composition\index{$End(X)$-composition} is associative, \ie
\newline $\gamma(\alpha;\gamma(\alpha^1;\alpha^1_1, \dots,
\alpha^1_{n_1}), \gamma(\alpha^2;\alpha^2_1, \dots, \alpha^2_{n_2}),
\dots, \gamma(\alpha^q;\alpha^q_1, \dots, \alpha^q_{n_q}))$ is the
\\ same as $\gamma(\gamma(\alpha; \alpha^1, \dots, \alpha^q),
\alpha_1^1, \dots \alpha_{n_1}^1, \alpha_1^2, \dots, \alpha_{n_2}^2,
\dots, \alpha_1^q, \dots, \alpha_{n_q}^q)$.
\item
$End(X)$-composition is unital\index{$End(X)$-composition}, \ie for
$\alpha \in End(X)(w;w_1, \dots, w_q)$ we have
$\gamma(\alpha;1_{w_1}, \dots, 1_{w_q})=\alpha=\gamma(1_w;\alpha)$.
\item
$End(X)$-functoriality is functorial\index{$End(X)$-functoriality},
\ie for functions
\newline $\xymatrix@1{\{1, \dots, p\} \ar[r]^{\iota} & \{1, \dots,
q\} \ar[r]^{\theta} & \{1, \dots, r\}}$ the composition
$$\xymatrix{
End(X)(w;w_{\theta \iota (1)}, \dots, w_{\theta \iota (p)})
\ar[r]^{()^{\iota}} & End(X)(w;w_{\theta (1)}, \dots, w_{\theta
(q)})}$$
$$\xymatrix@1{\ar[r]^-{()^{\theta}} & End(X)(w;w_{1}, \dots, w_{r})}$$
is the same as $$\xymatrix@1{ End(X)(w;w_{\theta \iota (1)}, \dots,
w_{\theta \iota (p)}) \ar[r]^-{()^{\theta \circ \iota}} &
End(X)(w_1, \dots, w_r)}$$ and for the identity $id_q:\{1, \dots,
q\} \rightarrow \{1, \dots, q\}$ the map \newline
$()^{id_q}:End(X)(w;w_1, \dots, w_q) \rightarrow End(X)(w;w_1,
\dots, w_q)$ is the identity.
\item
The $End(X)$-compositions\index{$End(X)$-composition} $\gamma$ are
equivariant\index{equivariant} with respect to
$End(X)$-\\functoriality\index{$End(X)$-functoriality} in the sense
that if $\iota:\{1, \dots,p\} \rightarrow \{1, \dots, q\}$ is a
function, $\alpha \in End(X)(w;w_{\iota (1)}, \dots, w_{\iota(p)})$,
and $\alpha_{\ell} \in End(X)(w_{\ell};w_{\ell 1}, \dots, w_{\ell
p_{\ell}})$ for $\ell = 1, \dots, q$ then
$$\gamma(\alpha^{\iota};\alpha_1, \dots, \alpha_q) = \gamma(\alpha;
\alpha_{\iota (1)}, \dots, \alpha_{\iota (p)})^{\bar{\iota}},$$
where $\bar{\iota}:\{1,2, \dots, p_{\iota (1)} + \cdots + p_{\iota
(p)} \} \rightarrow \{1,2, \dots, p_1 + \cdots + p_q\}$ is the
function obtained by parsing the sequence $1,2, \dots, p_1+ \cdots
+p_q$ into consecutive blocks $B_1, \dots, B_q$ of lengths $p_1,
\dots, p_q$ and then writing them in the order $B_{\iota (1)},
\dots, B_{\iota (p)}$ as in Example \ref{Endexample1}.
\item
The $End(X)$-compositions\index{$End(X)$-composition} $\gamma$ are
equivariant\index{equivariant} with respect to
\newline $End(X)$-functoriality\index{$End(X)$-functoriality} in the sense that if $\alpha \in
End(X)(w;w_1, \dots, w_q)$, \newline $\alpha_{\ell} \in
End(X)(w_{\ell};w_{\ell \iota_{\ell}(1)}, \dots, w_{\ell
\iota_{\ell}(p_{\ell}')})$, and $\iota_{\ell}:\{1, \dots, p_{\ell}'
\} \rightarrow \{1, \dots, p_{\ell}\}$ are functions for $\ell=1,
\dots, q$ then
$$\gamma(\alpha;(\alpha_1)^{\iota_1}, \dots, (\alpha_q)^{\iota_q})=
\gamma(\alpha; \alpha_1, \dots, \alpha_q)^{\iota_1+ \cdots \iota_q}
$$ where $\iota_1+ \cdots + \iota_q:\{1, \dots, p_1'+ \cdots +
p_q'\} \rightarrow \{1, \dots, p_1 + \dots + p_q\}$ is the function
obtained by placing $\iota_1, \dots, \iota_q$ side by side.
\item
$End(I)$-functoriality\index{$End(I)$-functoriality} is functorial,
\ie \newline for functions $\xymatrix@1{\{1, \dots,n\} \ar[r]^f &
\{1, \dots, m\} \ar[r]^g & \{1, \dots, \ell \}}$ and words \newline
 $w, w_1, \dots, w_q \in Mor_{End(I)^k}(n,1)$
the composition
$$\xymatrix@1{End(X)(w;w_1, \dots, w_q) \ar[r]^-{()_f} & End(X)(w_f;
(w_1)_f, \dots, (w_q)_f)}$$
$$\xymatrix@1{\ar[r]^-{()_g} & End(X)((w_f)_g;((w_1)_f)_g,\dots,
((w_q)_f)_g )}$$ is the same as
$$\xymatrix@1{End(X)(w;w_1, \dots, w_q) \ar[r]^-{()_{g \circ f}} &
End(X)(w_{g \circ f};(w_1)_{g \circ f}, \dots, (w_q)_{g \circ f})}$$
and for the identity $id_n:\{1, \dots n\} \rightarrow \{1, \dots,
n\}$ the map
\newline
$()_{id_n}:End(X)(w;w_1, \dots, w_q) \rightarrow End(X)(w;w_1,
\dots, w_q)$ is the identity.
\item
$End(I)$-substitution\index{$End(I)$-substitution} is associative.
\newline Let $w, w_1, \dots, w_q \in Mor_{End(I)^k} (m,1)$, $t_i \in
End(I)(k_i)$ for $i=1, \dots, m$ and $s_{ij} \in End(I)(k_{ij})$ for
$1 \leq i \leq m$ and $1 \leq j \leq k_i$. Let $$\aligned
v&:=\gamma^{\times
k}(w;t_1^{\times k}, \dots, t_m^{\times k}) \\
v_{\ell}&:=\gamma^{\times k}(w_{\ell};t_1^{\times k}, \dots,
t_m^{\times k}) \\
u&:=\gamma^{\times k}(v;s_{11}^{\times k}, s_{12}^{\times k}, \dots,
s_{1k_1}^{\times k}, s_{21}^{\times k}, \dots, s_{31}^{\times k},
\dots, s_{m1}^{\times k}, \dots, s_{mk_m}^{\times k}) \\
u_{\ell}&:=\gamma^{\times k}(v_{\ell};s_{11}^{\times k},
s_{12}^{\times k}, \dots, s_{1k_1}^{\times k}, s_{21}^{\times k},
\dots, s_{31}^{\times k}, \dots, s_{m1}^{\times k}, \dots,
s_{mk_m}^{\times k})
\endaligned$$ for $\ell=1, \dots, q$. Then the composition
$$\xymatrix@1@C=5pc{End(X)(w;w_1, \dots, w_q) \ar[r]^{(t_1, \dots, t_m)^*}
& End(X)(v;v_1, \dots, v_q)}$$
$$\xymatrix@1@C=5pc{\ar[r]^-{(s_{11}, \dots, s_{mk_m})^*} & End(X)(u;u_1, \dots
u_q)}$$ is the same as
$$\xymatrix@1@C=5pc{End(X)(w;w_1, \dots, w_q) \ar[r]^{(r_1, \dots, r_m)^*}
& End(X)(u;u_1, \dots, u_q)}$$ where $r_i=\gamma^{\times
k}(t_i^{\times k};s_{i1}^{\times k}, s_{i2}^{\times k}, \dots,
s_{ik_i}^{\times k})=\gamma_{End(I)}(t_i;s_{i1}, s_{i2}, \dots,
s_{ik_i})^{\times k}$ for $i=1, \dots, m$. Note that
$u=\\\gamma^{\times k}(w;\gamma^{\times k}(t_1^{\times
k};s_{11}^{\times k}, s_{12}^{\times k}, \dots, s_{1k_1}^{\times
k}), \dots, \gamma^{\times k}(t_m^{\times k};s_{m1}^{\times k},
s_{m2}^{\times k}, \dots, s_{mk_m}^{\times k}))$.
\item
$End(I)$-substitution\index{$End(I)$-substitution} is unital.
\newline For the unit\index{unit} $1 \in End(I)(1)$ of the theory $End(I)$ and
\newline $w, w_1, \dots, w_q \in Mor_{End(I)^k}(m,1)$ the function
$$(1, \dots, 1)^*:End(X)(w;w_1, \dots, w_q) \rightarrow End(X)(w;w_1, \dots,
w_q)$$ is the identity.
\item
$End(X)$-composition\index{$End(X)$-composition} is
$End(I)$-equivariant\index{equivariant}. \newline If $f:\{1, \dots,
m\} \rightarrow \{1, \dots, \ell\}$ is a function, $w,w_i, w_{ij}
\in Mor_{End(I)^k}(m,1)$, $\alpha \in End(X)(w;w_1, \dots w_q)$, and
$\alpha_j \in End(X)(w_j;w_{j1}, \dots, w_{jp_j})$ for $j=1, \dots,
q$, then
$$\gamma(\alpha_f;(\alpha_1)_f, \dots, (\alpha_q)_f)=\gamma(\alpha;\alpha_1,
\dots, \alpha_q)_f.$$
\item
$End(X)$-functoriality\index{$End(X)$-functoriality} and
$End(I)$-functoriality\index{$End(I)$-functoriality} commute.
\newline For functions $\iota:\{1, \dots, p\} \rightarrow \{1,
\dots, q\}$ and $f:\{1, \dots, m\} \rightarrow \{1, \dots, \ell\}$
and morphism $\alpha \in End(X)(w;w_{\iota (1)}, \dots, w_{\iota
(p)})$ we have $(\alpha^{\iota})_f=(\alpha_f)^{\iota}$.
\item
$End(X)$-functoriality\index{$End(X)$-functoriality} and
$End(I)$-substitution\index{$End(I)$-substitution} commute.
\newline The diagram
$$\xymatrix@R=3pc@C=3pc{End(X)(w;w_{\iota (1)}, \dots, w_{\iota (p)}) \ar[r]^-{()^{\iota}}
\ar[d]_{(u_1, \dots, u_m)^*} & End(X)(w;w_1, \dots, w_q)
\ar[d]^{(u_1, \dots, u_m)^*}
\\ End(X)(v;v_{\iota (1)}, \dots, v_{\iota (p)}) \ar[r]_-{()^{\iota}}
& End(X)(v;v_1, \dots, v_q) }$$ commutes.
\item
$End(I)$-functoriality\index{$End(I)$-functoriality} and
$End(I)$-substitution\index{$End(I)$-substitution} commute, in the
sense that if $f_i:\{1, \dots, k_i\} \rightarrow \{1, \dots, k_i'\}$
are functions and $u_i \in End(I)(k_i)$ for $i=1, \dots, m$ and $w,
w_1, \dots, w_q \in End(I)^k(m)$, then the diagram below commutes.

\begingroup
\vspace{-2\abovedisplayskip} \small
$$\xymatrix@R=3pc@C=.05pc{\hspace{1.5cm} End(X)(w;w_1, \dots, w_q) \ar[r]^-{(u_1, \dots, u_m)^*}
\ar[dr]_{((u_1)_{f_1}, \dots, (u_m)_{f_m})^* \hspace{3em}} &
End(X)(v;v_1, \dots, v_q) \ar[d]^{()_{f_1+\cdots + f_m}}
\\ & End(X)(v_{f_1+\cdots + f_m};(v_1)_{f_1+\cdots + f_m}, \dots,
(v_q)_{f_1+\cdots + f_m})}$$
\endgroup
\noindent Note that
$$ \aligned \gamma^{\times k}(w;(u_1)_{f_1}^{\times k}, \dots,
(u_m)_{f_m}^{\times k}) &=
\gamma^{\times k}(w;u_1, \dots, u_m)_{f_1+ \cdots + f_m}  \\
& = v_{f_1 + \cdots + f_m}.
\endaligned$$
\item
$End(I)$-functoriality\index{$End(I)$-functoriality} and
$End(I)$-substitution\index{$End(I)$-substitution} commute, in the
sense that if $f:\{1, \dots, m\} \rightarrow \{1, \dots, \ell\}$ is
a function and $u_i \in End(I)(k_i)$ for $i=1, \dots, \ell$, then
the diagram
$$\xymatrix@R=3pc@C=3pc{End(X)(w;w_1, \dots, w_q) \ar[r]^-{()_f}
\ar[d]_{(u_{f1}, \dots, u_{fm})^*} & End(X)(w_f;(w_1)_f, \dots,
(w_q)_f) \ar[d]^{(u_1, \dots, u_{\ell})^*}
\\ End(X)(v;v_1, \dots, v_q) \ar[r]_-{()_{\bar{f}}}
& End(X)(v_{\bar{f}};(v_1)_{\bar{f}}, \dots, (v_q)_{\bar{f}}) }$$
commutes, where $v=\gamma^{\times k}(w;u_{f1}, \dots, u_{fm})$ and
$v_{\bar{f}}=\gamma^{\times k}(w_f;u_1, \dots, u_{\ell})$ \etc
\item
$End(I)$-substitution\index{$End(I)$-substitution} and
$End(I)$-composition\index{$End(I)$-composition} commute.
\newline Let $w, w_i, w_{ij} \in Mor_{End(I)^k}(m,1)$ and $u_i \in
End(I)(k_i)$ for $i=1, \dots, m$. Let $\alpha \in End(X)(w;w_1,
\dots, w_q)$, $\alpha_{\ell} \in End(X)(w_{\ell};w_{\ell 1}, \dots,
w_{\ell p_{\ell}})$ for $\ell=1, \dots, q$ and $\beta:=(u_1, \dots,
u_m)^* \alpha$ \etc Then
$$(u_1, \dots, u_m)^* \gamma(\alpha;\alpha_1, \dots, \alpha_q)=
\gamma(\beta;\beta_1, \dots, \beta_q).$$
\end{enumerate}
\end{thm}

This concludes our motivational discussion of the 2-theory $End(X)$
fibered over the theory $End(I)$ for a 2-functor $X:I^2 \rightarrow
Cat$. Next we turn to the general discussion.
\index{2-theory!endomorphism 2-theory|)}\index{theory!endomorphism
2-theory|)}\index{2-theory!$(End(X),End(I))$|)}\index{$(End(X),End(I))$|)}\index{$End(X)$|)}
\index{endomorphism 2-theory|)}

\section{2-Theories and Algebras over 2-Theories}
A general 2-theory has all of the properties described in the
example above. P. Hu and I. Kriz introduce the notion of a 2-theory
in \cite{hu} as follows.

\begin{defn}
A {\it 2-theory $\Theta$ fibered over the theory
$T$}\index{2-theory!2-theory fibered over a
theory|textbf}\index{theory!2-theory fibered over a
theory|seeonly{2-theory}}, written $(\Theta,T)$ for short, is a
natural number $k$, a theory $T$, and a
contravariant\index{contravariant} functor $\Theta:T \rightarrow
Cat$ from the category $T$ to the 2-category $Cat$ of small
categories such that
\begin{itemize}
\item
$Obj \hspace{1mm}  \Theta(m) = \coprod_{n \geq 0} Mor_{T^k}(m,n)$
for all $m \in \mathbb{N}$, where $T^k$ is the theory with the same
objects as $T$, but with $Mor_{T^k}(m,n)=Mor_{T}(m,n)^k$
\item
If $w_1, \dots, w_n \in Mor_{T^k}(m,1)$, then the word in
$Mor_{T^k}(m,n)$ with which the $n$-tuple $w_1, \dots, w_n$ is
identified is the product in $\Theta(m)$ of $w_1, \dots, w_n$
\item
For $w \in Mor_T(m,n)$ the functor $\Theta(w):\Theta(n) \rightarrow
\Theta(m)$ is $\Theta(w)(v)=v \circ w^{\times k}$ on objects $v \in
Mor_{T^k}(n,j)$.
\end{itemize}
For objects $w_1, \dots, w_n, w \in Mor_{T^k}(m,1) \subseteq Obj
\hspace{1mm} \Theta(m)$ we set
$$\Theta(w;w_1, \dots, w_n):=Mor_{\Theta(m)}(\prod_{i=1}^n w_i,w).$$
\end{defn}

The second condition explains the choice of notation $\prod_{i=1}^n
w_i$. Given a 2-theory such as this, it has operations and relations
as in Theorem \ref{2-theoryoperations}. Vice-a-versa, given sets
$\Theta(w;w_1, \dots, w_n):=Mor_{\Theta(m)}(\prod_{i=1}^n w_i,w)$
with operations and relations as in Theorems
\ref{2-theoryoperations} and \ref{2-theoryrelations} we get a
2-theory. We refer to these operations and relations as the {\it
operations and relations of 2-theories}\index{operations and
relations of 2-theories}. Recall that a pseudo algebra $I$ over a
theory $T$ is a category such that for every word $w \in T(n)$ we
have a functor $\Phi_n(w):I^n \rightarrow I$. Moreover, for every
operation of theories (composition, substitution, and identity) we
have a coherence iso and for every relation of theories we have a
coherence diagram\index{operations and relations of theories}. A
pseudo $(\Theta, T)$-algebra can be defined analogously.

\begin{defn}
Let $(\Theta,T)$ be a 2-theory. A {\it pseudo $(\Theta,T)$-algebra
over $I^k$} consists of the following data\index{algebra!pseudo
algebra|textbf}\index{algebra!pseudo $(\Theta,T)$-algebra over
$I^k$}\index{algebra!pseudo algebra over a
2-theory|textbf}\index{algebra!pseudo algebra over
$(\Theta,T)$|textbf}:
\begin{itemize}
\item
a small pseudo $T$-algebra $I$ with structure maps {\small
$\Phi:T(n) \rightarrow Functors(I^n,I)$}
\item
a strict 2-functor $X:I^k \rightarrow Cat$
\item
set maps {\small $\phi:\Theta(w;w_1, \dots, w_n) \rightarrow
End(X)(\Phi(w);\Phi(w_1),
 \dots, \Phi(w_n))$}, where $\Phi(w)$ means to apply $\Phi$ to each component
of $w$ to make $I^k$ into the product pseudo $T$-algebra of $k$
copies of $I$
\item
a coherence iso modification\index{coherence
isomorphism|textbf}\index{coherence iso modification} for each
operation of 2-theories and these coherence iso
modifications\index{coherence iso modification} satisfy coherence
diagrams indexed\index{indexed} by the relations of
2-theories\index{operations and relations of 2-theories}.
\end{itemize}
\end{defn}

A morphism of pseudo $(\Theta,T)$-algebras over $I^k$ is similar to
a morphism of pseudo $T$-algebras.

\begin{defn}
Let $X,Y:I^k \rightarrow Cat$ be pseudo $(\Theta,T)$-algebras over
$I^k$. Then a {\it morphism $H:X \rightarrow Y$ of pseudo
$(\Theta,T)$-algebras over $I^k$}\index{morphism!morphism of pseudo
$(\Theta,T)$-algebras|textbf} is a strict 2-natural transformation
$H:X \Rightarrow Y$ with coherence iso modifications\index{coherence
isomorphism|textbf}\index{coherence iso modification}
$\rho_{\alpha}$ indexed\index{indexed} by elements $\alpha \in
\Theta(w;w_1, \dots, w_n)$, where $w,w_1, \dots, w_n \in Obj
\hspace{1mm} \Theta(m)$.
$$\xymatrix@R=4pc@R=3pc{X \circ \Phi(w_1) \circ d^m \times \cdots \times
X \circ \Phi(w_n) \circ d^m \ar@{=>}[rr]^-{\phi_X(\alpha)}
\ar@<11ex>@{=>}[d]^{H * i_{\Phi(w_n)} * i_{d^m}}
 \ar@<-11ex>@{=>}[d]_{H * i_{\Phi(w_1)} * i_{d^m}} \ar@{}[d]|{\cdots}
& & X \circ \Phi(w) \circ d^m \ar@{=>}[d]^{H * i_{\Phi(w)}* i_{d^m}}
\ar@{~>}[ld]_{\rho_{\alpha}}
\\ Y \circ \Phi(w_1) \circ d^m \times \cdots \times
Y \circ \Phi(w_n) \circ d^m \ar@{=>}[rr]_-{\phi_Y(\alpha)} & & Y
\circ \Phi(w) \circ d^m }$$ The coherence iso
modification\index{coherence iso modification} $\rho_{\alpha}$ is
required to commute with all coherence iso modifications of the
pseudo algebra structure.
\end{defn}

The 2-cells of pseudo $(\Theta,T)$-algebras over $I^k$ are also
similar to the 2-cells of pseudo $T$-algebras.

\begin{defn}
Let $G,H:X \rightarrow Y$ be morphisms of pseudo
$(\Theta,T)$-algebras over $I^k$. Then a {\it
2-cell}\index{2-cell!2-cell in the 2-category of pseudo
$(\Theta,T)$-algebras|textbf} $\sigma: G \Rightarrow H$ is a
modification which commutes with the coherence iso modifications
$\rho^G$ and $\rho^H$ appropriately.
\end{defn}

\begin{thm}
The pseudo $(\Theta,T)$-algebras over $I^k$ form a 2-category.
\end{thm}
\begin{pf}
Routine.
\end{pf}

\section{The Algebraic Structure of Rigged Surfaces} \label{sec:rigged}
The purpose of this section is to introduce the category of rigged
surfaces\index{rigged surface|(textbf} as an example of a pseudo
algebra\index{algebra!pseudo algebra!examples of pseudo algebras}
over a 2-theory fibered over a theory\index{2-theory!2-theory
fibered over a theory} and to describe its stack\index{stack}
structure. This approach was introduced in \cite{hu} by P. Hu and I.
Kriz. In their terminology, a smooth, compact, not necessarily
connected, 2-dimensional manifold\index{manifold} $x$ with a complex
structure\index{complex structure} is called a {\it rigged surface}
if each boundary component $k$ comes equipped with a
parametrization\index{parametrization}
diffeomorphism\index{diffeomorphism} $f_k:S^1 \rightarrow k$ which
is analytic with respect to the complex structure on $x$, \ie the
diffeomorphism\index{diffeomorphism} $f_k$ extends to a
holomorphic\index{holomorphic} map when we go into local
coordinates. A boundary component\index{boundary components} $k$ is
called {\it inbound}\index{inbound|textbf} or {\it
outbound}\index{outbound|textbf} depending on the
orientation\index{orientation} of its
parametrization\index{parametrization} $f_k$ with respect to the
orientation\index{orientation} on $k$ induced by the complex
structure\index{complex structure}. The convention is to call the
identity parametrization\index{parametrization} of the boundary of
the unit disk {\it inbound}\index{inbound|textbf}. A {\it
morphism}\index{morphism!morphism of rigged surfaces|textbf} of
rigged surfaces is a holomorphic\index{holomorphic}
diffeomorphism\index{diffeomorphism} which preserves the boundary
parametrizations.\index{boundary
parametrization}\index{parametrization}

The structure of the category of rigged surfaces has the following
features, which were studied in \cite{hu}. For finite sets $a$ and
$b$, let $Obj \hspace{1mm}  X_{a,b}$ denote the set of rigged
surfaces $x$ equipped with a bijection between the inbound boundary
components\index{boundary components} of $x$ and $a$ as well as a
bijection between the outbound boundary components of $x$ and $b$.
For $x,y \in Obj \hspace{1mm} X_{a,b}$, let $Mor_{X_{a,b}}(x,y)$ be
the morphisms of rigged surfaces which preserve the bijections with
$a$ and $b$. For finite sets $a,b,c,$ and $d$ we can take the
disjoint union\index{disjoint union|textbf} of any two rigged
surfaces $x \in Obj \hspace{1mm}  X_{a,b}$ and $y \in Obj
\hspace{1mm} X_{c,d}$ and the result is an element of $Obj
\hspace{1mm}  X_{a \coprod c, b \coprod d}$. We can apply this
process to morphisms as well, and we get a functor $\coprod: X_{a,b}
\times X_{c,d} \rightarrow X_{a \coprod c, b \coprod d}$ called {\it
disjoint union}\index{disjoint union|textbf}. Note that this functor
is indexed\index{indexed|textbf} by the finite sets $a,b,c,$ and
$d$. For finite sets $a,b,$ and $c$ we also have a {\it gluing
functor}\index{gluing|textbf} $\check{?}:X_{a \coprod c, b \coprod
c} \rightarrow X_{a,b}$ which identifies an inbound boundary
component $k$ with an outbound\index{outbound} boundary component
$k'$ according to $f_k(z) \sim f_{k'}(z)$ for all $z \in S^1$
whenever $k$ and $k'$ are labelled\index{label} by the same element
of $c$. This gluing\index{gluing|textbf} functor is also
indexed\index{indexed|textbf} by the finite sets $a,b,$ and $c$.
There is also a {\it unit $0$} in $X_{0,0}$ given by the empty
set\index{empty set}. These disjoint union\index{disjoint
union|textbf} functors, gluing\index{gluing|textbf} functors, and
unit\index{unit|textbf} along with their coherence
isos\index{coherence isomorphism} and coherence
diagrams\index{coherence diagram} give the category of rigged
surfaces the structure of a {\it pseudo algebra\index{algebra!pseudo
algebra} over the 2-theory of commutative monoids with
cancellation}. More precisely, if $I$ denotes the category of finite
sets and bijections, then the assignment $(a,b) \mapsto X_{a,b}$
defines a strict 2-functor $X:I^2 \rightarrow Cat$ which is a pseudo
algebra over the 2-theory which we now describe.

We define the 2-theory $(\Theta,T)$ of {\it commutative monoids with
cancellation} as follows\index{2-theory!2-theory of commutative
monoids with cancellation|textbf}\index{commutative monoid with
cancellation|textbf}. Let $T$ be the theory of commutative
monoids\index{theory!theory of commutative
monoids}\index{commutative monoid!theory of commutative monoids} and
let $+:2 \rightarrow 1$ and $0:0 \rightarrow 1$ be the usual words
in the theory of commutative monoids. Let $k=2$. The 2-theory
$\Theta$ is generated by three words: {\it addition}
$+$\index{addition|textbf}, {\it
cancellation}\index{cancellation|textbf} $\check{?}$, and {\it
unit}\index{unit|textbf} $0$. These are described in terms of a
general algebra $X:I^2 \rightarrow Sets$ over $(\Theta, T)$ as
follows. Note that $+$ and $0$ have two meanings.
\index{$+$}\index{$\check{?}$}\index{0 \protect\unskip}
$$+:X_{a,b} \times X_{c,d} \rightarrow X_{a+c, b +d} $$
$$\check{?}:X_{a+c,b+c} \rightarrow X_{a,b} $$
$$0 \in X_{0,0}  $$

These generating words\index{generating words}\index{word!generating
words} must satisfy the following axioms.
\begin{enumerate}
\item
The word $+$ is {\it commutative}\index{commutative|textbf}. $$
\xymatrix@R=3pc@C=3pc{X_{a,b} \times X_{c,d} \ar[r]^{+} \ar[d] &
X_{a+c, b +d} \ar@{=}[d]
\\ X_{c,d} \times X_{a,b} \ar[r]_{+} & X_{c+a,d+b}}$$
\item
The word $+$ is {\it associative}\index{associative|textbf}.
$$\xymatrix@C=5pc@R=3pc{ (X_{a,b} \times X_{c,d}) \times X_{e,f} \ar[d]
\ar[r]^{+ \times 1_{X_{e,f}}} & X_{a+c,b+d} \times X_{e,f} \ar[d]^+
\\ X_{a,b} \times (X_{c,d} \times X_{e,f})
\ar[d]_{1_{X_{a,b}}\times +} & X_{(a+c)+e,(b+d)+f} \ar@{=}[d]
\\ X_{a,b} \times X_{c+e,d+f} \ar[r]_+ & X_{a+(c+e),b+(d+f)}
 }$$
\item
The word $+$ has {\it unit}\index{unit|textbf} $0 \in X_{0,0}$.
$$\xymatrix@R=3pc@C=3pc{X_{a,b} \times \{0\} \ar[r]^{ +  } \ar[dr]_{pr_1} &
X_{a+0,b+0} \ar@{=}[d]
\\ & X_{a,b}} $$
\item
The word $\check{?}$ is {\it transitive}\index{transitive|textbf}.
$$\xymatrix@C=3pc@R=3pc{X_{(a+c)+d, (b+c)+d} \ar@{=}[d] \ar[r]^-{\check{?}} & X_{a+c,b+c}
\ar[d]^{\check{?}} \\
X_{a+(c+d),b+(c+d)} \ar[r]_-{\check{?}} & X_{a,b} }$$
\item
The word $\check{?}$ {\it distributes}\index{distributivity|textbf}
over the word $+$.
$$\xymatrix@R=3pc@C=3pc{X_{a+c,b+c} \times X_{e,f} \ar[r]^{+} \ar[dd]_{\check{?} \times
1_{X_{e,f}}} & X_{(a+c)+e,(b+c)+f} \ar@{=}[d] \\ &
X_{(a+e)+c,(b+f)+c} \ar[d]^{\check{?}}
\\ X_{a,b} \times X_{e,f}
\ar[r]_{+} & X_{a+e,b+f} }$$
\item
Trivial cancellation\index{trivial
cancellation|textbf}\index{cancellation!trivial cancellation|textbf}
is trivial.
$$\xymatrix@1@C=3pc@R=3pc{X_{a+0,b+0} \ar[r]^-{\check{?}} \ar@{=}[dr] & X_{a,b}
\ar[d]^{1_{X_{a,b}}} \\ & X_{a,b}}$$
\end{enumerate}

The category of rigged surfaces forms a pseudo
algebra\index{algebra!pseudo algebra!examples of pseudo algebras}
over this 2-theory of commutative monoids with
cancellation\index{2-theory!2-theory of commutative monoids with
cancellation}. The category $I$ of finite sets and bijections
equipped with the operation $\coprod$ is a pseudo
algebra\index{algebra!pseudo algebra!examples of pseudo algebras}
over the theory $T$ of commutative monoids\index{theory!theory of
commutative monoids}\index{commutative monoid!theory of commutative
monoids}. The pseudo algebra structure on $X:I^2 \rightarrow Cat$ is
given by assigning a fixed choice of $\coprod$ to $+$,
gluing\index{gluing} of manifolds to $\check{?}$, and the empty set
to 0. This defines the structure maps $\Theta(w;w_1, \dots, w_n)
\rightarrow End(X)(w;w_1, \dots, w_n)$.

In \cite{hu} and \cite{hu1} the algebraic structure of
holomorphic\index{holomorphic}\index{rigged surface!holomorphic
families of rigged surfaces|textbf} families of rigged surfaces is
captured by a stack\index{algebra!pseudo algebra!stacks of pseudo
algebras}\index{stack!stack of pseudo algebras|textbf} of pseudo
algebras over the 2-theory of commutative monoids with
cancellation\index{cancellation}\index{SPCMC|textbf}\index{2-theory!2-theory
of commutative monoids with cancellation}, which is also called a
stack of lax\index{lax} commutative monoids with cancellation
(SLCMC)\index{cancellation}\index{SLCMC|textbf}\index{stack!stack of
lax commutative monoids with cancellation}. We describe this stack
now. Let $\mathcal{B}$ be the category of finite dimensional complex
manifolds\index{manifold!complex manifold} with morphisms
holomorphic\index{holomorphic} maps. A collection $\{B_i \rightarrow
B\}_i$ of (open) holomorphic\index{holomorphic} embeddings are a
cover\index{Grothendieck cover} if their combined image covers $B$.
This makes $\mathcal{B}$ into a Grothendieck
topology\index{Grothendieck topology}. For any finite dimensional
complex manifold\index{manifold!complex manifold}\index{complex
manifold} $B$ let $I^B$ denote the category of covering
spaces\index{covering space} of $B$ with finite fibers and morphisms
given by isomorphisms of covering spaces. The category $I^B$ is a
pseudo commutative monoid\index{commutative monoid!pseudo
commutative monoid} under $\coprod$. Let $s$ and $t$ be objects of
$\mathcal{B}$. Define $X^B_{s,t}$ as the category of
holomorphic\index{holomorphic} families\index{rigged
surface!holomorphic families of rigged surfaces|textbf} of rigged
surfaces over $B$ with inbound\index{inbound} boundary
components\index{boundary components} labelled\index{label} by the
covering space\index{covering space} $s$ of $B$ and
outbound\index{outbound} boundary components labelled\index{label}
by the covering space $t$ of $B$. Such a holomorphic\index{rigged
surface!holomorphic families of rigged
surfaces|textbf}\index{holomorphic}\index{holomorphic families of
rigged surfaces|textbf} family $x$ is by definition a complex
manifold\index{manifold!complex manifold}\index{complex manifold}
$x$ with analytic boundary and a transverse holomorphic map $p:x
\rightarrow B$ such that $x_b=p^{-1}(b)$ is a rigged surface for all
$b \in B$. Moreover, the boundary parametrizations\index{boundary
parametrization}\index{parametrization} of $p^{-1}(b)$ vary
holomorphically\index{holomorphic} with $b$ in the precise sense on
page 330 of \cite{hu}. To say that the inbound\index{inbound}
boundary components\index{boundary components} of $x$ are
labelled\index{label} by the covering space\index{covering space}
$s$ means that for each $b \in B$ the rigged surface $x_b$ is
equipped with a bijection between its inbound\index{inbound}
boundary components and the fiber of $s$ over $b$. The explanation
for the covering space $t$ labelling\index{label} the
outbound\index{outbound} boundary components\index{boundary
components} is similar. With these definitions as well as disjoint
union\index{disjoint union}, gluing\index{gluing}, and empty
set\index{empty set}, the functor $X^B:(I^B)^2 \rightarrow Cat$ is a
pseudo algebra\index{algebra!pseudo algebra!examples of pseudo
algebras} over the 2-theory of commutative monoids with
cancellation\index{cancellation}\index{2-theory!2-theory of
commutative monoids with cancellation}.

Let $\mathcal{C}$ denote the 2-category of pseudo
algebras\index{algebra!pseudo algebra} over the 2-theory of
commutative monoids with
cancellation\index{cancellation}\index{2-theory!2-theory of
commutative monoids with cancellation}. This 2-category admits
bilimits\index{bilimit}, which we prove in a special case in the
next section. Define a contravariant pseudo
functor\index{contravariant}\index{functor!pseudo functor}
$G:\mathcal{B} \rightarrow \mathcal{C}$ by taking a finite
dimensional complex manifold\index{manifold!complex
manifold}\index{complex manifold} $B$ to the pseudo
algebra\index{algebra!pseudo algebra} $X^B$ over the 2-theory of
commutative monoids with cancellation\index{cancellation} with
underlying pseudo commutative monoid $I^B$. Then $G$ takes
Grothendieck covers\index{Grothendieck cover} to bilimits because it
does so on the underlying categories comprising the pseudo
algebras\index{algebra!pseudo algebra}. Hence $G$ is a
stack\index{stack}. It is in this sense that the category of rigged
surfaces forms a stack.\index{rigged surface|)}

\section{Weighted Pseudo Limits of Pseudo $(\Theta,T)$-Algebras}
The 2-category of pseudo $(\Theta,T)$-algebras admits
weighted\index{limit!weighted pseudo limit}\index{weighted} pseudo
limits, just like the 2-category of pseudo $T$-algebras. In the
following theorem we prove this for pseudo $(\Theta,T)$-algebras
with fixed underlying pseudo $T$-algebra $I^k$. The proof can be
modified to the general case of pseudo $(\Theta, T)$-algebras with
different underlying pseudo $T$-algebras by taking the pseudo limit
of the underlying pseudo $T$-algebras as well.

\begin{thm}
Let $\mathcal{J}$ be a 1-category and $\mathcal{C}$ the 2-category
of pseudo $(\Theta,T)$-algebras over $I^k$. Let $F:\mathcal{J}
\rightarrow \mathcal{C}$ be a pseudo functor. Then $F$ admits a
pseudo limit $(X,\pi)$ in $\mathcal{C}$, where $\pi:\Delta_X
\Rightarrow F$ is a universal pseudo cone.\index{limit!pseudo
limit|textbf}\index{algebra!pseudo algebra!pseudo limits of pseudo
algebras|textbf}
\end{thm}
\begin{pf}
Let $\gamma$ and $\delta$ be the 2-cells in $\mathcal{C}$ which make
$F$ into a pseudo functor. For each $j \in Obj \hspace{1mm}
\mathcal{J}$, let $X^j:I^k \rightarrow Cat$ be the strict 2-functor
belonging to the pseudo $(\Theta,T)$-algebra $Fj$. Then for each
fixed object $i \in Obj \hspace{1mm} I^k$ and each object $j \in
\mathcal{J}$ we have a category $X^j_i$. For each morphism $f:j
\rightarrow m$ in $\mathcal{J}$, the map $Ff:X^j \Rightarrow X^m$ is
a strict 2-natural transformation which gives us a functor
$(Ff)_i:X^j_i \rightarrow X^m_i$ for each $i \in Obj \hspace{1mm}
I^k$. Thus for fixed $i$ we have a pseudo functor $F_i:\mathcal{J}
\rightarrow Cat$ defined by $j \mapsto X^j_i$ and $f \mapsto
(Ff)_i$. The coherence isos of $F_i$ are the coherence iso
modifications of $F$ evaluated at $i$.

Let $X_i:=PseudoCone(\mathbf{1},F_i)$, where $\mathbf{1}$ is the
terminal object\index{terminal object} in the category of small
categories. Then it is known from Chapter \ref{sec:laxlimitsinCat}
that $X_i$ is the pseudo limit of $F_i$ in $Cat$. Proceeding
analogously on morphisms of $I^k$, we obtain a strict 2-functor
$X:I^k \rightarrow Cat$ defined by $i \mapsto X_i$. More precisely,
if $h:i_1 \rightarrow i_2$ is a morphism in $I^k$ and $\eta \in Obj
\hspace{1mm} X_{i_1}$, then $X_h(\eta)(j):=X^j_h(\eta(j))$ for $j
\in Obj \hspace{1mm}  \mathcal{J}$.

A more conceptual way to view the construction of the strict
2-functor $X:I^k \rightarrow Cat$ is the following. For $i \in I^k$,
let $F_i:\mathcal{J} \rightarrow Cat$ be the pseudo functor from
above. For a morphism $h:i_1 \rightarrow i_2$ in $I^k$, let
$F_h:F_{i_1} \Rightarrow F_{i_2}$ be the pseudo natural
transformation given by $F_h(j):= X_h^j$. The pseudo natural
transformation $F_h$ is actually strictly 2-natural because $Ff:X^j
\Rightarrow X^m$ is a strict 2-natural transformation for each $f:j
\rightarrow m$ in $\mathcal{J}$. Thus $i \mapsto F_i$ and $h \mapsto
F_h$ define a strict functor $I^k \rightarrow Functors(\mathcal{J},
Cat)$. Now recall that $PseudoCone(\mathbf{1}, -)$ is a covariant
functor from $Functors(\mathcal{J},Cat)$ to $Cat$. The composition
$$\xymatrix@1@C=7pc{ I^k \ar[r] & Functors(\mathcal{J}, Cat)
\ar[r]^-{PseudoCone(\mathbf{1},-)} & Cat }$$ is $X:I^k \rightarrow
Cat$.

We claim that this 2-functor $X:I^k \rightarrow Cat$ has the
structure of a pseudo $(\Theta, T)$-algebra. The argument is like
Lemma \ref{limitisalgebra}, although the coherences need some care.
First we define maps $\phi:\Theta(w;w_1, \dots, w_n) \rightarrow
End(X)(\Phi(w);\Phi(w_1), \dots, \Phi(w_n))$, where $w_1, \dots,
w_n, w \in Mor_{T^k}(m,1)$. Let $\alpha \in \Theta(w;w_1, \dots,
w_n)$. We need to define a natural transformation
$$\phi(\alpha): X \circ \Phi(w_1) \circ d^m \times \cdots \times
X \circ \Phi(w_n) \circ d^m \Rightarrow X \circ \Phi(w) \circ d^m$$
``componentwise,'' where $d^m:I^m \rightarrow (I^m)^k$ is the
diagonal functor. Let $$\xymatrix@1{\phi_j:\Theta(w;w_1, \dots, w_n)
\ar[r] & End(X^j)(\Phi(w);\Phi(w_1), \dots, \Phi(w_n))}$$ be the
maps that make $X^j:I^k \rightarrow Cat$ into a pseudo
$(\Theta,T)$-algebra for each $j \in Obj \hspace{1mm} \mathcal{J}$.
Let $i \in I^m$. We define a functor
$$\xymatrix@1{(\phi(\alpha))_i: X_{\Phi(w_1) \circ d^m (i)} \times \cdots \times
X_{\Phi(w_n) \circ d^m (i)} \ar[r] & X_{\Phi(w) \circ d^m (i)}}$$
and show that $i \mapsto (\phi(\alpha))_i$ is natural. Recall that
objects of $$X_{\Phi(w_{\ell}) \circ d^m (i)}=PseudoCone(\mathbf{1},
F_{\Phi(w_{\ell}) \circ d^m (i)})$$ can be identified with a subset
of
$$\{(a_j)_j \times (\varepsilon_f)_f \in
\prod_{j \in Obj \hspace{1mm}  \mathcal{J}} Obj \hspace{1mm}
X^j_{\Phi(w_{\ell}) \circ d^m (i)} \times \prod_{f \in Mor
\hspace{1mm} \mathcal{J}} Mor \hspace{1mm}  X^{Tf}_{\Phi(w_{\ell})
\circ d^m (i)} \vert $$ $$
 \varepsilon_f:(Ff)_{\Phi(w_{\ell}) \circ d^m (i)}(a_{Sf})\rightarrow a_{Tf}
\text{ is iso for all } f \in Mor \hspace{1mm}  \mathcal{J} \}$$ by
Remark \ref{PL1}. A similar statement holds for morphisms according
to \ref{PL2}. Let $\eta^{\ell} = (a^{\ell}_j)_j \times
(\varepsilon^{\ell}_f)_f \in Obj \hspace{1mm}  X_{\Phi(w_{\ell})
\circ d^m (i)}$ and $(\xi^{\ell}_j)_j \in Mor \hspace{1mm}
 X_{\Phi(w_{\ell}) \circ d^m (i)}$ for $1 \leq \ell \leq n$.
Define $$a_j:=(\phi_j(\alpha))_i(a_j^1, \dots, a_j^n)$$ and
$$\varepsilon_f:=(\phi_{Tf}(\alpha))_i (\varepsilon^1_f, \dots,
\varepsilon^n_f) \circ (\rho^{Ff}_{\alpha})_i(a^1_{Sf}, \dots,
a^n_{Sf}).$$ Note that
$$\xymatrix@1{(\rho^{Ff}_{\alpha})_i(a^1_{Sf}, \dots, a^n_{Sf}):
(Ff)_{\Phi(w_{\ell}) \circ d^m (i)}(\phi_{Sf}(\alpha))_i(a^1_{Sf},
\dots, a^n_{Sf}) \ar[r] &} $$
$$(\phi_{Tf}(\alpha))_i((Ff)_{\Phi(w_{\ell}) \circ d^m (i)}
(a^1_{Sf}),\dots, (Ff)_{\Phi(w_{\ell}) \circ d^m (i)}(a^n_{Sf}))$$
and the composition in the definition of $\varepsilon_f$ makes
sense. Also define $$\xi_j:=(\phi_j(\alpha))_i(\xi^1_j, \dots,
\xi^n_j).$$ Then $\phi(\alpha)$ is defined ``componentwise'' by
$$(\phi(\alpha))_i(\eta^1, \dots, \eta^n):=(a_j)_j \times (\varepsilon_f) $$
and
$$(\phi(\alpha))_i((\xi^1_j)_j, \dots, (\xi^n_j)_j):=(\xi_j)_j.$$
By an argument similar to the proof of Lemma \ref{limitisalgebra},
these images are actually in $X_{\Phi(w) \circ d^m (i)}$. Next note
that $i \mapsto (\phi(\alpha))_i$ is natural because $i \mapsto
(\phi_j(\alpha))_i$ is natural for all $j \in Obj \hspace{1mm}
\mathcal{J}$, \ie $i \mapsto (\phi(\alpha))_i$ is natural in each
``coordinate'' and is therefore natural. Hence we have constructed
set maps  $\phi:\Theta(w;w_1, \dots, w_n) \rightarrow
End(X)(\Phi(w);\Phi(w_1), \dots,
 \Phi(w_n))$.

We define the coherence iso modifications\index{coherence iso
modification} for $\phi$ to be those modifications which have the
coherence iso modifications for $\phi_j$ in the $j$-th coordinate.
For example, we define the identity modification $I_w:1_{\Phi(w)}
\rightsquigarrow \phi(1_w)$ by
$$ I_w((a_j)_j \times (\varepsilon_f)_f) := (I_w^j(a_j))_j$$
for $i \in I^m$ and $(a_j)_j \times (\varepsilon_f)_f \in X_{\Phi(w)
\circ d^m (i)}$. The arrow $I_w((a_j)_j \times (\varepsilon_f)_f)$
is an arrow
 in the category $X_{\Phi(w) \circ d^m (i)}$ by an argument like the
proof of Lemma \ref{limitisalgebra}. Similarly, we can show that
these assignments are modifications and that the coherence diagrams
are satisfied because everything is done componentwise. Hence $X:I^k
\rightarrow Cat$ has the structure of a pseudo $(\Theta,T)$-algebra.

Next we need a universal pseudo cone $\pi:\Delta_X \Rightarrow F$,
where $\Delta_X:\mathcal{J} \rightarrow \mathcal{C}$ is the constant
functor which takes everything to $X$. Define a natural
transformation $\pi_j:X \Rightarrow X^j$ by letting $\pi_j(i):X_i
\Rightarrow X^j_i$ be the projection. The natural transformation
$\pi_j$ commutes with the $(\Theta,T)$ structure maps, and so
$\pi_j$ is a morphism of pseudo $(\Theta,T)$-algebras by taking the
coherence iso modifications\index{coherence iso modification} to be
trivial. The assignment $j \mapsto \pi_j$ is pseudo natural with
coherence 2-cell $\tau_{j,m}(f):Ff \circ \pi_j \Rightarrow \pi_m$
for each $f:j \rightarrow m$ in $\mathcal{J}$ as in the 1-theory
case. A similar argument to the 1-theory case shows that
$\tau_{j,m}(f)$ is a 2-cell in $\mathcal{C}$. Hence, we have a
pseudo natural transformation $\pi:\Delta_X \Rightarrow F$. We can
prove the universality of $\pi$ by applying the argument in the
lemmas leading up to Theorem \ref{isoofcategories} to $X_i
\rightarrow X_i^j$ for each fixed $i \in Obj \hspace{1mm}  I^k$ and
then passing to functors $I^k \rightarrow Cat$. We must of course
take the coherence isos into consideration.

We conclude that $(X,\pi)$ is a pseudo limit of the pseudo functor
$F:\mathcal{J} \rightarrow \mathcal{C}$.
\end{pf}

\begin{thm} \label{ThetaTlimits}
The 2-category of pseudo $(\Theta,T)$-algebras over $I^k$ admits
pseudo limits.\index{algebra!pseudo algebra!pseudo limits of pseudo
algebras|textbf}\index{limit!pseudo limit|textbf}
\end{thm}
\begin{pf}
This follows immediately from the previous theorem.
\end{pf}

\begin{lem} \label{cotensorproductsThetaT}
The 2-category $\mathcal{C}$ of pseudo $(\Theta,T)$-algebras admits
cotensor products\index{cotensor product|textbf}.
\end{lem}
\begin{pf}
Let $J \in Obj \hspace{1mm}  Cat$ and let $F:I^k \rightarrow Cat$ be
a pseudo $(\Theta,T)$-algebra. Define a strict 2-functor $P:I^k
\rightarrow Cat$ by $P_i:=(F_i)^J$, which is the 1-category of
1-functors $J \rightarrow F_i$. We claim that $P$ has the structure
of a pseudo $(\Theta,T)$-algebra. This structure is obtained by
doing the operations pointwise. Let $\phi: \Theta(w;w_1, \dots, w_n)
\rightarrow End(F)(\Phi(w);\Phi(w_1), \dots, \Phi(w_n))$ denote the
maps which make $F$ into a pseudo $(\Theta,T)$-algebra. Then define
$$\phi^P :\Theta(w;w_1, \dots, w_n) \rightarrow
End(P)(\Phi(w);\Phi(w_1), \dots, \Phi(w_n))$$
$$\phi^P(\alpha)_i(\eta^1, \dots,
\eta^n)(j):=\phi(\alpha)_i(\eta^1(j), \dots, \eta^n(j))$$ for
functors $\eta^{\ell}:J \rightarrow X_{\Phi(w_{\ell}) \circ d^m(i)}$
with $1 \leq \ell \leq n$. Coherence isos can also be defined in
this manner. Then the coherence diagrams commute because they
commute pointwise. Hence $P$ is a pseudo $(\Theta,T)$-algebra.

A proof similar to the proof of Lemma \ref{cotensorproductsT} shows
that $P$ is the cotensor product of $J$ and $F$. We must apply the
argument for $F$ in Lemma \ref{cotensorproductsT} to each $F_i$ for
$i \in Obj \hspace{1mm} I^k$.
\end{pf}

\begin{thm}
The 2-category $\mathcal{C}$ of pseudo $(\Theta,T)$-algebras admits
weighted pseudo limits.\index{algebra!pseudo algebra!pseudo limits
of pseudo algebras|textbf}\index{limit!weighted pseudo
limit|textbf}\index{weighted|textbf}
\end{thm}
\begin{pf}
By Theorem \ref{ThetaTlimits} it admits pseudo limits, and hence it
admits pseudo equalizers\index{equalizer!pseudo equalizer}. The
2-category $\mathcal{C}$ obviously admits products. By Lemma
\ref{cotensorproductsThetaT} it admits cotensor products. Hence by
Theorem \ref{streetpseudo} it admits weighted pseudo limits.
\end{pf}

\begin{thm}
The 2-category $\mathcal{C}$ of pseudo $(\Theta,T)$-algebras admits
weighted bilimits.\index{algebra!pseudo algebra!bilimits of pseudo
algebras}\index{bilimit|textbf}\index{bilimit!weighted
bilimit|textbf}\index{weighted|textbf}
\end{thm}
\begin{pf}
The 2-category $\mathcal{C}$ admits weighted pseudo limits, so it
also admits weighted bilimits.
\end{pf}

\index{product|see{2-product, bitensor product, cotensor product,
tensor product}} \index{stack!stack of lax commutative monoids with
cancellation|see{SLCMC}} \index{coequalizer|see{bicoequalizer}}
\index{limit|see{bilimit}}\index{colimit|see{bicolimit}}\index{pseudo|see{algebra,
colimit, commutative monoid, commutative monoid with cancellation,
cone, equalizer, equivalence, functor, isomorphism, limit, module,
natural transformation }} \index{lax|see{algebra, colimit, functor,
limit}}
\index{functor|see{2-functor}}\index{pullback|see{2-pullback}}
\index{ring|see{commutative semi-ring}}
\index{weighted|see{bicolimit, bilimit, colimit,
limit}}\index{adjoint|see{biadjoint,
quasiadjoint}}\index{adjunction|see{biadjunction, quasiadjunction}}

\backmatter

\providecommand{\bysame}{\leavevmode\hbox
to3em{\hrulefill}\thinspace}
\providecommand{\MR}{\relax\ifhmode\unskip\space\fi MR }
\providecommand{\MRhref}[2]{%
  \href{http://www.ams.org/mathscinet-getitem?mr=#1}{#2}
} \providecommand{\href}[2]{#2}


\begin{theindex}

  \item $(End(X),End(I))$, \textbf{147--154}
  \item $+$, 156
  \item $Alg'$, 113, 129
  \item $Cat$, 9, 138
  \item $Cat_0$, 52
  \item $End(I)$, 147
  \item $End(I)$-composition, 154
  \item $End(I)$-functoriality, 150, 152--154
  \item $End(I)$-substitution, 150, 153, 154
  \item $End(X)$, \textbf{40--41}, \textbf{57}, \textbf{147--154}
  \item $End(X)$-composition, 150, 152, 153
  \item $End(X)$-functoriality, 150, 152, 153
  \item $FiniteSets$, 65
  \item $Graph'$, 114, 129, 130
  \item $Hom$, 18
  \item $I$, 61, 66
  \item $Lex$, 15
  \item $Psd$, 19
  \item $PseudoCone$, 14
  \item $R_{G'}$, 114, 129--130
  \item $\Gamma$, 45
  \item $\check{?}$, 156
  \item $c$, 61, 66
  \item $s$, 61, 66
  \item 0 \protect\unskip, 156
  \item 0-cell, 9
  \item 1-cell, 9
  \item 2-adjoint, 127
    \subitem left 2-adjoint, 127
  \item 2-category, 3, 9
  \item 2-cell, 2, 3, 9
    \subitem 2-cell in the 2-category of pseudo $(\Theta,T)$-algebras,
        \textbf{155}
    \subitem 2-cell in the 2-category of pseudo $T$-algebras,
        \textbf{66}
  \item 2-colimit, 15
  \item 2-equalizer, 20
  \item 2-functor
    \subitem diagonal 2-functor, 112
    \subitem forgetful 2-functor, viii, 2, \textbf{113}, 117,
        \textbf{120}, 121, 126, 127
  \item 2-monad, viii, 3, 52, 61, 66, \textbf{70, 71}, 113
  \item 2-product, 19, 20, 37, 56, 58, 66, 80
  \item 2-pullback, 20
  \item 2-theory, viii, 1, 147
    \subitem $(End(X),End(I))$, \textbf{147--154}
    \subitem 2-theory fibered over a theory, \textbf{154}, 156
    \subitem 2-theory of commutative monoids with cancellation, viii, 1,
        5, 7, \textbf{156}, 158, 159
    \subitem endomorphism 2-theory, \textbf{147--154}

  \indexspace

  \item abelian groups, 113
    \subitem theory of abelian groups, 113
  \item addition, \textbf{156}
  \item adjoint, \see{biadjoint, quasiadjoint}{162}
    \subitem lax adjoint, 2, \textbf{86}
    \subitem pseudo adjoint, 2
  \item adjunction, 1, 81, \see{biadjunction, quasiadjunction}{162}
    \subitem soft adjunction, 81
  \item admits bilimits, \textbf{18}
  \item admits pseudo limits, \textbf{18}
  \item algebra, viii, 1, 39, 51
    \subitem $C$-algebra, 52, \textbf{70}
    \subitem $T$-algebra, \textbf{51}, 52
    \subitem $\mathbf{T}$-algebra, \textbf{56}
    \subitem $\mathcal{T}$-algebra, \textbf{58}, \textbf{68}
    \subitem algebra over a theory, 1, 2, \textbf{51}
    \subitem algebra over a theory enriched in groupoids, 39,
        \textbf{58}
    \subitem algebra over a theory on a set of objects, \textbf{56}
    \subitem algebra over the theory of theories, 56
    \subitem categorical $T$-algebra, 51
    \subitem free pseudo $T$-algebra on a pseudo $S$-algebra,
        \textbf{114--115}, 120, 129
    \subitem functorial $T$-algebra, 51
    \subitem lax algebra, 2, 61
    \subitem lax algebra over a 2-theory, 147
    \subitem pseudo $(\Theta,T)$-algebra over $I^k$, 155
    \subitem pseudo $T$-algebra, \textbf{61}, \textbf{68}
    \subitem pseudo algebra, viii, 1--3, 5, 7, 15, 39, 56, \textbf{61},
        63, 113, 115, 126, 130, 133, 145, 147, \textbf{155},
        156, 159
      \subsubitem bicolimits of pseudo algebras, 129, 136
      \subsubitem bilimits of pseudo algebras, \textbf{80}, 137, 161
      \subsubitem examples of pseudo algebras, 6, 65, 156, 158
      \subsubitem pseudo limits of pseudo algebras, \textbf{73},
        \textbf{80}, \textbf{159}, \textbf{161}
      \subsubitem stacks of pseudo algebras, 137, 146, 158
    \subitem pseudo algebra over $(\Theta,T)$, \textbf{155}
    \subitem pseudo algebra over $T$, \textbf{61}
    \subitem pseudo algebra over a 2-theory, 3, 147, \textbf{155}
    \subitem strict algebra, 3, 7, 52, 113
  \item arity, \textbf{39}
  \item arrow, 21
    \subitem biuniversal arrow, 2, 81, \textbf{83}, 91, 93, 96,
        \textbf{99}, \textbf{110}, 117, 121, 129
    \subitem universal arrow, \textbf{81}
  \item associative, \textbf{157}

  \indexspace

  \item basis for a Grothendieck topology, \textbf{137}
  \item biadjoint, viii, 2, 3, 81, \textbf{86}, 112
    \subitem left biadjoint, viii, 2, \textbf{86}, \textbf{99}, 113,
        \textbf{120}, 121, 126
    \subitem right biadjoint, \textbf{86}
  \item biadjunction, 81, \textbf{85}, 93, \textbf{110}
  \item bicategory, 2
  \item bicoequalizer, \textbf{20}, 136
  \item bicolimit, viii, 2, 3, 15, 18, 112, \textbf{129}, 131, 136
    \subitem example of bicolimit, 15
    \subitem weighted bicolimit, viii, 3, 20, \textbf{20}, \textbf{29},
        129, \textbf{136}
  \item bicoproduct, \textbf{20}, 136
  \item bilifting, 81
  \item bilimit, viii, 1--3, 6, 14, \textbf{14}, \textbf{38},
        \textbf{80}, 112, 139, 140, \textbf{145}, 146, 159,
        \textbf{161}
    \subitem bilimit of a diagram, 18, 139
    \subitem conical bilimit, 14
    \subitem example of bilimit, 15
    \subitem indexed bilimit, 14, 20
    \subitem weighted bilimit, 14, 20, \textbf{38}, \textbf{80},
        \textbf{161}
  \item birepresentation, 20, 81
  \item bitensor product, 3, \textbf{20}, 129, \textbf{132}, 136
  \item biuniversal arrow, 2, 81, \textbf{83}, 91, 93, 96, \textbf{99},
        \textbf{110}, 117, 121, 129
  \item boundary components, viii, 1, 156, 158
  \item boundary parametrization, 1, 156, 158

  \indexspace

  \item cancellation, viii, 1, 5--7, \textbf{156}, 158, 159
    \subitem trivial cancellation, \textbf{158}
  \item category of descent data, \textbf{141}
  \item category of small categories, 9
  \item central extension, 6
  \item chiral, 6
  \item cocycle condition, \textbf{141, 142}, 144
  \item coequalizer, 28, \see{bicoequalizer}{162}
    \subitem pseudo coequalizer, 28
  \item coherence 2-cell, \textbf{10}, 11--13, 64, 65
  \item coherence diagram, viii, 1, 2, 5, \textbf{10}, \textbf{12},
        \textbf{61--64}, 156
  \item coherence iso modification, 155, 160, 161
  \item coherence isomorphism, viii, 1, 2, 5, 16, 17, 27,
        \textbf{61--64}, 76, \textbf{155}, 156
  \item colimit, 1, 14, 15, \see{bicolimit}{162}
    \subitem lax colimit, 2
    \subitem pseudo colimit, 2, \textbf{14}, 15, 18, \textbf{21},
        \textbf{28, 29}, 31, 129
    \subitem weighted pseudo colimit, 2, \textbf{28, 29}
  \item commutative, \textbf{157}
  \item commutative monoid, 5
    \subitem pseudo commutative monoid, \textbf{65}, 158
      \subsubitem examples of pseudo commutative monoids, 5, 65
    \subitem theory of commutative monoids, 1, 65, 126, 127, 156, 158
  \item commutative monoid with cancellation, \textbf{5}, \textbf{156}
    \subitem pseudo commutative monoid with cancellation, 1, 6
      \subsubitem examples of pseudo commutative monoids with cancellation,
        5, 6
  \item commutative semi-ring
    \subitem pseudo commutative semi-ring, 6, \textbf{65}
      \subsubitem example of pseudo commutative semi-ring, 65
    \subitem theory of commutative semi-rings, 65
  \item complex manifold, 158, 159
  \item complex structure, viii, 1, 156
  \item composition, 10, 40, \textbf{41}, 55, \textbf{56}, 147, 148
  \item composition axiom, 11, 12
  \item composition of morphisms of pseudo $T$-algebras, 65
  \item composition of pseudo functors, \textbf{11}
  \item cone, 14
    \subitem pseudo cone, \textbf{14}
  \item conformal field theory, viii, 1--3, \textbf{6}, 5--7, 14
  \item congruence, \textbf{21}, 22, 24--26, 52, 114--116, 130,
        133, 134
  \item conical, 14, \textbf{19}
  \item continuous map, 9
  \item contravariant, 1, 6, 137, 138, \textbf{140}, 141, 142, 145, 146,
        148, 150, 151, 154, 159
  \item cotensor product, 2, 3, \textbf{19}, 20, \textbf{37}, 73,
        \textbf{79}, 80, \textbf{161}
  \item covering space, 7, 158
    \subitem stack of covering spaces, 7

  \indexspace

  \item descent data, viii, \textbf{141}
  \item descent object, 140
  \item diagonal 2-functor, 112
  \item diffeomorphism, 1, 156
  \item direct sum, 65
  \item directed graph, \textbf{21}, 52, 113, 114, 129, 130, 133, 139
  \item discrete category, 9
  \item disjoint union, viii, 1, 5, 6, 21, 29, 126, \textbf{156}, 158
  \item distributivity, viii, \textbf{158}

  \indexspace

  \item elementary 2-cell, 67, 68
  \item elliptic cohomology, 5
  \item empty set, 156, 158
  \item endomorphism 2-theory, \textbf{147--154}
  \item endomorphism theory, \textbf{40, 41}, \textbf{54}
  \item equalizer, 138
    \subitem pseudo equalizer, 19, 20, 37, 80, 161
  \item equivalence, \textbf{84}
    \subitem pseudo natural equivalence, \textbf{111}
  \item equivariant, 42, 152, 153
  \item exact, 138
    \subitem left exact, 15

  \indexspace

  \item factorizing 2-cell, 119
  \item finitary monad, 39
  \item forgetful 2-functor, viii, 2, \textbf{113}, 117, \textbf{120},
        121, 126, 127
  \item forgetful functor, 47, 113, 114
  \item free category, 21, 114
  \item free finitary monad, 39
  \item free functor, 113
  \item free groupoid, 21
  \item free pseudo $T$-algebra on a pseudo $S$-algebra,
        \textbf{114--115}, 120, 129
  \item free theory, 39, \textbf{47}, \textbf{56}, 59, 66, 113, 114,
        129, 130, 133
    \subitem free theory functor, \textbf{47}
  \item Freyd's Adjoint Functor Theorem, 114, 129
  \item functor, 10, \see{2-functor}{162}
    \subitem forgetful functor, 113, 114
    \subitem free functor, 113
    \subitem free theory functor, \textbf{47}
    \subitem lax functor, 2, \textbf{10}
    \subitem pseudo functor, 1--3, 6, \textbf{10}, 137, 141, 142,
        145, 146, 159
      \subsubitem composition of pseudo functors, \textbf{11}
  \item functorial, 42

  \indexspace

  \item generating words, 157
  \item Giraud stack, \textbf{142}
  \item gluing, viii, 1, 5, 6, \textbf{156}, 158
  \item graph, 21
    \subitem directed graph, \textbf{21}
  \item Grothendieck cover, 1, 6, \textbf{137}, 138, 140, 142, 145,
        158, 159
  \item Grothendieck site, 1, \textbf{137}
  \item Grothendieck topology, viii, \textbf{137}, 142, 145, 146, 158
  \item group, 39, 51
    \subitem theory of groups, 39, \textbf{51}
  \item group homomorphism, 51
  \item groupoid, 9, 29, 37--39, 56--59, 61, 66, 68

  \indexspace

  \item Hilbert space, 6
  \item Hilbert tensor product, 6
  \item holomorphic, 1, 6, 137, 156, 158
  \item holomorphic families of rigged surfaces, 1, 137, \textbf{158}
  \item homomorphism of bicategories, 81
  \item homotopy, 9
  \item horizontal composition, \textbf{10}, 13

  \indexspace

  \item inbound, 5, \textbf{156}, 158
  \item index, 14
  \item indexed, 1, 3, 6, 61, 155, \textbf{156}
  \item indexing, 1
  \item indexing category, 18, 21
  \item indices, 6
  \item initial object, 15, 45
  \item iso, 2, 82
  \item isomorphism, 84
    \subitem pseudo isomorphism, \textbf{84}, 110
    \subitem pseudo natural pseudo isomorphism, \textbf{111}

  \indexspace

  \item label, 5, \textbf{6}, 156, 158
  \item Lawvere, viii
  \item Lawvere theory, 1, 6, \textbf{39}, 56
    \subitem enriched Lawvere theory, \textbf{56}
  \item lax, \textbf{2}, \textbf{5}, 10, 14, 147, 158,
        \see{algebra, colimit, functor, limit}{162}
  \item LCMC, 6
  \item left exact, 15
  \item limit, 1, 14, \see{bilimit}{162}
    \subitem indexed pseudo limit, 3, \textbf{19}
    \subitem lax limit, 2, \textbf{14}
    \subitem pseudo limit, viii, 2, 3, \textbf{14}, \textbf{31},
        \textbf{37}, \textbf{73}, \textbf{159}, \textbf{161}
      \subsubitem example of pseudo limit, 15
    \subitem weighted pseudo limit, viii, 2, 3, \textbf{19},
        \textbf{37}, \textbf{71}, 73, \textbf{80}, 159,
        \textbf{161}
  \item line bundle, 6

  \indexspace

  \item manifold, viii, 1, 156
    \subitem complex manifold, 6, 158, 159
  \item many sorted theory, \textbf{54}
  \item modification, \textbf{13}
  \item modular functor, \textbf{6}
    \subitem one dimensional modular functor, 6
  \item module
    \subitem pseudo module, 6
  \item monad, 39, 52, 56
    \subitem finitary monad, 39
    \subitem free finitary monad, 39
  \item Moonshine, 5
  \item morphism, viii, 3
    \subitem lax morphism, 2
    \subitem morphism of $T$-algebras, \textbf{51}
    \subitem morphism of descent objects, \textbf{141}
    \subitem morphism of pseudo $(\Theta,T)$-algebras, \textbf{155}
    \subitem morphism of pseudo $T$-algebras, \textbf{64}, 121
    \subitem morphism of rigged surfaces, \textbf{156}
    \subitem morphism of theories, 2, \textbf{47}, \textbf{49, 50}, 56
    \subitem morphism of theories enriched in groupoids, \textbf{58}
    \subitem morphism of theories on a set of objects, \textbf{56}
    \subitem pseudo morphism of pseudo $T$-algebras, \textbf{64}
    \subitem pseudo morphism of theories, \textbf{63}

  \indexspace

  \item natural transformation, 3
    \subitem pseudo natural transformation, 3, 10, \textbf{11}

  \indexspace

  \item object, 21
  \item object with descent data, \textbf{141}
  \item operations and relations of 2-theories, 155
  \item operations and relations of theories, 155
  \item operations of 2-theories, 149, 150
  \item operations of theories, 45, 61
  \item orientation, 1, 156
  \item outbound, 5, 156, \textbf{156}, 158

  \indexspace

  \item parametrization, 1, 156, 158
  \item path category, 21
  \item path integral, 5
  \item product, 40,
        \see{2-product, bitensor product, cotensor product, tensor product}{162}
  \item pseudo, \textbf{5},
        \see{algebra, colimit, commutative monoid, commutative monoid with cancellation, cone, equalizer, equivalence, functor, isomorphism, limit, module, natural transformation }{162}
  \item pullback, 137, 138, \see{2-pullback}{162}

  \indexspace

  \item quantum field theory, 5
  \item quasiadjoint, 113
  \item quasiadjunction, 3, 81
    \subitem transcendental quasiadjunction, 81
  \item quasicolimit, 3
  \item quasilimit, 3
  \item quotient category, 21

  \indexspace

  \item relation, 64
  \item relations of 2-theories, 151
  \item relations of theories, 45, 61
  \item rigged surface, viii, 1, 3, 5--7, 137, 147, \textbf{156--159}
    \subitem holomorphic families of rigged surfaces, 1, 137,
        \textbf{158}
  \item ring, 113, \see{commutative semi-ring}{162}
    \subitem theory of rings, 113

  \indexspace

  \item section, 6
  \item Segal, Graeme, viii, 1, 5
  \item semi-ring, 6, \see{commutative semi-ring}{65}
  \item sheaf, 6, \textbf{137}
  \item SLCMC, \textbf{5}, 6, \textbf{158}
    \subitem central extension of SLCMC's, 6
    \subitem examples of SLCMC's, 6
    \subitem morphism of SLCMC's, 6
  \item soft adjunction, 81
  \item source, \textbf{21}
  \item SPCMC, \textbf{158}
  \item stability axiom, \textbf{137}
  \item stack, viii, 1--3, 5--7, 14, 137, 138, \textbf{140},
        \textbf{142}, \textbf{145}, 156, 159
    \subitem Giraud stack, \textbf{142}
    \subitem stack of categories, \textbf{140}
    \subitem stack of covering spaces, 7
    \subitem stack of lax commutative monoids with cancellation, 5, 158,
        \see{SLCMC}{162}
    \subitem stack of pseudo algebras, \textbf{158}
  \item state space, \textbf{6}
  \item string theory, 5
  \item substituted word, 40
  \item substitution, \textbf{41}, 54, \textbf{56}, 147, 148
  \item substitution maps, 40
  \item substitution morphism, 64
  \item symmetric monoidal category, 1, 65

  \indexspace

  \item target, \textbf{21}
  \item tensor product, 6, 19, 28, 65
  \item terminal object, 10, 14, 22, 39, 40, 45, 54, 73, 147, 159
  \item theory, viii, 1, 2, \textbf{39}, 40, \textbf{41}, \textbf{45},
        \textbf{47}, \textbf{51}, 54, 56, 146, 147
    \subitem 2-theory fibered over a theory, {\it see} 2-theory
    \subitem endomorphism 2-theory, \textbf{147--154}
    \subitem endomorphism theory, \textbf{40, 41}, \textbf{54}
    \subitem endomorphism theory enriched in groupoids, \textbf{57}
    \subitem enriched theory, \textbf{56}
    \subitem free theory, 39, \textbf{56}, 59, 66, 113, 129, 130, 133
    \subitem many sorted theory, \textbf{39}, \textbf{54}
    \subitem theory enriched in categories, 71
    \subitem theory enriched in groupoids, \textbf{56--58}, 59, 61, 66,
        68
    \subitem theory indexed over a theory, 1
    \subitem theory of abelian groups, 113
    \subitem theory of commutative monoids, 1, 65, 126, 127, 156, 158
    \subitem theory of commutative semi-rings, 65
    \subitem theory of groups, 39, \textbf{51}
    \subitem theory of rings, 113
    \subitem theory of theories, \textbf{56}, 64, 66, 68, 151
    \subitem theory on a set of objects, \textbf{54}, 56
    \subitem trivial theory, 126, 127
  \item topological space, 9
  \item trace class, 6
  \item trace map, \textbf{6}
  \item transcendental quasiadjunction, 81
  \item transitive, \textbf{157}
  \item transitivity axiom, \textbf{137}
  \item trivial cancellation, \textbf{158}
  \item trivial theory, 126
  \item tuple, 40

  \indexspace

  \item unit, viii, 5, \textbf{19, 20}, 42, \textbf{56}, 147, 148, 150,
        153, \textbf{156, 157}
  \item unit axiom, 11, 12
  \item universal arrow, \textbf{81}

  \indexspace

  \item vector space, 6, 65
  \item vertical composition, \textbf{9}, 13
  \item vertical identity, \textbf{9}

  \indexspace

  \item weight, 14
  \item weighted, viii, 2, 3, 9, 14, \textbf{19, 20}, 21,
        \textbf{28, 29}, 31, \textbf{37, 38}, \textbf{71}, 73,
        \textbf{80}, 129, \textbf{136}, 159, \textbf{161},
        \see{bicolimit, bilimit, colimit, limit}{162}
  \item word, \textbf{39}
    \subitem generating words, 157
    \subitem substituted word, 40, 148

  \indexspace

  \item Yoneda's Lemma for bicategories, 81

\end{theindex}

\end{document}